\documentclass[11pt]{article}
\usepackage{mysettings}

\usepackage{cleveref}
\usepackage{comment}

\usepackage{wrapfig}

\theoremstyle{plain}
\newtheorem{theorem}{Theorem}
\newtheorem{lemma}{Lemma}
\newtheorem{proposition}{Proposition}
\newtheorem{corollary}{Corollary}
\newtheorem{assumption}{Assumption}

\theoremstyle{definition}
\newtheorem{definition}[theorem]{Definition}

\theoremstyle{remark}
\newtheorem{remark}[theorem]{Remark}

\usepackage{bm}

\newcommand{\kminInt}{K_{\min}^{\mathcal{M}_{1}}}
\newcommand{\kmaxInt}{K_{\max}^{\mathcal{M}_{1}}}
\newcommand{\kmagInt}{\widetilde{K}_{\mathcal{M}_{1}}}

\newcommand{\kminExt}{K_{\min}^{\mathcal{M}_{2}}}
\newcommand{\kmaxExt}{K_{\max}^{\mathcal{M}_{2}}}
\newcommand{\kmagExt}{\widetilde{K}_{\mathcal{M}_{2}}}

\newcommand{\fundBound}{\kappa}

\newcommand{\ROI}{S_{\operatorname{reg}}}
\newcommand{\ROIdiam}{D_{\operatorname{reg}}}

\newcommand{\SCf}{m_{f}}
\newcommand{\smoothf}{L_{f}}

\newcommand{\compScale}{L_{\operatorname{cmp}}}

\newcommand{\goodProb}{p_{\operatorname{good}}}

\newcommand{\radiusScaling}{\alpha_{t+1}}

\newcommand{\ROIraucheC}{C_{\log}}
\newcommand{\tubeExpansionConstant}[1]{L_{#1}^{\kmagInt}}

\newcommand{\mm}[1]{\textcolor{purple}{#1}}
\usepackage{graphicx} 
\usepackage[labelfont=bf]{caption}
\usepackage{amsmath, amssymb, tabularx, caption, float, pdflscape}
\usepackage{booktabs,tabularx,longtable,array}
\usepackage{placeins} 
\hypersetup{
    linkcolor=black,
    colorlinks=true,
    citecolor=blue,
    urlcolor=violet
}
\title{Computing Fr\'{e}chet Means on Unknown Submanifolds of Riemannian Manifolds}
\author{
    Ian McPherson \thanks{Department of Applied Mathematics and Statistics, Johns Hopkins University, Baltimore, MD 21218, USA.}
    \and Mauro Maggioni \footnotemark[1]
    \thanks{Department of Mathematics, Johns Hopkins University, Baltimore, MD 21218, USA.}
}
\date{}
\begin{document}
\marginsize{1in}{1in}{1in}{1in}
\maketitle
\renewcommand{\thefootnote}{\relax}
\footnotetext{\textbf{Funding:} Maggioni is partially supported by AFOSR awards FA9550-21-1-0317 and FA9550-23-1-0445.}
\renewcommand{\thefootnote}{\arabic{footnote}}
\begin{abstract}
    Empirical averages are fundamental statistical primitives. On nonlinear spaces, their generalization, the Fr\'{e}chet mean, generally lacks a closed-form solution, and its computation requires solving an optimization problem. We study the computation of intrinsic Fr\'{e}chet means for data supported on an unknown compact $d$-dimensional Riemannian manifold $\mathcal{M}_{1}$ embedded in an ambient metric space $(X,d_{\operatorname{obs}})$. Such settings arise when observations are better represented in non-Euclidean spaces, like probability distributions in Wasserstein space. We assume access only to a local sampler on $\mathcal{M}_{1}$ and the observation metric $d_{\operatorname{obs}}$, without knowledge of its intrinsic geometry. We propose a graph-based iterative algorithm that adaptively samples the manifold and approximates geodesic interpolation. Under mild assumptions, we establish an $\mathcal{O}(1/T)$ expected objective gap guarantee using at most $\widetilde{\mathcal{O}}(T^{4d+3})$ sampled points. As a byproduct of our analysis we extend the metric equivalence guarantees of ISOMAP from Euclidean to Riemannian ambient spaces, which may be of independent interest. Numerical experiments demonstrate the effectiveness of our method for computing intrinsic averages on unknown data manifolds.
\end{abstract}

\section{Introduction}
Averaging is among the most fundamental operations in statistics. It reduces a collection of observations to a representative value and underlies estimation, aggregation, and noise reduction. In $\mathbb{R}^{D}$,
the sample mean admits two equivalent descriptions: 
the arithmetic average, 
and 
the 
minimizer of the average squared distance to the sample points. The former 
explicitly depends on the linearity of $\mathbb{R}^{D}$,
while the latter 
only on the underlying distance. 
Thus, the variational formulation
generalizes means
to general metric spaces, albeit without a closed-form expression. 

Intrinsic averaging is particularly relevant in high-dimensional data analysis. The \textit{manifold hypothesis} posits that high-dimensional data often vary along or concentrate near a manifold of substantially lower intrinsic dimension. When the ambient space is $\mathbb{R}^{D}$, the arithmetic mean may lie \textit{off} this manifold. More generally, averages formed using only the ambient geometry ignore intrinsic structure, resulting in representative averages \textit{off} the submanifold. This is worrisome for scientific inference, where these representatives should reflect real phenomena.

Motivated by this perspective, let $\mathcal M_1$ be an \textit{unknown} Riemannian submanifold of an ambient Riemannian manifold $\mathcal M_2$. Given observations $z_1,\ldots,z_N\in\mathcal M_1$, let
    $P_N:=\frac{1}{N}\sum_{r=1}^N\delta_{z_r}$
denote their empirical measure. We consider the empirical Fr\'echet problem
\begin{equation}\label{eqn:intro-frechet-problem-statement}
    \widehat{x}_N
    \in
    \operatorname*{argmin}_{x\in\mathcal M_1}
    f_{P_N}(x),
    \qquad \text{where} \quad
    f_{P_N}(x)
    :=
    \frac{1}{2N}
    \sum_{r=1}^N d_{\mathcal M_1}^2(x,z_r).
\end{equation}

Our difficulty is that $\mathcal{M}_{1}$ is itself unknown:
we only have access to a local sampling mechanism for which realized samples are interpreted as points within $\mathcal{M}_{2}$. This setup occurs within molecular simulation, where
low-dimensional structure is inferred from sampled
configuration and explored through simulation-driven sampling
\cite{bello2023staying,trstanova2020local}. This general setting poses a computational barrier, as standard first-order Riemannian optimization packages may not be used to solve \eqref{eqn:intro-frechet-problem-statement}, since we have no access to $\mathcal{M}_{1}$.
Despite this, we ask:

\begin{tcolorbox}[myboxstyle]
\centering
\textit{
Can we iteratively compute an empirical Fr\'echet mean on an unknown embedded submanifold with access to only a local sampler and ambient distances? If so, can we provide a provable end-to-end sampling and convergence guarantee?
}
\end{tcolorbox}

Under explicit geometric and
localization assumptions, we answer this in the affirmative. Our starting point is geodesic interpolation: at each iteration, sample a reference point, and interpolate a fraction along the connecting geodesic segment. This procedure is well suited for 
our setting for two reasons:
\textit{(i)} 
it is a case of
Riemannian stochastic gradient descent applied to $f_{P_{N}}$; and
\textit{(ii)} it only requires geodesic interpolation, which can be approximated with
graphs. 
The former enables 
convergence rates while the latter enables statistical guarantees.

We adhere to the following algorithmic principle: only approximate the geometry required at the current optimization step.
Rather than using a uniformly fine global graph representation of $\mathcal M_1$, our method locally refines the graph around \textit{regions of interest}, 
increasing the resolution as the algorithm progresses. We formalize this local information by an idealized sampling scheme for theoretical analysis, 
relaxing to an implementable approximation thereof for numerics.

\subsection{Contributions and Outline}
We summarize our main contribution: providing an algorithm with theoretical guarantees.

\begin{enumerate}[leftmargin=0.5cm]
    \item[] (\textbf{Graph-based computation on an unknown submanifold}) We introduce an adaptive graph-interpolation method for computing empirical Fr\'{e}chet means on unknown manifolds. Under geometric and localization assumptions, the final iterate satisfies
    \begin{equation}\label{eqn:main-contribution}
        \mathbb{E}\left[f_{P_{N}}(x) - f_{P_{N}}(\widehat{x}_{N})\right] \leq \frac{C}{\gamma + T} = \mathcal{O}(T^{-1}),
    \end{equation}
    where $C$ is a geometric constant independent of $T$, preserving the rate of 
    geodesic interpolation. 
    Sufficient per-iteration and cumulative sample sizes scale as  $\widetilde{\mathcal{O}}(t^{4d+2})$ and $\widetilde{\mathcal{O}}(T^{4d+3}),$
    respectively, up to geometric constants, where $d = \dim(\mathcal{M}_{1})$. See Theorem~\ref{thm:informal-main-thm} for an informal statement, and Theorem~\ref{thm:main-conv-thm} for the full formal statement.
\end{enumerate}

\noindent As two additional contributions of independent interest, we establish two generalizations.

\begin{enumerate}[leftmargin=0.5cm]
    \item[] (\textbf{ISOMAP-type metric recovery for Riemannian ambient $\mathcal{M}_{2}$})
    We establish a metric-recovery theorem for an embedded submanifold $\mathcal{M}_{1} \subseteq \mathcal{M}_{2}$ when $\mathcal{M}_{2}$ is Riemannian. Theorem~\ref{thm:deterministic-metric-equivalence-for-ISOMAP} gives deterministic conditions for when a graph weighted by ambient distances satisfies 
    \[
    (1-\lambda_{1})d_{\mathcal{M}_{1}}(x,y) \leq d_{G}(x,y) \leq (1+\lambda_{2})d_{\mathcal{M}_{1}}(x,y).
    \]
    The key step is establishing a curved-ambient chord--arc estimate via comparison arguments. 
    Corollary~\ref{cor:probabilistic-isomap} recovers a high-probability guarantee using a covering argument. These results capture \cite{bernstein2000graph} as a special case when $\mathcal{M}_{2} = \mathbb{R}^{D}$. See Section~\ref{sec:generalized-isomap} for details.
    
    \item[] (\textbf{Biased Riemannian SGD})
    We establish an $\mathcal{O}(t^{-1})$ convergence guarantee for Riemannian SGD with biased gradient estimators with controlled bias. 
    This provides a Riemannian analogue of classical biased-SGD analyses in Euclidean optimization and complements existing Riemannian stochastic-approximation results.
    See Section~\ref{sec:biased-sgd-conv-theorem} for details.
    
\end{enumerate}

\paragraph{Outline of Paper.} Section 2 formalizes the empirical Fr\'echet mean problem,
introduces 
our graph-based interpolation algorithm, and presents an informal version of our main result. Then we describe the implementable version used in our numerics, and provide a literature review. Section 3 develops the convergence analysis: Section 3.2 establishes a biased Riemannian SGD result, Section 3.3 proves ISOMAP-type metric-recovery guarantees, and Sections 3.4-3.5 combine these ingredients to obtain the main convergence and sampling guarantees. Section 4 contains the numerical experiments, applying the algorithm to submanifolds of $\mathbb{R}^{D}$ and $\mathcal{P}_{2,\operatorname{ac}}(\mathbb{R}^{D})$.

\section{Problem Formulation and Our Algorithmic Approach}
\subsection{Empirical Fr\'{e}chet Means and Geodesic Interpolation}\label{sec:problem-statement-and-algorithm}
\begin{algorithm}[t]
\caption{Geodesic Interpolation}\label{alg:geodesic-interpolation}
\begin{algorithmic}
\State \textbf{Input:} References $(z_{r})_{r=1}^{N}$, initialization $\widetilde{x}_{0}$, stepsizes $(\eta_{t})_{t\geq0}$.\vspace{.1cm}
\For{$t=0,1,2,\dots$} 
\State \textbf{Uniformly Pick Reference Point:} Draw $I_{t+1}\sim\mathrm{Unif}([N])$ and set
    $y_{t+1}:=z_{I_{t+1}}$;\vspace{.1cm}
\State \textbf{Update by Interpolation:} $\widetilde{x}_{t+1} \leftarrow \exp_{\widetilde{x}_{t}}(\eta_{t}\log_{\widetilde{x}_{t}}(y_{t+1}))$.\vspace{.1cm}
\EndFor
\end{algorithmic}
\end{algorithm}

We aim to solve the statistical/optimization problem of \eqref{eqn:intro-frechet-problem-statement}. 
Whenever $\widehat{x}_{N}$ is used below, we assume that the argmin 
is nonempty and fix a minimizer; later assumptions ensure this.

When $\mathcal{M}_{1} = \mathbb{R}^{d}$, the minimizer admits the familiar unique closed-form expression
\begin{equation} \label{eqn:closed-form-Euclidean-average}
    \widehat{x}_{N} = \frac{1}{N} \sum_{r=1}^{N} z_{r}.
\end{equation}
The explicit formula, derived from first-order optimality conditions,
relies fundamentally on the linearity of $\mathbb{R}^{d}$. In contrast, no closed form solution is available on nonlinear spaces,
whereas the variation definition \eqref{eqn:intro-frechet-problem-statement} remains meaningful.
The solution is commonly called a \textit{Fr\'echet mean} or \textit{metric barycenter}. Hadamard spaces provide a canonical setting in which the solution is unique; in general, existence and uniqueness require additional hypotheses.

To solve \eqref{eqn:intro-frechet-problem-statement} when $\mathcal{M}_{1} = \mathbb{R}^{d}$, recall the algorithm of recursive averaging.
Let $(I_{t})_{t \geq 1}$ be i.i.d. random variables with $I_{t} \sim \operatorname{Unif}[N]$, set $y_{t} := z_{I_{t}}$, and define $\widetilde{x}_{t} := \frac{1}{t}\sum_{i=1}^{t} y_{i}, $ for $t \geq 1$. Then,
\begin{equation}\label{eqn:recursive-averaging}
    \widetilde{x}_{t+1} = \frac{t}{t+1}\cdot \frac{1}{t}\sum_{i=1}^{t} y_{i} + \frac{1}{t+1}y_{t+1} = \frac{t}{t+1}\widetilde{x}_{t} + \frac{1}{t+1}y_{t+1}. 
\end{equation}
By the strong law of large numbers, $\widetilde{x}_{t} \rightarrow \widehat{x}_{N}$ almost surely. This update also has a geometric interpretation: \eqref{eqn:recursive-averaging} interpolates 
between the current average and a newly sampled point. 

Intuitively, one generalizes the update \eqref{eqn:recursive-averaging} to Riemannian manifolds by interpolating on the respective geodesic segments.
That is, 
whenever $y_{t+1}$ lies in a normal neighborhood of $\widetilde{x}_{t}$ define 
\begin{equation} \label{eqn:geodesic-interpolation-intro}
    \widetilde{x}_{t+1} = \exp_{\widetilde{x}_{t}}^{\mathcal{M}_{1}} \left(\frac{1}{t+1} \log_{\widetilde{x}_{t}}^{\mathcal{M}_{1}}(y_{t+1})\right).
\end{equation}
We refer to \eqref{eqn:geodesic-interpolation-intro} as \textit{geodesic interpolation} (e.g., Algorithm~\ref{alg:geodesic-interpolation}). 

Although motivated by statistical origins, Algorithm~\ref{alg:geodesic-interpolation}
has another optimization interpretation
which we use to provide convergence guarantees. When the logarithmic map is well-defined,
\begin{equation*}
    \nabla_{\mathcal{M}_{1}}\left(\frac{1}{2}d_{\mathcal{M}_{1}}^{2}(\cdot,z)\right)(x) = -\log_{x}^{\mathcal{M}_{1}}(z).
\end{equation*}
Let
$\mathcal{F}_{t} := \sigma(x_{0}, I_{1}, \dots, I_{t})$. Since $\widetilde{x}_{t}$ is $\mathcal{F}_{t}$-measurable and $I_{t+1}$ is independent of $\mathcal{F}_{t}$ we have
\begin{equation*}
    \mathbb{E}[-\log_{x_{t}}^{\mathcal{M}_{1}}(y_{t+1}) \mid \mathcal{F}_{t}] = -\frac{1}{N} \sum_{r=1}^{N}\log_{x_{t}}^{\mathcal{M}_{1}}(z_{r}) = \nabla_{\mathcal{M}_{1}} f_{P_{N}}(x_{t}).
\end{equation*}
Thus, the update \eqref{eqn:geodesic-interpolation-intro}
is Riemannian stochastic gradient descent (SGD) with the unbiased gradient estimator $g_{t} := - \log_{\widetilde{x}_{t}}^{\mathcal{M}_{1}}(y_{t+1})$ and stepsize $\eta_{t} = 1/(t+1)$.

\subsection{Unknown Submanifolds and Graph-Based Geodesic Interpolation}
\label{subsec:unknown-submanifold-graph-interpolation}

\begin{figure}[t]
    \centering
    \includegraphics[width=0.9\linewidth]{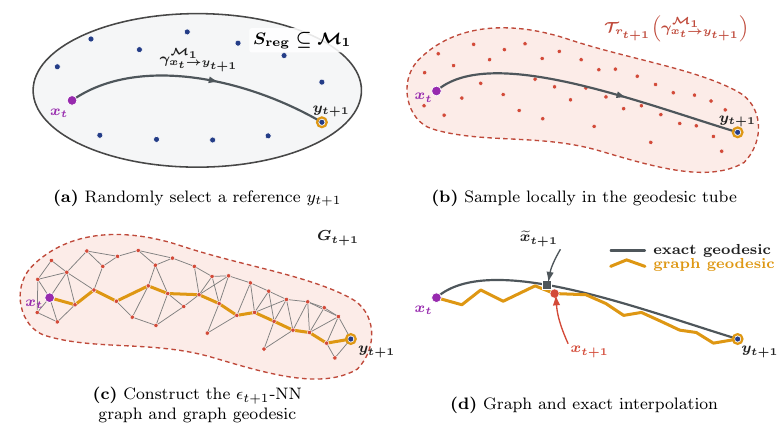}
    \caption{Graph Approximation of Geodesic Interpolation Step with Tube Sampler}
    \label{fig:idealized-step}
\end{figure}
Let $(\mathcal{M}_{1},g_{1}) \hookrightarrow$ 
$(\mathcal{M}_{2},g_{2})$, with $g_{1}$ induced by $g_2$, and intrinsic and ambient distances by $d_{\mathcal{M}_1}$ and $d_{\mathcal{M}_{2}}$, respectively. 
Given $(z_{r})_{r=1}^{N} \subseteq \mathcal{M}_{1}$, Algorithm~\ref{alg:geodesic-interpolation} would suffice to find $\widehat{x}_{N}$ if $\mathcal{M}_{1}$ were known.
However, update \eqref{eqn:geodesic-interpolation-intro} cannot be directly implemented as $\mathcal{M}_{1}$ is unknown.
The saving grace is that we can interpret Algorithm~\ref{alg:geodesic-interpolation} purely through the lens of geodesics: for $\eta_{t} \in [0,1]$, 
the update 
interpolates an $\eta_{t}$-fraction along the intrinsic geodesic from $x_{t}$ to chosen reference $y_{t+1}$. Thus, we ask: \textit{can we emulate these updates via graph interpolation}? 

To do so, at first glance one may 
demand a fine representation of the full region of interest. 
We avoid this for two reasons: $i)$ fine representations of untraversed regions 
are wasteful; $ii)$ an initial fine construction is both prohibitively expensive statistically and computationally. In contrast, we only use \textit{fine enough, local} graph representations $G_{t}$, where graph fidelity is proportional to the expected objective gap.
Intuitively, the larger the objective gap, the further the iterate is from the minimizer, and thus coarser approximations still provide useful directional information. 
As the gap decreases, we require finer representations for useful information.

At iteration $t+1$, we assume access to local sampling oracles on $\mathcal{M}_{1}$. Informally, given a current iterate $x_t$, a sampled reference point $y_{t+1}$, and a scale parameter $\widetilde{r}_{t+1}$, the sampler produces points of $\mathcal{M}_1$ within a region -- defined by local radius $\widetilde{r}_{t+1}$ -- near the interpolating geodesic $\gamma_{x_{t} \rightarrow y_{t+1}}^{\mathcal{M}_{1}}$. Explicitly, let $\ROI \subseteq \mathcal{M}_{1}$ denote the region of interest containing the reference points and potential new samples. We consider
\begin{equation}\label{eqn:geodesic-tube}
    S_{\widetilde{r}_{t+1}}(\gamma_{x_{t} \rightarrow y_{t+1}}^{\mathcal{M}_{1}}) := \ROI \cap T_{\widetilde{r}_{t+1}}(\gamma_{x_{t} \rightarrow y_{t+1}}^{\mathcal{M}_{1}}),
\end{equation}
where $T_{\widetilde{r}_{t+1}}(\gamma_{x_{t} \rightarrow y_{t+1}}^{\mathcal{M}_{1}}) := \{p \in \mathcal{M}_{1} \colon \operatorname{dist}_{\mathcal{M}_{1}}(p, \gamma_{x_{t} \rightarrow y_{t+1}}^{\mathcal{M}_{1}}) \leq \widetilde{r}_{t+1}\}$. On each localized tube, let the sampling oracle be $\mu_{\gamma, r}$, the normalized Riemannian volume:
\begin{equation}\label{eqn:normalized-Riemannian-volume}
    \mu_{\gamma,r}(A) := \frac{\operatorname{vol}_{\mathcal{M}_{1}}(A \cap \ROI \cap T_{r}(\gamma))}{\operatorname{vol}_{\mathcal{M}_{1}}(\ROI \cap T_{r}(\gamma))}.
\end{equation}
We emphasize this
is an idealized sampling oracle, and 
the cost of realizing such an oracle depends inversely on the tube radius $\widetilde{r}_{t+1}$; the smaller $\widetilde{r}_{t+1}$, the larger the cost.

\begin{algorithm}[t]
\caption{Graph-Based Geodesic Interpolation}
\label{algo:graph-algorithm}
\begin{algorithmic}
\Require References $(z_r)_{r=1}^N$, initialization $x_0$,
stepsizes $(\eta_t)_{t=0}^{T-1}$, refinement rules
$(\mathsf{Refine}_{t})_{t=1}^{T}$.
\For{$t=0,\ldots,T-1$}
    \State \textbf{Uniformly Pick Reference Point:} Draw $I_{t+1}\sim\mathrm{Unif}([N])$ and set
    $y_{t+1}:=z_{I_{t+1}}$; \vspace{.1cm}
    \State \textbf{Graph Refinement:} 
        $G_{t+1}
        \leftarrow
        \mathsf{Refine}_{t+1}(x_t,y_{t+1};G_t);$ \Comment{Refine by Alg.~\ref{algo:tube-sample-oracle} or~\ref{algo:ball-refinement}} \vspace{.1cm}
    \If{$X_{t}$ and $y_{t+1}$ are disconnected in $G_{t+1}$}
    \State $x_{t+1} \gets x_{t}$.
    \Else
    \State \textbf{Update by Graph Interpolation:} Compute a shortest path $\gamma^{G_{t+1}}_{x_t\to y_{t+1}}$, then set
\begin{equation*}
    x_{t+1} := v_{j_{t}}^{(t)}, \qquad \text{with} \quad j_t:=\operatorname*{argmin}_{0\le j\le m_t}
\left|d_{G_{t+1}}(x_t,v_j^{(t)})
-\eta_t d_{G_{t+1}}(x_t,y_{t+1})\right|
\end{equation*}
\State \qquad and $\{v_{j}^{(t)}\}_{j=0}^{m_{t}}$ the vertices of $\gamma_{x_{t} \rightarrow y_{t+1}}^{G_{t+1}}$. 
\EndIf
\EndFor
\State \Return $x_T$.
\end{algorithmic}
\end{algorithm}
The resulting update is the graph-based analogue of Algorithm~\ref{alg:geodesic-interpolation}. At iteration $t+1$, draw a reference point $y_{t+1} = z_{I_{t+1}}$ with $I_{t+1} \sim \operatorname{Unif}[N]$. If $y_{t+1} = x_{t}$, set $x_{t+1} = x_{t}$ and take $G_{t+1} = \{x_{t}\}$; no samples are needed. Otherwise, perform the following:
create a local graph representation $G_{t+1}$ by first locally sampling in a region around $\gamma_{x_{t} \rightarrow y_{t+1}}^{\mathcal{M}_{1}}$ (e.g. by Algorithm~\ref{algo:tube-sample-oracle}), then construct an $\epsilon_{t+1}$ nearest-neighbors graph $G_{t+1}$ using these samples with $x_{t}$ and $y_{t+1}$. Lastly, compute a graph-shortest path $\gamma_{x_t\to y_{t+1}}^{G_{t+1}}$ from $x_t$ to $y_{t+1}$, which exists with high probability as established in Section~\ref{sec:idealized-convegence-analysis}, with ordered vertices $v_0^{(t)}=x_t,\, v_1^{(t)},\ldots, v_{m_t}^{(t)}=y_{t+1}$. Update the iterate to $x_{t+1}$ by moving an $\eta_{t}$-fraction, measured by graph distance, along the path, choosing the smaller index when ties occur:
\begin{equation}\label{eqn:graph-update}
\begin{aligned}
x_{t+1} &:= v_{j_{t}}^{(t)}, \qquad \text{with} \quad j_t:=\operatorname*{argmin}_{0\le j\le m_t}
\left|d_{G_{t+1}}(x_t,v_j^{(t)})
-\eta_t d_{G_{t+1}}(x_t,y_{t+1})\right|.
\end{aligned}
\end{equation}
In words, this update mirrors ideal geodesic interpolation by replacing the intrinsic geodesic with a local graph approximation thereof. 
The full algorithm is summarized by Algorithm~\ref{algo:graph-algorithm}.

\begin{theorem}[\textbf{Informal}]\label{thm:informal-main-thm} Suppose $\mathcal{M}_{1} \subseteq \mathcal{M}_{2}$ with $\operatorname{dim}(\mathcal{M}_{1}) = d$, and assume access to Algorithm~\ref{algo:tube-sample-oracle}. Suppose we are given a sample $\{z_{r}\}_{r=1}^{N} \subseteq \mathcal{M}_{1}$, and define
\begin{equation}\label{eqn:objective-functional-informal}
    f_{P_{N}}(x) := \frac{1}{2N}\sum_{r=1}^{N} d_{\mathcal{M}_{1}}^{2}(x, z_{r}).
\end{equation}

Under suitable technical assumptions, run Algorithm~\ref{algo:graph-algorithm} for $T \geq 1$ iterations with step-size schedule $(\eta_{t})_{t\geq 0}$ and local scale parameters $(r_{t})_{t\geq 0}$, where $\eta_{t} = \mathcal{O}(t^{-1})$ and $r_{t} = \mathcal{O}(t^{-2})$. Let Algorithm~\ref{algo:tube-sample-oracle} sample $n_{t}$ points with
\begin{equation}\label{eqn:statistical-cost-per-iterate-informal}
     (\mathrm{Local} \text{ } \mathrm{Sampling} \text{ } \mathrm{Cost}) \qquad n_{t} \gtrsim t^{4d+2}\left(\log(t^{4d+2}) + \log(T^{2})\right)
\end{equation}
at each iteration. Then, the local $\epsilon_{t}$-neighbors graph $G_{t}$ construction costs at most
\begin{equation}\label{eqn:computational-cost-per-iterate-informal}
    (\mathrm{Local} \text{ } \mathrm{Construction} \text{ } \mathrm{Cost}) \qquad \operatorname{cost}^{\operatorname{dist}}_{t} \leq \binom{n_{t} + 2}{2} = \mathcal{O}(n_{t}^{2}) 
\end{equation}
pairwise-distance oracle calls per iteration. The final iterate $x_{T}$ produced by Algorithm~\ref{algo:graph-algorithm} satisfies
\begin{equation}\label{eqn:guarantee-informal}
    \mathbb{E}_{\operatorname{alg}}\left[f_{P_{N}}(x_{T}) - f_{P_{N}}(\widehat{x}_{N})\right] \leq \mathcal{O}\left(T^{-1}\right),
\end{equation}
where $\widehat{x}_{N}$ is the minimizer of \eqref{eqn:objective-functional-informal} within $\ROI$.
\end{theorem}

\begin{algorithm}[t]
\caption{Idealized Tube-Based Graph Refinement}
\label{algo:tube-sample-oracle}
\begin{algorithmic}
\Require Current iterate $x_t$, reference $y_{t+1}$,
sampling radius $\widetilde{r}_{t+1}$, connectivity scale $\epsilon_{t+1}$,
and sample size $n_{t+1}$.

\State \textbf{Sample from Geodesic Tube:}
Draw
\[
    U_{t+1}
    :=
    \{u_k\}_{k=1}^{n_{t+1}}
    \overset{\mathrm{i.i.d.}}{\sim} 
    \mu_{\gamma_{x_{t} \rightarrow y_{t+1}}^{\mathcal{M}_{1}}, \widetilde{r}_{t+1}},
\]
with $\mu_{\gamma,r}$ defined by \eqref{eqn:normalized-Riemannian-volume}.

\State \textbf{Construct Local Graph:}
Let $G_{t+1}$ be the ambient-distance-weighted
$\epsilon_{t+1}$-neighborhood graph on
    $U_{t+1}\cup\{x_t,y_{t+1}\},$
with edge weights weighted by ambient distance weights;

\State \Return $G_{t+1}$.
\end{algorithmic}
\end{algorithm}
The informal theorem highlights the two main features of the method. First, despite no access to the intrinsic geometry of $\mathcal{M}_{1}$, Algorithm~\ref{algo:graph-algorithm} achieves an $\mathcal{O}(T^{-1})$ expected objective-gap rate. This matches the minimax-optimal iteration dependence for stochastic first-order methods on strongly convex Euclidean objectives \cite{lower-bound-SGD}. The graph approximation therefore introduces no loss in the outer iteration rate, provided the local graph representations are accurate enough.

Second, the cost of achieving this rate is statistical but local. Rather than constructing a fine global graph approximation of $\mathcal{M}_{1}$, the algorithm constructs a \textit{fine enough, local} approximation around $\gamma_{x_{t} \rightarrow y_{t+1}}^{\mathcal{M}_{1}}$
The resulting per-iteration sampling cost grows polynomially with $t$ with exponential dependence on the intrinsic dimension $d$, ignoring logarithmic costs for probabilistically guaranteeing $T$ good graph constructions, all up to geometric constants of $\mathcal{M}_{1}$.

\subsection{Replacing the Idealized Tube Sampler}\label{sec:practical-graph-algo-description}
The refinement step (e.g. Algorithm~\ref{algo:tube-sample-oracle}) of the previous algorithm has two practical limitations: $i)$ the tube sampler requires explicit knowledge of $\gamma_{x\rightarrow y}^{\mathcal{M}_{1}}$, the object the graph-geodesics approximate; $ii)$ the refinement constructs a \textit{fresh} local graph at every iteration, wasting
previously sampled points. We address these limitations by approximating the sampling mechanism of Algorithm~\ref{algo:tube-sample-oracle} with a local sampler requiring only knowledge of previously sampled points.
To explain, consider the uniform ball sampler: given $v\in\mathcal{M}_{1}$, radius $r>0$, and sample count $b$, the oracle returns $b$ samples from
    $\operatorname{Unif}\left(B_{\mathcal{M}_{1}}(v,r)\right).$ We note other samplers may be used.

\begin{figure}[t]
    \centering
    \includegraphics[width=0.9\linewidth]{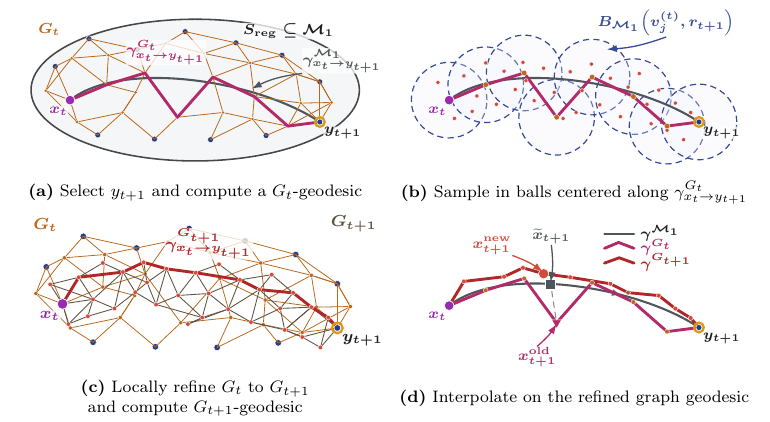}
    \caption{Graph Approximation of Geodesic Interpolation Step with Ball Sampler}
    \label{fig:graph-approximation-with-ball-sampler}
\end{figure}

Suppose that after iteration $t$ we are given a graph representation $G_t$ and the current iterate $x_t$. At iteration $t+1$, we begin by drawing 
    $y_{t+1}=z_{I_{t+1}},
    $ where $
    I_{t+1}\sim\operatorname{Unif}[N],$
as before. We differ in the graph refinement step. Assuming that $x_t$ and $y_{t+1}$ are connected in $G_t$, compute a graph-shortest path $\gamma_{x_t\rightarrow y_{t+1}}^{G_t} = (v_0^{(t)},v_1^{(t)},\ldots,v_{m_t}^{(t)}),$ with $v_{0}^{(t)} = x_{t}$ and $v_{m_t}^{(t)}=y_{t+1}$ 
which provides a coarse approximation of $\gamma_{x_t\rightarrow y_{t+1}}^{\mathcal{M}_{1}}$.

For each vertex $v_j^{(t)}$, call the ball sampler with radius $r_{t+1}$ and per-center sample count $b_{t+1}$. Denote the resulting samples by
\[
    U_{t+1}
    :=
    \bigcup_{j=0}^{m_t}
    \left\{
        u_{j,k}^{(t+1)}
    \right\}_{k=1}^{b_{t+1}},
    \quad \text{with }
    u_{j,k}^{(t+1)}
    \sim
    \operatorname{Unif}
    \left(
        B_{\mathcal{M}_{1}}
        \left(v_j^{(t)},r_{t+1}\right)
    \right).
\]
Panels (a)--(b) of Figure~\ref{fig:graph-approximation-with-ball-sampler} depict the coarse graph path and corresponding local sampling region.

Then, we refine $G_{t}$ \textit{locally} with these new samples.
To \textit{locally} refine, we select \textit{nearby} $v \in G_{t}$ near $\gamma_{x_{t} \rightarrow y_{t+1}}^{G_{t}}$ -- distance measured by $d_{G_{t}}$ -- and remove preexisting edges where both endpoints lie in this set. Then, we construct an $\epsilon_{t+1}$-neighborhood graph on the selected vertices together with $U_{t+1}$ using ambient distance $d_{\mathcal{M}_{2}}$. We update $G_{t}$ to $G_{t+1}$ by adding the new sampled vertices to the vertex set, and all edges of the $\epsilon_{t+1}$ nearest neighbors graph to the remaining edges of $G_{t}$. 
The resulting graph refinement is visualized in panel (c) of Figure~\ref{fig:graph-approximation-with-ball-sampler}, and formally defined by Algorithm~\ref{algo:ball-refinement}. Finally, we compute a graph-shortest path $\gamma_{x_t\rightarrow y_{t+1}}^{G_{t+1}}$ in the refined graph and set $x_{t+1}$ according to \eqref{eqn:graph-update} depicted in panel (d) of Figure~\ref{fig:graph-approximation-with-ball-sampler}. 

The radius $r_{t+1}$ controls a localization tradeoff. Larger radii enable more discrepancy between the coarse graph path and the unknown intrinsic geodesic, but generate samples further away from the desired region. A smaller radius concentrates the refinement near the coarse path, but is appropriate only when $G_t$ is already sufficiently accurate. As the graph refines, previously sampled vertices may reduce the number of additional samples required in later iterations. 

Relative to the idealized method, 
additional approximation arises from centering the local sampling regions along
$\gamma_{x_t\rightarrow y_{t+1}}^{G_t}$
rather than
$\gamma_{x_t\rightarrow y_{t+1}}^{\mathcal{M}_{1}}$,
and 
reusing the graph and its samples. An analysis of these effects would require controlling the coverage produced by the graph-centered balls and the dependence introduced by adaptive graph refinement. We leave this extension for future work. We use Algorithm~\ref{algo:graph-algorithm} with Algorithm~\ref{algo:ball-refinement}
for the numerical experiments in Section~\ref{sec:numerical-experiments}.

\begin{algorithm}[t]
\caption{Ball-Based Persistent Graph Refinement}
\label{algo:ball-refinement}
\begin{algorithmic}
\Require Current iterate $x_t$, reference $y_{t+1}$,
current graph $G_t$, ball radius $r_{t+1}$,
per-center sample size $b_{t+1}$,
rewiring radius $\tau_{t+1}$,
and connectivity scale $\epsilon_{t+1}$.
\If{$x_{t}$ and $y_{t+1}$ are disconnected in $G_{t}$}
\State $G_{t+1} \gets G_{t}$.
\Else{}
\State \textbf{Compute Coarse Graph Geodesic:}
Compute a shortest path
\[
    \Gamma_t
    :=
    \gamma_{x_t\rightarrow y_{t+1}}^{G_t}
    =
    (v_0^{(t)},\ldots,v_{m_t}^{(t)}).
\]

\State \textbf{Sample along Coarse Geodesic:}
For each $j=0,\ldots,m_t$, draw
\[
    \{u_{j,k}\}_{k=1}^{b_{t+1}}
    \overset{\mathrm{i.i.d.}}{\sim}
    \operatorname{Unif}
    B_{\mathcal M_1}
    (v_j^{(t)},r_{t+1}),
\]
and set
    $U_{t+1}
    :=
    \bigcup_{j=0}^{m_t}
    \{u_{j,k}\}_{k=1}^{b_{t+1}}$;

\State \textbf{Select Local Rewiring Region:}
With $v_{j}^{(t)} \in \Gamma_{t}$, set
\[
    V_{t+1}^{\mathrm{loc}}
    :=
    \left\{
        u\in V(G_t):
        \min_{0\le j\le m_t}
        d_{G_t}(u,v_j^{(t)})
        \le \tau_{t+1}
    \right\}, \quad 
    E_{t+1}^{\mathrm{loc}}
    :=
    \left\{
        \{u,v\}\in E(G_t):
        u,v\in V_{t+1}^{\mathrm{loc}}
    \right\}.
\]

\State \textbf{Local Refinement:}
Let $\widetilde G_{t+1}$ be the
$\epsilon_{t+1}$-neighborhood graph on
    $V_{t+1}^{\mathrm{loc}}\cup U_{t+1}$
with ambient distance weights;

\State \textbf{Update Persistent Graph:}
Set
\[
    V(G_{t+1})
    :=
    V(G_t)\cup U_{t+1} \quad \text{and} \quad 
    E(G_{t+1})
    :=
    \bigl(E(G_t)\setminus E_{t+1}^{\mathrm{loc}}\bigr)
    \cup E(\widetilde G_{t+1}).
\]

\State \Return $G_{t+1}$.
\EndIf
\end{algorithmic}
\end{algorithm}

\subsection{Related Work}

\paragraph{Fr\'echet means: computation and statistics.}
Classical geometric and algorithmic treatments of the Riemannian center of mass include \cite{Kar77,ATV13}, while stochastic geodesic recursions for computing $p$-means on manifolds were studied in \cite{ADPY12}. Of particular interest, recent work focuses on computing Fr\'echet means in Wasserstein space \((\mathcal P_{2,\operatorname{ac}}(\mathbb{R}^{d}), W_{2})\).
Existence and uniqueness of such means have been characterized in \cite{AC11}. The closest direct predecessor of our stochastic-interpolation viewpoint is \cite{CMRS20}, analyzing first-order methods for Bures--Wasserstein barycenters; the online Wasserstein recursion of \cite{BFRT25} is similarly based on repeated interpolation toward sampled measures. The Bures--Wasserstein manifold is a geodesically convex submanifold of \(\mathcal P_{2,\operatorname{ac}}(\mathbb{R}^{d})\) with respect to the ambient metric $W_{2}$. In contrast, we allow $\mathcal{M}_{1}$ to \textit{not} be geodesically convex with respect to the ambient metric. The analysis of \cite{CMRS20} obtains global rates through a Polyak--{\L}ojasiewicz inequality in this specific Bures--Wasserstein setting, whereas our local convexity and smoothness constants are induced by  localization based on geometric properties of the given $\mathcal{M}_{1}$. We impose these 
conditions to ensure existence and uniqueness of the target Fr\'echet mean. Additionally, both of these prior algorithmic approaches assume access to the relevant geometric structures -- intrinsic geodesics or tangent vectors (i.e. transport maps) -- while we explicitly assume no direct access. A complementary recent method, GEORCE-FM \cite{RHM25}, jointly optimizes discretized geodesics and a Fr\'echet mean. Although they approximate both geodesics and the Fr\'echet mean, they require a known local chart and access to the metric tensor of a known manifold, whereas we explicitly do not.

Graphs have been used for Fr\'echet-type estimation on unknown manifolds due to their geometric expressivity. 
\cite{CM12} estimates geodesic distances by neighborhood-graph shortest paths and uses a modified ISOMAP representation to construct a manifold mean for functional data in $L^{2}$. More recently, \cite{GM26} combines a global ISOMAP representation with Fr\'echet regression to map Euclidean coordinates back to a general metric space, with an emphasis on submanifolds of $\mathcal{P}_{2,\operatorname{ac}}(\mathbb{R}^{d})$. These approaches both heavily rely on isometric to $\mathbb R^{d}$ assumptions; the former asking for the intrinsic geodesic geometry to be Euclidean and the latter for the ambient restricted metric itself to be Euclidean. In contrast, we ask for neither the intrinsic nor ambient spaces to be Euclidean.
\cite{TGTP20} extends ISOMAP convergence to noisy observations lying near an unknown manifold. 
Their analysis permits vanishing off-manifold perturbations and supports subsequent inference on Fr\'echet means, but specializes to one-dimensional curve segments embedded in $\mathbb R^{d}$, where intrinsic distance reduces to arc length and the Fr\'chet mean can be computed by averaging arc-length coordinates. In contrast, our metric-recovery results apply to embedded submanifolds of arbitrary intrinsic dimension in a Riemannian ambient space,
although our present analysis assumes that graph vertices lie exactly on the submanifold. 

\paragraph{Manifold learning in $\mathbb R^D$ and beyond.} The \textit{manifold hypothesis} -- high-dimensional data concentrates around/on low-dimensional surfaces -- has given rise to various algorithmic approaches for exploiting this mathematical structure over the past three decades. Of note, ISOMAP estimates intrinsic geodesic distances by shortest paths in a neighborhood graph \cite{TdSL00}; quantitative shortest-path approximation results include \cite{ACLG19}. These results provide non-asymptotic sampling guarantees to ensure metric equivalence with high probability, but require the ambient space to be $\mathbb R^d$. We directly extend the results of \cite{TdSL00}, allowing for replacement of $\mathbb R^D$ by a Riemannian manifold $\mathcal M_2$ with bounded sectional curvature, capturing the classical result as a special case. A motivating line of research centers on submanifolds of $\mathcal{P}_{2,\operatorname{ac}}(\mathbb R^D)$. WASSMAP \cite{HHK23} applies the ISOMAP construction with pairwise Wasserstein distances, and \cite{HMST25} develops graph-based recovery of Wasserstein submanifolds with quantitative graph-scale and sampling conditions and asymptotic recovery guarantees. These works provide asymptotic guarantees for when  $\mathcal{P}_{2,\operatorname{ac}}(\mathbb{R}^{D})$ is the ambient space. Related spectral work establishes graph-Laplacian convergence for manifold samples observed through a general metric-space ambient \cite{XS26}, while \cite{ORST25} proves variational convergence of graph $p$-Dirichlet energies for Laplace learning in Wasserstein space. The metric-recovery result in Section~\ref{sec:generalized-isomap} differs in both scope and use. It gives a local, two-sided multiplicative comparison for ambient-distance-weighted graphs on an embedded submanifold $\mathcal M_1 \subset \mathcal M_2$, where $\mathcal M_2$ is Riemannian, with admissible scales controlled by curvature, the second fundamental form, and injectivity radii. Extending our geometric analysis to the Wasserstein setting -- afforded by second-order Otto calculus -- is a natural direction for future work.

\paragraph{Riemannian optimization on known manifolds.}
Riemannian optimization treats the feasible manifold as part of the problem specification, allowing structural constraints and symmetries -- such as orthogonality, fixed rank, or positive definiteness -- to be encoded intrinsically. 
Foundational treatments developed Riemannian analogues of gradient and conjugate-gradient methods, Newton and quasi-Newton schemes, and trust-region algorithms using tangent-space models, retractions, and vector transports \cite{EAS98,AMS08,Bou23}. For smooth objectives on general manifolds, modern complexity analyses typically target approximate first- or second-order stationarity
\cite{BAC19}. When geodesic convexity is available,
global objective-gap guarantees have been established for deterministic and stochastic gradient, subgradient, and proximal-point methods \cite{ZS16,BFM17}. Nonsmooth Riemannian optimization includes early subgradient and proximal-point schemes \cite{FO98,FO02}, proximal-gradient methods for composite objectives \cite{HW22}, and bundle methods for black-box nonsmooth objectives \cite{HNP23,BHJ24,DGM26}. While these Riemannian optimization methods have vastly expanded in the previous three decades, they assume an explicit manifold model and access to the geometric (exact or inexact) primitives. We differ as the feasible submanifold is unknown, so updates must be constructed only from ambient information, and analytical tools can only be invoked implicitly.

A foundational treatment of Riemannian stochastic gradient descent is \cite{Bon13}; finite-time analyses for geodesically convex objectives and averaged stochastic methods appear in \cite{ZS16,TFBJ18}. With exact geodesics, Algorithm~\ref{alg:geodesic-interpolation} is an instance of Riemannian SGD with an unbiased stochastic gradient estimator, whereas replacing the geodesic by a graph path produces an implicit estimator that is generally biased. Euclidean descent analyses based on conditional alignment of the mean search direction and affine control of its conditional second moment are developed in \cite{bottou2018optimization}, while broader classes of biased gradient estimators are organized and compared in \cite{demidovich2023guidezoobiasedsgd}. Our analysis gives a tangent-space counterpart of the former descent framework: conditional alignment and second-moment control are imposed in $T_{x_t}\mathcal M_1$, and local strong geodesic convexity and smoothness yield an $\mathcal O(t^{-1})$ last-iterate expected objective-gap bound under decaying stepsizes.

\paragraph{Riemannian optimization without explicit manifold access.}
A recent line couples manifold learning with optimization when the feasible geometry is available only through samples or a learned representation. Shustin, Avron, and Sober \cite{SAS22} consider a smooth submanifold of Euclidean space and use Manifold--MLS to reconstruct the tangent spaces, retractions, gradients, and vector transports needed by approximate Riemannian methods from noiseless samples. Diepeveen and Weber \cite{DW25} take a different route: they assume access to a learned pullback geometry on Euclidean data space, replace the Levi--Civita connection by an iso-connection, and develop first-order methods for Euclidean-convex objectives and iso-barycentres. Their iso-barycentre is defined through an iso-logarithmic zero condition and need not coincide with the ordinary Fr\'echet minimizer. Hauptmann, Korolev, and Thorpe \cite{HKT26} work directly with a point cloud sampled from a compact positive-reach submanifold of Euclidean space; for a smooth ambient energy, they estimate tangent projections by local PCA and apply a learned normal correction, obtaining convergence to a neighborhood of a local minimizer whose radius vanishes with the time-step and manifold-estimation errors.
The present method is parallel to these approaches. It is specialized to the empirical Fr\'echet functional and reconstructs neither a complete differential-geometric optimization oracle nor a global learned metric; it approximates only the intrinsic geodesics required during a run of the geodesic interpolation algorithm. This specialization permits a general Riemannian ambient space and yields an explicit coupling between graph resolution and the outer $\mathcal O(T^{-1})$ rate, at the cost of the idealized localized tube-sampling oracle and the localization assumptions of Section~\ref{sec:assumptions}. Addressing a more general class of objectives with non-Euclidean ambient spaces presents a natural direction for future work.

\section{Convergence Guarantee}\label{sec:idealized-convegence-analysis}

\setlength{\abovecaptionskip}{4pt}
\setlength{\belowcaptionskip}{0pt}
\begin{table}[ht]
\centering
\footnotesize
\renewcommand{\arraystretch}{0.95}
\begin{tabularx}{\textwidth}{@{}>{\raggedright\arraybackslash}p{0.28\textwidth}X@{}}
\toprule
\textbf{Symbol} & \textbf{Meaning} \\
\midrule
\((\mathcal M_1,g_1)\hookrightarrow(\mathcal M_2,g_2)\) & Unknown embedded data manifold and ambient Riemannian manifold\\
\(d_{\mathcal M_1},\ d_{\mathcal M_2}\) & Intrinsic distance on \(\mathcal M_1\) and ambient distance on \(\mathcal M_2\). \\
\(T_x\mathcal M,\ \langle\cdot,\cdot\rangle_x,\ \lVert\cdot\rVert_x\) & Tangent space at \(x\), with the inner product and norm. \\
\(\exp_x^{\mathcal M_j},\ \log_x^{\mathcal M_j}\) & Exponential map and its local inverse on \(\mathcal M_j\), for \(j\in\{1,2\}\). \\
\(\gamma_{x\to y}^{\mathcal M_1}\) & A minimizing intrinsic geodesic segment from \(x\) to \(y\) in \(\mathcal M_1\). \\
\(\Pi,\ \kappa\) & Second fundamental form of \(\mathcal M_1\subset\mathcal M_2\) and its uniform bound. \\
\(K_{\min}^{\mathcal M_j},\ K_{\max}^{\mathcal M_j},\ \widetilde K_{\mathcal M_j}\) & Sectional-curvature bounds, and \(\widetilde K_{\mathcal M_j}=\max\{|K_{\min}^{\mathcal M_j}|,|K_{\max}^{\mathcal M_j}|\}\). \\
\(S_{\mathrm{reg}},\ D_{\mathrm{reg}}\) & 
Region of interest and its intrinsic diameter. \\
\(r_{\mathrm{inj}}^{\mathcal M_j}(S)\) & Uniform injectivity radius of \(\mathcal M_j\) over \(S\subseteq\mathcal M_1\). \\
\(\rho_S(s)\) & Largest $d_{\mathcal{M}_{1}}(x,y)$ for \(x,y \in S\) with ambient distance at most \(s\). \\
\bottomrule
\end{tabularx}
\caption{Geometry notation.}
\label{tab:notation-main-geometry}
\end{table}

\begin{table}[ht]
\centering
\footnotesize
\renewcommand{\arraystretch}{0.95}
\begin{tabularx}{\textwidth}{@{}>{\raggedright\arraybackslash}p{0.28\textwidth}X@{}}
\toprule
\textbf{Symbol} & \textbf{Meaning} \\
\midrule
\(f_{P_N}\) & Empirical Fréchet objective \(f_{P_N}(x)=(2N)^{-1}\sum_{r=1}^N d_{\mathcal M_1}^2(x,z_r)\). \\
\(\widehat x_N\) & A minimizer of \(f_{P_N}\) over \(\mathcal M_1\). \\
\(I_{t+1},\ y_{t+1}=z_{I_{t+1}}\) & Uniformly sampled reference index and respective reference point. \\
\(\eta_t\) & Interpolation fraction 
used in the update from \(t\) to \(t+1\). \\
\(\widetilde x_{t+1}\), \(x_{t+1}\) & Exact intrinsic interpolation and graph interpolation between \(x_t\) and \(y_{t+1}\). \\
\(G_t,\ V(G_t),\ E(G_t),\ d_{G_t}\) & Current weighted graph, vertex and edge sets, and shortest-path metric. \\
\(\gamma_{x_t\to y_{t+1}}^{G_{t+1}}=(v_0^{(t)},\ldots,v_{m_t}^{(t)})\) & Refeined-graph shortest path used for the \((t+1)\)-st interpolation update. \\
\(\widetilde g_t=-\log_{x_t}^{\mathcal M_1}(y_{t+1})\) & Exact stochastic gradient estimator for Algorithm~\ref{alg:geodesic-interpolation}. \\
\(g_t=-\eta_t^{-1}\log_{x_t}^{\mathcal M_1}(x_{t+1})\) & Implicit graph-induced stochastic gradient estimator for Algorithm~\ref{algo:graph-algorithm}. \\
\(m_f,\ L_f\) & Local strong-convexity and smoothness constants for \(f_{P_N}\) on \(S_{\mathrm{reg}}\). \\
\bottomrule
\end{tabularx}
\caption{Optimization and algorithm notation.}
\label{tab:notation-main-optimization}
\end{table}

\begin{table}[ht]
\centering
\footnotesize
\renewcommand{\arraystretch}{0.95}
\begin{tabularx}{\textwidth}{@{}>{\raggedright\arraybackslash}p{0.28\textwidth}X@{}}
\toprule
\textbf{Symbol} & \textbf{Meaning} \\
\midrule
\(\xi_t,\ \mathcal F_t\) & Local-sampler randomness at iteration \(t\), filtration encoding until iteration \(t\). \\
\(\alpha_{t+1},\ r_{t+1}=\alpha_{t+1}\ell_{t+1}\) & Relative and absolute sample tube radius, where \(\ell_{t+1}=d_{\mathcal M_1}(x_t,y_{t+1})\). \\
\(\lambda_1^{(t)},\ \lambda_2^{(t)},\ \lambda_t\) & Lower and upper graph-metric distortion tolerances; on \(E_t\) both are \(\lambda_t\). \\
\(\epsilon_t\) & Ambient graph-connectivity scale at iteration \(t\). \\
\(\delta_t\) & Intrinsic fill/covering radius for sampling guarantee. \\
\(T_t,\ \widetilde r_t,\ \widetilde{T}_t\) & Geodesic tube, inflated radius, and inflated tube. \\
\(\mathfrak S_{\alpha_t,\epsilon_t},\ \mu_{\widetilde{\mathcal T}},\ \mathbf q^{(t)}(\delta)\) & Inflated-tube family restricted to $\ROI$, and sampling measures\\
\(E_t\) & Event that the local graph at iteration \(t\) satisfies local metric equivalence. \\
\(\Omega_t=\bigcap_{s=1}^tE_s\) & Good-run event that metric equivalence has held through iteration \(t\). \\
\(p_t,\ q_t\) & Conditional 
event probabilities \(\mathbb P(E_{t+1}\mid\mathcal F_t)\) and \(\mathbb P(E_{t+1}^c\mid\mathcal F_t)\). \\
\(\mathfrak{e}_{t+1}=\epsilon_{t+1}/[2\eta_t(1-\lambda_{t+1})]\) & Contribution of nearest-vertex graph-interpolation rounding. \\
\bottomrule
\end{tabularx}
\caption{Statistics, probability, and sampling notation.}
\label{tab:notation-main-probability}
\end{table}

Within this section, we make precise the interpretation of the graph-based algorithm as 
Riemannian SGD with a biased gradient estimator, bias controlled by graph approximation quality. We collect the notations used for geometry, graphs and optimization, and probability and statistics in Tables~\ref{tab:notation-main-geometry},~\ref{tab:notation-main-optimization}, and~\ref{tab:notation-main-probability}. A full notation list is provided in Appendix~\ref{sec:full-notation-appendix}.

The convergence analysis has three components. In Section~\ref{sec:biased-sgd-conv-theorem}, we establish a Riemannian SGD convergence result for biased estimators whose conditional means remain aligned with the Riemannian gradient and whose conditional second moments are controlled. In Section~\ref{sec:generalized-isomap}, we prove high-probability metric-equivalence guarantees for local $\epsilon$-neighborhood graphs when $\mathcal{M}_{1}$ is embedded in an ambient Riemannian manifold $\mathcal{M}_{2}$. In Section~\ref{sec:tube-containment-and-gradient-estimator}, we show how these graph guarantees control the discrepancy between graph interpolation and exact geodesic interpolation, and hence the bias of the implicit gradient estimator. Together, these ingredients yield the convergence guarantee for the graph-based algorithm in Section~\ref{sec:main-result-theorem}.

\subsection{Assumptions}\label{sec:assumptions}

We begin with the geometry of the intrinsic manifold $\mathcal{M}_{1}$ and the ambient space $\mathcal{M}_{2}$.

\begin{assumption}(\textbf{Intrinsic Geometry and Observed Distances})
\label{ass:embedded-submanifold}
Let $(\mathcal M_1,g_1)$ be a connected, compact
$C^3$ Riemannian manifold without boundary, with
$d=\dim(\mathcal M_1) < +\infty$. Let
\begin{equation}\label{eqn:intrinsic-sec-curv-bound}
-\infty<\kminInt
\leq \operatorname{sec}_{\mathcal M_1}
\leq \kmaxInt<+\infty,
\end{equation}
and set
$\kmagInt:=\max\{|\kminInt|,|\kmaxInt|\}$.
Let $d_{\mathrm{obs}}$ be a metric on $\mathcal M_1$.
Assume that at least one of the following holds.
\begin{enumerate}[label=(\alph*)]
\item \textbf{Riemannian ambient geometry.}
The manifold $\mathcal M_1$ is $C^3$ embedded in a Riemannian
manifold $(\mathcal M_2,g_2)$, with $g_1$ induced by $g_2$ and
$d_{\mathrm{obs}}=d_{\mathcal M_2}$ on $\mathcal M_1\times\mathcal M_1$.
Moreover,
\begin{equation}\label{eqn:extrinsic-sec-curv-bound}
-\infty<\kminExt
\leq \operatorname{sec}_{\mathcal M_2}
\leq \kmaxExt<+\infty.
\end{equation}
Set $\kmagExt:=\max\{|\kminExt|,|\kmaxExt|\}$.

\item \textbf{Local metric compatibility.}
There exist $K_{\mathrm{met}}\geq0$ and $r_{\mathrm{met}}>0$
such that, whenever $x,y\in\mathcal M_1$ and
$d_{\mathcal M_1}(x,y)\leq r_{\mathrm{met}}$,
\begin{equation}\label{eqn:abstract-local-metric-compatibility}
\begin{aligned}
d_{\mathrm{obs}}(x,y)
&\leq d_{\mathcal M_1}(x,y),\\
d_{\mathcal M_1}(x,y)^2-d_{\mathrm{obs}}(x,y)^2
&\leq K_{\mathrm{met}}d_{\mathcal M_1}(x,y)^4.
\end{aligned}
\end{equation}
\end{enumerate}
\end{assumption}

Under Assumption~\ref{ass:embedded-submanifold}\textup{(a)} we note that assuming \eqref{eqn:intrinsic-sec-curv-bound} is unnecessary. This is since
the second fundamental form $\Pi$ is continuous, so compactness
gives a constant $0<\fundBound<+\infty$ such that
\begin{equation}\label{eqn:second-fund-bound}
\|\Pi_p(v,v)\|_p\leq\fundBound,
\qquad
p\in\mathcal M_1,\quad
v\in T_p\mathcal M_1,\quad \text{and} \quad \|v\|_p=1.
\end{equation}
The Gauss equation permits the choice
$\kmagInt\leq\kmagExt+2\fundBound^2$ in this branch.
In the case where we do not have control on the ambient space, we include \eqref{eqn:intrinsic-sec-curv-bound} with Assumption~\ref{ass:embedded-submanifold}\textup{(b)}. This second branch is a consequence of the metric assumptions of \cite{XS26}, in the absence of ambient Riemannian structure; namely, for when $\mathcal{M}_{1} \subseteq \mathcal{P}_{2,\operatorname{ac}}(\mathbb{R}^{d})$. We emphasize that Assumption~\ref{ass:embedded-submanifold}\textup{(a)} implies Assumption~\ref{ass:embedded-submanifold}\textup{(b)} with informative constants depending on the geometry of $\mathcal{M}_{1}$ and $\mathcal{M}_{2}$.

We turn to our data localization assumption. We impose intrinsic localization for the Fr\'echet and interpolation
arguments and, under
Assumption~\ref{ass:embedded-submanifold}\textup{(a)},
ambient localization for the chord--arc comparisons. The purpose is threefold: $i)$ sufficient condition for existence and local uniqueness of $\widehat{x}_{N}$; $ii)$ existence and uniqueness of interpolating $\mathcal{M}_{1}$-geodesics and $\mathcal{M}_{2}$-geodesics; $iii)$ enabling for use of comparison theorems for quantitative control of discrepancies between $d_{\mathcal{M}_{1}}$ and $d_{\mathcal{M}_{2}}$.
For $S\subseteq\mathcal M_1$, define the uniform injectivity radius
\[
r_{\operatorname{inj}}^{\mathcal M_j}(S)
:=\inf_{p\in S}\operatorname{inj}_{\mathcal M_j}(p),
\]
for $j=1$ and, under branch \textup{(a)}, also for $j=2$.
We use the convention $1/\sqrt0=+\infty$.

\begin{assumption}(\textbf{Localization})
\label{ass:region-of-interest}
Suppose Assumption~\ref{ass:embedded-submanifold} holds. Let $\ROI\subseteq\mathcal{M}_{1}$ be a compact, geodesically convex, full dimensional region with $\operatorname{vol}_{\mathcal{M}_{1}}(\ROI) \in (0, \infty)$. Let $(z_{r})_{r=1}^{N}$ and all samples generated by Algorithm~\ref{algo:tube-sample-oracle} be contained in $\ROI$. Let
\begin{equation}\label{eqn:diameter-of-region-of-interest}
    \ROIdiam =
    \sup_{x,y\in\ROI}d_{\mathcal{M}_{1}}(x,y) 
     <
    \min\left\{
       r_{\operatorname{int}}, r_{\operatorname{amb}}
    \right\},
\end{equation}
where
\[
r_{\operatorname{int}}
:=
\frac12\min\left\{
\frac{\pi}{\sqrt{(\kmaxInt)_+}},
r_{\operatorname{inj}}^{\mathcal M_1}(\ROI)
\right\},
\]
and, for the selected branch of
Assumption~\ref{ass:embedded-submanifold},
\[
r_{\operatorname{amb}}
:=
\begin{cases}
\displaystyle
\min\left\{
\frac{\pi}{\fundBound},
\frac{\pi}{\sqrt{(\kmaxExt)_+}},
r_{\operatorname{inj}}^{\mathcal M_2}(\ROI)
\right\},
& \text{for Assumption~\ref{ass:embedded-submanifold}\textup{(a)}},\\[1em]
+\infty,
& \text{for Assumption~\ref{ass:embedded-submanifold}\textup{(b)}}.
\end{cases}
\]
\end{assumption}

The reasons for Assumption~\ref{ass:region-of-interest} are threefold. First, the intrinsic
curvature, diameter, and injectivity-radius restrictions -- encoded by $r_{\operatorname{int}}$ -- 
keep us in the standard regime for existence and uniqueness of the empirical
Fr\'echet mean; see, e.g., Afsari~\cite{afsari2011riemannian}. Second, they ensure
that the minimizing geodesics, logarithmic maps, Hessian bounds, and
comparison arguments used below are well defined on
$S_{\mathrm{reg}}$. 
Note that all iterates remain in $\ROI$ as the sampled points do.
Lastly, 
for the
normalized-volume sampler used by Algorithm~\ref{algo:tube-sample-oracle}, these localization conditions imply a lower bound on the probability of sampling within small intrinsic balls, uniform with respect to all admissible sampling regions. This enables later derived statistical guarantees.
We make this precise in Lemma~\ref{lem:sampling-mass-lower-bound} in Appendix~\ref{sec:proof-of-sampling-mass-lower-bound}.

\subsection{Biased Riemannian Stochastic Gradient Descent (SGD)}\label{sec:biased-sgd-conv-theorem}

\begin{algorithm}[H]
\caption{Riemannian Stochastic Gradient Descent}\label{algo:stochastic-gradient-descent}
\begin{algorithmic}
\State \textbf{Input:} Initialization $x_{0}$, stepsize schedule $(\eta_{t})_{t\geq0}$ \vspace{.1cm}
\While{Convergence Criterion Not Met}
\State \textbf{Generate Gradient Estimate:} $g_{t}\in T_{x_{t}}\mathcal{M}$.
\State \textbf{Update via Exponential Map:} $x_{t+1}\leftarrow\exp_{x_{t}}(-\eta_{t}g_{t})$.
\EndWhile
\end{algorithmic}
\end{algorithm}

Consider minimizing $f\colon\mathcal{M}\to\mathbb{R}$ on a Riemannian manifold $\mathcal{M}$. Algorithm~\ref{algo:stochastic-gradient-descent} generates
\begin{equation}\label{eqn:riemannian-sgd-update}
    x_{t+1}=\exp_{x_{t}}(-\eta_{t}g_{t}),
\end{equation}
where $g_{t}$ estimates $\nabla_{\mathcal{M}}f(x_{t})$. To analyze this algorithm, we work on a \textit{nice} region containing the entire trajectory and every geodesic step used in the analysis, formalized in the following.

\begin{assumption}\label{ass:function-assumptions-biased-SGD}
Let $(x_{t})_{t\geq0}$ be generated by Algorithm~\ref{algo:stochastic-gradient-descent}. Let $U\subseteq\mathcal{M}$ be a strongly geodesically convex neighborhood containing all iterates and a stationary minimizer $x_{\ast}\in U$ of $f$ over $U$, with $f_{\ast}:=f(x_{\ast})$ and $\nabla_{\mathcal{M}}f(x_{\ast})=0$. Assume that $f$ is continuously differentiable on $U$ and satisfies the following conditions.
\begin{enumerate}
    \item (\textbf{Local Strong Convexity}) There exists $\SCf>0$ such that, for all $x,y\in U$,
    \begin{equation}\label{eqn:local-strong-convexity}
        f(y)\geq f(x)
        +\langle\nabla_{\mathcal{M}}f(x),\log_{x}(y)\rangle_{x}
        +\frac{\SCf}{2}d_{\mathcal{M}}(x,y)^{2}.
    \end{equation}

    \item (\textbf{Local Geodesic Smoothness}) There exists $\smoothf>0$ such that, for all $x,y\in U$,
    \begin{equation}\label{eqn:local-smoothness}
        f(y)\leq f(x)
        +\langle\nabla_{\mathcal{M}}f(x),\log_{x}(y)\rangle_{x}
        +\frac{\smoothf}{2}d_{\mathcal{M}}(x,y)^{2}.
    \end{equation}
\end{enumerate}
\end{assumption}

\begin{remark}
Assumption~\ref{ass:function-assumptions-biased-SGD} implies the Polyak--\L{}ojasiewicz lower bound
\begin{equation}\label{eqn:local-pl-inequality}
    2\SCf\bigl(f(x)-f_{\ast}\bigr)
    \leq \|\nabla_{\mathcal{M}}f(x)\|_{x}^{2},
    \qquad \text{for } x\in U.
\end{equation}
\end{remark}

For the Fr\'echet mean functional $f_{P_N}$, this is a consequence of Assumptions~\ref{ass:embedded-submanifold} and~\ref{ass:region-of-interest}; proof deferred to Appendix~\ref{sec:proof-of-functional-SC-smooth-parameters}.

\begin{lemma}\label{lem:functional-SC-smooth-parameters}
Assume Assumptions~\ref{ass:embedded-submanifold} and~\ref{ass:region-of-interest} 
hold. Then the Fr\'{e}chet mean functional defined in \eqref{eqn:intro-frechet-problem-statement} is $\SCf$-strongly convex and $\smoothf$-smooth on $\ROI\subseteq\mathcal{M}_{1}$, where
\begin{equation*}
    \SCf:=
    \begin{cases}
        1, & \kmaxInt\leq0,\\
        \ROIdiam\sqrt{\kmaxInt}\cot(\sqrt{\kmaxInt}\ROIdiam), & \kmaxInt>0,
    \end{cases}
\end{equation*}
and
\begin{equation*}
    \smoothf:=
    \begin{cases}
        1, & \kminInt\geq0,\\
        \ROIdiam\sqrt{-\kminInt}\coth(\sqrt{-\kminInt}\ROIdiam), & \kminInt<0.
    \end{cases}
\end{equation*}
In particular, $0<\SCf\leq\smoothf<+\infty$ under \eqref{eqn:diameter-of-region-of-interest}.
\end{lemma}

Under Assumption~\ref{ass:function-assumptions-biased-SGD}, the standard unbiased analysis applies when $\mathbb{E}[g_{t}\mid\mathcal{F}_{t}]=\nabla_{\mathcal{M}}f(x_{t})$. We allow bias, retaining only the alignment and second-moment conditions needed for descent.

\begin{assumption}\label{assump:gradient-estimator-assumptions}
Let $(\mathcal{F}_{t})_{t\geq0}$ be a filtration such that $x_{t}$ is $\mathcal{F}_{t}$-measurable, and let $\eta_{t}\geq0$ be deterministic or $\mathcal{F}_{t}$-measurable. Assume the following for every $t\geq0$.
\begin{enumerate}
    \item (\textbf{Well-defined normal steps}) The random vector $g_{t}\in T_{x_{t}}\mathcal{M}$ is conditionally square-integrable, and $-\eta_{t}g_{t}$ lies almost surely in the injectivity domain of $x_{t}$. The radial geodesic $s\mapsto\exp_{x_{t}}(-s\eta_{t}g_{t})$, $s\in[0,1]$, remains in $U$. Consequently,
    \[
        \log_{x_{t}}(x_{t+1})=-\eta_{t}g_{t},
        \quad \text{and} \quad
        d_{\mathcal{M}}(x_{t},x_{t+1})=\eta_{t}\|g_{t}\|_{x_{t}}.
    \]

    \item (\textbf{Estimator Alignment}) There exists $c_{1}>0$ such that
    \begin{equation}\label{eqn:ass-alignment-inequality}
        \left\langle\nabla_{\mathcal{M}}f(x_{t}),
        \mathbb{E}[g_{t}\mid\mathcal{F}_{t}]\right\rangle_{x_{t}}
        \geq c_{1}\|\nabla_{\mathcal{M}}f(x_{t})\|_{x_{t}}^{2}.
    \end{equation}

    \item (\textbf{Controlled Second Moment}) There exist $\sigma_{0}^{2},\sigma_{1}^{2}\geq0$ such that
    \begin{equation}\label{eqn:gradient-estimator-second-moment-bound}
        \mathbb{E}[\|g_{t}\|_{x_{t}}^{2}\mid\mathcal{F}_{t}]
        \leq \sigma_{0}^{2}
        +\sigma_{1}^{2}\|\nabla_{\mathcal{M}}f(x_{t})\|_{x_{t}}^{2}.
    \end{equation}
\end{enumerate}
\end{assumption}

The alignment condition \eqref{eqn:ass-alignment-inequality} requires the conditional mean to retain a uniform component in the gradient direction, while \eqref{eqn:gradient-estimator-second-moment-bound} controls the quadratic term in the one-step bound. This second-moment condition controls the magnitude of the estimator. As a consequence, we can control objective gap improvement achieved by update \eqref{eqn:riemannian-sgd-update}, in expectation.

\begin{lemma}\label{lem:expected-one-step-improvement}
Assume Assumptions~\ref{ass:function-assumptions-biased-SGD} and~\ref{assump:gradient-estimator-assumptions} hold. If $\eta_{t}\smoothf\sigma_{1}^{2}\leq2c_{1},$ 
then
\begin{equation}\label{eqn:expected-one-step-improvement}
    \mathbb{E}[f(x_{t+1})\mid\mathcal{F}_{t}]-f(x_{t})
    \leq
    -2\eta_{t}\SCf\left(c_{1}-\frac{\eta_{t}\smoothf\sigma_{1}^{2}}{2}\right)
    [f(x_{t})-f_{\ast}]
    +\frac{\eta_{t}^{2}\smoothf\sigma_{0}^{2}}{2}.
\end{equation}
\end{lemma}

\begin{proof}
By local smoothness, the normal-step identities, and \eqref{eqn:riemannian-sgd-update},
\begin{align*}
    f(x_{t+1})-f(x_{t})
    &\leq
    \langle\nabla_{\mathcal{M}}f(x_{t}),\log_{x_{t}}(x_{t+1})\rangle_{x_{t}}
    +\frac{\smoothf}{2}d_{\mathcal{M}}(x_{t+1},x_{t})^{2}\\
    &=-\eta_{t}\langle\nabla_{\mathcal{M}}f(x_{t}),g_{t}\rangle_{x_{t}}
    +\frac{\eta_{t}^{2}\smoothf}{2}\|g_{t}\|_{x_{t}}^{2}.
\end{align*}
Taking conditional expectations and applying \eqref{eqn:ass-alignment-inequality} and \eqref{eqn:gradient-estimator-second-moment-bound} gives
\begin{align*}
    \mathbb{E}[f(x_{t+1})\mid\mathcal{F}_{t}]-f(x_{t})
    &\leq
    -\eta_{t}\left(c_{1}-\frac{\eta_{t}\smoothf\sigma_{1}^{2}}{2}\right)
    \|\nabla_{\mathcal{M}}f(x_{t})\|_{x_{t}}^{2}
    +\frac{\eta_{t}^{2}\smoothf\sigma_{0}^{2}}{2}.
\end{align*}
The coefficient of the squared gradient is nonnegative as $\eta_{t} L_{f} \sigma_{1}^{2} \leq 2c_{1}$. So, \eqref{eqn:local-pl-inequality} yields \eqref{eqn:expected-one-step-improvement}.
\end{proof}

Iterating the one-step estimate gives the following rates. The proof is deferred to \Cref{sec:biased-sgd-convergence-proof}.

\begin{theorem}\label{thm:biased-sgd-convergence}
Suppose Assumptions~\ref{ass:function-assumptions-biased-SGD} and~\ref{assump:gradient-estimator-assumptions} hold, and let $(x_{t})_{t\geq0}$ be generated by Algorithm~\ref{algo:stochastic-gradient-descent}.
\begin{enumerate}
    \item (\textbf{Constant Stepsize}) If $0<\eta_{t}=\bar{\eta}
        <
        \min\big\{
            \frac{c_{1}}{\smoothf\sigma_{1}^{2}},
            \frac{1}{\SCf c_{1}}
        \big\}$  for all $t$, with the convention $c_{1}/(\smoothf\sigma_{1}^{2})=+\infty$ when $\sigma_{1}^{2}=0$, then, for every $t\geq0$,
    \begin{equation}\label{eqn:constant-stepsize-convergence-guarantee}
        \mathbb{E}[f(x_{t})]-f_{\ast}
        \leq
        \frac{\bar{\eta}\smoothf\sigma_{0}^{2}}{2\SCf c_{1}}
        +(1-\bar{\eta}\SCf c_{1})^{t}
        \left(
            f(x_{0})-f_{\ast}
            -\frac{\bar{\eta}\smoothf\sigma_{0}^{2}}{2\SCf c_{1}}
        \right).
    \end{equation}

    \item (\textbf{Decaying Stepsize}) Let $\beta>1/(\SCf c_{1})$ and set
    \[
        \eta_{t}:=\frac{\beta}{\gamma+t}.
    \]
    If $\gamma>0$ is chosen so that $
        \eta_{0}<
        \min\big\{
            \frac{c_{1}}{\smoothf\sigma_{1}^{2}},
            \frac{1}{\SCf c_{1}}
        \big\},$
    then, for every $t\geq0$,
    \begin{equation}\label{eqn:decaying-stepsize-convergence-guarantee}
        \mathbb{E}[f(x_{t})]-f_{\ast}
        \leq \frac{\nu}{\gamma+t},
    \end{equation}
    where
    \begin{equation*}
        \nu:=\max\left\{\gamma(f(x_{0})-f_{\ast}),
            \frac{\beta^{2}\smoothf\sigma_{0}^{2}}{2(\beta\SCf c_{1}-1)}
        \right\}.
    \end{equation*}
\end{enumerate}
\end{theorem}

\subsection{ISOMAP for Non-Euclidean Data}
\label{sec:generalized-isomap}

Classical ISOMAP \cite{bernstein2000graph} gives conditions under which graph geodesics approximate intrinsic geodesics when $\mathcal{M}_{1}\subseteq\mathbb{R}^{d}$. We extend the metric-equivalence argument to when $\mathcal{M}_{1}\subseteq\mathcal{M}_{2}$. The following helps compare short intrinsic segments with ambient geodesic chords.

\begin{definition}
Let $S\subseteq\mathcal{M}_{1}$ be nonempty. For $s\geq0$, define
\begin{equation}\label{eqn:localized-transfer-extrinsic-to-intrinsic}
\rho_S(s)
:=
\sup\left\{
d_{\mathcal M_1}(x,y):
x,y\in S,\ d_{\mathrm{obs}}(x,y)\leq s
\right\}.
\end{equation}
\end{definition}

Under Assumption~\ref{ass:embedded-submanifold}\textup{(a)},
for $L>0$, define the ambient Jacobi comparison factor
\begin{equation*}
    C_{J}(L)
    :=
    \max\left\{
        1,
        \frac{A_{\kminExt}(L)}{B_{\kmaxExt}(L)}
    \right\},
\end{equation*}
where $A_{\kminExt}$ and $B_{\kmaxExt}$ are defined in
\eqref{eqn:constants-for-Jacobi-control}. This technical factor accounts for the ambient-curvature contribution in a short chord-arc estimate.

We now give deterministic conditions under which an ambient-weighted graph satisfies a two-sided metric comparison. The proof is deferred to Appendix~\ref{sec:proof-of-deterministic-metric-equivalence-for-isomap}.

\begin{theorem}\label{thm:deterministic-metric-equivalence-for-ISOMAP}
Assume Assumption~\ref{ass:embedded-submanifold}\textup{(a)} holds. Let $\{x_{i}\}_{i=1}^{n}\subseteq S\subseteq\mathcal{M}_{1},$ and $G$ be a graph such that every edge $xy$ has weight $d_{\mathcal{M}_{2}}(x,y)$. Fix
$\lambda_{1},\lambda_{2}\in(0,1)$. Let $\compScale>0$ satisfy
\[
    \compScale\leq\frac{\pi}{\fundBound},
    \quad \text{and} \quad
    \compScale<\frac{\pi}{\sqrt{\kmaxExt}}
    \quad\text{if }\kmaxExt>0.
\]
Assume that $s_{0}>0$ satisfies the following conditions.
\begin{enumerate}
    \item (\textbf{Localization for Comparisons})
        $\rho_{S}(s_{0})\leq\compScale;$

    \item (\textbf{Short-Edge Distortion})
        $\rho_{S}(s_{0})
        \leq
        \left(
            \frac{\fundBound^{2}}{24}
            +\frac{5\kmagExt C_{J}(\compScale)}{18}
        \right)^{-1/2}
        \sqrt{\lambda_{1}};$

    \item (\textbf{Ambient Normal Neighborhood})
        $\rho_{S}(s_{0})
        <
        r_{\operatorname{inj}}^{\mathcal{M}_{2}}(S)
        :=
        \inf_{p\in S}\operatorname{inj}_{\mathcal{M}_{2}}(p).$
\end{enumerate}
Let $\epsilon_{\min},\epsilon_{\max},\delta>0$ satisfy
\begin{enumerate}
    \setcounter{enumi}{3}
    \item (\textbf{Connectivity and Covering})
        $\epsilon_{\min}\leq\epsilon_{\max}\leq s_{0},
        $ and $
        \delta\leq\frac{\lambda_{2}\epsilon_{\min}}{4}.$
\end{enumerate}
Assume further that the following hold for the graph construction and data sampling:
\begin{enumerate}
    \setcounter{enumi}{4}
    \item (\textbf{Graph Construction})
    $G$ contains every edge $xy$ satisfying
    $d_{\mathcal{M}_{2}}(x,y)\leq\epsilon_{\min}$ and contains no edge
    $xy$ satisfying
    $d_{\mathcal{M}_{2}}(x,y)>\epsilon_{\max}$.

    \item (\textbf{Data Sampling})
    The data set $\{x_{i}\}_{i=1}^{n}$ is a $\delta$-cover of $S$ with
    respect to $d_{\mathcal{M}_{1}}$.

    \item (\textbf{Geodesic Convexity of $S$})
    Every pair of points in $S$ is joined by a minimizing
    $\mathcal{M}_{1}$-geodesic contained in $S$.
\end{enumerate}
Then
\begin{equation}\label{eqn:metric-equivalence}
    (1-\lambda_{1})d_{\mathcal{M}_{1}}(x,y)
    \leq d_{G}(x,y)
    \leq(1+\lambda_{2})d_{\mathcal{M}_{1}}(x,y)
\end{equation}
for every $x,y\in\{x_{i}\}_{i=1}^{n}$.
\end{theorem}

The lower bound in \eqref{eqn:metric-equivalence} is edgewise. The first
three conditions ensure that every graph edge lies in a regime where its
ambient length is at least $(1-\lambda_{1})$ times the corresponding
intrinsic distance. Summing along graph paths then rules out ambient
shortcuts. The normal-neighborhood condition is imposed on
$\rho_{S}(s_{0})$, rather than on $s_{0}$, because the comparison argument
follows the intrinsic geodesic joining the edge endpoints.

The upper bound follows by discretizing a minimizing intrinsic geodesic in
$S$. The $\delta$-cover provides nearby vertices, and the graph contains
the resulting short edges because
$\delta\leq\lambda_{2}\epsilon_{\min}/4$. Thus the exact graph and
covering scales are
\begin{equation}\label{eqn:non-interpretable-scales}
    \epsilon_{\max}\leq s_{0},
    \quad \text{and} \quad
    \delta\leq\frac{\lambda_{2}\epsilon_{\min}}{4}.
\end{equation}
When $\mathcal{M}_{2}=\mathbb{R}^{d}$, the ambient-curvature term vanishes
and the theorem reduces to the classical Euclidean chord--arc argument.

Without a controlled Riemannian ambient space, Assumption~\ref{ass:embedded-submanifold}\textup{(b)} gives the same outcomes.

\begin{proposition}
\label{prop:abstract-metric-isomap}
Assume Assumption~\ref{ass:embedded-submanifold}\textup{(b)}.
Let $\lambda_1,\lambda_2\in(0,1)$,
$0<\epsilon_{\min}\leq\epsilon_{\max}$, and $\delta>0$.
Let $V(G)=\{x_i\}_{i=1}^n\subseteq S\subseteq\mathcal M_1$,
and give every edge $uv$ weight $d_{\mathrm{obs}}(u,v)$.
Assume the graph-construction, data-sampling, and
geodesic-convexity conditions of
Theorem~\ref{thm:deterministic-metric-equivalence-for-ISOMAP},
with $d_{\mathcal M_2}$ replaced by $d_{\mathrm{obs}}$.
If
\[
\delta\leq\frac{\lambda_2\epsilon_{\min}}4,
\qquad
\rho_S(\epsilon_{\max})\leq r_{\mathrm{met}},
\quad \text{and} \quad
K_{\mathrm{met}}\rho_S(\epsilon_{\max})^2
\leq2\lambda_1-\lambda_1^2,
\]
then \eqref{eqn:metric-equivalence} holds for every pair of vertices.
\end{proposition}

\begin{proof}
Subdividing a minimizing intrinsic geodesic into segments of length
at most $r_{\mathrm{met}}$ and applying the triangle inequality
shows that $d_{\mathrm{obs}}\leq d_{\mathcal M_1}$ globally.
For every graph edge $uv$, put
$\ell=d_{\mathcal M_1}(u,v)\leq\rho_S(\epsilon_{\max})$.
Then
\[
d_{\mathrm{obs}}(u,v)
\geq\ell\sqrt{1-K_{\mathrm{met}}\ell^2}
\geq(1-\lambda_1)\ell.
\]
Summing along graph paths gives the lower bound.
The upper-bound proof of
Lemma~\ref{lem:global-upper-bound-on-graph-distance}
applies with $d_{\mathrm{obs}}$ in place of $d_{\mathcal M_2}$,
using $d_{\mathrm{obs}}\leq d_{\mathcal M_1}$.
\end{proof}

To obtain a high-probability statement, it remains to guarantee the covering condition. The proof of the following corollary is deferred to Appendix~\ref{sec:proof-of-probabilistic-isomap}.

\begin{corollary}\label{cor:probabilistic-isomap}
Assume the hypotheses of either
Theorem~\ref{thm:deterministic-metric-equivalence-for-ISOMAP}
or Proposition~\ref{prop:abstract-metric-isomap}
hold on a compact geodesically convex set
$S\subseteq\mathcal M_1$, except possibly the data-sampling condition.
Let
$\{x_{i}\}_{i=1}^{n}$ be i.i.d. samples from a probability measure
$\mu_{S}$ supported on $S$, and let $\zeta\in(0,1)$. Define
\begin{equation}
    q_{S}(\delta/2)
    :=
    \inf_{x\in S}
    \mu_{S}\bigl(
        B_{\mathcal{M}_{1}}(x,\delta/2)\cap S
    \bigr),
\end{equation}
and assume $q_{S}(\delta/2)>0$. If
\begin{equation}\label{eqn:probabilistic-isomap-sampling-bound}
    n
    \geq
    \frac{1}{q_{S}(\delta/2)}
    \left(
        \log\mathcal{N}(S,\delta/2)
        +\log\left(\frac{1}{\zeta}\right)
    \right),
\end{equation}
then, with probability at least $1-\zeta$,
\begin{equation}\label{eqn:probabilistic-isomap-metric-equivalence}
    (1-\lambda_{1})d_{\mathcal{M}_{1}}(x_{i},x_{j})
    \leq d_{G}(x_{i},x_{j})
    \leq(1+\lambda_{2})d_{\mathcal{M}_{1}}(x_{i},x_{j})
\end{equation}
for all $i,j\in\{1,\dots,n\}$.
\end{corollary}

Finally, the modulus $\rho_{S}$ can be replaced by an explicit local scale on compact sets.

\begin{lemma}\label{lem:local-intrinsic-extrinsic-compatibility}
Assume Assumption~\ref{ass:embedded-submanifold} holds,
and let $S\subseteq\mathcal M_1$ be nonempty and compact.
There exist constants $C_{\mathrm{loc},S}\geq1$ and
$s_{\mathrm{loc},S}>0$ such that
\[
d_{\mathcal M_1}(x,y)
\leq C_{\mathrm{loc},S}d_{\mathrm{obs}}(x,y)
\]
whenever $x,y\in S$ and
$d_{\mathrm{obs}}(x,y)\leq s_{\mathrm{loc},S}$.
In branch \textup{(b)}, one may take
$C_{\mathrm{loc},S}=\sqrt2$.
\end{lemma}

\begin{proof}
First suppose Assumption~\ref{ass:embedded-submanifold}\textup{(a)}
holds, so $d_{\mathrm{obs}}=d_{\mathcal M_2}$.
For each $p\in S$, choose an adapted coordinate neighborhood $U_{p}\subseteq\mathcal{M}_{2}$ in which $\mathcal{M}_{1}\cap U_{p}$ is represented by a coordinate subspace. After shrinking $U_{p}$, the ambient and intrinsic metric tensors are uniformly equivalent to the Euclidean metric on a convex coordinate neighborhood. Thus, there exist $C_{p}\geq1$ and $s_{p}>0$ such that
    $d_{\mathcal{M}_{1}}(x,y)
    \leq C_{p}d_{\mathcal{M}_{2}}(x,y)$
for $x,y\in S\cap U_{p}$ satisfying $d_{\mathcal{M}_{2}}(x,y)\leq s_{p}$. A finite subcover of the compact set $S$, together with a Lebesgue-number argument, yields uniform constants $C_{\mathrm{loc},S}$ and $s_{\mathrm{loc},S}$.

Under Assumption~\ref{ass:embedded-submanifold}\textup{(b)},
the local upper bound makes the identity
$(S,d_{\mathcal M_1})\to(S,d_{\mathrm{obs}})$ continuous.
Compactness implies that its inverse is uniformly continuous,
and hence $\rho_S(s)\to0$ as $s\downarrow0$.
Choose $0<r_*\leq r_{\mathrm{met}}$ with
$K_{\mathrm{met}}r_*^2\leq1/2$, and then choose
$s_{\mathrm{loc},S}>0$ so that
$\rho_S(s_{\mathrm{loc},S})\leq r_*$.
For $d_{\mathrm{obs}}(x,y)\leq s_{\mathrm{loc},S}$,
\[
d_{\mathrm{obs}}(x,y)^2
\geq(1-K_{\mathrm{met}}r_*^2)d_{\mathcal M_1}(x,y)^2
\geq\tfrac12d_{\mathcal M_1}(x,y)^2.
\]
Thus $C_{\mathrm{loc},S}=\sqrt2$ is admissible.
\end{proof}

Thus,
    $\rho_{S}(s_{0})
    \leq C_{\mathrm{loc},S}s_{0}$
whenever $s_{0}\leq s_{\mathrm{loc},S}$. Under Assumption~\ref{ass:embedded-submanifold}\textup{(a)},
the following is sufficient for the scale conditions of
Theorem~\ref{thm:deterministic-metric-equivalence-for-ISOMAP}:
\begin{equation}\label{eqn:simplified-assumptions-on-intrinsic-ambient-distances}
    s_{0}
    <
    \min\left\{
        s_{\mathrm{loc},S},
        \frac{\compScale}{C_{\mathrm{loc},S}},
        \frac{1}{C_{\mathrm{loc},S}}
        \left(
            \frac{\fundBound^{2}}{24}
            +\frac{5\kmagExt C_{J}(\compScale)}{18}
        \right)^{-1/2}
        \sqrt{\lambda_{1}},
        \frac{r_{\operatorname{inj}}^{\mathcal{M}_{2}}(S)}
        {C_{\mathrm{loc},S}}
    \right\}.
\end{equation}

Particularly, for Assumption~\ref{ass:embedded-submanifold}\textup{(a)}, for sufficiently small
$\lambda:=\lambda_1=\lambda_2$, the following scales may be chosen:
\begin{equation}\label{eqn:readable-scale-requirements}
    \epsilon_{\max}\leq s_{0}\lesssim\sqrt{\lambda},
    \quad \text{and} \quad
    \delta\leq\frac{\lambda\epsilon_{\min}}{4}
    \lesssim\lambda^{3/2},
\end{equation}
where the hidden constants depend only on the local metric-comparison, curvature, second-fundamental-form, and injectivity-radius bounds on $S$.

\subsection{Graph Representation and Interpolation as Biased Riemannian SGD} \label{sec:tube-containment-and-gradient-estimator}

This subsection formalizes the random local graph construction-interpolation process used by Algorithm~\ref{algo:graph-algorithm} with the tube sampler, Algorithm~\ref{algo:tube-sample-oracle}. We first define the filtration with respect to both algorithmic and sampling randomness. This enables us to define the per-iteration event for when local metric equivalence holds. We then use local metric equivalence to control the discrepancy between graph and exact interpolation, enabling interpretation of Algorithm~\ref{algo:graph-algorithm} as Riemannian SGD with structured estimators. We finish by providing sufficient sampling rates to recover these structured gradient estimators with high probability.

\paragraph{Random process.}
Fix a horizon $T\geq1$, and suppose $x_{0}$ and $G_{0}$ are given. At iteration $t+1$, where $0\leq t\leq T-1$, draw
\[
    I_{t+1}\sim\operatorname{Unif}([N])
\]
independently of the past, and set $
    y_{t+1}:=z_{I_{t+1}}.
$
Given a radius $r_{t+1}>0$, define the target tube
\begin{equation*}
    \mathcal{T}_{t+1}
    :=
    T_{r_{t+1}}
    \bigl(
        \gamma_{x_{t}\rightarrow y_{t+1}}^{\mathcal{M}_{1}}
    \bigr),
\end{equation*}
where the tube is defined in \eqref{eqn:geodesic-tube}. Let $\xi_{t+1}$ collect all randomness used by the local sampler conditional on the current state and the selected reference point. We encode the total randomness generated through iteration $t$ by
\begin{equation}\label{eqn:filtration-definition}
    \mathcal{F}_{t}
    :=
    \sigma(
        x_{0},
        I_{1},\xi_{1},
        \ldots,
        I_{t},\xi_{t}
    ),
    \qquad \text{for }
    0\leq t\leq T,
\end{equation}
with $\mathcal{F}_{0}:=\sigma(x_{0})$. This construction articulates the two-stage randomness per iteration: $i)$ random choice of target point; $ii)$ local sampling for graph approximations. We assume the graph construction, shortest-path selection, and interpolation are measurable and deterministic once $(x_{0},I_{1},\xi_{1},\ldots,I_{t},\xi_{t})$ are fixed.

To control the discrepancy between $\gamma_{x_{t} \rightarrow y_{t+1}}^{\mathcal{M}_{1}}$ and $\gamma_{x_{t} \rightarrow y_{t+1}}^{G_{t+1}}$, we must control vertices just outside of $\mathcal{T}_{t+1}$ by one admissible ambient edge-length. Set
$r_{t+1}:=\alpha_{t+1}d_{\mathcal M_1}(x_t,y_{t+1}),$ and 
$\widetilde r_{t+1}:=r_{t+1}
+\rho_{S_{\mathrm{reg}}}(\epsilon_{t+1}).$
The target tube, inflated tube, and sampling support are respectively
\begin{equation}\label{eqn:inflated-tube}
\begin{split}
\mathcal T_{t+1}
&:= \mathcal T_{r_{t+1}}\bigl(\gamma^{\mathcal M_1}_{x_t\to y_{t+1}}\bigr),
\\
\widetilde{\mathcal{T}}_{t+1}
&:= \mathcal T_{\widetilde r_{t+1}}
\bigl(\gamma^{\mathcal M_1}_{x_t\to y_{t+1}}\bigr),
\qquad
\widetilde S_{t+1}:=S_{\mathrm{reg}}\cap\widetilde{\mathcal{T}}_{t+1}.
\end{split}
\end{equation}
The graph is constructed on the new samples in $\widetilde S_{t+1}$
together with $x_t$ and $y_{t+1}$. Define event
\begin{equation}\label{eqn:defining-per-iterate-good-event}
\begin{split}
E_{t+1}:=\bigl\{\,
&(1-\lambda_{t+1})d_{\mathcal M_1}(u,v)
\le d_{G_{t+1}}(u,v),
\quad u,v\in V(G_{t+1}),\\
&d_{G_{t+1}}(a,v)\le
(1+\lambda_{t+1})d_{\mathcal M_1}(a,v), \quad a\in\{x_t,y_{t+1}\},\quad v\in V(G_{t+1})
\,\bigr\}.
\end{split}
\end{equation}
We define the event for local metric equivalence throughout the run inductively:
\begin{equation}\label{eqn:trajectory-of-good-events}
    \Omega_{0}:=\Omega,
    \quad \text{and} \quad
    \Omega_{t}:=\bigcap_{s=1}^{t}E_{s}.
\end{equation}

\paragraph{Metric equivalence and tube containment.}
Fix $0\leq t\leq T-1$ and suppose $0<\eta_{t}\leq1$. Define the exact
stochastic gradient estimator
\begin{equation}\label{eqn:unbiased-gradient-estimator}
    \widetilde g_{t}
    :=
    -\log_{x_{t}}^{\mathcal{M}_{1}}(y_{t+1}).
\end{equation}
Because $I_{t+1}$ is uniform and independent of $\mathcal{F}_{t}$, we have 
    $\mathbb{E}[\widetilde g_{t}\mid\mathcal{F}_{t}]
    =
    \nabla_{\mathcal{M}_{1}}f(x_{t}).$
The corresponding exact-geodesic update is
\[
    \widetilde x_{t+1}
    :=
    \exp_{x_{t}}^{\mathcal{M}_{1}}(-\eta_{t}\widetilde g_{t}).
\]
Let $x_{t+1}$ be the graph-based iterate by Algorithm~\ref{algo:graph-algorithm}. By Assumption~\ref{ass:region-of-interest}, the logarithm from $x_t$ to $x_{t+1}$ is unique, and the implicit gradient estimator is
\begin{equation}\label{eqn:biased-graph-gradient-estimator}
    g_{t}
    :=
    -\frac{1}{\eta_{t}}
    \log_{x_{t}}^{\mathcal{M}_{1}}(x_{t+1}).
\end{equation}
When $x_t=y_{t+1}$, we set
$x_{t+1}=\widetilde x_{t+1}=x_t$ and
$g_t=\widetilde g_t=0$. For the fresh-graph construction,
$G_{t+1}$ consists of the single vertex $x_t$, so $E_{t+1}$
holds without sampling. For what follows, we assume $y_{t+1} \neq x_{t}$. The following gives a sufficient condition for the nodes of a graph
geodesic segment to be within a geodesic tube of the true geodesic segment; proof deferred to
Appendix~\ref{sec:proof-of-conditions-on-metric-equivalence-for-tube-containment}.

\begin{lemma}\label{lem:conditions-on-metric-equivalence-for-tube-containment}
Assume Assumptions~\ref{ass:embedded-submanifold} and 
\ref{ass:region-of-interest} hold. Let $x_{t}$ and $y_{t+1}$ be given by Algorithm~\ref{algo:graph-algorithm}. Set
$\ell_{t+1}
    :=
    d_{\mathcal{M}_{1}}(x_{t},y_{t+1}) > 0,
    $ and $
    r_{t+1}
    :=
    \radiusScaling\ell_{t+1}$ for $\alpha_{t} \in (0,1]$. Define
\begin{equation}\label{eqn:function-for-tube-containment}
\widetilde{\Psi}(\ell,r)
:=
\begin{cases}
\displaystyle
\frac{2}{\sqrt{\kmaxInt}\ell}
\arccos\left(
    \cos\left(\frac{\sqrt{\kmaxInt}\ell}{2}\right)
    \cos\left(\sqrt{\kmaxInt}r\right)
\right),
& \kmaxInt>0,\\[1.2em]
\displaystyle
\sqrt{1+4\left(\frac{r}{\ell}\right)^{2}},
& \kmaxInt\leq0.
\end{cases}
\end{equation}
Let $G_{t+1}$ be the $\epsilon_{t+1}$-neighborhood graph constructed on the vertices in the inflated tube $\widetilde{\mathcal S}_{t+1}$ with points $x_{t}, y_{t+1}$. Suppose event $E_{t+1}$ holds. If
\begin{equation}\label{eqn:bound-on-metric-equivalence-parameters-for-tube-containment}
    \frac{1+\lambda_{2}^{(t+1)}}
    {1-\lambda_{1}^{(t+1)}}
    < \widetilde{\Psi}(\ell_{t+1}, \alpha_{t+1}\ell_{t+1}) = \widetilde{\Psi}(\ell_{t+1}, r_{t+1}),
\end{equation}
then every vertex of $\gamma_{x_{t}\rightarrow y_{t+1}}^{G_{t+1}}$ belongs to $\mathcal{T}_{t+1}$.
\end{lemma}

We next show that the sufficient condition really just depends on the radius scaling $\alpha$. The proof follows from direct computation, and is given in Appendix~\ref{sec:appendix-proof-of-comparison-with-quadratic} for completion.
\begin{lemma}\label{lem:comparison-with-quadratic}
    Assume Assumption~\ref{ass:region-of-interest} holds. Then, for all $\ell \in (0,\ROIdiam]$ and $\alpha \in (0,1]$ 
    \begin{equation} \label{eqn:comparison-controlling-tube-containment-function}
        c_{\Psi}\alpha^{2} \leq \widetilde{\Psi}(\ell, \alpha \ell) - 1 \leq 2 \alpha^{2},
    \end{equation}
    with
    \begin{equation}\label{eqn:explicit-constant-for-lower-bound-control-of-comparison}
    c_{\Psi} := \begin{cases} \displaystyle \frac{4}{1 + \sqrt{5}}\cos\left(\frac{\sqrt{\kmaxInt} \ROIdiam}{2}\right)\operatorname{sinc}\left(\sqrt{\kmaxInt} \ROIdiam\right), \qquad &\kmaxInt > 0, \vspace{0.2em} \\
    \displaystyle \frac{4}{1 + \sqrt{5}}, \qquad & \kmaxInt \leq 0.
        \end{cases}
    \end{equation}
\end{lemma}

Thus, we obtain a sufficient condition on $\lambda_{t+1} = \lambda_{1}^{(t+1)} = \lambda_{2}^{(t+1)}$ as a function of radius scaling -- regardless of geodesic -- to guarantee tube containment. Explicitly, suppose
\begin{equation}\label{eqn:sufficicient-metric-condition-for-tube-containment}
    \lambda_{t+1} \leq \frac{c_{\Psi}}{8} \alpha_{t+1}^{2}.
\end{equation}
Then, since $0 < c_{\Psi} \leq 4/(1+\sqrt{5}) < 2$, it follows that for all $\alpha_{t+1} \in (0,1]$ that $\lambda_{t+1} < 1/2$, and 
\begin{equation*}
    \frac{1+\lambda_{t+1}}{1-\lambda_{t+1}}
\le1+4\lambda_{t+1}
\le1+\frac{c_\Psi}{2}\alpha_{t+1}^2
<1+c_\Psi\alpha_{t+1}^2
\le\widetilde\Psi(\ell_{t+1},\alpha_{t+1}\ell_{t+1}).
\end{equation*}
where the first inequality follows from $\lambda_{t+1} \leq 1/2$, and the second from \eqref{eqn:sufficicient-metric-condition-for-tube-containment}, and the last from Lemma~\ref{lem:comparison-with-quadratic}. Then, by Lemma~\ref{lem:conditions-on-metric-equivalence-for-tube-containment} tube containment follows. 

Intuitively, graph interpolation and exact geodesic interpolation should produce similar iterates under tube containment.
The following formalizes this; proof deferred to Appendix~\ref{sec:proof-of-discrepancy-from-geodesic-interpolation}.

\begin{lemma}\label{lem:discrepancy-from-geodesic-interpolation}
Assume Assumptions~\ref{ass:embedded-submanifold}
and~\ref{ass:region-of-interest} hold.
Let $\eta_t,\alpha_{t+1}\in(0,1]$, and let
$x_t,y_{t+1},x_{t+1}$ be generated by Algorithm~\ref{algo:graph-algorithm}, and $\widetilde{x}_{t+1} = \exp_{x_{t}}(\eta_{t} \log_{x_{t}}(y_{t+1}))$. Suppose $E_{t+1}$ and \eqref{eqn:bound-on-metric-equivalence-parameters-for-tube-containment}
hold. 
Then
\begin{equation}\label{eqn:discrepancy-from-geodesic-interpolation}
\begin{aligned}
    d_{\mathcal{M}_{1}}(x_{t+1},\tilde x_{t+1})
    \leq
    \left(
        \frac{
            \lambda_{1}^{(t+1)}
            +
            \lambda_{2}^{(t+1)}
        }{
            1-\lambda_{1}^{(t+1)}
        }
        \eta_{t}
        +2\radiusScaling
    \right)
    d_{\mathcal{M}_{1}}(x_{t},y_{t+1}) + \frac{\epsilon_{t+1}}{2(1-\lambda_{1}^{(t+1)})}.
\end{aligned}
\end{equation}
\end{lemma}

By controlling $d_{\mathcal{M}_{1}}(x_{t+1}, \widetilde{x}_{t+1})$, we can control $\|g_{t} - \tilde{g}_{t}\|$. To do so, we use a Rauch-comparison argument, as in the proof of Lemma~4.5 of \cite{fetecau2022wellposedness}.
For $x,y,z \in \ROI$ and $q \in \gamma_{y \rightarrow z}^{\mathcal{M}_{1}}$,
\begin{equation}\label{eqn:uniform-rauche-comparison-constant}
    \|D_{q}\log_{x}\| \leq \ROIraucheC(\ROIdiam, \kmaxInt)
    :=
    \begin{cases}
        1,
        & \kmaxInt\leq0,\\[0.4em]
        \displaystyle
        \frac{\sqrt{\kmaxInt}\ROIdiam}
        {\sin(\sqrt{\kmaxInt}\ROIdiam)},
        & \kmaxInt>0.
    \end{cases}
\end{equation}
Integrating this bound along the minimizing geodesic from $y$ to $z$ gives
\begin{equation}\label{eqn:Rauch-comparison-estimate}
    \left\|
        \log_{x}^{\mathcal{M}_{1}}(y)
        -
        \log_{x}^{\mathcal{M}_{1}}(z)
    \right\|_{x}
    \leq
    \ROIraucheC(\ROIdiam, \kmaxInt)
    d_{\mathcal{M}_{1}}(y,z).
\end{equation}
Thus, after noting $\|g_{t} - \tilde{g}_t\| = \eta_{t}^{-1}\|\log_{x_{t}}(x_{t+1}) - \log_{x_{t}}(\widetilde{x}_{t+1})\|$, the iterate discrepancy is transferred by direct application of \eqref{eqn:Rauch-comparison-estimate} with Lemma~\ref{lem:discrepancy-from-geodesic-interpolation}; proof is omitted.
\begin{corollary}
\label{cor:implicit-gradient-discrepancy}
Suppose the hypotheses of
Lemma~\ref{lem:discrepancy-from-geodesic-interpolation} hold. Let $\widetilde{g}_{t}$ and $g_{t}$ be defined by \eqref{eqn:unbiased-gradient-estimator} and \eqref{eqn:biased-graph-gradient-estimator}, respectively, and $\ROIraucheC := \ROIraucheC(\ROIdiam, \kmaxInt)$ by \eqref{eqn:uniform-rauche-comparison-constant}. Then
\begin{align}
\|g_t-\widetilde g_t\|_{x_t}
\leq C_{\log}
\Bigg[
    \left(
    \frac{
    \lambda_1^{(t+1)}+\lambda_2^{(t+1)}
    }{
    1-\lambda_1^{(t+1)}
    }
    +
    \frac{2\alpha_{t+1}}{\eta_t}
    \right)
    d_{\mathcal M_1}(x_t,y_{t+1})
    +
    \frac{\epsilon_{t+1}}
    {2\eta_t(1-\lambda_1^{(t+1)})}
\Bigg].
\end{align}
\end{corollary}

To simplify later notation, define the following quantity, which captures the tangent-space discrepancy due to \textit{graph interpolation rounding}:
\begin{equation}\label{eqn:rounding-error-bound}
    \mathfrak{e}_{t+1} := \frac{\epsilon_{t+1}}{2\eta_{t}(1 - \lambda_1^{(t+1)})}.
\end{equation}

\paragraph{Conditional second-moment and alignment bounds.} 
The event $E_{t+1}$ enforces metric equivalence on $\widetilde{\mathcal{S}}_{t+1}$, resulting in a quantifiable discrepancy between graph geodesic interpolation and $\mathcal{M}_{1}$-geodesic interpolation. As a consequence, we obtain control on the second-moment of $g_{t}$, up to the factor $\mathbb{P}(E_{t+1} \mid \mathcal{F}_{t})$. Proof deferred to Appendix~\ref{sec:proof-of-metric-equivalence-gives-second-moment-bounds}.

\begin{lemma}
\label{lem:metric-equivalence-gives-second-moment-bounds}
Assume Assumptions~\ref{ass:embedded-submanifold} and~\ref{ass:region-of-interest} hold. Assume the tube-containment condition \eqref{eqn:bound-on-metric-equivalence-parameters-for-tube-containment} holds on $E_{t+1}$. Let $0<\eta_t\leq1$, $0\leq t\leq T-1$ be given, and $p_t:=\mathbb P(E_{t+1}\mid\mathcal F_t)>0$.
Define $C_{t+1}:=
(1+\lambda_2^{(t+1)})/(1-\lambda_1^{(t+1)})$, and $\mathfrak{e}_{t+1}$ by \eqref{eqn:rounding-error-bound}.
Then, for every $\vartheta>0$,
\begin{align}
\label{eqn:metric-equivalence-second-moment-bound-general}
 \mathbb E[\|g_t\|_{x_t}^2\mid\mathcal F_t,E_{t+1}]
 \leq
 \frac{(1+\vartheta)C_{t+1}^2}{p_t}
 \left(
 2f_{P_N}(\widehat x_N)
 +\frac{1}{m_f}\|\nabla_{\mathcal M_1}f_{P_N}(x_t)\|_{x_t}^2
 \right)
 +(1+\vartheta^{-1})\mathfrak e_{t+1}^2.
\end{align}
In particular, if
$\lambda_1^{(t+1)}=\lambda_2^{(t+1)}\leq1/2$, then taking
$\vartheta=1/3$ gives
\begin{equation}
\label{eqn:metric-equivalence-second-moment-bound}
 \mathbb E[\|g_t\|_{x_t}^2\mid\mathcal F_t,E_{t+1}]
 \leq
 \frac{1}{p_t}
 \left(
 24f_{P_N}(\widehat x_N)
 +\frac{12}{m_f}\|\nabla_{\mathcal M_1}f_{P_N}(x_t)\|_{x_t}^2
 \right)
 +4\mathfrak e_{t+1}^2.
\end{equation}
\end{lemma}

It remains to control alignment between $g_{t}$ and $\nabla_{\mathcal{M}_{1}}f(x_{t})$. 
This is the content of the following; proof is deferred to Appendix~\ref{sec:proof-of-tube-containment-to-gradient-estimator-control}.

\begin{lemma}
\label{lem:tube-containment-to-gradient-estimator-control}
Assume Assumptions~\ref{ass:embedded-submanifold} and~\ref{ass:region-of-interest} hold. Let $(x_{t})_{t=0}^{T}$ be generated by Algorithm~\ref{algo:graph-algorithm}. Assume the tube-containment condition \eqref{eqn:bound-on-metric-equivalence-parameters-for-tube-containment} holds on $E_{t+1}$. Fix $c_{1}\in(0,1)$ and $0\leq t\leq T-1$. Suppose $\mathbb{P}(E_{t+1} \mid \mathcal{F}_{t})>0$. Define $D_{t} := \max_{1\leq i\leq N} d_{\mathcal{M}_{1}}(x_{t},z_{i})$, $\mathfrak{e}_{t+1}$ by \eqref{eqn:rounding-error-bound}, $\ROIraucheC := \ROIraucheC(\ROIdiam,\kmaxInt)$ by \eqref{eqn:uniform-rauche-comparison-constant},
and $q_{t} := \mathbb{P}(E_{t+1}^{c} \mid \mathcal F_{t}).$
If
\begin{align}
\label{eqn:tube-radius-vs-distance-to-minimizer}
 &\ROIraucheC
 \left[
 \left(
 \frac{\lambda_1^{(t+1)}+\lambda_2^{(t+1)}}
 {1-\lambda_1^{(t+1)}}\eta_t
 +2\alpha_{t+1}
 \right)D_t
 + \eta_{t}\mathfrak{e}_{t+1}
 \right]
 +2q_{t}\eta_tD_t
 \leq
 (1-c_1)\eta_t
 \|\nabla_{\mathcal M_1}f_{P_N}(x_t)\|_{x_t},
\end{align}
then
\begin{equation}
\label{eqn:expected-gradient-points-in-right-direction}
 \left\langle
 \nabla_{\mathcal M_1}f_{P_N}(x_t),
 \mathbb E[g_t\mid\mathcal F_t,E_{t+1}]
 \right\rangle
 \geq
 c_1\|\nabla_{\mathcal M_1}f_{P_N}(x_t)\|^2.
\end{equation}
\end{lemma}

\paragraph{Uniform lower bound on
$\mathbb P(E_{t+1}\mid\mathcal F_t)$.} Lemmas~\ref{lem:metric-equivalence-gives-second-moment-bounds} and~\ref{lem:tube-containment-to-gradient-estimator-control} both depend on $\mathbb{P}(E_{t+1} \mid \mathcal{F}_{t})$, which depends on the specific sampling region $\widetilde{S}_{t+1}$. We now provide sample rates for \textit{all} potential admissible sampling regions. Define the following family:
for $\alpha\in(0,1]$ and $\epsilon>0$, let
\begin{equation}\label{eqn:sampling-family}
\mathfrak{S}_{\alpha,\epsilon}
:=
\left\{
S_{\mathrm{reg}}\cap
T_{\alpha d_{\mathcal M_1}(x,y)
+\rho_{S_{\mathrm{reg}}}(\epsilon)}
\bigl(\gamma^{\mathcal M_1}_{x\to y}\bigr)
:\ x,y\in S_{\mathrm{reg}},\ x\ne y
\right\}.
\end{equation}
For $S\in\mathfrak{S}_{\alpha_t,\epsilon_t}$, let $\mu_S$ denote
normalized Riemannian volume on $S$, and define
\begin{equation}\label{eqn:uniform-lower-bound-for-family}
q^{(t)}(s)
:=
\inf_{S\in\mathfrak{S}_{\alpha_t,\epsilon_t}}
\inf_{v\in S}
\mu_S\bigl(B_{\mathcal M_1}(v,s)\bigr),
\qquad \text{for }s>0.
\end{equation}
Lemma~\ref{lem:sampling-mass-lower-bound} yields a uniform lower bound at the sampling scales below, yielding the following sample complexities via covering number bounds; proof of both
are deferred to Appendix~\ref{sec:proof-of-sampling-to-bound-probability-of-bad-event}.

\begin{lemma}\label{lem:sampling-to-bound-probability-of-bad-event}
Assume Assumptions~\ref{ass:embedded-submanifold} and~\ref{ass:region-of-interest}
hold. Fix $1\le t\le T$. Let $0<\alpha_t\le1$,
$0<\lambda_t\le1/2$, and let $\epsilon_t,\delta_t>0$ satisfy
\[
\epsilon_t\le s_{\mathrm{loc},S_{\mathrm{reg}}},
\quad \text{and} \quad
\epsilon_t\le
\frac{\alpha_tD_{\mathrm{reg}}}{C_{\mathrm{loc},S_{\mathrm{reg}}}}.
\]
Given Assumption~\ref{ass:embedded-submanifold}\textup{(a)}, suppose the first four scale conditions of
Theorem~\ref{thm:deterministic-metric-equivalence-for-ISOMAP}
hold with $S=S_{\mathrm{reg}}$,
$\epsilon_{\min}=\epsilon_{\max}=\epsilon_t$, and
$\lambda_1=\lambda_2=\lambda_t$.
For Assumption~\ref{ass:embedded-submanifold}\textup{(b)}, instead suppose
\[
\delta_t\leq\frac{\lambda_t\epsilon_t}{4},
\qquad
\rho_{S_{\mathrm{reg}}}(\epsilon_t)\leq r_{\mathrm{met}},
\quad \text{and} \quad
K_{\mathrm{met}}\rho_{S_{\mathrm{reg}}}(\epsilon_t)^2
\leq2\lambda_t-\lambda_t^2.
\]
Set 
$\widetilde r_{\max}^{(t)}
:=\alpha_tD_{\mathrm{reg}}
+\rho_{S_{\mathrm{reg}}}(\epsilon_t),$
and, for $s>0$, define
\begin{equation}\label{eqn:covering-number-family}
\mathcal N_t(s):=
\left(
3+\frac{4L_{\widetilde r_{\max}^{(t)}}^{\widetilde K_{\mathcal M_1}}
(D_{\mathrm{reg}}+2\widetilde r_{\max}^{(t)})}{s}
\right)
\left(
1+\frac{8L_{\widetilde r_{\max}^{(t)}}^{\widetilde K_{\mathcal M_1}}
\widetilde r_{\max}^{(t)}}{s}
\right)^{d-1}.
\end{equation}
Let $\zeta_t\in(0,1)$. There exists $c_{\mathrm{mass}}(\ROI,\ROIdiam,\kmagInt) > 0$ such that following holds. Conditioned on
$\mathcal H_t:=\mathcal F_{t-1}\vee\sigma(I_t)$, if
\begin{equation}\label{eqn:explicit-sampling-budget}
n_t\ge
\left\lceil
\frac{2^{d}(D_{\mathrm{reg}}+\widetilde r_{\max}^{(t)})
(\widetilde r_{\max}^{(t)})^{d-1}}
{c_{\mathrm{mass}}(\ROI,\ROIdiam,\kmagInt)\delta_{t}^d}
\left(\log\mathcal N_t(\delta_{t}/2)+\log\frac1{\zeta_t}\right)
\right\rceil,
\end{equation}
samples are drawn with the uniform sampler on $\widetilde S_t$, then
\[
\mathbb P(E_t\mid\mathcal F_{t-1})\ge1-\zeta_t,
\quad \text{and} \quad
\mathbb P(E_t^c\mid\mathcal F_{t-1})\le\zeta_t.
\]
\end{lemma}

\subsection{Convergence Guarantee}\label{sec:main-result-theorem}
To emphasize the convergence rate and the hierarchy of algorithmic scales, we state the specialization of Theorem~\ref{thm:main-conv-thm} obtained by fixing the tunable constants as
\begin{equation*}
    c_{1} = \goodProb = \frac{1}{2}, \qquad \beta = \frac{8}{\SCf}, \quad \text{and} \quad \gamma = \frac{576 \smoothf}{\SCf^{2}}.
\end{equation*}
The full result and its proof, and its specialization to Theorem~\ref{thm:main-conv-thm-body}, are given in Appendix~\ref{sec:proof-of-main-conv-thm}.

\begin{theorem}
\label{thm:main-conv-thm-body} 
Suppose Assumptions~\ref{ass:embedded-submanifold} and~\ref{ass:region-of-interest} hold. Initialize
$x_0\in S_{\mathrm{reg}}$, and fix $T\geq 1$.
Let $m_f$ and $L_f$ be the strong-geodesic-convexity and
smoothness constants from Lemma~\ref{lem:functional-SC-smooth-parameters}. 
Set $\gamma
    :=
    576\smoothf/(\SCf)^2,$ 
and define
\begin{equation}
    \eta_t
    :=
    \frac{8}{m_f(\gamma+t)},
    \qquad
    0\leq t\leq T-1, 
    \quad \text{and} \quad \alpha_{t}
    :=
    \frac{8}{m_f(\gamma+t)^2},
    \qquad
    1\leq t\leq T.
    \label{eq:simplified-stepsize}
\end{equation}

There exist uniform constants
    $0<C_\lambda\leq 1,
    C_\epsilon>0,
    C_{\mathrm{conv}}>0,$
such that the following holds: For all $1\leq t \leq T$, let
\begin{equation}
    \lambda_t
    :=
    C_\lambda\alpha_t^2,
    \qquad
    \epsilon_t
    :=
    C_\epsilon\alpha_t,
    \quad \text{and} \quad
    \delta_t
    :=
    \frac{\lambda_t\epsilon_t}{4}.
    \label{eq:simplified-geometric-schedules}
\end{equation}

Run Algorithm~\ref{algo:graph-algorithm} for $T$ iterations, where Algorithm~\ref{algo:tube-sample-oracle} draws $n_{t}$ conditionally i.i.d. samples from normalized Riemannian volume on $\widetilde{S}_{t}$,
defined by \eqref{eqn:inflated-tube}. Let $c_{\mathrm{mass}}(\ROI,\ROIdiam,\kmagInt) > 0$ be the
constant in Lemma~\ref{lem:sampling-mass-lower-bound}. If
\begin{equation}
    n_t
    \geq
    \left\lceil
\frac{2^{d}(D_{\mathrm{reg}}+\widetilde r_{\max}^{(t)})
(\widetilde r_{\max}^{(t)})^{d-1}}
{c_{\mathrm{mass}}(\ROI,\ROIdiam,\kmagInt)\delta_{t}^d}
\left(\log\mathcal N_t(\delta_{t}/2)+\log\frac1{\zeta_T}\right)
\right\rceil, \quad 
    \zeta_T
    :=
    \min\left\{
        \frac{8}{m_fT(\gamma+T)},
        \frac12
    \right\},
    \label{eq:simplified-sampling-rate}
\end{equation}
where 
$\widetilde r_{\max}^{(t)}
:=\alpha_tD_{\mathrm{reg}}
+\rho_{S_{\mathrm{reg}}}(\epsilon_t),$ and $\mathcal{N}_{t}(\cdot)$ is defined by \eqref{eqn:covering-number-family}. Then the final iterate satisfies
\begin{equation}
    \mathbb E_{\mathrm{alg}}
    \left[
        f_{P_N}(x_T)
        -
        f_{P_N}(\widehat x_N)
    \right]
    \leq
    \frac{C_{\mathrm{conv}}}{\gamma+T}.
    \label{eq:simplified-convergence-rate}
\end{equation}
where
$C_{\mathrm{conv}}
=C_{\mathrm{conv}}(\SCf,\smoothf,\ROIdiam,\kminInt,\kmaxInt)$ and $\widehat{x}_{N}$ is the minimizer of \eqref{eqn:objective-functional-informal} in $\ROI$.
\end{theorem}
We emphasize Theorem~\ref{thm:main-conv-thm-body} suppresses only numerical
constants and parameter choices obscuring the principal rate. 
The graph-scale constants depend on
the intrinsic geometry and localization parameters together
with the constants of the selected metric-recovery branch; explicitly depending on the ambient geometry when Assumption~\ref{ass:embedded-submanifold}\textup{(a)} holds.
A few comments on this specialized theorem are in order.

First, the bound \eqref{eq:simplified-convergence-rate} has $\mathcal{O}(T^{-1})$ outer iteration dependence. This matches the minimax-optimal rate for stochastic first-order methods on strongly convex Euclidean objectives with unbiased gradient estimators \cite{lower-bound-SGD}. As a special case, this recovers the rate of geodesic interpolation. We emphasize that this outer rate is obtained without direct access to the intrinsic exponential map, logarithmic map, or tangent structure of $\mathcal{M}_{1}$. 

Second, the schedules separate the optimization and geometric approximation scales. At the update from \(x_t\) to \(x_{t+1}\),
\begin{equation}\label{eqn:stepsize-radius-scalings}
    \eta_t=\mathcal{O}((\gamma+t)^{-1}),
    \quad \text{and} \quad
    \alpha_{t+1}=\mathcal{O}((\gamma+t)^{-2}).
\end{equation}
The stepsize \(\eta_t\) is the standard strongly convex SGD scaling yielding an \(\mathcal{O}(T^{-1})\) objective-gap rate. By Lemma~\ref{lem:discrepancy-from-geodesic-interpolation}, the discrepancy between graph and exact geodesic interpolation satisfies
\[
d_{\mathcal{M}_1}(x_{t+1},\widetilde{x}_{t+1})
\lesssim
(\lambda_{t+1}\eta_t+\alpha_{t+1})
d_{\mathcal{M}_1}(x_t,y_{t+1}) + \epsilon_{t+1}.
\]
Because the implicit gradient estimator divides the update displacement by \(\eta_t\), the terms corresponding to tube radius and rounding error contribute respectively
$\alpha_{t+1}/\eta_{t}
=
\mathcal{O}((\gamma+t)^{-1})$ and $\epsilon_{t+1}/\eta_{t} = \mathcal{O}((\gamma+t)^{-1})$.
These error terms are small enough to preserve the required conditional
alignment when the objective gap is above a threshold used in the convergence argument of
Appendix~\ref{sec:proof-of-sgd-like-rate-on-good-event}. Below said threshold, a separate stability bound controls
the next objective gap. Together, these bounds recover the
$O(T^{-1})$ convergence rate.

Third, given the tube radius scaling, the following scalings directly follow:
\begin{equation} \label{eqn:graph-scalings}
    \lambda_{t} = \mathcal{O}((\gamma+ t)^{-4}), \qquad \epsilon_{t} = \mathcal{O}((\gamma+t)^{-2}), \quad \text{and} \quad \delta_{t} = \mathcal{O}((\gamma +t)^{-6}).
\end{equation}
The scaling $\lambda_{t} = \mathcal{O}(\alpha_{t}^{2})$ ensures the
containment condition of Lemma~\ref{lem:conditions-on-metric-equivalence-for-tube-containment} holds. Together with the connectivity and rounding scales, this gives the estimator-discrepancy bound used to drive objective gap improvements in Appendix~\ref{sec:proof-of-sgd-like-rate-on-good-event}.
The schedules of the graph construction parameters $\epsilon_t$
and $\delta_t$ follow from
Theorem~\ref{thm:deterministic-metric-equivalence-for-ISOMAP}
or Proposition~\ref{prop:abstract-metric-isomap}, together with
Lemma~\ref{lem:local-intrinsic-extrinsic-compatibility}.
Intuitively, smaller radii result in smaller geodesic tubes, which require finer graph representations.

Lastly, with the scalings of \eqref{eqn:stepsize-radius-scalings} and \eqref{eqn:graph-scalings}, and Lemma~\ref{lem:local-intrinsic-extrinsic-compatibility} in mind, we have that
\begin{equation*}
    \frac{2^{d}(D_{\mathrm{reg}}+\widetilde r_{\max}^{(t)})
(\widetilde r_{\max}^{(t)})^{d-1}}
{c_{\mathrm{mass}}(\ROI,\ROIdiam,\kmagInt)\delta_{t}^d}\lesssim \frac{(\widetilde{r}^{(t)}_{\max})^{d-1}}{\delta_{t}^{d}} \lesssim \frac{(\alpha_{t} + \epsilon_{t+1}\mm{\epsilon_t?})^{d-1}}{\delta_{t}^{d}} = \mathcal{O}((\gamma +t)^{4d+2}),
\end{equation*}
with $\mathcal{N}_{t}(\delta_{t}/2)$ having the same complexity by an analogous argument. Thus, for the uniform localized sampler, this objective gap guarantee is
obtained with sufficient per-iteration and cumulative sampling costs
of $\widetilde{\mathcal O}(t^{4d+2})$ and $\widetilde{\mathcal O}(T^{4d+3})$, suppressing the poly-log terms and $\gamma$.

\section{Numerical Experiments}
\label{sec:numerical-experiments}
We present numerics of Algorithm~\ref{algo:graph-algorithm} with graph refinement Algorithm~\ref{algo:ball-refinement} in various settings. Sections~\ref{sec:numerics-sphere} and~\ref{sec:numerics-half-moon} provide two toy examples where $\mathcal{M}_{1}$ is embedded in $\mathbb{R}^{d}$ and $\mathcal{P}_{2,\operatorname{ac}}(\mathbb{R}^{d})$ respectively. In Section~\ref{sec:k-means-MNIST} we use the proposed algorithm as a primitive for \textit{Lloyd's algorithm} for $k$-means clustering on the MNIST dataset. The code for these experiments is available at 
\begin{center}
\textcolor{myPurple}{https://github.com/mcphersonianoliver/graph-frechet}.
\end{center}

\paragraph{Implementation Details.} All experiments were run on a desktop computer running Windows~11 with a
six-core Intel Core i5-11400F processor at 2.60\,GHz (12 logical processors)
and 16\,GB of RAM. The implementation was written in Python~3.14.2.
We used NumPy for numerical computation and pseudorandom sampling,
NetworkX for weighted graph construction and shortest-path computations,
and SciPy together with Python Optimal Transport (POT) for the Wasserstein
experiments. Figures were generated with Matplotlib. 

Every stochastic experiment used an explicit NumPy random generator with a
recorded seed. Exact sphere geometry and half-moon parameter information were
used only for sampling and evaluation. Graph construction, shortest-path
selection, interpolation, and Fréchet updates used only the ambient-distance
and local-sampling interfaces described above.

\paragraph{Riemannian Structure of $(\mathbb{S}^{2},d_{\mathbb{S}_{2}})$.} 
We regard the unit sphere
\[
    \mathbb{S}^{2}
    := \{x\in\mathbb{R}^{3}:\|x\|_2=1\}
\]
as an embedded Riemannian submanifold of $\mathbb{R}^{3}$, metric
induced by the Euclidean inner product. Thus
\[
    T_x\mathbb{S}^{2}
    = \{v\in\mathbb{R}^{3}:\langle x,v\rangle=0\},
    \quad \text{and} \quad
    d_{\mathbb{S}^{2}}(x,y)
    = \arccos\!\left(\langle x,y\rangle\right).
\]
For $x,y\in\mathbb{S}^{2}$ with
$\theta=d_{\mathbb{S}^{2}}(x,y)\in(0,\pi)$, the minimizing geodesic is
\[
    \gamma_{x\to y}(s)
    =
    \frac{\sin((1-s)\theta)}{\sin\theta}x
    +
    \frac{\sin(s\theta)}{\sin\theta}y,
    \quad \text{for } s\in[0,1].
\]

\paragraph{(Formal) Riemannian Structure of
$(\mathcal{P}_{2,\operatorname{ac}}(\mathbb{R}^{d}),W_{2})$.} 
For $\mu\in\mathcal{P}_{2,\operatorname{ac}}(\mathbb{R}^{d})$, its tangent space is identified with
\[
    T_\mu\mathcal{P}_2(\mathbb{R}^{d})
    =
    \overline{\{\nabla\varphi:
        \varphi\in C_c^\infty(\mathbb{R}^{d})\}}^{\,L^2(\mu)},
\]
equipped with
\[
    \langle v,w\rangle_\mu
    := \int_{\mathbb{R}^{d}}
        \langle v(x),w(x)\rangle\,d\mu(x).
\]
If $\nu\in\mathcal{P}_2(\mathbb{R}^{d})$ and
$T_{\mu \rightarrow \nu}$ denotes the optimal transport map from $\mu$ to $\nu$, then
\[
    W_2^2(\mu,\nu)
    =
    \int_{\mathbb{R}^{d}}
        \|T_{\mu \rightarrow \nu}(x)-x\|_2^2\,d\mu(x),
    \quad \text{and} \quad
    \log_\mu(\nu)=T_{\mu \rightarrow \nu}-\operatorname{Id}.
\]
The constant-speed Wasserstein geodesic from $\mu$ to $\nu$ is the
displacement interpolation
\[
    \mu_s
    =
    \big((1-s)\operatorname{Id}+sT_{\mu \rightarrow \nu}\big)_{\#}\mu,
    \quad \text{for } s\in[0,1].
\]
Accordingly, in the Wasserstein experiments the ambient edge weight between
two sampled measures is their $W_2$ distance. The particular numerical
representation of the measures and computation of these distances are
specified in the corresponding experimental sections.

\subsection{Controlled graph interpolation on the sphere}
\label{sec:numerics-sphere}
\begin{figure}[t]
    \centering
    \begin{subfigure}[t]{0.49\textwidth}
        \centering
        \includegraphics[width=\linewidth]{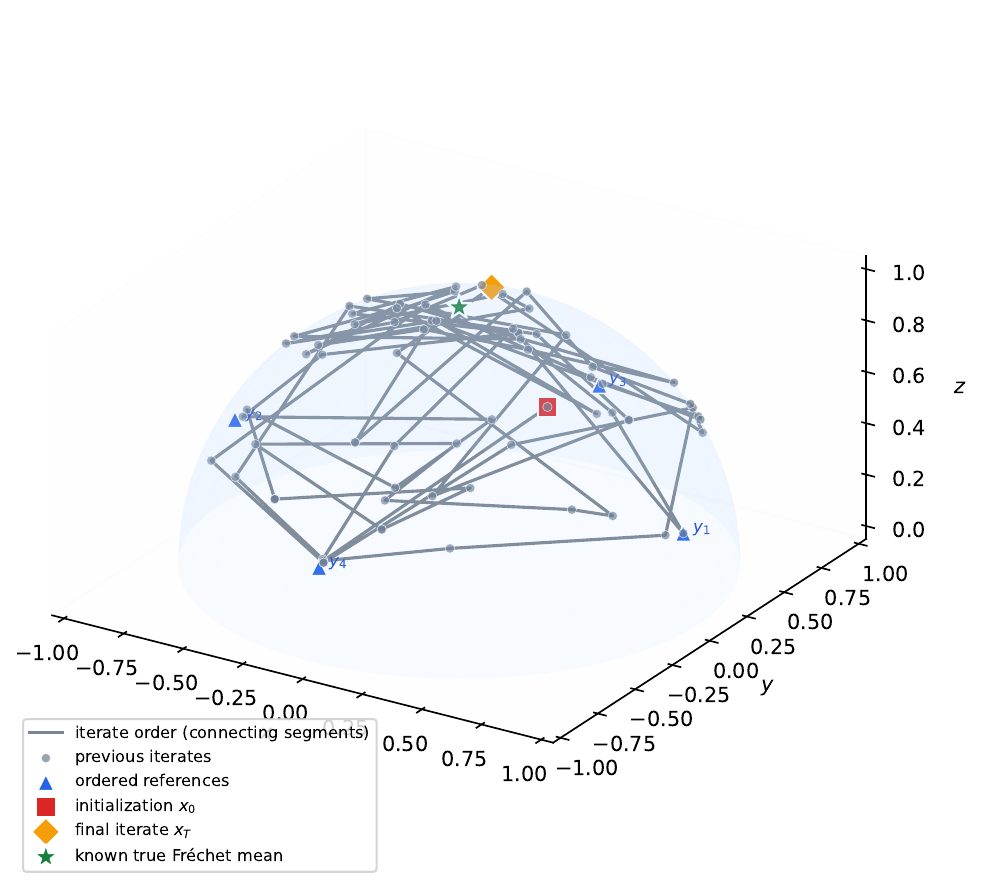}
        \caption{Optimization trajectory.}
        \label{fig:sphere-full-run}
    \end{subfigure}
    \hfill
    \begin{subfigure}[t]{0.49\textwidth}
        \centering
        \includegraphics[width=\linewidth]{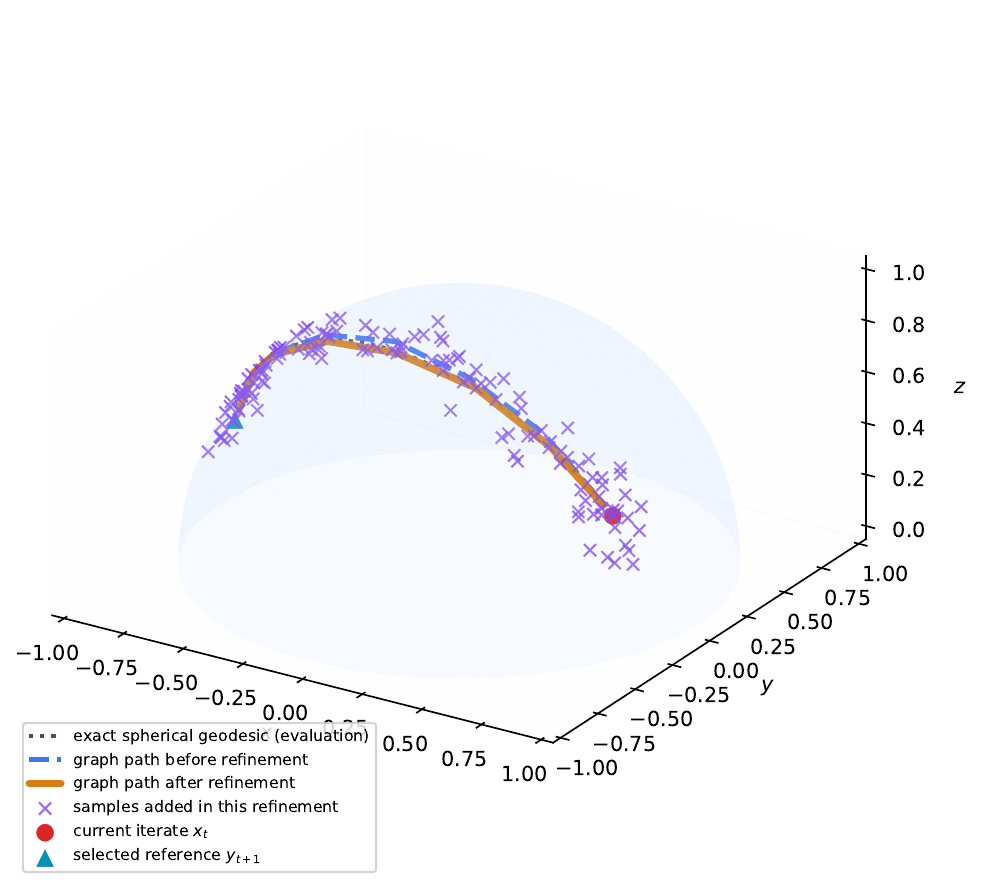}
        \caption{Representative local refinement step.}
        \label{fig:sphere-refinement-step}
    \end{subfigure}

    \caption{
    \textbf{a)} Representative 
    optimization trajectory, together with the initialization, four reference points, final iterate, and known Fr\'echet mean. 
    \textbf{b)} Representative update for which local refinement changes the selected graph shortest path. 
    }
    \label{fig:sphere-qualitative}
\end{figure}

Consider the unit sphere
    $ \mathbb{S}^{2}\subseteq \mathbb{R}^{3}$
equipped with Riemannian metric induced by $\mathbb{R}^{3}$, providing a controlled setting where the intrinsic geometry is analytically available.

Let $\theta=7\pi/18$. We consider the four reference points $\{z_{i}\}_{i=1}^{4}$ defined respectively by
\begin{equation}
\label{eq:sphere-reference-points}
\begin{aligned}
    (\sin\theta,0,\cos\theta), \qquad 
    (-\sin\theta,0,\cos\theta), \qquad
    (0,\sin\theta,\cos\theta), 
    \qquad \text{and} \qquad (0,-\sin\theta,\cos\theta).
\end{aligned}
\end{equation}
Writing $e_3=(0,0,1)$ and $v_0=(0.55,-0.25,0)\in T_{e_3}\mathbb S^2$, we initialize at
    $x_0
    =
    \exp^{\mathbb S^2}_{e_3}(v_0)$
By symmetry of the reference configuration, together with localization to the corresponding strongly convex spherical region, the Fr\'echet mean is the north pole,
\begin{equation}
\label{eq:sphere-known-mean}
    \widehat{x}_{4}
    =
    (0,0,1).
\end{equation}

We compare three interpolation procedures. The first is exact geodesic interpolation, Algorithm~\ref{alg:geodesic-interpolation}, which is implementable as geometry of $\mathbb{S}^{2}$ is known. This serves as an oracle optimization baseline. The second is our graph-algorithm with access to a tube sampler, i.e., Algorithm~\ref{algo:graph-algorithm} with the graph refinement with tube sampler Algorithm~\ref{algo:tube-sample-oracle}. We note that the sampling procedure uses intrinsic information, i.e., $\gamma_{x_{t} \rightarrow y_{t+1}}^{\mathbb{S}^{2}}$, that would be unavailable on an unknown submanifold. The numerical tube construction retains previously sampled points; its purpose here is to provide an idealized sampler-oracle baseline to compare with the following procedure.
The last procedure is the implementable persistent ball-sampler construction, namely Algorithm~\ref{algo:graph-algorithm} with graph refinement with ball sampler Algorithm~\ref{algo:ball-refinement}, described in Section~\ref{sec:practical-graph-algo-description}. The overarching goals of the experiment are therefore twofold: $(i)$ to determine whether graph-based interpolation accurately approximates Algorithm~\ref{alg:geodesic-interpolation}, and $(ii)$ to determine whether the implementable ball-sampler refinement reproduces the behavior of the idealized oracle-tube sampler.

\paragraph{Implementation Details.}
We perform $T=100$ outer iterations and report quantitative results across $N_{\operatorname{alg}}=5$ paired runs. Within each run, the exact, oracle-tube, and ball-based procedures use the same initialization, reference sequence, and stepsize schedule. The initial graph is initialized with $150$ points sampled uniformly with respect to surface area in the spherical cap centered at $(0,0,1)$, 
together with $x_0$ and the four reference points.

Graph edges are weighted by $d_{\mathbb R^3}$. Schedules for graph-connectivity, interpolation stepsize, and sampling are respectively given by 
\begin{equation}
    \epsilon_t
    =
    \frac{1}{\sqrt{1+t/10}}, \qquad \eta_t
    =
    \min\left\{\frac{1}{2}, \frac{4}{\sqrt{t+8}} \right\}, 
    \qquad \text{and} \qquad \alpha_t = \frac{0.25}{\sqrt{1+t/10}}.
\end{equation}
Let $L^{G}_{t}$ denote the length of $\gamma_{x_{t} \rightarrow y_{t+1}}^{G_{t}}$ pre-refinement. Then, the local sampling radius is
    $r_t
    =
    \alpha_t L^{G}_{t}.$
We apply a less aggressive schedule than suggested by Section~\ref{sec:idealized-convegence-analysis} as we use ball-samplers instead of tube-samplers. 
This schedule shrinks the sampling region at a slower rate, allowing for the identification of the ideal tube sample region. We draw $15$ samples per selected refinement center and use at most $30$ centers along the guide path.

\paragraph{Discussion.}
\begin{figure}[t]
    \centering
    \includegraphics[width=0.6\linewidth]{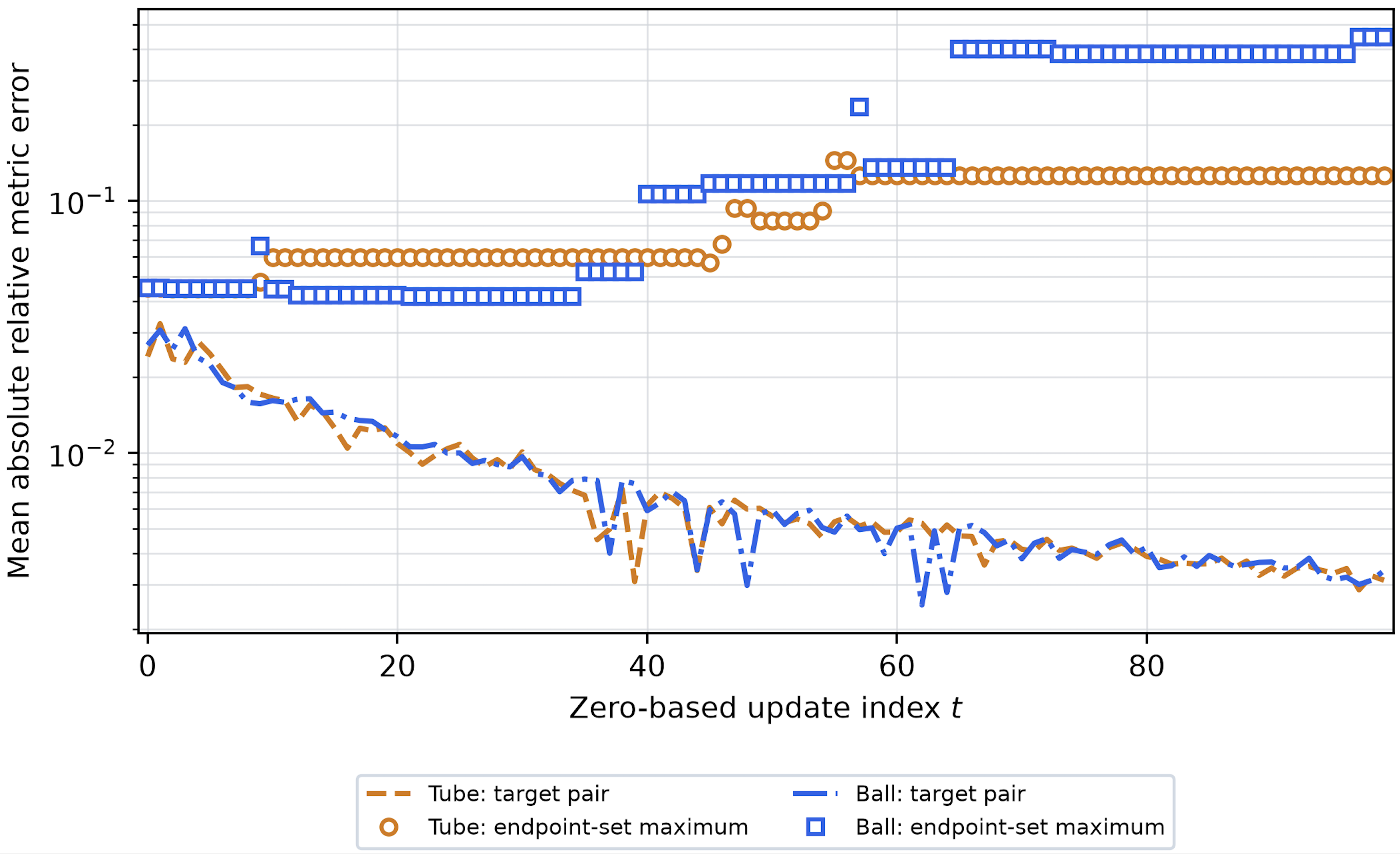}
    \caption{
    Plot of \(e_{\mathrm{target}}(t)\) and \(e_{\max}(t)\), averaged over \(N_{\operatorname{alg}}=5\) runs, with \(e_{\max}\) evaluated on the original/coarse endpoints.
    }
    \label{fig:sphere-quantitative}
\end{figure}
Figure~\ref{fig:sphere-qualitative} illustrates qualitatively what occurs during an optimization run. Figure~\ref{fig:sphere-full-run} shows a representative complete optimization trajectory, together with the reference configuration and known solution, while Figure~\ref{fig:sphere-refinement-step} isolates the mechanism at a single update. The pre-refinement graph geodesic provides the centers around which new samples are generated, and the post-refinement graph geodesic more closely approximates the true underlying geodesics, as desired. We move to quantify these observations.

Let $m\in\{\mathrm{tube},\mathrm{ball}\}$ index the two graph procedures, $j\in\{1,\ldots,N_{\operatorname{alg}}\}$ the paired runs, and $G_{t+1,j}^{m}$ represent the graph used for the $t+1$th graph interpolation. To test the intuition that the refined graph geodesic has less metric distortion than the graph approximation across the region of interest, we define
\begin{equation}
\label{eq:sphere-target-error}
    e_{\mathrm{target},j}^{(m)}(t)
    :=
    \left|
    \frac{
        d_{G_{t+1,j}^{(m)}}
        \left(
            x_{t,j}^{(m)},y_{t+1,j}
        \right)
    }{
        d_{\mathbb S^2}
        \left(
            x_{t,j}^{(m)},y_{t+1,j}
        \right)
    }
    -1
    \right|.
\end{equation}
This evaluates metric fidelity for the geodesic of interest. We compare this with the worst-case graph distortion for pairs of vertices of the initial graph:
\begin{equation}
\label{eq:sphere-max-error}
    e_{\max,j}^{(m)}(t)
    :=
    \max_{\{u,v\}\in\mathcal V(G_{0}^{(m)})}
    \left|
    \frac{
        d_{G_{t+1,j}^{(m)}}(u,v)
    }{
        d_{\mathbb S^2}(u,v)
    }
    -1
    \right|.
\end{equation}
To account for the randomness due to the sequence of target references and sampling, we average these two quantities with respect to the number of algorithm runs:
\begin{equation}\label{eq:average-quantities}
    \overline e_{\mathrm{target}}^{(m)}(t)
    :=
    \frac{1}{N_{\operatorname{alg}}}
    \sum_{j=1}^{N_{\operatorname{alg}}}
    e_{\mathrm{target},j}^{(m)}(t), \quad \text{and} \quad \overline e_{\mathrm{max}}^{(m)}(t) := \frac{1}{N_{\operatorname{alg}}}
    \sum_{j=1}^{N_{\operatorname{alg}}}
    e_{\mathrm{max}}^{(m)}(t)
\end{equation}

Figure~\ref{fig:sphere-quantitative} plots these quantities, reflecting two key observations: $i)$ the metric equivalence distortion for the target geodesic $\gamma_{x_{t} \rightarrow y_{t+1}}^{\mathcal{M}_{1}}$ is much lower than the global worst-case distortion, regardless of the sampling mechanism; $ii)$ the target metric equivalence behaves similarly for both sampling methods and decreases over the optimization run. The separation between the target metric equivalence and the global worst-case reflects the intended adaptive construction. As $\epsilon_{t}$ decreases, repeatedly traversed regions receive more samples, resulting in finer resolution, while untraversed regions do not improve in accuracy.

\subsection{Graph interpolation on a submanifold of Wasserstein space}
\begin{figure}[t]
    \centering
    \begin{subfigure}[t]{0.31\textwidth}
        \centering
        \includegraphics{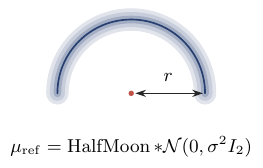} 
        \caption{Reference measure.}
        \label{fig:halfmoon-reference}
    \end{subfigure}
    \begin{subfigure}[t]{0.31\textwidth}
        \centering
        \includegraphics{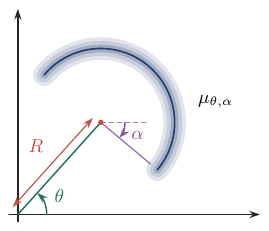}
        \caption{Representative $\mu_{\theta,\alpha} \in \mathcal{M}_{1}$.}
        \label{fig:halfmoon-single-point}
    \end{subfigure}
    \begin{subfigure}[t]{0.31\textwidth}
        \centering
        \includegraphics{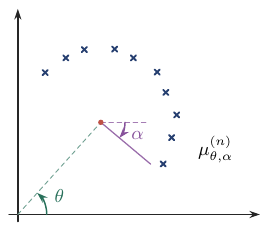}
        \caption{Discretized version of $\mu_{\theta,\alpha}$.}
        \label{fig:halfmoon-discretization}
    \end{subfigure}
    \caption{Construction of the half-moon family.}
    \label{fig:half-moon-manifold-parametrization}
\end{figure}

We next consider the infinite-dimensional ambient space
    $\mathcal{M}_{2}
    =
    \bigl(\mathcal{P}_{2,\mathrm{ac}}(\mathbb{R}^{2}),W_{2}\bigr).$
Following \cite{HMST25}, let $\mathcal{M}_{1}$ be defined by a reference measure acted upon by two types of rotations:
\[
    \mathcal{M}_{1} = \left\{
        \mu_{\theta,\alpha}
        :
        (\theta,\alpha)
        \in
        \mathbb{R} \times \mathbb{R}
    \right\}.
\]
The reference measure is a half-moon convolved with a Gaussian, with the local center a fixed radius from the origin. The parameters $\theta$ and $\alpha$ correspond to two types of rotations acting on $\mathbb{R}^{2}$, the former about the origin while the latter about the local center. We visualize this construction in \Cref{fig:half-moon-manifold-parametrization}.

In order to have a unique Fr\'echet mean, we restrict our attention to
\[
    \ROI := \{\mu_{\theta,\alpha} \colon (\theta,\alpha) \in [0,\pi] \times [-\pi/2, \pi/2]\}.
\]

We consider four reference measures $\{z_i\}_{i=1}^{4}\subset\mathcal{M}_{1}$, with parameter coordinates
\begin{equation}
\label{eq:half-moon-reference-points}
\begin{aligned}
    \left(\frac{\pi}{4},-\frac{\pi}{4}\right), \qquad
    \left(\frac{3\pi}{4},-\frac{\pi}{4}\right), \qquad
    \left(\frac{\pi}{4},\frac{\pi}{4}\right), \quad \text{and} \quad
    \left(\frac{3\pi}{4},\frac{\pi}{4}\right).
\end{aligned}
\end{equation}
Motivated by the symmetry of the reference configuration
in parameter space, we take
\begin{equation}
\label{eq:half-moon-known-mean}
    \widehat{x}_{4}
    =
    \mu_{\pi/2,0}
\end{equation}
as a reference target for numerical comparisons.
\begin{figure}[t]
    \centering
    \includegraphics[scale=0.5]{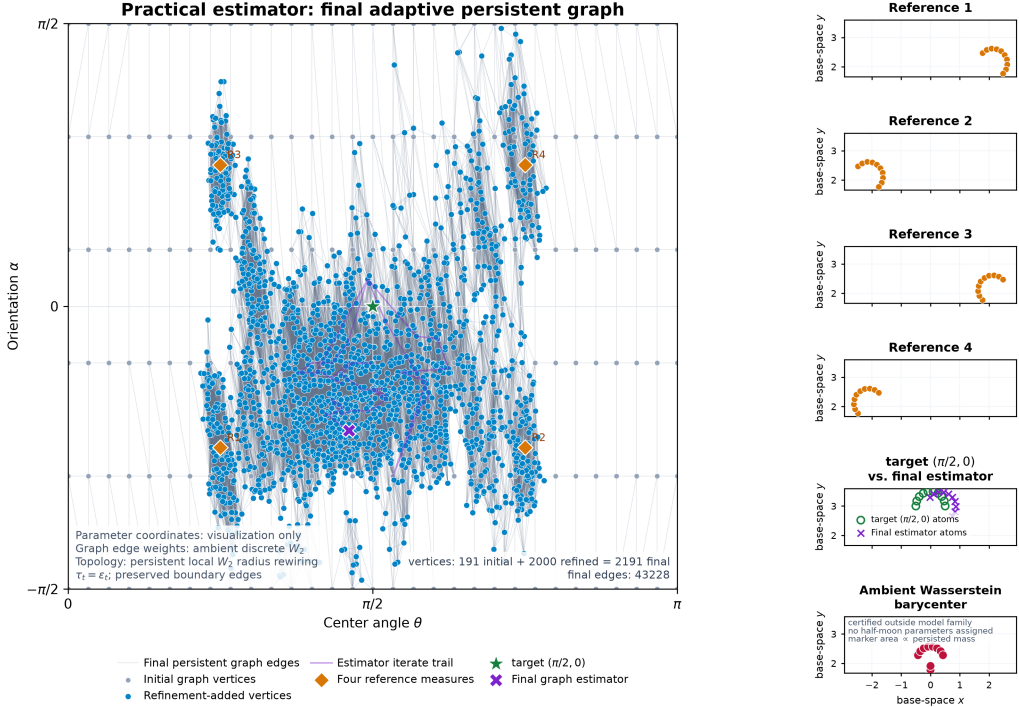}
    \caption{Final adaptive persistent graph in parameter
space, the four reference measures, the symmetry-based
reference target and final graph estimator, and the
ambient Wasserstein barycenter comparison.}
    \label{fig:half-moon-qualitative}
\end{figure}
\label{sec:numerics-half-moon}
\paragraph{Implementation Details.} We apply Algorithm~\ref{algo:graph-algorithm} with the local refinement mechanism of Algorithm~\ref{algo:ball-refinement} for $T = 100$ graph-interpolation updates for multiple optimization trajectories. The initial graph $G_{0}$ contains 191 vertices distributed across the parameter domain at roughly a $6:1$ ratio between $\theta$ and $\phi$. This ratio is chosen to roughly account for the anisotropy of the parameters $\theta$ and $\phi$ when realized as $W_{2}$ perturbations. At each update, all vertices are retained and at most 10 new half-moon samples are drawn around at most 5 refinement centers along a coarse graph geodesic. Vertices of the graph are measures in $\mathcal{M}_{1}$ and graph edges are weighted by the ambient distance $W_{2}$.

Because the intrinsic distance $d_{\mathcal{M}_{1}}$ is not available in
closed form, we separately construct a fixed high-resolution \(241 \times 41\) reference
graph on the half-moon family and use its shortest-path metric as a numerical
proxy for $d_{\mathcal{M}_{1}}$. This reference graph is frozen before
evaluation and is never exposed to the graph-based estimator. It plays the
same diagnostic role that the analytically known spherical distance plays in
Section~\ref{sec:numerics-sphere}.

\paragraph{Discussion.}
Figure~\ref{fig:half-moon-qualitative} illustrates the adaptive graph
construction in a nonlinear ambient space. In the parameter-space
visualization, newly sampled vertices concentrate around the regions
explored by the graph interpolation procedure, while other portions of
the domain retain the coarser initial representation, reflecting the intended local refinement mechanism.

The comparison with the symmetry-based reference target
$x_{\mathrm{sym}}=\mu_{\pi/2,0}$ illustrates the
finite-budget behavior of the graph interpolation
procedure.
Since each update
selects a graph vertex, the estimator remains within the discretized
half-moon family by construction. Quantitative
diagnostics of the graph geometry mirroring those of Figure~\ref{fig:sphere-quantitative}, evaluated using the frozen
high-resolution reference graph described above, are deferred to
Appendix~\ref{sec:additional-numerical-results} as they echo the same findings: our graph method \textit{locally} improves manifold representation.

The ambient Wasserstein barycenter comparison highlights a different
point. Although all four reference measures belong to the half-moon
family, the displayed ambient barycenter no longer remains in $\mathcal{M}_{1}$. 
This highlights how ambient averaging can leave the data manifold when the manifold is non-geodesically convex with respect to the ambient geometry. This approach is in contrast to prior Barycenter work that considers geodesically convex submanifolds of $\mathcal{P}_{2,\operatorname{ac}}(\mathbb{R}^{d})$, specifically the Bures-Wasserstein submanifold; see \cite{CMRS20}.
\begin{figure}[t]
    \centering
    \includegraphics[scale=0.55]{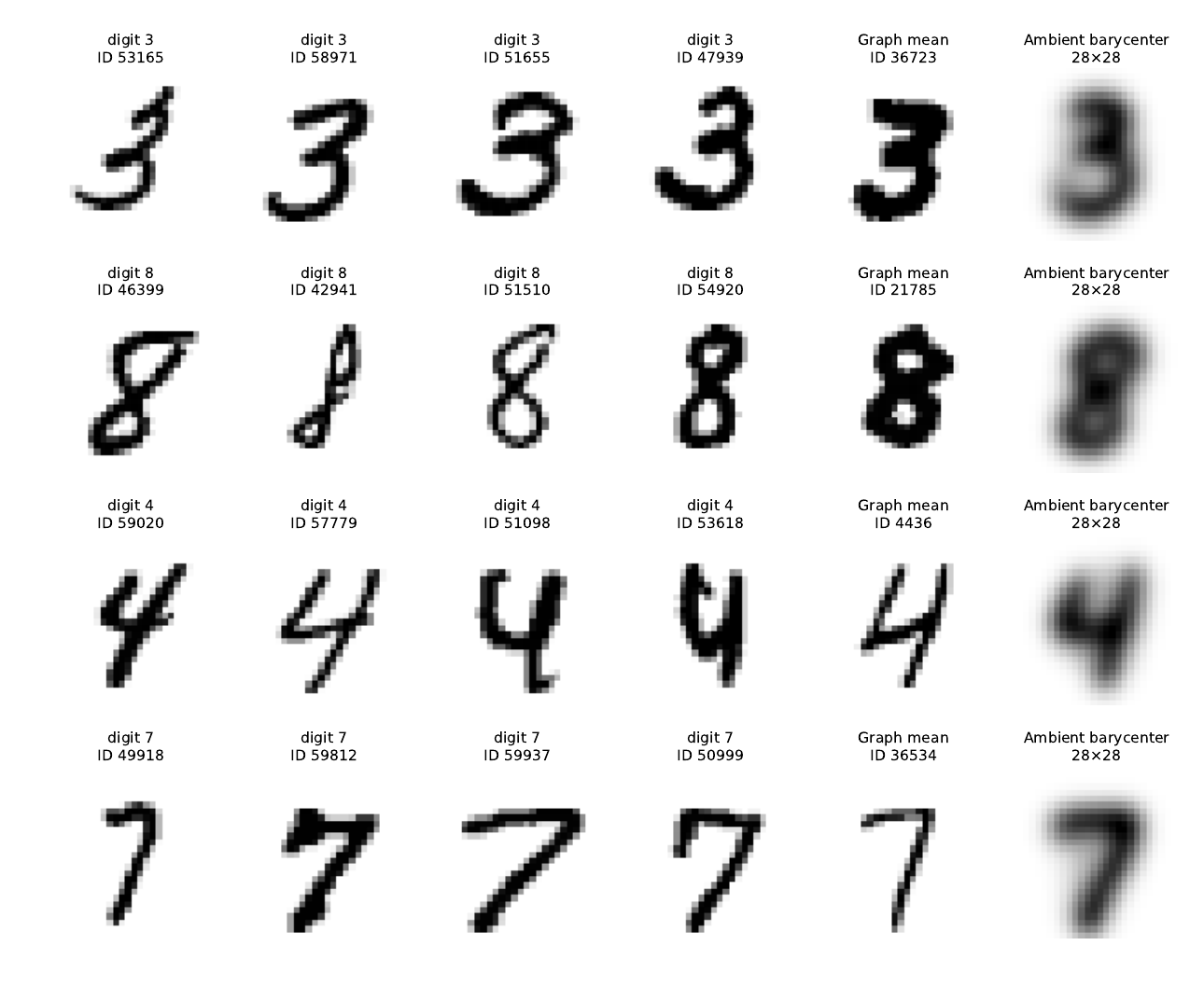}
    \caption{Single-digit averaging with $N=32$ observations per row.
    Four of the reference observations are shown, followed by the graph-based and ambient Fr\'echet representatives.
    }
    \label{fig:mnist-single-digit-means}
\end{figure}

\subsection{Graph-based averaging and clustering on MNIST}

\label{sec:k-means-MNIST}
Lastly, we consider the MNIST digits database \cite{lecun1998gradient}, interpreting each image as
a probability measure on the pixel space. We first examine representatives produced by
Fr\'echet-mean computations on collections from
one digit class, then from two classes. We then
use the same averaging primitive within Lloyd's algorithm for $k$-means
clustering. The first experiments illustrate the averaging operation
itself; the last highlights its use as a primitive.

\paragraph{Implementation Details.}
Distances for graph construction are computed by first compressing the $28 \times 28$ images to $14 \times 14$ by averaging and normalization, then computing the $W_{2}$ distance on image pairs. Graph representatives are displayed using their original $28\times28$ pixels.

To implement local sampling, we organize a finite reservoir
$\mathcal{R}$ of MNIST images into a cover tree~\cite{BKL06},
using the exact $W_2$ distance between their normalized
$14\times14$ representations. We emphasize that this is a cheaper numerical surrogate than $W_{2}$ on the original images; other cheap surrogates may be considered. A cover tree is a hierarchical
metric index whose levels correspond to geometrically decreasing
scales. These levels satisfy nesting, covering, and separation
conditions: points present at a coarse level remain present at
finer levels, each point has a nearby parent at the preceding
level, and distinct points at the same level are separated at
the corresponding scale. The construction requires only
pairwise metric evaluations, without access to intrinsic
coordinates or geodesics. The tree is constructed in advance
and reused across optimization runs.

Given a refinement center $x$ and radius $r$, a range query
identifies the reservoir points in the ambient metric ball,
\[
    \mathcal{R}_{r}(x)
    :=
    \{\nu\in\mathcal{R}:W_2(x,\nu)\leq r\}.
\]
The tree accelerates this query by pruning subtrees using the
triangle inequality. Specifically, if a node represents a
measure $v$ and $R_v$ bounds the distance from $v$ to every
measure in its subtree, then
    $W_2(x,v)>r+R_v$
implies that the entire subtree lies outside the query ball.

After excluding previously revealed or otherwise ineligible
reservoir points, we sample uniformly without replacement
from the remaining candidates, up to the requested sample
budget. Newly revealed samples are added to the persistent
$W_2$-weighted graph used for shortest-path interpolation.
The cover tree and this graph serve distinct purposes: the
former indexes ambient neighborhoods, while the latter
approximates intrinsic geometry. 
\begin{figure}[t]
    \centering
    \includegraphics[scale=0.58]{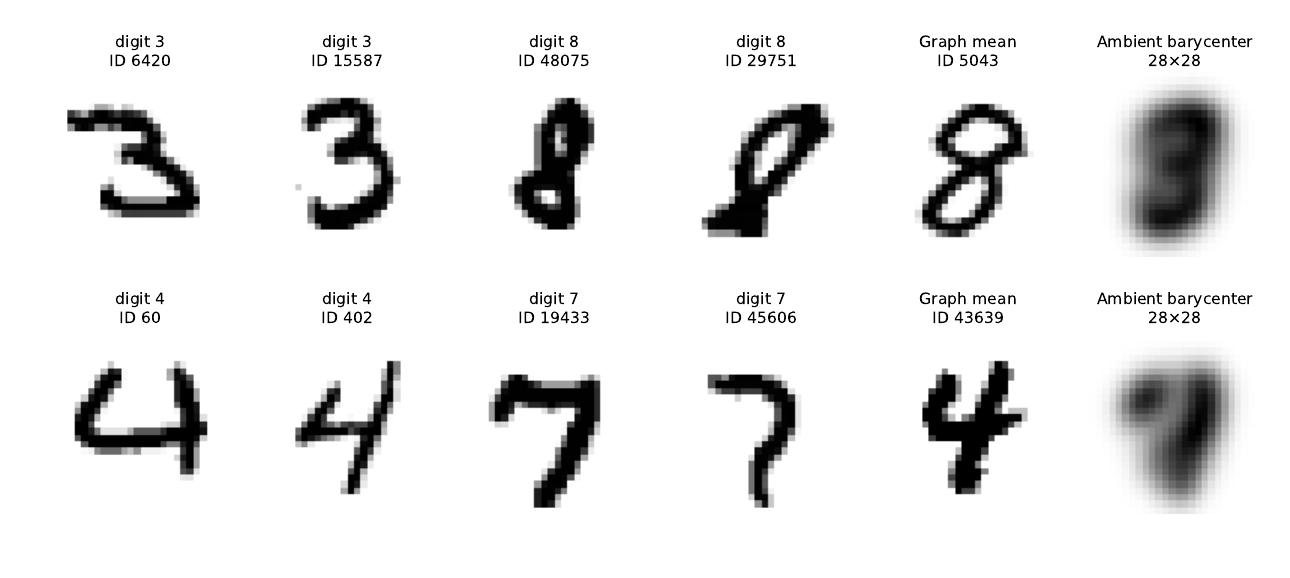}
    \caption{A single mean of each two-class collection: $4/7$ and $3/8$,
    with $16$ reference observations per digit. Two observations from
    each digit are displayed, outputs use all $32$ observations.
    }
    \label{fig:mnist-two-class-means}
\end{figure}

\paragraph{Graph-based averages of MNIST digits.}
For a fixed reference collection $\{z_i\}_{i=1}^{N}$, we apply
Algorithm~\ref{algo:graph-algorithm} with local refinement by
Algorithm~\ref{algo:ball-refinement} and compare its returned representative with
an ambient Wasserstein barycenter of the same collection. In these
averaging figures, the ambient barycenters were recomputed on the
normalized original $28\times28$ images for display and evaluation using entropic OT distances. 

We consider two subcollections of MNIST, namely the collection of all $3$s and $8$s and the collection of all $4$s and $7$s. The local sampler will be restricted to these subcollections, which play the role of an unknown submanifold of $\mathcal{P}_{2}(\mathbb{R}^{28 \times 28})$. We first consider averaging $N=32$ observations from a single digit class to check if the subclasses are geodesically separable. 
Figure~\ref{fig:mnist-single-digit-means} shows examples for each of the digits
$3$, $4$, $7$, and $8$. Each row displays four reference observations, the
graph-based and ambient representatives; both outputs
use the complete collection of $32$ observations. For each digit, the representative average is both out of the original reference sample and remains within the digit class, although the sampler has access to the other digits in the subcollection. We note that the ambient space representative lies outside of the sub-manifold, and loses image fidelity.

We next average collections containing an equal amount of observations from each of
two digit classes, the subcollections $4/7$ and $3/8$.
Figure~\ref{fig:mnist-two-class-means} shows two observations per digit
and the two resulting representatives. Each method computes
\emph{one mean of the entire $32$-observation collection}. For each collection, the computed graph average lies outside the reference collection. In particular, the average computed by the ambient method not only blurs the signal but also produces objects that do not correspond to any numerical digit. This is especially clear for the average of the $4$s and $7$s subcollection.
Additional MNIST averaging examples are provided in Appendix~\ref{sec:additional-numerical-results}. We note that all of the produced centers for both the single-digit and mixed-digit cohorts are out of the reference sample.

\begin{algorithm}[t]
\caption{Graph-based Lloyd's algorithm}
\label{algo:graph-lloyd}
\begin{algorithmic}[1]
\Require Observations $\{z_i\}_{i=1}^{N}$, initial graph $G_0$, and initial representatives $c_j^{(0)}\in\mathcal V(G_0)$.
\Require Local sampler $\mathsf S$, inner budgets $T_{\ell,j}$,
 outer budget $L_{\max}$.
\State $G\gets G_0$; $c_j\gets c_j^{(0)}$ for $j=1,\ldots,k$.
\For{$\ell=0,\ldots,L_{\max}-1$}
    \State $a_i\gets\operatorname*{arg\,min}_{1\leq j\leq k}
    d_G(z_i,c_j)^2$ for $i=1,\ldots,N$
    \State $C_j\gets\{i:a_i=j\}$ for $j=1,\ldots,k$.
    \For{$j=1,\ldots,k$}
        \If{$C_j\neq\varnothing$}
            \State $(c_j^+,G)\gets\mathsf{GraphMean}
            (c_j,\{z_i:i\in C_j\},G,\mathsf S,T_{\ell,j})$.
        \Else
            \State Choose $c_j^+\in\mathcal V(G)$ using a prespecified empty-cluster rule.
        \EndIf
    \EndFor
    \State $c_j\gets c_j^+$ for $j=1,\ldots,k$.
    \If{the prescribed stopping criterion is satisfied}
        \State \textbf{break}
    \EndIf
\EndFor
\State Recompute $a_i$ and $C_j$ using the final representatives and graph.
\State \Return $\{c_j\}_{j=1}^{k}$, $\{C_j\}_{j=1}^{k}$, and $G$.
\end{algorithmic}
\end{algorithm}

\paragraph{Graph-based Lloyd's algorithm}
We use the graph-based Fr\'echet mean as a representative-update primitive in Lloyd's algorithm for $k$-means clustering. 
Each observation is assigned to its nearest current representative in
graph distance. The observations assigned to a cluster then become the
reference set for a graph-based Fr\'echet-mean computation. For notation, write $\mathsf{GraphMean}(c,Z,G,\mathsf S,T)$ for a call to
Algorithm~\ref{algo:graph-algorithm}, initialized at $c$, with reference collection $Z$, graph $G$, local sampler $\mathsf S$, and $T$ inner updates. Algorithm~\ref{algo:graph-lloyd} gives a schematic outer iteration.

We consider clustering into two groups for the digit pairs $4/7$ and $3/8$, comparing Graph
Lloyd with Ambient Lloyd: Lloyd's algorithm with the ambient averaging primitive. 
For each digit pair, we fix a collection of $N=32$
observations and run each method
$N_{\operatorname{alg}}=100$ times, varying the random seeds
that determine initialization and, for Graph Lloyd,
local sampling. We group the runs into $20$ disjoint
batches of five and retain one run per method from each
batch using the selection rule specified below.
Thus, the aggregate summaries contain $20$ retained
outcomes per method and digit pair, obtained from
$100$ individual runs.

\begin{figure}[t]
    \centering
    \vspace{-1em}
    \includegraphics[scale=0.5]{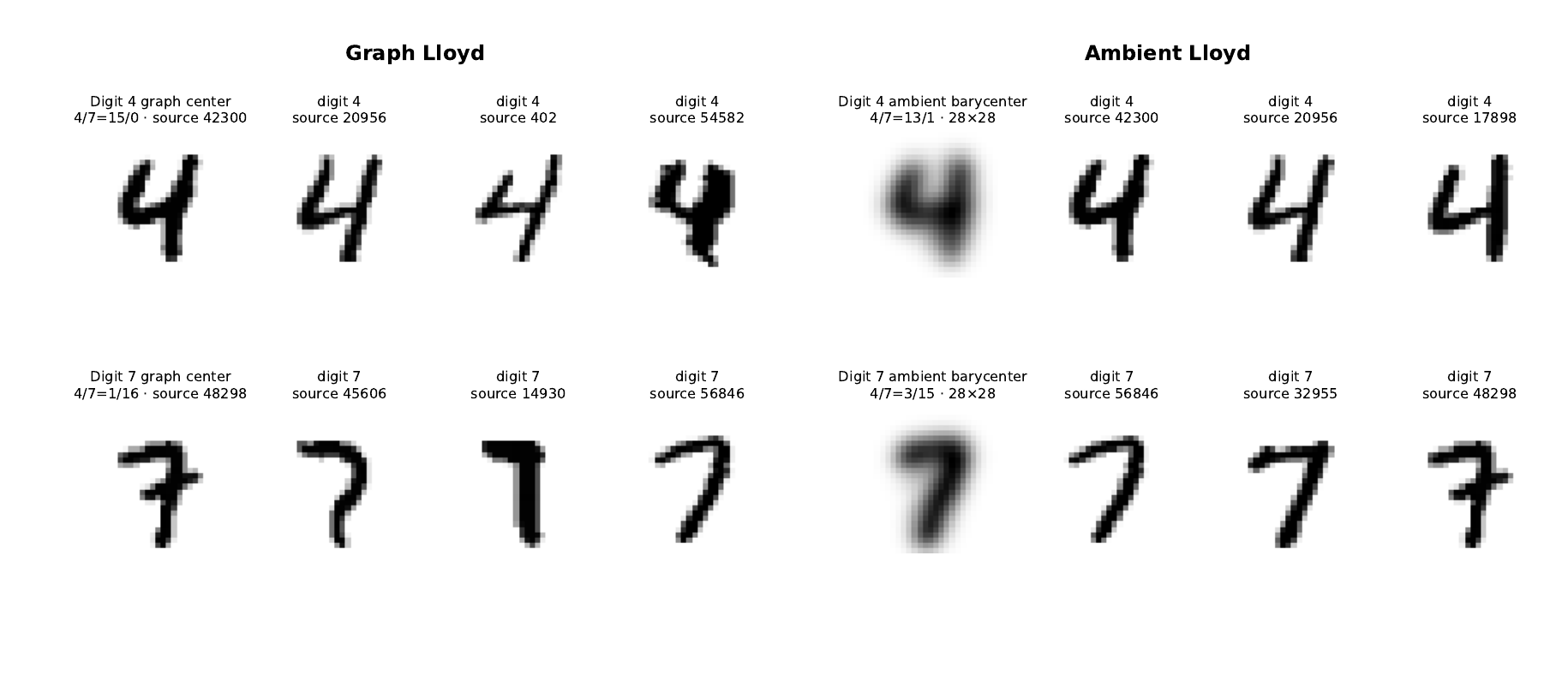}
    \vspace{-2em}
    \caption{Representative $4/7$ Lloyd example; the graph method attains $96.875\%$ accuracy (31/32) while the ambient
    attains $87.5\%$ accuracy (28/32). Counts give the numbers of $4$s and $7$s in each assigned cluster. 
    }
    \label{fig:mnist-representatives}
\end{figure}

\paragraph{Discussion.}
Figure~\ref{fig:mnist-representatives} displays a representative trial for $4$s and $7s$ in which both methods attain a high classification accuracy, which is computed post-algorithm run. We display the resulting representatives, together with their nearest assigned observations. Both methods exhibit better performance than random assignment, highlighting that
using Wasserstein geometry provides useful information about the space of digits. For each
method we highlight the two final representative centers, together with the displayed
nearest assigned observations. However, our method provides centers that lie directly on the submanifold while the ambient Lloyd method provides off-manifold averages.
Additionally, we maintain similar performance outcomes to the ambient space computations despite compressing our images for the use of cheap numerical distance surrogates. We note that
the trial selected was to exhibit both the use of these primitives for classification and the generated representatives; the role is illustrative rather than
an estimate of typical performance. Classification for the full sample of this run is deferred to Appendix~\ref{sec:additional-numerical-results}.

\begin{figure}[t]
    \centering
    \includegraphics[width=0.67\linewidth]{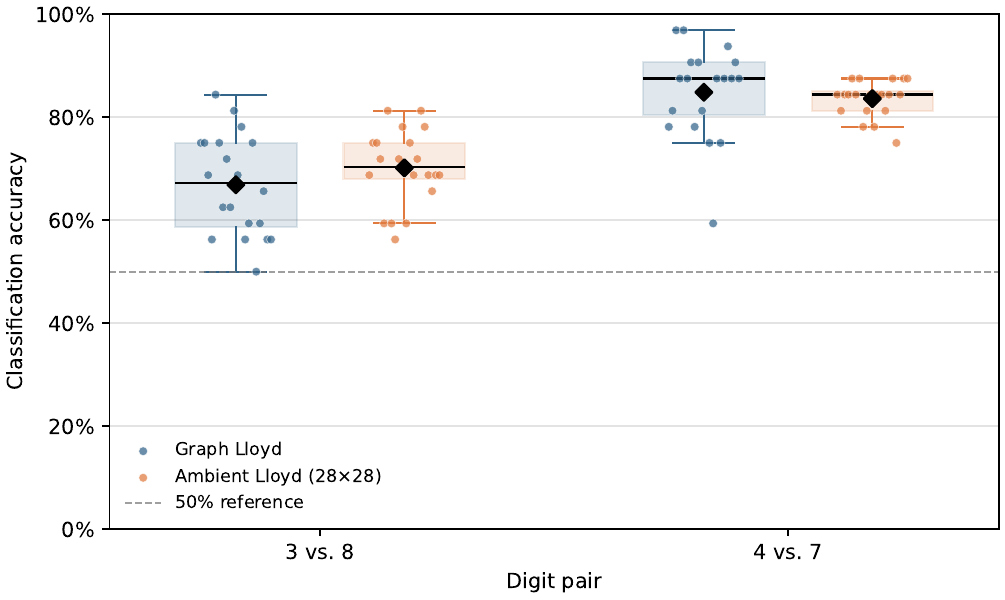}
    \caption{Oracle best-of-five final clustering accuracy
for Graph Lloyd and Ambient Lloyd.
For each method and digit pair, $100$ runs on the same
fixed collection of $32$ observations are grouped into
$20$ batches of five. Each point is the highest final
label-matched accuracy in one batch.
The summaries describe algorithmic variability on
these fixed observations. }
    \label{fig:mnist-accuracy}
\end{figure}

To quantify performance, consider the following metric. Let $a_i\in\{1,2\}$ denote the final cluster assignment and
$y_i\in\{1,2\}$ encode the corresponding digit label. We report 
\begin{equation}
\label{eq:mnist-clustering-accuracy}
    \operatorname{Acc}
    :=\max_{\pi\in\mathfrak S_2}
    \frac{1}{N}\sum_{i=1}^{N}
    \mathbf 1_{\{\pi(a_i)=y_i\}},
\end{equation}
which measures \textit{clustering
accuracy} after optimal label matching.
For method $m$, let
$\operatorname{Acc}_{b,r}^{(m)}$ denote the clustering accuracy
\eqref{eq:mnist-clustering-accuracy}
of run $r\in\{1,\ldots,5\}$ in batch $b$.
We report
\begin{equation}\label{eq:mnist-selected-accuracy}
A_b^{(m)}
:=
\max_{1\leq r\leq5}
\operatorname{Acc}_{b,r}^{(m)},
\qquad \text{for }b=1,\ldots,20.
\end{equation}
This is an oracle best-of-five summary: the ground-truth
digit labels are used to select the reported outcome
after the runs have completed. The results are reported in
Figure~\ref{fig:mnist-accuracy}.

We note that both methods have a harder time accurately clustering
the samples for the subcollection of $3$s and $8$s than $4$s and $7$s. We posit that this reflects the intuition that digits $3$ and $8$ are more similar in appearance than
the digits $4$ and $7$. Graph Lloyd exhibits greater
variability across trials than that of Ambient Lloyd consistent with the
randomness introduced by the graph-based averaging procedure. We note that both procedures
outperform random assignment, indicating how the Wasserstein geometry does capture likeness of the digits.

We remark that digit-label agreement, however, need not fully reflect the geometric structure
identified by either clustering algorithm. The clustering may also reflect morphological variation within each
digit class, such as differences in stroke thickness and the spatial
distribution of mass. Thus, a partition that does not align perfectly with the
digit labels may still capture meaningful structure.
Overall, these results indicate that Graph Lloyd can serve as a useful
primitive for the \textit{coarse downstream task of $k$-means clustering}.

\section{Conclusion}\label{sec:conclusion}

We studied the computation of empirical Fr\'echet means when the data is supported on an unknown Riemannian submanifold $\mathcal M_1\subseteq \mathcal M_2$ and the intrinsic geometry of $\mathcal M_1$ is unavailable. Rather than attempting to reconstruct the entire manifold, our approach approximates only the intrinsic geodesic interpolation required at the current iteration, leading to a graph-based averaging procedure where intrinsic geodesics are approximated by shortest paths in locally sampled graphs.

The analysis connects geometric reconstruction with stochastic optimization. We established a two-sided metric-recovery result for ambient-distance-weighted graphs on submanifolds of a general Riemannian ambient manifold, extending the classical ISOMAP-type comparison beyond Euclidean ambient spaces. These metric estimates yield quantitative control of graph interpolation relative to exact intrinsic interpolation. Interpreting the resulting update as an implicit biased Riemannian stochastic gradient step then allows the geometric approximation error to be propagated into an optimization guarantee. Under the geometric and localization assumptions of Section~\ref{sec:assumptions}, appropriately refined local graphs preserve the $O(T^{-1})$ expected objective-gap dependence of the corresponding exact-geodesic stochastic method. Additionally, this outer iteration guarantee is obtained with sufficient per-iteration and cumulative sampling costs of $\widetilde O(t^{4d+2})$ and $\widetilde O(T^{4d+3})$, respectively, where $d=\dim(\mathcal M_1)$.

The numerical experiments of Section~\ref{sec:numerical-experiments} examine the implementable Algorithm~\ref{algo:graph-algorithm} using the refinement procedure Algorithm~\ref{algo:ball-refinement} in three settings. On the sphere, where intrinsic geometry is available for evaluation, the experiments illustrate how local refinement improves metric accuracy along the paths used for optimization without uniformly refining the entire sampled region. The Wasserstein half-moon example highlights the distinction between averaging within a structured family and computing an unconstrained ambient barycenter. Finally, the MNIST experiments demonstrate the use of the graph-based mean both as a standalone averaging procedure and as a representative-update primitive within Lloyd's algorithm for clustering.

Several questions remain open. Most directly, the convergence analysis uses fresh local graphs sampled around the unknown intrinsic geodesic, whereas the practical procedure retains previously sampled points and places new samples in local balls centered along the current graph shortest path. Establishing comparable guarantees for this adaptive persistent construction would require controlling both graph-centered coverage and the dependence introduced by reusing the evolving graph. 
Finally, verifying
Assumption~\ref{ass:embedded-submanifold}\textup{(b)}
for the Wasserstein families considered numerically remains
a separate step for placing a larger class of examples within this
formal theoretical framework.

More broadly, the results suggest a principle for optimization on unknown geometries: it is unnecessary to learn a complete geometric model before optimization begins. For the Fr\'echet problem studied here, reconstructing sufficiently accurate local geodesic information only where it is needed is enough to recover a first-order optimization guarantee. Understanding how far this principle extends beyond Fr\'echet averaging is a natural direction for future work.

\section*{Acknowledgements}
The authors are grateful to Christian Kuemmerle, Matteo Levi and Marie Jose Saad for early discussions and experiments on an initial version of the project. 

\section*{Declaration of Generative AI Use}
Prior to the release of GPT-5.6, the authors had developed a preliminary manuscript containing the principal mathematical ideas, results, and their proofs, together with preliminary numerical experiments involving submanifolds of Wasserstein space. During subsequent revisions, OpenAI ChatGPT was used to assist in checking mathematical proofs, identifying and removing unnecessary assumptions, and improving grammar and exposition. Additionally, OpenAI Codex was used to generalize an earlier version of the numerical codebase to accommodate arbitrary ambient spaces. All AI-assisted mathematical and computational modifications were reviewed by the authors, who take full responsibility for the correctness of the mathematical arguments, numerical implementations, and results presented in this manuscript.
\bibliographystyle{alpha}
\bibliography{refs}

\begin{thebibliography}{BRGGK23}

\bibitem[ABRW12]{lower-bound-SGD}
Alekh Agarwal, Peter~L. Bartlett, Pradeep Ravikumar, and Martin~J. Wainwright.
\newblock Information-theoretic lower bounds on the oracle complexity of
  stochastic convex optimization.
\newblock {\em IEEE Transactions on Information Theory}, 58(5):3235--3249,
  2012.

\bibitem[AC11]{AC11}
Martial Agueh and Guillaume Carlier.
\newblock Barycenters in the {W}asserstein space.
\newblock {\em SIAM Journal on Mathematical Analysis}, 43(2):904--924, 2011.

\bibitem[ACG19]{ACLG19}
Ery Arias-Castro and Thibaut~Le Gouic.
\newblock Unconstrained and curvature-constrained shortest-path distances and
  their approximation.
\newblock {\em Discrete \& Computational Geometry}, 62(1):1--28, 2019.

\bibitem[ADPY12]{ADPY12}
Marc Arnaudon, Cl{\'e}ment Dombry, Anthony Phan, and Le~Yang.
\newblock Stochastic algorithms for computing means of probability measures.
\newblock {\em Stochastic Processes and their Applications}, 122(4):1437--1455,
  2012.

\bibitem[Afs11]{afsari2011riemannian}
Bijan Afsari.
\newblock Riemannian {$L^p$} center of mass: Existence, uniqueness, and
  convexity.
\newblock {\em Proceedings of the American Mathematical Society},
  139(2):655--673, 2011.

\bibitem[AMS08]{AMS08}
P.-A. Absil, Robert Mahony, and Rodolphe Sepulchre.
\newblock {\em Optimization Algorithms on Matrix Manifolds}.
\newblock Princeton University Press, 2008.

\bibitem[ATV13]{ATV13}
Bijan Afsari, Roberto Tron, and Ren{\'e} Vidal.
\newblock On the convergence of gradient descent for finding the {R}iemannian
  center of mass.
\newblock {\em SIAM Journal on Control and Optimization}, 51(3):2230--2260,
  2013.

\bibitem[BAC19]{BAC19}
Nicolas Boumal, P.-A. Absil, and Coralia Cartis.
\newblock Global rates of convergence for nonconvex optimization on manifolds.
\newblock {\em IMA Journal of Numerical Analysis}, 39(1):1--33, 2019.

\bibitem[BCN18]{bottou2018optimization}
L{\'e}on Bottou, Frank~E. Curtis, and Jorge Nocedal.
\newblock Optimization methods for large-scale machine learning.
\newblock {\em SIAM Review}, 60(2):223--311, 2018.

\bibitem[BdLT00]{bernstein2000graph}
Mira Bernstein, Vin {de Silva}, John~C. Langford, and Joshua~B. Tenenbaum.
\newblock Graph approximations to geodesics on embedded manifolds.
\newblock Technical report, Department of Psychology, Stanford University,
  2000.

\bibitem[BFM17]{BFM17}
Glaydston~C. Bento, Orizon~P. Ferreira, and Jefferson~G. Melo.
\newblock Iteration-complexity of gradient, subgradient and proximal point
  methods on {R}iemannian manifolds.
\newblock {\em Journal of Optimization Theory and Applications},
  173(2):548--562, 2017.

\bibitem[BFRT25]{BFRT25}
Julio Backhoff, Joaquin Fontbona, Gonzalo Rios, and Felipe Tobar.
\newblock Stochastic gradient descent for barycenters in {W}asserstein space.
\newblock {\em Journal of Applied Probability}, 62(1):15--43, 2025.

\bibitem[BHJ24]{BHJ24}
Ronny Bergmann, Roland Herzog, and Hajg Jasa.
\newblock The {R}iemannian convex bundle method, 2024.

\bibitem[BKL06]{BKL06}
Alina Beygelzimer, Sham Kakade, and John Langford.
\newblock Cover trees for nearest neighbor.
\newblock In {\em Proceedings of the 23rd International Conference on Machine
  Learning}, ICML '06, pages 97--104. Association for Computing Machinery,
  2006.

\bibitem[Bon13]{Bon13}
Silv{\`e}re Bonnabel.
\newblock Stochastic gradient descent on {R}iemannian manifolds.
\newblock {\em IEEE Transactions on Automatic Control}, 58(9):2217--2229, 2013.

\bibitem[Bou23]{Bou23}
Nicolas Boumal.
\newblock {\em An Introduction to Optimization on Smooth Manifolds}.
\newblock Cambridge University Press, 2023.

\bibitem[BRGGK23]{bello2023staying}
Juan~M. Bello-Rivas, Anastasia Georgiou, John Guckenheimer, and Ioannis~G.
  Kevrekidis.
\newblock Staying the course: iteratively locating equilibria of dynamical
  systems on {Riemannian} manifolds defined by point-clouds.
\newblock {\em Journal of Mathematical Chemistry}, 61:600--629, 2023.

\bibitem[CM12]{CM12}
Dong Chen and Hans-Georg M{\"u}ller.
\newblock Nonlinear manifold representations for functional data.
\newblock {\em The Annals of Statistics}, 40(1):1--29, 2012.

\bibitem[CMRS20]{CMRS20}
Sinho Chewi, Tyler Maunu, Philippe Rigollet, and Austin~J. Stromme.
\newblock Gradient descent algorithms for {B}ures--{W}asserstein barycenters.
\newblock In Jacob Abernethy and Shivani Agarwal, editors, {\em Proceedings of
  the Thirty-Third Conference on Learning Theory}, volume 125 of {\em
  Proceedings of Machine Learning Research}, pages 1276--1304. PMLR, 2020.

\bibitem[DGM26]{DGM26}
Mateo D{\'i}az, Benjamin Grimmer, and Ian McPherson.
\newblock Nonsmooth {R}iemannian optimization with inexact manifold primitives
  via bundle methods, 2026.

\bibitem[DMSR23]{demidovich2023guidezoobiasedsgd}
Yury Demidovich, Grigory Malinovsky, Igor Sokolov, and Peter Richtárik.
\newblock A guide through the zoo of biased sgd, 2023.

\bibitem[DW25]{DW25}
Willem Diepeveen and Melanie Weber.
\newblock Iso-riemannian optimization on learned data manifolds, 2025.

\bibitem[EAS98]{EAS98}
Alan Edelman, Tom{\'a}s~A. Arias, and Steven~T. Smith.
\newblock The geometry of algorithms with orthogonality constraints.
\newblock {\em SIAM Journal on Matrix Analysis and Applications},
  20(2):303--353, 1998.

\bibitem[FO98]{FO98}
O.~P. Ferreira and P.~R. Oliveira.
\newblock Subgradient algorithm on {R}iemannian manifolds.
\newblock {\em Journal of Optimization Theory and Applications}, 97(1):93--104,
  1998.

\bibitem[FO02]{FO02}
O.~P. Ferreira and P.~R. Oliveira.
\newblock Proximal point algorithm on {R}iemannian manifolds.
\newblock {\em Optimization}, 51(2):257--270, 2002.

\bibitem[FP22]{fetecau2022wellposedness}
Razvan~C. Fetecau and Francesco~S. Patacchini.
\newblock Well-posedness of an interaction model on riemannian manifolds.
\newblock {\em Communications on Pure and Applied Analysis}, 21(11):3559--3585,
  2022.

\bibitem[GM26]{GM26}
{\'A}lvaro Gajardo and Hans-Georg M{\"u}ller.
\newblock Global {Fr{\'e}chet} manifold learning for random objects, with
  application to low-dimensional {W}asserstein representations of
  distributional data.
\newblock {\em Journal of Machine Learning Research}, 27(61):1--49, 2026.

\bibitem[HHK23]{HHK23}
Keaton Hamm, Nick Henscheid, and Shujie Kang.
\newblock Wassmap: {W}asserstein isometric mapping for image manifold learning.
\newblock {\em SIAM Journal on Mathematics of Data Science}, 5(2):475--501,
  2023.

\bibitem[HKT26]{HKT26}
Andreas Hauptmann, Yury Korolev, and Matthew Thorpe.
\newblock Gradient descent on point clouds and applications in learned operator
  correction, 2026.

\bibitem[HMNP23]{HNP23}
Najmeh Hoseini~Monjezi, Soghra Nobakhtian, and Mohamad~Reza Pouryayevali.
\newblock A proximal bundle algorithm for nonsmooth optimization on
  {R}iemannian manifolds.
\newblock {\em IMA Journal of Numerical Analysis}, 43(1):293--325, 2023.

\bibitem[HMST25]{HMST25}
Keaton Hamm, Caroline Moosm{\"u}ller, Bernhard Schmitzer, and Matthew Thorpe.
\newblock Manifold learning in {W}asserstein space.
\newblock {\em SIAM Journal on Mathematical Analysis}, 57(3):2983--3029, 2025.

\bibitem[HW22]{HW22}
Wen Huang and Ke~Wei.
\newblock {R}iemannian proximal gradient methods.
\newblock {\em Mathematical Programming}, 194(1--2):371--413, 2022.

\bibitem[Jos17]{Jost2017}
J{\"u}rgen Jost.
\newblock {\em Riemannian Geometry and Geometric Analysis}.
\newblock Universitext. Springer, 7 edition, 2017.

\bibitem[Kar70]{Karcher1970BergerCurvatureTensorEstimates}
Hermann Karcher.
\newblock A short proof of berger's curvature tensor estimates.
\newblock {\em Proceedings of the American Mathematical Society}, 26:642--644,
  1970.

\bibitem[Kar77]{Kar77}
Hermann Karcher.
\newblock Riemannian center of mass and mollifier smoothing.
\newblock {\em Communications on Pure and Applied Mathematics}, 30(5):509--541,
  1977.

\bibitem[LBBH98]{lecun1998gradient}
Yann LeCun, L{\'e}on Bottou, Yoshua Bengio, and Patrick Haffner.
\newblock Gradient-based learning applied to document recognition.
\newblock {\em Proceedings of the IEEE}, 86(11):2278--2324, 1998.

\bibitem[Lee18]{lee2018introduction}
John~M Lee.
\newblock {\em Introduction to Riemannian Manifolds}, volume 176.
\newblock Springer, 2018.

\bibitem[ORST25]{ORST25}
Mary Chriselda~Antony Oliver, Michael Roberts, Carola-Bibiane Sch{\"o}nlieb,
  and Matthew Thorpe.
\newblock Laplace learning in {W}asserstein space, 2025.

\bibitem[RHM25]{RHM25}
Frederik~M{\"o}bius Rygaard, S{\o}ren Hauberg, and Steen Markvorsen.
\newblock Simultaneous optimization of geodesics and {Fr{\'e}chet} means, 2025.

\bibitem[SAS22]{SAS22}
Boris Shustin, Haim Avron, and Barak Sober.
\newblock Manifold free {R}iemannian optimization, 2022.

\bibitem[TdSL00]{TdSL00}
Joshua~B. Tenenbaum, Vin de~Silva, and John~C. Langford.
\newblock A global geometric framework for nonlinear dimensionality reduction.
\newblock {\em Science}, 290(5500):2319--2323, 2000.

\bibitem[TFBJ18]{TFBJ18}
Nilesh Tripuraneni, Nicolas Flammarion, Francis Bach, and Michael~I. Jordan.
\newblock Averaging stochastic gradient descent on {R}iemannian manifolds.
\newblock In S{\'e}bastien Bubeck, Vianney Perchet, and Philippe Rigollet,
  editors, {\em Proceedings of the 31st Conference on Learning Theory},
  volume~75 of {\em Proceedings of Machine Learning Research}, pages 650--687.
  PMLR, 2018.

\bibitem[TGTP20]{TGTP20}
Michael~W. Trosset, Mingyue Gao, Minh Tang, and Carey~E. Priebe.
\newblock Learning 1-dimensional submanifolds for subsequent inference on
  random dot product graphs, 2020.

\bibitem[TLL20]{trstanova2020local}
Zofia Trstanova, Ben Leimkuhler, and Tony Leli{\`e}vre.
\newblock Local and global perspectives on diffusion maps in the analysis of
  molecular systems.
\newblock {\em Proceedings of the Royal Society A: Mathematical, Physical and
  Engineering Sciences}, 476(2233):20190036, 2020.

\bibitem[XS26]{XS26}
Liane Xu and Amit Singer.
\newblock Manifold learning in metric spaces.
\newblock {\em Applied and Computational Harmonic Analysis}, 80:101813, 2026.

\bibitem[ZS16]{ZS16}
Hongyi Zhang and Suvrit Sra.
\newblock First-order methods for geodesically convex optimization.
\newblock In Vitaly Feldman, Alexander Rakhlin, and Ohad Shamir, editors, {\em
  Proceedings of the 29th Conference on Learning Theory}, volume~49 of {\em
  Proceedings of Machine Learning Research}, pages 1617--1638. PMLR, 2016.

\end{thebibliography}
\hypertarget{toc}{}
\tableofcontents
\begin{appendix}
    
\section{Full Tables of Notation}\label{sec:full-notation-appendix}

Within this section, we collect notation used throughout the rest of the Appendix.
\begingroup
\small
\setlength{\LTleft}{0pt}
\setlength{\LTright}{0pt}
\renewcommand{\arraystretch}{1.12}
\begin{longtable}{@{}>{\raggedright\arraybackslash}p{0.24\textwidth}>{\raggedright\arraybackslash}p{0.72\textwidth}@{}}
\toprule
\textbf{Symbol} & \textbf{Meaning} \\
\midrule
\endfirsthead
\multicolumn{2}{l}{\small\itshape Table continued from previous page}\\
\toprule
\textbf{Symbol} & \textbf{Meaning} \\
\midrule
\endhead
\midrule
\multicolumn{2}{r}{\small\itshape Continued on next page}\\
\endfoot
\bottomrule
\caption{Manifolds, basic geometry, and localization.}\label{tab:notation-app-geometry}\vspace{-1em}
\endlastfoot
\(\nabla_{\mathcal M}f,\ \operatorname{Hess}f\) & Riemannian gradient and Hessian of \(f\). \\
\(\operatorname{sec}_{\mathcal M}\) & Sectional-curvature function of \(\mathcal M\). \\
\(\gamma_{x\to y}^{\mathcal M}\) & A selected minimizing geodesic segment from \(x\) to \(y\) in \(\mathcal M\). \\
\(B_{\mathcal M}(x,r)\) & Intrinsic ball of radius \(r\) centered at \(x\) in \(\mathcal M\). \\
\(\operatorname{dist}_{\mathcal M}(p,A)\) & Intrinsic distance from \(p\) to a subset \(A\subseteq\mathcal M\). \\
\(\mathcal T_r(\gamma)\) & Closed radius-\(r\) geodesic tube \(\{p:\operatorname{dist}_{\mathcal M_1}(p,\gamma)\le r\}\). \\
\(\Pi\) & Second fundamental form of the embedding \(\mathcal M_1\subset\mathcal M_2\). \\
\(\kappa\) & Uniform bound satisfying \(\lVert\Pi_p(v,v)\rVert_p\le\kappa\) for unit \(v\in T_p\mathcal M_1\). \\
\(K_{\min}^{\mathcal M_j},\ K_{\max}^{\mathcal M_j}\) & Global lower and upper sectional-curvature bounds for \(\mathcal M_j\), \(j\in\{1,2\}\). \\
\(\widetilde K_{\mathcal M_j}\) & Absolute curvature bound \(\max\{|K_{\min}^{\mathcal M_j}|,|K_{\max}^{\mathcal M_j}|\}\). \\
\(r_{\mathrm{inj}}^{\mathcal M_j}(S)\) & Uniform injectivity radius \(\inf_{p\in S}\operatorname{inj}_{\mathcal M_j}(p)\). \\
\(S_{\mathrm{reg}}\) & Region of interest containing references, iterates, and local samples \\
\(D_{\mathrm{reg}}\) & Intrinsic diameter \(\sup_{x,y\in S_{\mathrm{reg}}}d_{\mathcal M_1}(x,y)\). \\
\(\rho_S(s)\) & Ambient-to-intrinsic modulus \(\sup\{d_{\mathcal M_1}(x,y):x,y\in S,\ d_{\mathcal M_2}(x,y)\le s\}\). \\
\(C_{\mathrm{loc},S},\ s_{\mathrm{loc},S}\) & \(d_{\mathcal M_1}(x,y)\le C_{\mathrm{loc},S}d_{\mathcal M_2}(x,y)\) at ambient scales at most \(s_{\mathrm{loc},S}\). \\

\end{longtable}
\endgroup

\begingroup
\small
\setlength{\LTleft}{0pt}
\setlength{\LTright}{0pt}
\renewcommand{\arraystretch}{1.12}
\begin{longtable}{@{}>{\raggedright\arraybackslash}p{0.24\textwidth}>{\raggedright\arraybackslash}p{0.72\textwidth}@{}}
\toprule
\textbf{Symbol} & \textbf{Meaning} \\
\midrule
\endfirsthead
\multicolumn{2}{l}{\small\itshape Table continued from previous page}\\
\toprule
\textbf{Symbol} & \textbf{Meaning} \\
\midrule
\endhead
\midrule
\multicolumn{2}{r}{\small\itshape Continued on next page}\\
\endfoot
\bottomrule
\caption{Fréchet objective and optimization.}\label{tab:notation-app-optimization}\vspace{-1 em}
\endlastfoot
\(z_r\) & The \(r\)-th fixed reference point in \(\mathcal M_1\). \\
\(P_N\) & Empirical measure \(N^{-1}\sum_{r=1}^N\delta_{z_r}\). \\
\(f_{P_N}\) & Empirical Fréchet functional \((2N)^{-1}\sum_{r=1}^N d_{\mathcal M_1}^2(\cdot,z_r)\). \\
\(\widehat x_N\) & A selected minimizer of \(f_{P_N}\) over \(\mathcal M_1\). \\
\(\varepsilon\) & Target objective accuracy in the definition of an \(\varepsilon\)-minimizer. \\
\(y_t=z_{I_t}\) & Randomly selected reference at iteration \(t\) with $I_{t} \overset{i.i.d.}{\sim} \mathrm{Unif}[N]$. \\
\(\eta_t\) & Interpolation fraction and SGD stepsize for the update from \(t\) to \(t+1\). \\
\(\widetilde g_t\) & Gradient estimator associated with Algorithm~\ref{alg:geodesic-interpolation}, see \eqref{eqn:unbiased-gradient-estimator}.\\
\(g_t\) & Graph estimator associated with Algorithm~\ref{algo:graph-algorithm}, see \eqref{eqn:biased-graph-gradient-estimator}.\\ 
\(D_t\) & Maximum reference distance \(\max_{1\le i\le N}d_{\mathcal M_1}(x_t,z_i)\).\\
\(m_f,\ L_f\) & Strong-convexity and smoothness constants for $f$, see Lemma~\ref{lem:functional-SC-smooth-parameters}.\\
\(\beta,\ \gamma\) & Numerator and offset in the decaying stepsize \(\eta_t=\beta/(\gamma+t)\).\\
\(\nu\) & Constant controlling the \(1/(\gamma+t)\) objective-gap envelope. \\
\(\widehat t=\gamma+t\) & Shifted iteration counter used in Appendix F.2. \\
\end{longtable}
\endgroup

\begingroup
\small
\setlength{\LTleft}{0pt}
\setlength{\LTright}{0pt}
\renewcommand{\arraystretch}{1.12}
\begin{longtable}{@{}>{\raggedright\arraybackslash}p{0.24\textwidth}>{\raggedright\arraybackslash}p{0.72\textwidth}@{}}
\toprule
\textbf{Symbol} & \textbf{Meaning} \\
\midrule
\endfirsthead
\multicolumn{2}{l}{\small\itshape Table continued from previous page}\\
\toprule
\textbf{Symbol} & \textbf{Meaning} \\
\midrule
\endhead
\midrule
\multicolumn{2}{r}{\small\itshape Continued on next page}\\
\endfoot
\bottomrule
\caption{Graphs, shortest paths, and interpolation.}\label{tab:notation-app-graphs}\vspace{-1em}
\endlastfoot
\(G,\ G_t\) & Ambient-distance-weighted graph and the graph available at iteration \(t\). \\
\(V(G),\ E(G)\) & Vertex and edge sets of a graph \(G\). \\
\(d_G\) & Weighted shortest-path metric induced by \(G\). \\
\(\gamma_{x\to y}^{G}\) & A selected shortest path from \(x\) to \(y\) in \(G\). \\
\(v_j^{(t)},\ m_t\) & Ordered vertices and final vertex index of a graph geodesic.\\
\(\lambda_1^{(t)},\ \lambda_2^{(t)},\ \lambda_t\) & Iteration-dependent distortions; on event \(E_t\) we have $\lambda_{1}^{(t)} = \lambda_{2}^{(t)} = \lambda_{t}$.\\
\(L_{\mathrm{cmp}}\) & Maximum intrinsic segment length allowed for Theorem~\ref{thm:deterministic-metric-equivalence-for-ISOMAP}. \\
\(s_0\) & Maximum ambient edge scale admitted for Theorem~\ref{thm:deterministic-metric-equivalence-for-ISOMAP}.\\
\(\epsilon_t\) & Iteration-dependent graph connectivity scale.\\ 
\(\delta_{t}\) & Intrinsic covering radius discretizing a set. \\
\end{longtable}
\endgroup

\begingroup
\small
\setlength{\LTleft}{0pt}
\setlength{\LTright}{0pt}
\renewcommand{\arraystretch}{1.12}
\begin{longtable}{@{}>{\raggedright\arraybackslash}p{0.24\textwidth}>{\raggedright\arraybackslash}p{0.72\textwidth}@{}}
\toprule
\textbf{Symbol} & \textbf{Meaning} \\
\midrule
\endfirsthead
\multicolumn{2}{l}{\small\itshape Table continued from previous page}\\
\toprule
\textbf{Symbol} & \textbf{Meaning} \\
\midrule
\endhead
\midrule
\multicolumn{2}{r}{\small\itshape Continued on next page}\\
\endfoot
\bottomrule
\caption{Sampling, covering, and probability.}\label{tab:notation-app-probability}\vspace{-1em}
\endlastfoot
\(r_{t}, \alpha_t\) & Radius and relative sampling tube-radius, $r_{t} = \alpha_{t}d_{\mathcal{M}_{1}}(x_{t}, y_{t+1})$. \\
\(\widetilde r_{t+1}=r_{t+1}+\rho_{S_{\mathrm{reg}}}(\epsilon_{t+1})\) & Inflated radius allowing one admissible ambient edge. \\
\(\mathcal T_{t+1},\ \widetilde{\mathcal T}_{t+1}\) & Target tube \(\mathcal T_{r_{t+1}}(\gamma_{x_t\to y_{t+1}}^{\mathcal M_1})\) and its inflated counterpart. \\
\(\mathfrak S_{\alpha_t,\epsilon_t}\) & Family of all admissible inflated tubes at iteration \(t\). \\
\(\xi_t\) & Local sampler randomness at iteration \(t\), conditional on $x_{t}$ and $y_{t+1}$. \\
\(\mathcal F_t\) & Filtration \(\sigma(x_0,I_1,\xi_1,\ldots,I_t,\xi_t)\) for the full graph-based process. \\
\(\mathcal H_t\) & Intermediate filtration \(\mathcal F_{t-1}\vee\sigma(I_t)\), isolating sampling randomness. \\
\(\widetilde r_{\max}^{(t)}\) & Largest admissible inflated radius at iteration \(t\): \(\alpha_t\ell_{\max}+\rho_{S_{\mathrm{reg}}}(\epsilon_t)\). \\
\(E_t\) & Event that graph on \(\widetilde{S}_t\) satisfies local metric equivalence; see \eqref{eqn:defining-per-iterate-good-event}.\\ 
\(\Omega_t\) & Cumulative good event \(\bigcap_{s=1}^tE_s\), with \(\Omega_0=\Omega\). \\
\(p_t,\ q_t\) & Conditional probabilities \(p_t=\mathbb P(E_{t+1}\mid\mathcal F_t)\) and \(q_t=\mathbb P(E_{t+1}^c\mid\mathcal F_t)\). \\
\(\zeta,\ \zeta_t,\ \zeta_T\) & Failure tolerances of metric equivalence for Theorem~\ref{thm:main-conv-thm}, 
and 
Theorem~\ref{thm:main-conv-thm-body}. \\
\(p_{\mathrm{good}}\) & Uniform target lower bound on good sample event $E_{t}$. \\
\end{longtable}
\endgroup

\begingroup
\small
\setlength{\LTleft}{0pt}
\setlength{\LTright}{0pt}
\renewcommand{\arraystretch}{1.12}
\begin{longtable}{@{}>{\raggedright\arraybackslash}p{0.24\textwidth}>{\raggedright\arraybackslash}p{0.72\textwidth}@{}}
\toprule
\textbf{Symbol} & \textbf{Meaning} \\
\midrule
\endfirsthead
\multicolumn{2}{l}{\small\itshape Table continued from previous page}\\
\toprule
\textbf{Symbol} & \textbf{Meaning} \\
\midrule
\endhead
\midrule
\multicolumn{2}{r}{\small\itshape Continued on next page}\\
\endfoot
\bottomrule
\caption{Comparison functions and auxiliary geometric notation.}\label{tab:notation-app-comparison}
\endlastfoot
\(\operatorname{ct}_K(t)\) & Model cotangent function for Hessian-comparison defined at \eqref{eqn:comparison-function-for-Hessian-comparison}.\\
\(\widetilde\Psi(\ell,r)\) & Curvature-dependent function for tube containment, see \eqref{eqn:function-for-tube-containment}.\\
\(C_{\log}(D_{\mathrm{reg}},K_{\max}^{\mathcal M_1})\) & Uniform Lipschitz bound for \(y\mapsto\log_x^{\mathcal M_1}(y)\) on \(S_{\mathrm{reg}}\), see \eqref{eqn:uniform-rauche-comparison-constant}.\\
\(s_c(t)\) & Model Jacobi function for Jacobi field comparisons, see \eqref{eqn:comparison-function}. \\
\(A_{K_{\min}^{\mathcal M_2}}(L),\ B_{K_{\max}^{\mathcal M_2}}(L)\) & Normalized hyperbolic and spherical Jacobi comparison factors, see \eqref{eqn:constants-for-Jacobi-control}. \\
\(C_J(L)\) & Ambient Jacobi comparison factor \(\max\{1,A_{K_{\min}^{\mathcal M_2}}(L)/B_{K_{\max}^{\mathcal M_2}}(L)\}\). \\
\(\widehat\gamma(L)\) & Radially transported surrogate \(\int_0^LP_{0\leftarrow1}^{\sigma_s}[\gamma'(s)]\,ds\). \\
\(L_{\widetilde K}(r)\) & Jacobi-based Lipschitz factor \(\max\{\cosh(\sqrt{\widetilde K}r),\sinh(\sqrt{\widetilde K}r)/(\sqrt{\widetilde K}r)\}\). \\
\(C_{\widetilde K}(\tau),\ S_{\widetilde K}(\tau)\) & Hyperbolic comparison functions used in Lemma 40, see \eqref{eqn:hyperbolic-comparison-funtions}.\\
\end{longtable}
\endgroup

\section{Functional Properties Proof: Lemma~\ref{lem:functional-SC-smooth-parameters}}\label{sec:proof-of-functional-SC-smooth-parameters} 
We first record the required Hessian bounds for the squared-distance function.

\begin{lemma}\label{lem:sc-and-smoothness-distance-squared}
Let $(\mathcal{M},g)$ be a Riemannian manifold satisfying $K_{\min}\leq \operatorname{sec}_{\mathcal{M}}\leq K_{\max}$. Let $S\subseteq\mathcal{M}$ be geodesically convex, let $y\in S$, and suppose that $S$ is contained in a normal neighborhood $U_y$ of $y$. Define
\[
D_y:=\sup_{x\in S}d_{\mathcal{M}}(x,y)<+\infty,
\qquad
\phi_y(x):=\frac{1}{2}d_{\mathcal{M}}^2(x,y).
\]
If $K_{\max}>0$, assume in addition that $D_y<\frac{\pi}{2\sqrt{K_{\max}}}$. Then $\phi_y\in C^2(U_y)$ and, for every $x\in S$,
\begin{equation}
m_yg_x \preceq \operatorname{Hess}\phi_y(x)\preceq L_yg_x,
\end{equation}
where
\[
m_y:=
\begin{cases}
1, & K_{\max}\leq 0,\\[2mm]
D_y\sqrt{K_{\max}}
\cot\!\bigl(\sqrt{K_{\max}}D_y\bigr), & K_{\max}>0,
\end{cases}
\]
and
\[
L_y:=
\begin{cases}
1, & K_{\min}\geq 0,\\[2mm]
D_y\sqrt{-K_{\min}}
\coth\!\bigl(\sqrt{-K_{\min}}D_y\bigr), & K_{\min}<0.
\end{cases}
\]
When $D_y=0$, the expressions above are interpreted by continuous extension.
\end{lemma}

\begin{proof}
Let
$r(x):=d_{\mathcal{M}}(x,y).$
Since $U_y$ is a normal neighborhood of $y$, the squared-distance
function $\varphi_y$ is twice continuously differentiable on $U_y$. On $U_y\setminus\{y\}$, by standard radial
identities 
\begin{equation} \label{eqn:standard-radial-identities}
    \nabla\phi_y=r\nabla r=-\log_x(y),
\qquad
\operatorname{Hess}\phi_y
=dr\otimes dr+r\,\operatorname{Hess}r.
\end{equation}

Define
\begin{equation} \label{eqn:comparison-function-for-Hessian-comparison}
\operatorname{ct}_K(t):=
\begin{cases}
\sqrt{K}\cot(\sqrt{K}t), & K>0,\\[1mm]
1/t, & K=0,\\[1mm]
\sqrt{-K}\coth(\sqrt{-K}t), & K<0.
\end{cases}
\end{equation}
By the Hessian comparison theorem
(e.g., Theorem~11.7 in \cite{lee2018introduction}), on
$U_y\setminus\{y\}$,
\begin{equation}\label{eqn:Hessian-comparison-for-radial-distance}
\operatorname{ct}_{K_{\max}}(r)
\bigl(g-dr\otimes dr\bigr)
\leq \operatorname{Hess}r
\leq
\operatorname{ct}_{K_{\min}}(r)
\bigl(g-dr\otimes dr\bigr).
\end{equation}

Fix $x\in S\setminus\{y\}$, let $u:=\nabla r(x)$, and decompose arbitrary $v\in T_x\mathcal{M}$ as
\[
v=v^\parallel+v^\perp,
\qquad
v^\parallel:=\langle v,u\rangle u,
\qquad
v^\perp\perp u.
\]
Using this decomposition of $v$ with \eqref{eqn:standard-radial-identities} and $\operatorname{Hess}r(u,\cdot)=0$, we have
\[
\operatorname{Hess}\phi_y(x)[v,v]
=
\|v^\parallel\|^2
+r\,\operatorname{Hess}r(x)[v^\perp,v^\perp].
\]
Therefore, by application of \eqref{eqn:Hessian-comparison-for-radial-distance}
\[
\|v^\parallel\|^2
+r\operatorname{ct}_{K_{\max}}(r)\|v^\perp\|^2
\leq
\operatorname{Hess}\phi_y(x)[v,v]
\leq
\|v^\parallel\|^2
+r\operatorname{ct}_{K_{\min}}(r)\|v^\perp\|^2.
\]
Thus, the eigenvalues of Hess$\phi_{y}$ for the radial direction is $1$, while for the transverse directions lie within $[r\operatorname{ct}_{K_{\max}}(r), r \operatorname{ct}_{K_{\min}}(r)]$.

The function $t\mapsto t\operatorname{ct}_K(t)$ is decreasing on
$(0,\pi/\sqrt{K})$ when $K>0$, identically equal to $1$ when $K=0$,
and increasing on $(0,+\infty)$ when $K<0$. Since $r\leq D_y$, it follows
that
\[
r\operatorname{ct}_{K_{\max}}(r)\geq m_y,
\qquad \text{and} \qquad
r\operatorname{ct}_{K_{\min}}(r)\leq L_y.
\]
Moreover, notice that $m_y\leq 1\leq L_y$. Therefore, all eigenvalues of Hess$\phi_{y}$ lie within $[m_{y}, L_{y}]$, i.e.,
\[
m_y\|v\|^2
\leq
\operatorname{Hess}\phi_y(x)[v,v]
\leq
L_y\|v\|^2.
\]
At $x=y$, normal coordinates centered at $y$ give $\operatorname{Hess}\phi_y(y)=g_y:$ the same inequalities hold.
\end{proof}

\begin{proof}[Proof of Lemma~\ref{lem:functional-SC-smooth-parameters}]
Let $\{z_r\}_{r=1}^N\subseteq S_{\mathrm{reg}}$ denote the reference
points. Then
\[
f(x)
=
\frac{1}{2N}\sum_{r=1}^N
d_{\mathcal{M}_{1}}^2(x,z_r)
=
\frac{1}{N}\sum_{r=1}^N\phi_{z_r}(x).
\]
For each $r$, define $D_{z_r}:= \sup_{x\in S_{\mathrm{reg}}}d_{\mathcal{M}_{1}}(x,z_r).$ Since $z_r\in S_{\mathrm{reg}}$, $D_{z_r}\leq D_{\mathrm{reg}}.$
Under Assumptions~\ref{ass:embedded-submanifold}
and~\ref{ass:region-of-interest}, with the strict intrinsic locality
conditions stated above, Lemma~\ref{lem:sc-and-smoothness-distance-squared}
applies to each $\phi_{z_r}$. The monotonicity of
$t\mapsto t\operatorname{ct}_K(t)$ and the inequality
$D_{z_r}\leq D_{\mathrm{reg}}$ give
$m_{f} \leq m_{z_r} $ and $
L_{z_r}\leq \smoothf $
for every $r\in\{1,\ldots,N\}$. Consequently,
\[
\SCf g_x
\preceq
\operatorname{Hess}f(x)
=
\frac{1}{N}\sum_{r=1}^N
\operatorname{Hess}\phi_{z_r}(x)
\preceq
\smoothf g_x
\]
for every $x\in S_{\mathrm{reg}}$. The claimed strong geodesic convexity and geodesic smoothness then follow from the standard equivalence between uniform Hessian bounds and the corresponding first-order inequalities on a geodesically convex set.
\end{proof}

\section{Biased SGD Convergence Proof: Theorem~\ref{thm:biased-sgd-convergence}}\label{sec:biased-sgd-convergence-proof}
For each $t \geq 0$, define the random objective gap
\begin{equation*}
    \Delta_{t} := f(x_{t}) - f^{\star}.
\end{equation*}
When the expected gap is needed, we write
\begin{equation*}
    \overline{\Delta}_{t} := \mathbb{E}[\Delta_{t}].
\end{equation*}
\begin{proof}[Proof of Theorem~\ref{thm:biased-sgd-convergence}]
For both stepsize schedules, $\eta_t\leq\eta_0 \leq \frac{c_1}{\smoothf \sigma_1^2},$ with the convention that $c_1/(\smoothf \sigma_1^2)=+\infty$ when $\sigma_1^2=0$. For every $t \geq 0$, Lemma~\ref{lem:expected-one-step-improvement}
gives
\begin{align*}
\mathbb{E}[f(x_{t+1})\mid\mathcal{F}_t]-f(x_t)
&\leq
-2\eta_t\SCf 
\left(
c_1-\frac{\eta_t\smoothf \sigma_1^2}{2}
\right)
[f(x_t)-f_\ast]
+\frac{\eta_t^2\smoothf \sigma_0^2}{2}\\
&\leq
-\eta_t\SCf c_1[f(x_t)-f_\ast]
+\frac{\eta_t^2\smoothf \sigma_0^2}{2}.
\end{align*}
Rearrangement, subtracting $f_{\ast}$, and taking total expectations yields
\begin{equation}\label{eqn:decaying-stepsize-pre-induction}
\overline{\Delta}_{t+1}
\leq
(1-\eta_t\SCf c_1)\overline{\Delta}_t
+\frac{\eta_t^2\smoothf \sigma_0^2}{2}.
\end{equation}

\paragraph{Constant Stepsize.} Subtracting $\bar{\eta}L_{f} \sigma_{0}^{2}/(2m_{f} c_{1})$ from both sides of
\eqref{eqn:decaying-stepsize-pre-induction} gives
\begin{align}
\overline{\Delta}_{t+1}- \frac{\bar{\eta}L_{f} \sigma_{0}^{2}}{2m_{f} c_{1}} & \leq (1-\bar{\eta} m_{f} c_{1})\overline{\Delta}_{t} + \frac{\bar{\eta}^{2}L_{f} \sigma_{0}^{2}}{2} - \frac{\bar{\eta} L_{f} \sigma_{0}^{2}}{2 m_{f} c_{1}} \nonumber\\
&=
(1 - \bar{\eta} m_{f} c_{1})\left(\overline{\Delta}_t-\frac{\bar{\eta} L_{f} \sigma_{0}^{2}}{2m_{f} c_{1}}\right).
\label{eqn:recursion-for-constant-stepsize}
\end{align}
Iterating \eqref{eqn:recursion-for-constant-stepsize} from $0$ to
$t-1$ and recalling $\overline{\Delta}_{t} = \mathbb{E}[f(x_{t})] - f_{\ast}$ recovers \eqref{eqn:constant-stepsize-convergence-guarantee}. Moreover, noting that $(1 - \bar{\eta} m_{f} c_{1}) \in [0,1)$ since $0 \leq \bar{\eta} < 1/(m_{f}c_{1})$, establishes the limit of the recursion. 

\paragraph{Decaying Stepsize.} Let $a = \beta m_{f} c_{1}$ and $b = \beta^{2} L_{f} \sigma_{0}^{2}/2$. Since $\beta>1/(\SCf c_1)$, we have $a>1$. Also, note that
$\eta_0\leq1/(\SCf c_1)$ is equivalent to $\gamma\geq a$, and thus $0 \leq 1 - a/(\gamma + t) < 1$ for every $t \geq 0$. 
We show \eqref{eqn:decaying-stepsize-convergence-guarantee} holds for every $t\geq0$ by induction. The case $t=0$ follows from
$\nu\geq\gamma\overline{\Delta}_0$.

Suppose the bound holds for some $t\geq0$, and set $\widehat{t} = \gamma + t$. Then, 
\begin{align*}
\overline{\Delta}_{t+1} \leq
\left(1-\frac{a}{\widehat{t}}\right)\frac{\nu}{\widehat{t}}+\frac{b}{\widehat{t}^2} =
\frac{\nu(\widehat{t}-1)}{\widehat{t}^2} - \frac{(a-1)\nu-b}{\widehat{t}^2} \leq
\frac{\nu(\widehat{t}-1)}{\widehat{t}^2} \leq
\frac{\nu}{\widehat{t}+1}.
\end{align*}
The first inequality follows from \eqref{eqn:decaying-stepsize-pre-induction} and the induction hypothesis, the second since $\nu \geq \frac{\beta^{2}L_{f} \sigma_{0}^{2}}{2(\beta \SCf  c_{1} - 1)} = \frac{b}{a-1}$, and the last from $(\widehat{t}-1)(\widehat{t}+1)\leq \widehat{t}^2.$
\end{proof}

\section{Generalized ISOMAP Guarantees}
\subsection{Deterministic Metric Equivalence Proof: Theorem~\ref{thm:deterministic-metric-equivalence-for-ISOMAP}}\label{sec:proof-of-deterministic-metric-equivalence-for-isomap}
We prove the upper and lower bounds separately. The upper bound is the usual discretization argument for ISOMAP; we include the proof for completion.

\begin{lemma}[Adapted from Theorem~2 of \cite{bernstein2000graph}]\label{lem:global-upper-bound-on-graph-distance}
    Assume that $S\subseteq\mathcal{M}_{1}$ is geodesically convex, that $\{x_i\}_{i=1}^{n}$ is a $\delta$-cover of $S$ with respect to $d_{\mathcal{M}_{1}}$, and that $G$ contains every edge $uv$ satisfying $d_{\mathcal{M}_{2}}(u,v)\leq\epsilon_{\min}$. If $\lambda_{2}\in(0,1)$ and
    \begin{equation*}
        \delta\leq \frac{\lambda_{2}\epsilon_{\min}}{4},
    \end{equation*}
    then, for all $x,y\in\{x_i\}_{i=1}^{n}$,
    \begin{equation*}
        d_G(x,y)\leq (1+\lambda_{2})d_{\mathcal{M}_{1}}(x,y).
    \end{equation*}
\end{lemma}

\begin{proof}
    Fix $x,y\in\{x_i\}_{i=1}^{n}$ and write $L=d_{\mathcal{M}_{1}}(x,y)$. The claim is immediate when $L=0$. Otherwise, let $\gamma\colon[0,L]\to S$ be a unit-speed minimizing $\mathcal{M}_{1}$-geodesic from $x$ to $y$. Since $\delta\leq\lambda_{2}\epsilon_{\min}/4$ and $\lambda_{2}<1$, we have $\epsilon_{\min}-2\delta>0$. Set
    \begin{equation*}
        m:=\left\lceil \frac{L}{\epsilon_{\min}-2\delta}\right\rceil,
        \qquad
        p_j:=\gamma\left(\frac{jL}{m}\right),
        \qquad j=0,\dots,m.
    \end{equation*}
    Set $q_0=x$ and $q_m=y$. For each $j\in\{1,\dots,m-1\}$, choose a data point $q_j$ with $d_{\mathcal{M}_{1}}(q_j,p_j)\leq\delta$. Omitting repeated consecutive vertices if necessary, we obtain since $d_{\mathcal{M}_{2}} \leq d_{\mathcal{M}_{1}}$:
    \begin{align*}
        d_{\mathcal{M}_{2}}(q_{j-1},q_j)
        \leq d_{\mathcal{M}_{1}}(q_{j-1},p_{j-1})
        +d_{\mathcal{M}_{1}}(p_{j-1},p_j)
        +d_{\mathcal{M}_{1}}(p_j,q_j) 
        \leq 2\delta+\frac{L}{m}
        \leq \epsilon_{\min}.
    \end{align*}
    Hence every consecutive pair is joined by an edge of $G$. The resulting graph path satisfies
    \begin{align*}
        d_G(x,y)
        &\leq \sum_{j=1}^{m}d_{\mathcal{M}_{2}}(q_{j-1},q_j) 
        \leq L+2\delta(m-1) 
        \leq \left(1+\frac{2\delta}{\epsilon_{\min}-2\delta}\right)L
        \leq (1+\lambda_{2})L,
    \end{align*}
    where the last inequality follows from $\delta\leq\lambda_{2}\epsilon_{\min}/4$ and $\lambda_{2}\in(0,1)$.
\end{proof}

To prove the lower bound, we use the following local chord--arc estimate, whose proof is deferred to \Cref{sec:proof-of-local-chord-arc-lower-bound}.

\begin{lemma}\label{lem:local-chord-arc-lower-bound}
    Assume Assumption~\ref{ass:embedded-submanifold}\textup{(a)} holds. Let $\gamma\colon[0,L]\to\mathcal{M}_{1}$ be a unit-speed minimizing $\mathcal{M}_{1}$-geodesic from $x=\gamma(0)$ to $y=\gamma(L)$, viewed as a curve in $\mathcal{M}_{2}$. Assume that $\gamma([0,L])$ is contained in a normal neighborhood of $x$ in $\mathcal{M}_{2}$, that $L\leq\pi/\kappa$, and, if $\kmaxExt>0$, that $L<\pi/\sqrt{\kmaxExt}$. Then
    \begin{equation*}
        \left(1-\left[\frac{\kappa^{2}}{24}
        +\frac{5\kmagExt C_J(L)}{18}\right]L^{2}\right)L
        \leq d_{\mathcal{M}_{2}}(x,y),
 \qquad \text{where} \quad 
        C_J(L):=\max\left\{1,\frac{A_{\kminExt}(L)}{B_{\kmaxExt}(L)}\right\},
    \end{equation*}
    with
    \begin{align*}
        A_{\kminExt}(L)
        &:={}
        \begin{cases}
            1, & \kminExt\geq 0,\\
            \frac{\sinh(\sqrt{-\kminExt}L)}{\sqrt{-\kminExt}L}, & \kminExt<0,
        \end{cases} \quad \text{and} \quad
        B_{\kmaxExt}(L)
        &:=
        \begin{cases}
            1, & \kmaxExt\leq 0,\\
            \frac{\sin(\sqrt{\kmaxExt}L)}{\sqrt{\kmaxExt}L}, & \kmaxExt>0.
        \end{cases}
    \end{align*}
\end{lemma}

\begin{proof}[Proof of Theorem~\ref{thm:deterministic-metric-equivalence-for-ISOMAP}]
    The upper bound follows from Lemma~\ref{lem:global-upper-bound-on-graph-distance}. We prove the lower bound.

    Fix vertices $x,y$ of $G$, and let $x=x_0,x_1,\dots,x_m=y$ be the vertices of a shortest path in $G$. Set $\ell_i:=d_{\mathcal{M}_{1}}(x_i,x_{i+1}),$ for $i=0,\dots,m-1.$ Each graph edge has ambient length at most $\epsilon_{\max}\leq s_0$, and hence
    \begin{equation*}
        \ell_i
        \leq \rho_S(\epsilon_{\max})
        \leq \rho_S(s_0)
        \leq \compScale.
    \end{equation*}
    Let $\gamma_i\colon[0,\ell_i]\to\mathcal{M}_{1}$ be a unit-speed minimizing geodesic from $x_i$ to $x_{i+1}$. For every $s\in[0,\ell_i]$,
    \begin{equation*}
        d_{\mathcal{M}_{2}}(x_i,\gamma_i(s))
        \leq d_{\mathcal{M}_{1}}(x_i,\gamma_i(s))
        =s
        \leq \ell_i
        \leq \rho_S(s_0)
        <r_{\mathrm{inj}}^{\mathcal{M}_{2}}(S).
    \end{equation*}
    Thus $\gamma_i([0,\ell_i])$ lies in a normal neighborhood of $x_i$. The remaining length restrictions in Lemma~\ref{lem:local-chord-arc-lower-bound} follow from $\ell_i\leq\compScale$ and the assumptions on $\compScale$.

    On the admissible interval, $A_{\kminExt}$ is nondecreasing and $B_{\kmaxExt}$ is positive and nonincreasing; hence $C_J$ is nondecreasing. Therefore,
    \begin{align*}
        \left[\frac{\kappa^{2}}{24}
        +\frac{5\kmagExt C_J(\ell_i)}{18}\right]\ell_i^{2}
        \leq
        \left[\frac{\kappa^{2}}{24}
        +\frac{5\kmagExt C_J(\compScale)}{18}\right]\rho_S(s_0)^{2} 
        \leq \lambda_1,
    \end{align*}
    where the last inequality is the short-edge distortion condition in Theorem~\ref{thm:deterministic-metric-equivalence-for-ISOMAP}. By Lemma~\ref{lem:local-chord-arc-lower-bound}, we have
    \begin{equation}\label{eqn:used-local-lower-bound}
        (1-\lambda_1)d_{\mathcal{M}_{1}}(x_i,x_{i+1})
        \leq d_{\mathcal{M}_{2}}(x_i,x_{i+1}),
    \end{equation}
    for every $i$. Summing over the shortest path and applying the triangle inequality in $\mathcal{M}_{1}$ yields
    \begin{align*}
        d_G(x,y)
        &=\sum_{i=0}^{m-1}d_{\mathcal{M}_{2}}(x_i,x_{i+1}) 
        \geq (1-\lambda_1)\sum_{i=0}^{m-1}d_{\mathcal{M}_{1}}(x_i,x_{i+1}) 
        \geq (1-\lambda_1)d_{\mathcal{M}_{1}}(x,y).
    \end{align*}
    This proves the lower bound and completes the proof.
\end{proof}

\subsubsection{Proof of Lemma~\ref{lem:local-chord-arc-lower-bound}}\label{sec:proof-of-local-chord-arc-lower-bound}
We proceed in two steps. First, we establish a chord--arc inequality for the Cartan development of the $\mathcal{M}_{1}$-geodesic in the tangent space $T_x\mathcal{M}_{2}$. The proof is deferred to \Cref{sec:proof-of-cartan-development-chord}.

\begin{lemma}\label{lem:cartan-development-chord}
    Assume Assumption~\ref{ass:embedded-submanifold}\textup{(a)} holds. Let $\gamma\colon[0,L]\to\mathcal{M}_{1}$ be a unit-speed $\mathcal{M}_{1}$-geodesic with $\gamma(0)=x$, viewed as a curve in $\mathcal{M}_{2}$. Define its Cartan development $\tilde\gamma\colon[0,L]\to T_x\mathcal{M}_{2}$ by
    \begin{equation}\label{eqn:Cartan-development}
        \tilde\gamma(s):=\int_0^s P_{0\hookleftarrow \tau}^{\gamma}[\gamma'(\tau)]\,d\tau,
    \end{equation}
    where $P_{0\hookleftarrow \tau}^{\gamma}\colon T_{\gamma(\tau)}\mathcal{M}_{2}\to T_x\mathcal{M}_{2}$ denotes parallel transport along $\gamma$. If $L\leq\pi/\kappa$, then
    \begin{equation*}
        \frac{2}{\kappa}\sin\left(\frac{\kappa L}{2}\right)
        \leq \|\tilde\gamma(L)\|_{T_x\mathcal{M}_{2}}.
    \end{equation*}
\end{lemma}

Second, we compare the endpoint of the Cartan development with the ambient logarithm. The proof is deferred to \Cref{sec:proof-of-discrepancy-from-cartan-development}.

\begin{lemma}\label{lem:discrepancy-from-cartan-development}
    Assume Assumption~\ref{ass:embedded-submanifold}\textup{(a)} holds. Let $\gamma\colon[0,L]\to\mathcal{M}_{2}$ be a smooth unit-speed curve from $x=\gamma(0)$ to $y=\gamma(L)$, and let $\tilde\gamma$ be defined by \eqref{eqn:Cartan-development}. Assume that $\gamma([0,L])$ is contained in a normal neighborhood of $x$ in $\mathcal{M}_{2}$ and, if $\kmaxExt>0$, that $L<\pi/\sqrt{\kmaxExt}$. Then
    \begin{equation*}
        \|\log_x^{\mathcal{M}_{2}}(y)-\tilde\gamma(L)\|_{T_x\mathcal{M}_{2}}
        \leq \frac{5\kmagExt C_J(L)}{18}L^3,
\qquad \text{where} \quad
        C_J(L):=\max\left\{1,\frac{A_{\kminExt}(L)}{B_{\kmaxExt}(L)}\right\},
    \end{equation*}
    with
    \begin{align*}
        A_{\kminExt}(L)
        :={}
        \begin{cases}
            1, & \kminExt\geq 0,\\
            \frac{\sinh(\sqrt{-\kminExt}L)}{\sqrt{-\kminExt}L}, & \kminExt<0,
        \end{cases} \quad \text{and} \quad
        B_{\kmaxExt}(L)
        :={}
        \begin{cases}
            1, & \kmaxExt\leq 0,\\
            \frac{\sin(\sqrt{\kmaxExt}L)}{\sqrt{\kmaxExt}L}, & \kmaxExt>0.
        \end{cases}
    \end{align*}
\end{lemma}

If $\mathcal{M}_{2}$ is a Hilbert space, then $\kmagExt=0$, the discrepancy vanishes, and the classical Euclidean argument of \cite{bernstein2000graph} is recovered.

\begin{proof}[Proof of Lemma~\ref{lem:local-chord-arc-lower-bound}]
    Since $y$ lies in a normal neighborhood of $x$,
    \begin{align*}
        d_{\mathcal{M}_{2}}(x,y)
        &=\|\log_x^{\mathcal{M}_{2}}(y)\|_{T_x\mathcal{M}_{2}} \\
        &\geq \|\tilde\gamma(L)\|_{T_x\mathcal{M}_{2}}
        -\|\log_x^{\mathcal{M}_{2}}(y)-\tilde\gamma(L)\|_{T_x\mathcal{M}_{2}} \\
        &\geq \frac{2}{\kappa}\sin\left(\frac{\kappa L}{2}\right)
        -\frac{5\kmagExt C_J(L)}{18}L^3
    \end{align*}
    where the second inequality follows from Lemma~\ref{lem:cartan-development-chord} and~\ref{lem:discrepancy-from-cartan-development}. Applying $\sin u\geq u-u^3/6$ with $u=\kappa L/2$ gives
    \begin{equation*}
        d_{\mathcal{M}_{2}}(x,y)
        \geq \left(1-\left[\frac{\kappa^2}{24}
        +\frac{5\kmagExt C_J(L)}{18}\right]L^2\right)L.
    \end{equation*}
\end{proof}

\subsubsection{Chord--Arc on Cartan Development: Proof of \texorpdfstring{Lemma~\ref{lem:cartan-development-chord}}{Lemma: Cartan Development Chord}}\label{sec:proof-of-cartan-development-chord}
\begin{proof}
    We have $\tilde\gamma'(s)=P_{0\hookleftarrow s}^{\gamma}[\gamma'(s)]$ and $\|\tilde\gamma'(s)\|_{T_x\mathcal{M}_{2}} = 1,$ by \eqref{eqn:Cartan-development} and the isometry of parallel transport.
    Differentiating a vector field after parallel transport gives $\tilde\gamma''(s)
        =P_{0\hookleftarrow s}^{\gamma}
        [\nabla_{\gamma'(s)}^{\mathcal{M}_{2}}\gamma'(s)].$
    Since $\gamma$ is a geodesic in $\mathcal{M}_{1}$, the Gauss formula yields $\nabla_{\gamma'}^{\mathcal{M}_{2}}\gamma'=\Pi^{\mathcal{M}_{1}}(\gamma',\gamma'),$ giving
    \begin{equation*}
        \|\widetilde\gamma''(s)\|_{T_x\mathcal M_2}
=\|\Pi_{\gamma(s)}(\gamma'(s),\gamma'(s))\|_{T_{\gamma(s)}\mathcal M_2}
\le\kappa.
    \end{equation*}
    The inequality follows from Assumption~\ref{ass:embedded-submanifold} implying \eqref{eqn:second-fund-bound}.

    The curve $s\mapsto\tilde\gamma'(s)$ lies on the unit sphere of $T_x\mathcal{M}_{2}$. Therefore, for $0\leq s\leq t\leq L$, 
    \begin{equation*}
        \arccos\left(\left\langle\tilde\gamma'(s),\tilde\gamma'(t)\right\rangle_{T_x\mathcal{M}_{2}}\right)
        \leq\int_s^t\|\tilde\gamma''(\tau)\|_{T_x\mathcal{M}_{2}}\,d\tau
        \leq\kappa(t-s).
    \end{equation*}
    Now, $\left\langle\tilde\gamma'(s),\tilde\gamma'(t)\right\rangle_{T_x\mathcal{M}_{2}} \geq\cos(\kappa|t-s|)$ as $\kappa L\leq\pi$ and cosine decreases on $[0,\pi]$, implying
    \begin{align*}
        \|\tilde\gamma(L)\|_{T_x\mathcal{M}_{2}}^2
        &=\int_0^L\int_0^L
        \left\langle\tilde\gamma'(s),\tilde\gamma'(t)\right\rangle_{T_x\mathcal{M}_{2}}\,ds\,dt \\
        &\geq\int_0^L\int_0^L\cos(\kappa|t-s|)\,ds\,dt \\
        &=\frac{4}{\kappa^2}\sin^2\left(\frac{\kappa L}{2}\right).
    \end{align*}
    Taking square roots completes the proof.
\end{proof}

\subsubsection{Discrepancy from Cartan Development: Proof of Lemma~\ref{lem:discrepancy-from-cartan-development}}\label{sec:proof-of-discrepancy-from-cartan-development}
\paragraph{Geometric preliminaries and notation.}
Throughout this subsection, $P_{a\hookleftarrow b}^{\alpha}$ denotes parallel transport along a curve $\alpha$ from parameter $b$ to parameter $a$. We use the curvature convention
\begin{equation*}
    R(U,V)W
    :=\nabla_U\nabla_VW-\nabla_V\nabla_UW-\nabla_{[U,V]}W,
\end{equation*}
so a Jacobi field $J$ along a geodesic with velocity $X$ satisfies $\nabla_X\nabla_XJ+R(J,X)X=0$.

Let $x=\gamma(0)$ and set
\begin{equation*}
    v(s):=\log_x^{\mathcal{M}_{2}}(\gamma(s)),
    \qquad
    \Sigma(r,s):=\exp_x^{\mathcal{M}_{2}}(rv(s)),
    \qquad (r,s)\in[0,1]\times[0,L].
\end{equation*}
For fixed $s$, write $\sigma_s(r):=\Sigma(r,s)$ and
\begin{equation*}
    a(s):=\|v(s)\|_{T_x\mathcal{M}_{2}}
    =d_{\mathcal{M}_{2}}(x,\gamma(s)).
\end{equation*}
Then $\sigma_s$ is a geodesic of constant speed $a(s)$, and $a(s)\leq s$ because $\gamma$ is unit-speed. Define the variation fields
\begin{equation*}
    X(r,s):=\partial_r\Sigma(r,s),
    \qquad
    J(r,s):=\partial_s\Sigma(r,s).
\end{equation*}
For each fixed $s$, $J(\cdot,s)$ is a Jacobi field along $\sigma_s$, and
\begin{equation*}
    J(0,s)=0,
    \qquad
    J(1,s)=\gamma'(s),
    \qquad
    \nabla_rJ(0,s)=\nabla_sX(0,s)=v'(s).
\end{equation*}
Moreover, $\|X(r,s)\|=a(s)\leq s$ for every $(r,s)$.

For $c\in\mathbb{R}$, define the comparison function
\begin{equation}\label{eqn:comparison-function}
    s_c(t):=
    \begin{cases}
        t, & c=0,\\
        \dfrac{1}{\sqrt{c}}\sin(\sqrt{c}t), & c>0,\\
        \dfrac{1}{\sqrt{-c}}\sinh(\sqrt{-c}t), & c<0.
    \end{cases}
\end{equation}

\paragraph{Auxiliary lemmas.}
We first record a curvature-tensor norm bound.

\begin{lemma}[Consequence of Berger's estimate; cf. \cite{Karcher1970BergerCurvatureTensorEstimates}]\label{lem:curvature-tensor-norm-bound}
    Suppose $K_{\min}\leq\operatorname{sec}_{\mathcal{M}}\leq K_{\max}$ and set $\widetilde K:=\max\{|K_{\min}|,|K_{\max}|\}$. Then, for all $A,B,C\in T_p\mathcal{M}$,
    \begin{equation}\label{eqn:adapted-berger-curvature-tensor-estimate}
        \|R(A,B)C\|_p
        \leq 2\widetilde K\|A\|_p\|B\|_p\|C\|_p.
    \end{equation}
\end{lemma}

\begin{proof}
    Let $\operatorname{Rm}_p(A,B,C,D):=\langle R(A,B)C,D\rangle_p$. Equation~9 of \cite{Karcher1970BergerCurvatureTensorEstimates} gives
    \begin{equation*}
        \|\operatorname{Rm}_p\|
        \leq\left(\widetilde K^2
        +\frac{25}{36}(K_{\max}-K_{\min})^2\right)^{1/2}
        \leq\frac{\sqrt{34}}{3}\widetilde K
        \leq 2\widetilde K.
    \end{equation*}
    Therefore,
    \begin{align*}
        \|R(A,B)C\|_p
        &=\sup_{\|D\|_p=1}|\operatorname{Rm}_p(A,B,C,D)| 
        \leq 2\widetilde K\|A\|_p\|B\|_p\|C\|_p.
    \end{align*}
\end{proof}

The next lemma controls the radial variation field.

\begin{lemma}\label{lem:Jacobi-field-control}
    Under the assumptions of Lemma~\ref{lem:discrepancy-from-cartan-development}, for all $(r,s) \in [0,1] \times [0,L]$
    \begin{equation}
        \begin{gathered}
            \|J(r,s)\|_{\Sigma(r,s)}\leq C_J(L)r,
            \qquad \text{where} \quad
            C_J(L):=\max\left\{1,
            \frac{A_{\kminExt}(L)}{B_{\kmaxExt}(L)}\right\}.
        \end{gathered}
    \end{equation}
    Here
    \begin{equation}\label{eqn:constants-for-Jacobi-control}
        \begin{aligned}
            A_{\kminExt}(L)
            :={}
            \begin{cases}
                1, & \kminExt\geq0,\\
                \frac{\sinh(\sqrt{-\kminExt}L)}{\sqrt{-\kminExt}L}, & \kminExt<0,
            \end{cases} \quad \text{and} \quad
            B_{\kmaxExt}(L)
            :={}
            \begin{cases}
                1, & \kmaxExt\leq0,\\
                \frac{\sin(\sqrt{\kmaxExt}L)}{\sqrt{\kmaxExt}L}, & \kmaxExt>0.
            \end{cases}
        \end{aligned}
    \end{equation}
\end{lemma}

\begin{proof}
    Fix $s\in[0,L]$ and abbreviate $a=a(s)$. If $a=0$, then $D\exp_x(0)=I$ gives $J(r,s)=r\gamma'(s)$, so the claim follows from $C_J(L)\geq1$. Assume henceforth that $a>0$.

    Reparametrize the radial geodesic by arclength: for $t \in [0,a]$
    \begin{equation*}
        \bar\sigma_s(t):=\sigma_s(t/a),
        \qquad \text{and} \qquad
        \bar J(t,s):=J(t/a,s).
    \end{equation*}
    Then $\bar\sigma_s$ is unit-speed, and
    \begin{equation}\label{eqn:endpoints-reparameterized-jacobi}
        \bar J(0,s)=0,
        \qquad
        \bar J(a,s)=\gamma'(s).
    \end{equation}
    Decompose $\bar J=\bar J^{\top}+\bar J^{\perp}$ and $\gamma'(s)=\gamma'(s)^{\top}+\gamma'(s)^{\perp}$ into components tangent and orthogonal to $\dot{\bar\sigma}_s$.

    \paragraph{Tangential Component to $\dot{\bar{\sigma}}_{s}$.} By the geodesic equation, the tangential component of a Jacobi field along a unit-speed geodesic is affine in arclength. Thus, for a scalar $\lambda_s$ satisfying
    \begin{equation}\label{eqn:tangential-curve-speed}
        \gamma'(s)^{\top}=\lambda_s\dot{\bar\sigma}_s(a),
    \end{equation}
    the endpoint conditions give
    \begin{equation*}
        \bar J^{\top}(t,s)
        =\frac{t}{a}P_{t\hookleftarrow a}^{\bar\sigma_s}[\gamma'(s)^{\top}].
    \end{equation*}
    Consequently,
    \begin{equation}\label{eqn:bound-on-tangential-jacobi}
        \|J^{\top}(r,s)\|
        =r\|\gamma'(s)^{\top}\|.
    \end{equation}

    \paragraph{Normal Component to $\dot{\bar{\sigma}}_{s}$.} By the Jacobi Comparison Theorem (Theorem 11.9) \cite{lee2018introduction},
    \begin{equation*}
        \|\bar J^{\perp}(t,s)\|
        \leq s_{\kminExt}(t)\|\nabla_t\bar J^{\perp}(0,s)\|,
        \quad \text{and} \quad 
        \|\bar J^{\perp}(a,s)\|
        \geq s_{\kmaxExt}(a)\|\nabla_t\bar J^{\perp}(0,s)\|.
    \end{equation*}
    The denominator is positive because $a\leq L<\pi/\sqrt{\kmaxExt}$ when $\kmaxExt>0$. Hence
    \begin{equation*}
        \|J^{\perp}(r,s)\|
        \leq \frac{s_{\kminExt}(ar)}{s_{\kmaxExt}(a)}
        \|\gamma'(s)^{\perp}\|.
    \end{equation*}
    Moreover,
    \begin{align*}
        \frac{s_{\kminExt}(ar)}{s_{\kmaxExt}(a)}
        &=r\frac{s_{\kminExt}(ar)}{ar}\frac{a}{s_{\kmaxExt}(a)} 
        \leq r\frac{A_{\kminExt}(L)}{B_{\kmaxExt}(L)}.
    \end{align*}
    Indeed, $s_{\kminExt}(z)/z$ is bounded above by $A_{\kminExt}(L)$ on $(0,L]$, while $s_{\kmaxExt}(z)/z$ is bounded below by $B_{\kmaxExt}(L)$ on the same interval. Thus
    \begin{equation}\label{eqn:bound-on-normal-jacobi}
        \|J^{\perp}(r,s)\|
        \leq r\frac{A_{\kminExt}(L)}{B_{\kmaxExt}(L)}
        \|\gamma'(s)^{\perp}\|.
    \end{equation}
    Combining \eqref{eqn:bound-on-tangential-jacobi} and \eqref{eqn:bound-on-normal-jacobi}, and using $\|\gamma'(s)\|=1$, gives
    \begin{equation*}
        \|J(r,s)\| \leq r \sqrt{\|\gamma'(s)^{\top}\|^{2} + \left(\frac{A_{\kminExt}(L)}{B_{\kmaxExt}(L)}\right)^{2} \|\gamma'(s)^{\perp}\|^{2}}
        \leq C_J(L)r.
    \end{equation*}
\end{proof}

\begin{proof}[Proof of Lemma~\ref{lem:discrepancy-from-cartan-development}]
Define the radially transported surrogate
\begin{equation*}
    \widehat\gamma(L)
    :=\int_0^L P_{0\hookleftarrow1}^{\sigma_s}[\gamma'(s)]\,ds.
\end{equation*}
We will establish the following two identities. First,
\begin{equation}\label{eqn:discrepancy-identity-as-integrating-curvature}
    \log_x^{\mathcal{M}_{2}}(y)-\widehat\gamma(L)
    =\int_0^L\int_0^1(1-r)
    P_{0\hookleftarrow r}^{\sigma_s}
    [R^{\mathcal{M}_{2}}(J(r,s),X(r,s))X(r,s)]\,dr\,ds.
\end{equation}
Second, define $H_s:=P_{0\hookleftarrow s}^{\gamma}\circ P_{1\hookleftarrow0}^{\sigma_s} \colon T_x\mathcal{M}_{2}\to T_x\mathcal{M}_{2}.$ Then, for every $u\in T_x\mathcal{M}_{2}$,
\begin{equation}\label{eqn:holonomy-identity-as-integrating-curvature}
    (H_s-I)u
    =\int_0^s\int_0^1
    P_{0\hookleftarrow t}^{\gamma}\circ P_{1\hookleftarrow r}^{\sigma_t}
    \left[R^{\mathcal{M}_{2}}(X(r,t),J(r,t))
    P_{r\hookleftarrow0}^{\sigma_t}u\right]dr\,dt.
\end{equation}

Before establishing \eqref{eqn:discrepancy-identity-as-integrating-curvature} and \eqref{eqn:holonomy-identity-as-integrating-curvature}, we highlight how the desired bound is recovered. Setting $v_s:=P_{0\hookleftarrow1}^{\sigma_s}[\gamma'(s)],$ we have that $\|v_s\|=1$ and
\begin{equation*}
    \tilde\gamma(L)-\widehat\gamma(L)
    =\int_0^L(H_s-I)v_s\,ds.
\end{equation*}
By \eqref{eqn:discrepancy-identity-as-integrating-curvature}, Lemma~\ref{lem:curvature-tensor-norm-bound}, Lemma~\ref{lem:Jacobi-field-control}, and $\|X(r,s)\|=a(s)\leq s$,
\begin{align*}
    \|\log_x^{\mathcal{M}_{2}}(y)-\widehat\gamma(L)\|
    &\leq 2\kmagExt C_J(L)
    \int_0^L s^2\int_0^1r(1-r)\,dr\,ds 
    =\frac{\kmagExt C_J(L)}{9}L^3.
\end{align*}
Similarly, \eqref{eqn:holonomy-identity-as-integrating-curvature} and the isometry of parallel transport give
\begin{align*}
    \|\tilde\gamma(L)-\widehat\gamma(L)\|
    &\leq 2\kmagExt C_J(L)
    \int_0^L\int_0^s\int_0^1 rt\,dr\,dt\,ds =\frac{\kmagExt C_J(L)}{6}L^3.
\end{align*}
Consequently, after a triangle inequality the stated bound is proven:
\begin{align*}
    \|\log_x^{\mathcal{M}_{2}}(y)-\tilde\gamma(L)\| \leq  \|\log_x^{\mathcal{M}_{2}}(y)-\widehat \gamma(L)\| +  \|\widehat \gamma(L)-\tilde\gamma(L)\|
    &\leq \frac{5\kmagExt C_J(L)}{18}L^3.
\end{align*}

\paragraph{Proof of Identity \eqref{eqn:discrepancy-identity-as-integrating-curvature}.} Fix $s\in[0,L]$ and define
\begin{equation*}
    U_s(r):=P_{0\hookleftarrow r}^{\sigma_s}[J(r,s)],
    \qquad
    F_s(r):=P_{0\hookleftarrow r}^{\sigma_s}
    [R^{\mathcal{M}_{2}}(J(r,s),X(r,s))X(r,s)].
\end{equation*}
Differentiation after parallel transport and the Jacobi equation give
\begin{equation*}
    U_s''(r)=-F_s(r),
    \qquad U_s(0)=0.
\end{equation*}
Thus
\begin{equation*}
    U_s(1)=U_s'(0)-\int_0^1(1-r)F_s(r)\,dr.
\end{equation*}
Now $U_s(1)=P_{0\hookleftarrow1}^{\sigma_s}[\gamma'(s)]$, while torsion-freeness gives
\begin{equation*}
    U_s'(0)=\nabla_rJ(0,s)=\nabla_sX(0,s)=v'(s) = \frac{d}{ds} \log_{x}^{\mathcal{M}_{2}}(\gamma(s)).
\end{equation*}
Therefore, by rearrangement
\begin{equation}\label{eqn:pointwise-identity-for-derivative-of-log}
    \frac{d}{ds}\log_x^{\mathcal{M}_{2}}(\gamma(s))
    =P_{0\hookleftarrow1}^{\sigma_s}[\gamma'(s)]
    +\int_0^1(1-r)P_{0\hookleftarrow r}^{\sigma_s}
    [R^{\mathcal{M}_{2}}(J(r,s),X(r,s))X(r,s)]\,dr.
\end{equation}
Integrating \eqref{eqn:pointwise-identity-for-derivative-of-log} over $s\in[0,L]$ establishes \eqref{eqn:discrepancy-identity-as-integrating-curvature}.

\paragraph{Proof of Identity \eqref{eqn:holonomy-identity-as-integrating-curvature}.} Fix $u\in T_x\mathcal{M}_{2}$ and define
\begin{equation*}
    V(r,t):=P_{r\hookleftarrow0}^{\sigma_t}u.
\end{equation*}
Then $\nabla_rV=0$ and $\nabla_tV(0,t)=0$. Since the coordinate vector fields of the variation commute, the curvature identity gives
\begin{equation*}
    \nabla_r\nabla_tV=R^{\mathcal{M}_{2}}(X,J)V.
\end{equation*}
Integrating along the radial geodesic yields
\begin{equation*}
    \nabla_tV(1,t)
    =\int_0^1P_{1\hookleftarrow r}^{\sigma_t}
    [R^{\mathcal{M}_{2}}(X(r,t),J(r,t))V(r,t)]\,dr.
\end{equation*}
Because differentiation commutes with transport back along $\gamma$,
\begin{equation*}
    \frac{d}{dt}(H_tu)
    =P_{0\hookleftarrow t}^{\gamma}[\nabla_tV(1,t)].
\end{equation*}
Integrating from $0$ to $s$ and using $H_0=I$ proves \eqref{eqn:holonomy-identity-as-integrating-curvature}.
\end{proof}

\subsection{Probabilistic ISOMAP Proof: Corollary~\ref{cor:probabilistic-isomap}}\label{sec:proof-of-probabilistic-isomap}
The only random hypothesis in Corollary~\ref{cor:probabilistic-isomap} is the $\delta$-covering condition. The following standard occupancy argument verifies it.

\begin{lemma}\label{lem:sampling-for-delta-net}
    Let $S\subseteq\mathcal{M}_{1}$ be compact, and let $\{X_i\}_{i=1}^{n}$ be i.i.d. samples from a Borel probability measure $\mu_S$ supported on $S$. For $r>0$, define
    \begin{equation}\label{eqn:min-local-sampling-mass}
        q_S(r):=\inf_{x\in S}
        \mu_S\bigl(B_{\mathcal{M}_{1}}(x,r)\cap S\bigr).
    \end{equation}
    Let $\delta>0$ and suppose $q_S(\delta/2)>0$. Then, for every $\zeta\in(0,1)$, the sample $\{X_i\}_{i=1}^{n}$ is a $\delta$-cover of $S$ with probability at least $1-\zeta$, provided
    \begin{equation}\label{eqn:covering-number-to-delta-sampling}
        n\geq\frac{1}{q_S(\delta/2)}
        \left(\log\mathcal{N}(S,\delta/2)
        +\log\left(\frac{1}{\zeta}\right)\right).
    \end{equation}
\end{lemma}

\begin{proof}
    Let $N:=\mathcal{N}(S,\delta/2)$ and choose a $\delta/2$-net $z_1,\dots,z_N\in S$. For each $\ell\in\{1,\dots,N\}$, let
    \begin{equation*}
        E_\ell:=\left\{\text{some }X_i\text{ belongs to }
        B_{\mathcal{M}_{1}}(z_\ell,\delta/2)\right\}.
    \end{equation*}
    Since every such ball has $\mu_S$-mass at least $q_S(\delta/2)$,
    \begin{equation*}
        \mathbb{P}(E_\ell^c)
        \leq (1-q_S(\delta/2))^n
        \leq \exp(-nq_S(\delta/2)).
    \end{equation*}
    Hence, by the union bound,
    \begin{equation*}
        \mathbb{P}\left(\bigcap_{\ell=1}^{N}E_\ell\right)
        \geq 1-N\exp(-nq_S(\delta/2)).
    \end{equation*}

    On $\bigcap_{\ell=1}^{N}E_\ell$, fix $x\in S$ and choose $z_\ell$ with $d_{\mathcal{M}_{1}}(x,z_\ell)\leq\delta/2$. The event $E_\ell$ supplies an $X_i$ with $d_{\mathcal{M}_{1}}(z_\ell,X_i)\leq\delta/2$, and therefore $d_{\mathcal{M}_{1}}(x,X_i)\leq\delta$. Thus the sample is a $\delta$-cover. Finally, \eqref{eqn:covering-number-to-delta-sampling} implies
    \begin{equation*}
        N\exp(-nq_S(\delta/2))\leq\zeta,
    \end{equation*}
    which proves the claim.
\end{proof}

\begin{proof}[Proof of Corollary~\ref{cor:probabilistic-isomap}]
By Lemma~\ref{lem:sampling-for-delta-net} and
\eqref{eqn:probabilistic-isomap-sampling-bound},
the random sample is a $\delta$-cover of $S$ with probability
at least $1-\zeta$.
On this event, apply
Theorem~\ref{thm:deterministic-metric-equivalence-for-ISOMAP}
when Assumption~\ref{ass:embedded-submanifold}\textup{(a)} holds, or
Proposition~\ref{prop:abstract-metric-isomap}
when Assumption~\ref{ass:embedded-submanifold}\textup{(b)} holds.
This gives \eqref{eqn:probabilistic-isomap-metric-equivalence}.
\end{proof}

\section{Graph Algorithm Discrepancy from Riemannian SGD}

The following two lemmas provide an approximate equivalence between distance traveled by geodesic interpolation on $\mathcal{M}_{1}$ and graph interpolation, in terms of the metric equivalence and rounding errors. The first quantifies the effect of nearest-vertex rounding.

\begin{lemma}
\label{lem:nearest-vertex-graph-interpolation}
Let $\gamma^G_{x\to y}=(v_0,\ldots,v_m),$ with $
  v_0=x $ and $ v_m=y,$
be a shortest path on a weighted graph $G$. Suppose every edge of this path has weight $0 \leq w \leq \epsilon$. For $\eta\in[0,1]$, let
\[
  x^+\in\arg\min_{0\leq j\leq m}
  \left|d_G(x,v_j)-\eta d_G(x,y)\right|.
\]
Then, 
\begin{equation}
\label{eq:nearest-vertex-rounding}
  \left|d_G(x,x^+)-\eta d_G(x,y)\right|
  \leq \frac{\epsilon}{2}.
\end{equation}
If $\eta=0$ or $\eta=1$, the
left-hand side is zero.
\end{lemma}

\begin{proof}
If $m=0$, then $x=y=x^+$ and the claim is immediate. Assume $m\geq1$, and write $w_j$ for the weight of the edge $v_{j-1}v_j$. Define the cumulative path lengths
\[
  s_0:=0,
  \qquad
  s_j:=\sum_{i=1}^j w_i,
  \quad 1\leq j\leq m.
\]
Every prefix $(v_0,\ldots,v_j)$ is a shortest path from $x$ to $v_j$.
Indeed, if a strictly shorter $x$--$v_j$ path existed, concatenating it with
the suffix $(v_j,\ldots,v_m)$ would produce an $x$--$y$ walk strictly
shorter than $\gamma^G_{x\to y}$; deleting cycles from that walk cannot
increase its length and would contradict the shortest-path property. Hence
  $s_j=d_G(x,v_j),0\leq j\leq m,$
and $s_m=d_G(x,y)$.

Set $a:=\eta s_m\in[0,s_m]$. If $a=s_0$ or $a=s_m$, the corresponding
endpoint has zero error, proving the endpoint assertion. Otherwise there is
an index $k\in\{0,\ldots,m-1\}$ such that
$s_k\leq a\leq s_{k+1}$. Therefore
\[
 \min\{|a-s_k|,|s_{k+1}-a|\}
 \leq \frac{s_{k+1}-s_k}{2}
 =\frac{w_{k+1}}{2}
 \leq \frac{\epsilon}{2}.
\]
The minimum over all path vertices is no larger than the minimum over these
two bracketing vertices. Every element of the argmin attains that same
minimum value, so the estimate also holds under ties.
\end{proof}

The second compares the distance traveled via geodesic interpolation versus graph interpolation, accounting for metric equivalence and rounding errors of the previous lemma. 

\begin{lemma}
\label{lem:metric-equivalence-consequence}
Let $\lambda_{1}, \lambda_{2} \in (0,1)$, and suppose that 
\[
(1-\lambda_1)d_{\mathcal M_1}(x_t,v)
\le d_{G_{t+1}}(x_t,v)
\le(1+\lambda_2)d_{\mathcal M_1}(x_t,v),
\qquad
v\in\{x_{t+1},y_{t+1}\}.
\]
Let $\eta_t\in[0,1]$, let $\widetilde x_{t+1}$ be exact intrinsic interpolation, and suppose $x_{t+1}$ is selected by the nearest-vertex rule on an $\epsilon_{t+1}$-edge-bounded shortest path. Then
\begin{align}
\label{eq:graph-travel-distance-rounding}
 \frac{(1-\lambda_1) d_{\mathcal{M}_{1}}(x_{t}, \widetilde{x}_{t+1}) - \epsilon_{t+1}/2}{1+\lambda_2}
 \leq d_{\mathcal M_1}(x_t,x_{t+1})
 \leq
 \frac{(1+\lambda_2)d_{\mathcal M_1}(x_t,\widetilde x_{t+1}) + \epsilon_{t+1}/2}{1-\lambda_1}.
\end{align}
\end{lemma}

\begin{proof}
Set
 $\omega_{t+1}:= d_{G_{t+1}}(x_t,x_{t+1}) -\eta_t d_{G_{t+1}}(x_t,y_{t+1}).$
Lemma~\ref{lem:nearest-vertex-graph-interpolation} gives
$|\omega_{t+1}|\leq\epsilon_{t+1}/2$. Then,
\begin{align*}
 (1-\lambda_1)d_{\mathcal M_1}(x_t,x_{t+1})
 &\leq d_{G_{t+1}}(x_t,x_{t+1}) \\
 &\leq \eta_t d_{G_{t+1}}(x_t,y_{t+1})+|\omega_{t+1}|\\
 &\leq \eta_t(1+\lambda_2)
 d_{\mathcal M_1}(x_t,y_{t+1})+\frac{\epsilon_{t+1}}2.
\end{align*}
The first inequality follows from metric equivalence, the second from triangle inequality, and the last from metric equivalence in conjunction with the bound on $|\omega_{t+1}|$. Dividing both sides by $(1-\lambda_{1})$ recovers the desired upper bound of \eqref{eq:graph-travel-distance-rounding}. The lower bound is entirely analogous.
\end{proof}

\subsection{Proof of Lemma~\ref{lem:conditions-on-metric-equivalence-for-tube-containment}}\label{sec:proof-of-conditions-on-metric-equivalence-for-tube-containment}

\begin{proof}
For notational simplicity, suppress the iteration indices and write
\begin{equation*}
    \gamma:=\gamma_{x\rightarrow y}^{\mathcal{M}_{1}},
    \qquad
    \ell:=d_{\mathcal{M}_{1}}(x,y),
    \qquad r := \alpha \ell, \qquad 
    \tilde r:=r+\rho_{\ROI}(\epsilon).
\end{equation*}
Let $\mathcal{P}_{x\rightarrow y}^{G}$ be a graph path containing a vertex outside
$T_{r}(\gamma)$. We prove that
\begin{equation}\label{eqn:inequality-for-tube-containment}
    (1+\lambda_{2})d_{\mathcal{M}_{1}}(x,y)
    <
    \ell(\mathcal{P}_{x\rightarrow y}^{G}).
\end{equation}
If so, $\mathcal{P}_{x \rightarrow y}^{G}$ is not a graph geodesic as metric equivalence gives
    $d_{G}(x,y)
    \leq
    (1+\lambda_{2})d_{\mathcal{M}_{1}}(x,y).$

Let $z$ be the first vertex of $\mathcal{P}_{x\rightarrow y}^{G}$ outside
$T_{r}(\gamma)$, and let $u$ be its predecessor. By the graph edge-length bound,
$d_{\mathcal{M}_{2}}(u,z)\leq\epsilon$. Since $u\in T_{r}(\gamma)$,
\begin{align*}
    \operatorname{dist}_{\mathcal{M}_{1}}(z,\gamma)
    \leq
    d_{\mathcal{M}_{1}}(z,u)
    +
    \operatorname{dist}_{\mathcal{M}_{1}}(u,\gamma) 
    \leq
    \rho_{\ROI}(\epsilon)+r
    =
    \tilde r.
\end{align*}
Thus $z\in T_{\tilde r}(\gamma)$.

Choose a minimizing $\mathcal{M}_{1}$-geodesic from $x$ to $z$, and let $q$ be
its first point on $\partial T_{r}(\gamma)$. Then
\begin{equation*}
    d_{\mathcal{M}_{1}}(x,z)
    =
    d_{\mathcal{M}_{1}}(x,q)
    +
    d_{\mathcal{M}_{1}}(q,z),
\end{equation*}
and the triangle inequality gives
\begin{equation*}
    d_{\mathcal{M}_{1}}(x,q)
    +
    d_{\mathcal{M}_{1}}(q,y)
    \leq
    d_{\mathcal{M}_{1}}(x,z)
    +
    d_{\mathcal{M}_{1}}(z,y).
\end{equation*}
Therefore,
\begin{align}
    \ell(\mathcal{P}_{x\rightarrow y}^{G})
    &=
    \ell(\mathcal{P}_{x\rightarrow z}^{G})
    +
    \ell(\mathcal{P}_{z\rightarrow y}^{G})
    \nonumber\\
    &\geq
    d_{G}(x,z)+d_{G}(z,y)
    \nonumber\\
    &\geq
    (1-\lambda_{1})
    \left(
        d_{\mathcal{M}_{1}}(x,z)
        +
        d_{\mathcal{M}_{1}}(z,y)
    \right)
    \nonumber\\
    &\geq
    (1-\lambda_{1})
    \left(
        d_{\mathcal{M}_{1}}(x,q)
        +
        d_{\mathcal{M}_{1}}(q,y)
    \right).
    \label{eqn:lower-bound-graph-length-with-tube-boundary-and-back}
\end{align}
The second inequality uses metric equivalence on
$T_{\tilde r}(\gamma)$.

Choose a metric projection $p\in\gamma$ of $q$ onto $\gamma$. Since
$q\in\partial T_{r}(\gamma)$, we have $d_{\mathcal{M}_{1}}(p,q)=r.$
Set
\begin{equation*}
    a:=d_{\mathcal{M}_{1}}(x,p)\in[0,\ell].
\end{equation*}
Then $d_{\mathcal{M}_{1}}(p,y)=\ell-a$. The first-variation condition for a
closest point on a geodesic segment shows that the relevant hinge angle at $p$
is at least $\pi/2$. If $p$ lies in the interior of $\gamma$, both hinge angles
are equal to $\pi/2$; if $p$ is an endpoint, the nondegenerate hinge angle is at
least $\pi/2$. Before continuing, we note that the localization conditions of Assumption~\ref{ass:region-of-interest} ensures that the following comparison triangles lie in the admissible normal-neighborhood and curvature regimes.

Fix $K>0$ such that
$\operatorname{sec}_{\mathcal{M}_{1}}\leq K$ and the spherical comparison
triangles below are defined. The upper-curvature hinge comparison and the
spherical law of cosines yield
\begin{align*}
    d_{\mathcal{M}_{1}}(x,q)
    &\geq
    \frac{1}{\sqrt K}
    \arccos\left(
        \cos(\sqrt K\,a)\cos(\sqrt K\,r)
    \right),\\
    d_{\mathcal{M}_{1}}(q,y)
    &\geq
    \frac{1}{\sqrt K}
    \arccos\left(
        \cos(\sqrt K(\ell-a))\cos(\sqrt K\,r)
    \right).
\end{align*}
Define
\begin{equation*}
    \phi_{K}(s)
    :=
    \frac{1}{\sqrt K}
    \arccos\left(
        \cos(\sqrt K\,r)\cos(\sqrt K\,s)
    \right)
\end{equation*}
and
\begin{equation*}
    F_{K}(a):=\phi_{K}(a)+\phi_{K}(\ell-a).
\end{equation*}
The function $F_{K}$ is symmetric about $\ell/2$. Moreover,
\begin{equation*}
    \phi_{K}'(s)
    =
    \frac{
        \cos(\sqrt K\,r)\sin(\sqrt K\,s)
    }{
        \sqrt{
            1-\cos^{2}(\sqrt K\,r)\cos^{2}(\sqrt K\,s)
        }
    }.
\end{equation*}
When $\sqrt K\,\ell<\pi$, for $0\leq a<\ell/2$ we have
\begin{equation*}
    \sin(\sqrt K\,a)
    <
    \sin(\sqrt K(\ell-a)).
\end{equation*}
Since
\begin{equation*}
    u\longmapsto
    \frac{\cos(\sqrt K\,r)u}
    {\sqrt{1-\cos^{2}(\sqrt K\,r)+\cos^{2}(\sqrt K\,r)u^{2}}}
\end{equation*}
is nondecreasing on $[0,1]$, it follows that $F_{K}'(a)\leq0$ on
$[0,\ell/2]$. By symmetry, $F_{K}$ attains its minimum at $\ell/2$. Hence
\begin{equation}\label{eqn:lower-bound-for-tube-containment}
    d_{\mathcal{M}_{1}}(x,q)
    +
    d_{\mathcal{M}_{1}}(q,y)
    \geq
    \frac{2}{\sqrt K}
    \arccos\left(
        \cos\left(\frac{\sqrt K\,\ell}{2}\right)
        \cos(\sqrt K\,r)
    \right).
\end{equation}

Suppose first that $\kmaxInt>0$. Taking $K=\kmaxInt$ in
\eqref{eqn:lower-bound-for-tube-containment} and using
\eqref{eqn:bound-on-metric-equivalence-parameters-for-tube-containment} gives
\begin{equation*}
    d_{\mathcal{M}_{1}}(x,q)
    +
    d_{\mathcal{M}_{1}}(q,y)
    >
    \frac{1+\lambda_{2}}{1-\lambda_{1}}\ell.
\end{equation*}

Suppose next that $\kmaxInt\leq0$. Since
$\operatorname{sec}_{\mathcal{M}_{1}}\leq K$ for every $K>0$, we may let
$K\downarrow0$ in \eqref{eqn:lower-bound-for-tube-containment}. Using
\begin{equation*}
    \lim_{K\downarrow0}
    \frac{2}{\sqrt K}
    \arccos\left(
        \cos\left(\frac{\sqrt K\,\ell}{2}\right)
        \cos(\sqrt K\,r)
    \right)
    =
    \sqrt{\ell^{2}+4r^{2}},
\end{equation*}
we obtain
\begin{equation}\label{eqn:negative-curvature-case-for-tube-containment}
    d_{\mathcal{M}_{1}}(x,q)
    +
    d_{\mathcal{M}_{1}}(q,y)
    \geq
    \sqrt{\ell^{2}+4r^{2}}
    =
    \ell\sqrt{1+4\left(\frac r\ell\right)^{2}}.
\end{equation}
The nonpositive-curvature branch of
\eqref{eqn:bound-on-metric-equivalence-parameters-for-tube-containment}
therefore again gives
\begin{equation*}
    d_{\mathcal{M}_{1}}(x,q)
    +
    d_{\mathcal{M}_{1}}(q,y)
    >
    \frac{1+\lambda_{2}}{1-\lambda_{1}}\ell.
\end{equation*}

In either case, with \eqref{eqn:lower-bound-graph-length-with-tube-boundary-and-back} in mind the lower bounds on $d_{\mathcal{M}_{1}}(x,q) + d_{\mathcal{M}_{1}}(q,y)$ implies
\begin{equation*}
    \ell(\mathcal{P}_{x\rightarrow y}^{G})
    >
    (1+\lambda_{2})\ell,
\end{equation*}
which is \eqref{eqn:inequality-for-tube-containment}. Thus no graph geodesic
from $x$ to $y$ can leave $T_{r}(\gamma)$ completing the proof.
\end{proof}

\subsection{Proof of Lemma~\ref{lem:comparison-with-quadratic}} \label{sec:appendix-proof-of-comparison-with-quadratic}
\begin{proof}
    First, when $\kmaxInt \leq 0$, we have
    \begin{equation}\label{eqn:useful-identity-for-comparison-control}
        \widetilde{\Psi}(\ell, \alpha \ell) - 1 = \sqrt{1 + 4\alpha^{2}} - 1 = \frac{4\alpha^{2}}{1 + \sqrt{1+4\alpha^{2}}}.
    \end{equation}
    Then, it immediately follows since $\alpha \in (0,1]$ that
    \begin{equation*}
        \frac{4\alpha^{2}}{1 + \sqrt{5}} \leq \widetilde{\Psi}(\ell, \alpha \ell) - 1 \leq 2\alpha^{2}.
    \end{equation*} 
    
    Second, when $\kmaxInt > 0$, set 
    \begin{equation*}
        k := \sqrt{\kmaxInt}, \quad u := \frac{k\ell}{2}, \quad \text{and} \quad U := \frac{k \ROIdiam}{2},
    \end{equation*}
    where $U < \pi/4$ by Assumption~\ref{ass:region-of-interest}. Let $u > 0$. We first note that since $\alpha \in (0,1]$, we have 
    \begin{equation}
        \sin(2\alpha u) \geq 2 \alpha u \operatorname{sinc}(2U), \quad \text{and} \quad \sin u \leq \sqrt{\sin^{2} u + \cos^{2} u \sin^{2}(2\alpha u)} \leq u\sqrt{1+4\alpha^{2}}.\label{eqn:before-derivative-control}
    \end{equation}
    Then, by direct differentiation we have the following
    \begin{equation*}
        \partial_\alpha \widetilde{\Psi}(\ell, \alpha \ell) = \frac{2\cos u \sin(2\alpha u)}{\sqrt{\sin^{2}u + \cos^{2} u \sin^{2}(2\alpha u)}}
    \end{equation*}
    and it follows directly from \eqref{eqn:before-derivative-control} that we can upper and lower bound this partial derivative by
    \begin{equation*}
         \frac{4\alpha \cos U \operatorname{sinc}(2U)}{\sqrt{1+4\alpha^{2}}} \leq \frac{2\cos u \sin(2\alpha u)}{\sqrt{\sin^{2}u + \cos^{2} u \sin^{2}(2\alpha u)}} \leq 4\alpha u \cot u \leq 4 \alpha.
    \end{equation*}
    Integrating from $0$ to $\alpha$, and using $\widetilde{\Psi}(\ell,0) = 1$, it follows that 
    \begin{equation*}
        \cos(U) \operatorname{sinc}(2U)(\sqrt{1 + 4\alpha^{2}}-1) \leq \widetilde{\Psi}(\ell, \alpha \ell) - 1 \leq 2\alpha^{2}.
    \end{equation*}
    We can further lower bound the lower bound by using the identity \eqref{eqn:useful-identity-for-comparison-control} and $\alpha \in (0, 1]$,
    \begin{equation*}
        \cos(U) \operatorname{sinc}(2U)  \frac{4\alpha^{2}}{1+\sqrt{5}} \leq \cos(U) \operatorname{sinc}(2U) (\sqrt{1 + 4\alpha^{2}} - 1) \leq \widetilde{\Psi}(\ell, \alpha \ell) - 1.
    \end{equation*}
\end{proof}

\subsection{Proof of Lemma~\ref{lem:discrepancy-from-geodesic-interpolation}}\label{sec:proof-of-discrepancy-from-geodesic-interpolation}

\begin{proof}
For brevity, set
\[
 \ell_{t+1}:=d_{\mathcal M_1}(x_t,y_{t+1}),
 \qquad
 \widetilde d:=d_{\mathcal M_1}(x_t,\widetilde x_{t+1})=\eta_t\ell_{t+1},
 \qquad
 r_{t+1}:=\alpha_{t+1}\ell_{t+1}.
\]
By tube containment, $x_{t+1}\in\mathcal T_{r_{t+1}}(\gamma^{\mathcal M_1}_{x_t\to y_{t+1}})$. Hence there exists $p\in\gamma^{\mathcal M_1}_{x_t\to y_{t+1}}$ such that $d_{\mathcal M_1}(x_{t+1},p)\leq r_{t+1}$. Put $s:=d_{\mathcal M_1}(x_t,p)$. Since $p$ and $\widetilde x_{t+1}$ lie on the same minimizing geodesic segment, $d_{\mathcal M_1}(p,\widetilde x_{t+1})=|s-\widetilde d|$, and therefore
\begin{equation}
\label{eq:rounding-interpolation-triangle}
 d_{\mathcal M_1}(x_{t+1},\widetilde x_{t+1})
 \leq r_{t+1}+|s-\widetilde d|.
\end{equation}

The upper estimate of \eqref{eq:graph-travel-distance-rounding} in Lemma~\ref{lem:metric-equivalence-consequence} gives
\begin{align*}
 s-\widetilde d
 &\leq d_{\mathcal M_1}(x_t,x_{t+1})
      +d_{\mathcal M_1}(x_{t+1},p)-\widetilde d\\
 &\leq
 \left(
 \frac{1+\lambda_2^{(t+1)}}{1-\lambda_1^{(t+1)}}-1
 \right)\widetilde d
 +r_{t+1}
 +\frac{\epsilon_{t+1}}{2(1-\lambda_1^{(t+1)})}\\
 &=
 \left(
 \frac{\lambda_1^{(t+1)}+\lambda_2^{(t+1)}}
 {1-\lambda_1^{(t+1)}}\eta_t
 +\alpha_{t+1}
 \right)\ell_{t+1}
 +\frac{\epsilon_{t+1}}{2(1-\lambda_1^{(t+1)})}.
\end{align*}
The lower estimate of \eqref{eq:graph-travel-distance-rounding} in Lemma~\ref{lem:metric-equivalence-consequence} gives
\begin{align*}
 \widetilde d-s
 &\leq \widetilde d-d_{\mathcal M_1}(x_t,x_{t+1})
      +d_{\mathcal M_1}(x_{t+1},p)\\
 &\leq
 \left(
 1-\frac{1-\lambda_1^{(t+1)}}{1+\lambda_2^{(t+1)}}
 \right)\widetilde d
 +r_{t+1}
 +\frac{\epsilon_{t+1}}{2(1+\lambda_2^{(t+1)})}\\
 &=
 \left(
 \frac{\lambda_1^{(t+1)}+\lambda_2^{(t+1)}}
 {1+\lambda_2^{(t+1)}}\eta_t
 +\alpha_{t+1}
 \right)\ell_{t+1}
 +\frac{\epsilon_{t+1}}{2(1+\lambda_2^{(t+1)})}.
\end{align*}
Because $1-\lambda_1^{(t+1)}<1+\lambda_2^{(t+1)}$, the first of these two upper bounds dominates the second. Thus
\[
 |s-\widetilde d|
 \leq
 \left(
 \frac{\lambda_1^{(t+1)}+\lambda_2^{(t+1)}}
 {1-\lambda_1^{(t+1)}}\eta_t
 +\alpha_{t+1}
 \right)\ell_{t+1}
 +\frac{\epsilon_{t+1}}{2(1-\lambda_1^{(t+1)})}.
\]
Combining this inequality with \eqref{eq:rounding-interpolation-triangle} and $r_{t+1}=\alpha_{t+1}\ell_{t+1}$ proves
\eqref{eqn:discrepancy-from-geodesic-interpolation}.
\end{proof}

\subsection{Proof of Lemma~\ref{lem:metric-equivalence-gives-second-moment-bounds}}\label{sec:proof-of-metric-equivalence-gives-second-moment-bounds}

\begin{proof}
On $E_{t+1}$, the upper estimate in
Lemma~\ref{lem:metric-equivalence-consequence} gives
\begin{align*}
 \|g_t\|
 =\eta_t^{-1}d_{\mathcal M_1}(x_t,x_{t+1})\leq
 C_{t+1}\eta_t^{-1}
 d_{\mathcal M_1}(x_t,\widetilde x_{t+1})
 +\frac{\epsilon_{t+1}}
 {2\eta_t(1-\lambda_1^{(t+1)})}
 =C_{t+1}\|\widetilde g_t\|+\mathfrak e_{t+1}.
\end{align*}
By Young's inequality, if $a,b\geq0$ and $\vartheta>0$, then
$(a+b)^2\leq(1+\vartheta)a^2+(1+\vartheta^{-1})b^2$. Hence on $E_{t+1}$,
\[
 \|g_t\|^2
 \leq
 (1+\vartheta)C_{t+1}^2\|\widetilde g_t\|^2
 +(1+\vartheta^{-1})\mathfrak e_{t+1}^2.
\]
Taking conditional expectations and using nonnegativity gives
\begin{align*}
 \mathbb E[\|\widetilde g_t\|_{x_t}^2
 \mid\mathcal F_t,E_{t+1}]
 =\frac{\mathbb E[\mathbf 1_{E_{t+1}}
 \|\widetilde g_t\|_{x_t}^2\mid\mathcal F_t]}{p_t}
 \leq\frac{\mathbb E[\|\widetilde g_t\|_{x_t}^2
 \mid\mathcal F_t]}{p_t}.
\end{align*}
Conditional on $\mathcal F_t$, the reference index is uniform, so
\[
 \mathbb E[\|\widetilde g_t\|_{x_t}^2\mid\mathcal F_t]
 =2f_{P_N}(x_t)
 \leq 2f_{P_N}(\widehat x_N)
 +\frac{1}{m_f}\|\nabla^{\mathcal M_1}f_{P_N}(x_t)\|_{x_t}^2,
\]
where the final inequality is the Polyak--\L{}ojasiewicz consequence of strong geodesic convexity with $f_{P_N}(x_{t}) = f_{P_N}(\widehat{x}_{N}) + f_{P_N}(x_{t}) - f_{P_N}(\widehat{x}_{N})$. This proves \eqref{eqn:metric-equivalence-second-moment-bound-general}. Lastly, if $\lambda_{1}^{(t+1)} = \lambda_{2}^{(t+1)} \leq \frac12$, then $C_{t+1}\leq3$. Taking $\vartheta=1/3$ gives $(1+\vartheta)C_{t+1}^2\leq12$ and $1+\vartheta^{-1}=4$, proving
\eqref{eqn:metric-equivalence-second-moment-bound}.
\end{proof}

\subsection{Proof of Lemma~\ref{lem:tube-containment-to-gradient-estimator-control}}\label{sec:proof-of-tube-containment-to-gradient-estimator-control}

\begin{proof}
By application of Jensen's followed by Corollary~\ref{cor:implicit-gradient-discrepancy}, we have
\begin{align}
 &\left\|\mathbb E[g_t-\widetilde g_t
 \mid\mathcal F_t,E_{t+1}]\right\|
 \leq
 \frac{\ROIraucheC}{\eta_t}
 \left[
 \left(
 \frac{\lambda_1^{(t+1)}+\lambda_2^{(t+1)}}
 {1-\lambda_1^{(t+1)}}\eta_t
 +2\alpha_{t+1}
 \right)D_t
 +\eta_{t}\mathfrak{e}_{t+1}
 \right]. \label{eqn:prelim-for-gradient-norm-discrep}
\end{align}

It remains to control the change in the ideal estimator caused by
conditioning on $E_{t+1}$. Note,
\begin{align*}
 \mathbb E[\widetilde g_t\mid\mathcal F_t,E_{t+1}]
 - \mathbb{E}[\widetilde g_t \mid \mathcal{F}_{t}] = 
 q_t\mathbb E[\widetilde g_t\mid\mathcal F_t,E_{t+1}]
 -\mathbb E[\mathbf 1_{E_{t+1}^c}\widetilde g_t\mid\mathcal F_t],
\end{align*}
as $\mathbb{E}[\widetilde g_t \mid \mathcal{F}_{t}] = \mathbb{E}[1_{E_{t+1}}\widetilde g_t \mid \mathcal{F}_{t}] + \mathbb{E}[1_{E_{t+1}^{c}}\widetilde g_t \mid \mathcal{F}_{t}]$. 
Since $\|\widetilde g_t\|_{x_t}\leq D_t$ almost surely, and
$\mathbb E[\widetilde g_t\mid\mathcal F_t]
=\nabla_{\mathcal M_1}f_{P_N}(x_t)$, it follows after taking norms and triangle inequality that
\begin{equation}\label{eqn:prelim-2-for-gradient-norm-discrep}
 \left\|\mathbb E[\widetilde g_t\mid\mathcal F_t,E_{t+1}]
 -\nabla_{\mathcal M_1}f_{P_N}(x_t)\right\|_{x_t}
 \leq2q_tD_t.
\end{equation}
Therefore, we have that 
\begin{align*}
    \|\mathbb E[g_t \mid\mathcal F_t,E_{t+1}]
 -\nabla^{\mathcal M_1}f_{P_N}(x_t)\|
 &\leq \frac{\ROIraucheC}{\eta_{t}}
 \left[
 \left(
 \frac{\lambda_1^{(t+1)}+\lambda_2^{(t+1)}}
 {1-\lambda_1^{(t+1)}}\eta_t
 +2\alpha_{t+1}
 \right)D_t
 + \eta_{t}\mathfrak{e}_{t+1}
 \right]
 +2q_{t}D_t\\
 & \leq (1-c_{1})\|\nabla_{\mathcal{M}_{1}} f_{P_{N}}(x_{t})\|,
\end{align*}
where the first inequality follows \eqref{eqn:prelim-for-gradient-norm-discrep} and \eqref{eqn:prelim-2-for-gradient-norm-discrep}, and the second from \eqref{eqn:tube-radius-vs-distance-to-minimizer}. 

The desired inequality \eqref{eqn:expected-gradient-points-in-right-direction} follows by application of  Cauchy--Schwarz and the above bound:
\begin{align*}
    &
    \left\langle
        \nabla_{\mathcal{M}_{1}}f(x_{t}),
        \mathbb{E}[
            g_{t}
            \mid
            \mathcal{F}_{t},E_{t+1}
        ]
    \right\rangle
    =
    \|\nabla_{\mathcal{M}_{1}}f(x_{t})\|^{2}
    +
    \left\langle
        \nabla_{\mathcal{M}_{1}}f(x_{t}),
        \mathbb{E}[
            g_{t}
            \mid
            \mathcal{F}_{t},E_{t+1}
        ]
        -
        \nabla_{\mathcal{M}_{1}}f(x_{t})
    \right\rangle\\
    &\qquad \qquad \geq
    \|\nabla_{\mathcal{M}_{1}}f(x_{t})\|^{2}
    -
    \|\nabla_{\mathcal{M}_{1}}f(x_{t})\|
    \left\|
        \mathbb{E}[
            g_{t}
            \mid
            \mathcal{F}_{t},E_{t+1}
        ]
        -
        \nabla_{\mathcal{M}_{1}}f(x_{t})
    \right\|\\
    &\qquad \qquad \geq
    c_{1}
    \|\nabla_{\mathcal{M}_{1}}f(x_{t})\|^{2},
\end{align*}
finishing the proof.
\end{proof}

\subsection{Proof of Lemma~\ref{lem:sampling-to-bound-probability-of-bad-event}}\label{sec:proof-of-sampling-to-bound-probability-of-bad-event}
We first establish a covering number bound for an arbitrary geodesic tube $T_{r}(\gamma_{x \rightarrow y}^{\mathcal{M}_{1}})$ where $\gamma_{x \rightarrow y}^{\mathcal{M}_{1}} \subseteq \ROI$.  Proof is deferred to Appendix~\ref{sec:proof-of-geodesic-tube-covering-number-bound}.

\begin{lemma}\label{lem:geodesic-tube-covering-number-bound}
Suppose Assumptions~\ref{ass:embedded-submanifold} and~\ref{ass:region-of-interest} hold. Let $x, y \in \ROI$, $x \neq y$, and $\gamma_{x\rightarrow y}^{\mathcal{M}_{1}} \subseteq \mathcal{M}_{1}$ with $\ell := d_{\mathcal{M}_{1}}(x,y)$. Let $r > 0$ be arbitrary and
\begin{equation}\label{eqn:lipschitz-constant-for-covering of a tube}
\tubeExpansionConstant{r}
    :=
    \max\left\{
        \cosh\left(\sqrt{\kmagInt}\,r\right),
        \frac{\sinh\left(\sqrt{\kmagInt}\,r\right)}
             {\sqrt{\kmagInt}\,r}
    \right\}
\end{equation}
The quotient is interpreted as $1$ when
$\widetilde K_{\mathcal M_1}r^2=0$. Then, for every \(\delta>0\),
\begin{equation}\label{lem:tube-covering-number}
    \mathcal{N}\left(T_{r}(\gamma_{x\rightarrow y}^{\mathcal{M}_{1}}), \delta\right) \leq \left(3 + \frac{2\tubeExpansionConstant{r}(\ell + 2r)}{\delta}\right)\left(1 + \frac{4\tubeExpansionConstant{r}r}{\delta}\right)^{d-1}.
\end{equation}
\end{lemma}

\begin{corollary}\label{cor:uniform-covering-number-for-tube-family}
Suppose Assumptions~\ref{ass:embedded-submanifold} and~\ref{ass:region-of-interest}
hold. Let $\alpha\in(0,1]$ and $\epsilon>0$ be fixed. Let $\mathfrak{S}_{\alpha,\epsilon}$ be defined by \eqref{eqn:sampling-family},
$L_r^{\widetilde K_{\mathcal M_1}}$ by \eqref{eqn:lipschitz-constant-for-covering of a tube}, and $\widetilde r_{\max}:=\alpha D_{\mathrm{reg}}
+\rho_{S_{\mathrm{reg}}}(\epsilon)$. Then, for every $\delta>0$,
\begin{equation}\label{us:eq-clipped-cover}
\sup_{S\in\mathfrak{S}_{\alpha,\epsilon}}\mathcal N(S,\delta)
\le
\left(
3+\frac{4L_{\widetilde r_{\max}}^{\widetilde K_{\mathcal M_1}}
(D_{\mathrm{reg}}+2\widetilde r_{\max})}{\delta}
\right)
\left(
1+\frac{8L_{\widetilde r_{\max}}^{\widetilde K_{\mathcal M_1}}
\widetilde r_{\max}}{\delta}
\right)^{d-1}.
\end{equation}
\end{corollary}

\begin{proof}
Let $S\in\mathfrak{S}_{\alpha,\epsilon}$ be arbitrary, and write
$S=S_{\mathrm{reg}}\cap T_{\widetilde r}(\gamma^{\mathcal M_1}_{x\to y})$,
with $\ell=d_{\mathcal M_1}(x,y)$. By
Assumption~\ref{ass:region-of-interest}, $\ell\le D_{\mathrm{reg}}$ and
$\widetilde r\le\widetilde r_{\max}$.

Take a $\delta/2$-net of
$T_{\widetilde r}(\gamma^{\mathcal M_1}_{x\to y})$.
Discard each ball not intersecting $S$, and, for every
remaining ball, choose a point in the intersection. By the triangle inequality, the remaining points form a $\delta$-net of $S$.
Thus,
\[
\mathcal N(S,\delta)
\le \mathcal N\bigl(
T_{\widetilde r}(\gamma^{\mathcal M_1}_{x\to y}),\delta/2\bigr).
\]
Applying Lemma~\ref{lem:geodesic-tube-covering-number-bound} at radius $\delta/2$, followed by
monotonicity of $L_r^{\widetilde K_{\mathcal M_1}}$ in $r$, gives
\[
\begin{aligned}
\mathcal N(S,\delta)
&\le
\left(3+\frac{4L_{\widetilde r}^{\widetilde K_{\mathcal M_1}}
(\ell+2\widetilde r)}{\delta}\right)
\left(1+\frac{8L_{\widetilde r}^{\widetilde K_{\mathcal M_1}}
\widetilde r}{\delta}\right)^{d-1}
\\
&\le
\left(3+\frac{4L_{\widetilde r_{\max}}^{\widetilde K_{\mathcal M_1}}
(D_{\mathrm{reg}}+2\widetilde r_{\max})}{\delta}\right)
\left(1+\frac{8L_{\widetilde r_{\max}}^{\widetilde K_{\mathcal M_1}}
\widetilde r_{\max}}{\delta}\right)^{d-1}.
\end{aligned}
\]
The right-hand side is independent of $S$, finishing the proof.
\end{proof}

In order to turn this covering number into a sampling rate, via Corollary~\ref{cor:probabilistic-isomap}, we require a uniform lower-mass regularity for the sampling measure used by the algorithm. The following derives this for the normalized volume measure; proof is deferred to Appendix~\ref{sec:proof-of-sampling-mass-lower-bound}.

\begin{lemma}\label{lem:sampling-mass-lower-bound}
Assume Assumptions~\ref{ass:embedded-submanifold} and~\ref{ass:region-of-interest}
hold. There exists
$c_{\mathrm{mass}}(\ROI, \ROIdiam, \kmagInt) >0$, such that the following holds.
For every $x,y\in S_{\mathrm{reg}}$, let
$\gamma:=\gamma^{\mathcal M_1}_{x\to y}$ and
$\ell:=d_{\mathcal M_1}(x,y)$, with $\gamma=\{x\}$ when $x=y$.
For $0<r\le2D_{\mathrm{reg}}$, define
\[
S_r(\gamma):=S_{\mathrm{reg}}\cap T_r(\gamma),
\qquad
\mu_{\gamma,r}(A):=
\frac{\operatorname{vol}_{\mathcal M_1}(A\cap S_r(\gamma))}
{\operatorname{vol}_{\mathcal M_1}(S_r(\gamma))}.
\]
Then $\mu_{\gamma,r}$ is well-defined and, for every $0 < \delta \leq \min\{r,\ROIdiam\}$,
\begin{equation}\label{eqn:lower-mass-bound-general}
\inf_{v\in S_r(\gamma)}
\mu_{\gamma,r}\bigl(B_{\mathcal M_1}(v,\delta)\bigr)
\ge
c_{\mathrm{mass}}(\ROI, \ROIdiam, \kmagInt)
\frac{\delta^{d}}
{(\ell+r)r^{d-1}},
\end{equation}
where
\[
c_{\mathrm{mass}}(\ROI, \ROIdiam, \kmagInt):=
\frac{\operatorname{vol}_{\mathcal M_1}(S_{\mathrm{reg}})}
{\left[
2D_{\mathrm{reg}}C_{\log}
L_{D_{\mathrm{reg}}}^{\widetilde K_{\mathcal M_1}}
L_{2D_{\mathrm{reg}}}^{\widetilde K_{\mathcal M_1}}
\right]^d\max\{\omega_{d-1},2\omega_d\}},
\]
with $\omega_{j}$ denoting the volume of the Euclidean unit ball in dimension $j$, and $C_{\operatorname{log}}$ defined by \eqref{eqn:uniform-rauche-comparison-constant}.
\end{lemma}

We combine the covering and lower-mass bounds, then verify that the
covering event implies the endpoint-rooted comparisons defining $E_t$.

\begin{proof}[Proof of Lemma~\ref{lem:sampling-to-bound-probability-of-bad-event}]
Write $D:=D_{\mathrm{reg}}$ and $s_t:=\delta_t/2$.
For every $0<\epsilon\le D$, we have
$\rho_{S_{\mathrm{reg}}}(\epsilon)\ge\epsilon$.

Let $S\in\mathfrak S_{\alpha_t,\epsilon_t}$, with centerline
length $\ell$ and radius
$\widetilde r=\alpha_t\ell+\rho_{S_{\mathrm{reg}}}(\epsilon_t)$.
The scale restrictions and Lemma~\ref{lem:local-intrinsic-extrinsic-compatibility} give
\[
s_t\le\epsilon_t\le\widetilde r
\le\widetilde r_{\max}^{(t)}
\le\alpha_tD+C_{\mathrm{loc},S_{\mathrm{reg}}}\epsilon_t
\le2\alpha_tD\le2D,
\qquad s_t\le D.
\]
Thus Lemma~\ref{lem:sampling-mass-lower-bound} applies at radius $s_t$, and gives
\begin{equation}\label{us:eq-family-mass-small}
[q^{(t)}(s_t)]^{-1}
\le\frac{(D+\widetilde r_{\max}^{(t)})
(\widetilde r_{\max}^{(t)})^{d-1}}{c_{\mathrm{mass}}s_t^d}.
\end{equation}
Corollary~\ref{cor:uniform-covering-number-for-tube-family} also gives
$\mathcal N(S,s_t)\le\mathcal N_t(s_t)$.

We next justify the metric comparisons needed on a covering event.
Fix endpoints $x,y\in S_{\mathrm{reg}}$ with $x\ne y$, set
$\gamma:=\gamma^{\mathcal M_1}_{x\to y}$, and let
$S:=S_{\mathrm{reg}}\cap T_{\widetilde r}(\gamma)$.
For $a\in\{x,y\}$, define
\[
\varphi_a(w):=\tfrac12d_{\mathcal M_1}(a,w)^2,
\qquad
F_\tau^a(w):=\exp_a^{\mathcal M_1}
(\tau\log_a^{\mathcal M_1}(w)),\quad 0\le\tau\le1.
\]
Lemma~\ref{lem:sc-and-smoothness-distance-squared} and Assumption~\ref{ass:region-of-interest} give
$\operatorname{Hess}\varphi_a\succeq m_fg_1$ on $S_{\mathrm{reg}}$.
The curve $s\mapsto F_{e^{-s}}^a(w)$ is its negative gradient flow
and remains in $S_{\mathrm{reg}}$ by geodesic convexity.
For a minimizing path $\sigma$ in $S_{\mathrm{reg}}$, the
variational field $J:=\partial_uF_{e^{-s}}^a(\sigma(u))$ satisfies
\[
\frac{d}{ds}\|J\|^2
=-2\operatorname{Hess}\varphi_a[J,J]
\le-2m_f\|J\|^2.
\]
Integrating this inequality and then the length of the varied path gives
\[
d_{\mathcal M_1}(F_\tau^a(v),F_\tau^a(w))
\le\tau^{m_f}d_{\mathcal M_1}(v,w),\qquad 0<\tau\le1.
\]
For $v\in S$, choose a closest point $p\in\gamma$.
Since $F_\tau^a(p)\in\gamma$, it follows that
\[
\operatorname{dist}_{\mathcal M_1}(F_\tau^a(v),\gamma)
\le\tau^{m_f}d_{\mathcal M_1}(v,p)
\le\widetilde r.
\]
Therefore the minimizing geodesic from $a$ to $v$ stays in $S$;
the case $\tau=0$ follows by continuity.

Now let $\mathcal V\subseteq S$ be a $\delta_t$-cover containing
$x,y$, and construct its $\epsilon_t$-neighborhood graph $G$.
The edgewise lower-bound argument of
Theorem~\ref{thm:deterministic-metric-equivalence-for-ISOMAP}
when Assumption~\ref{ass:embedded-submanifold}\textup{(a)} holds, or
Proposition~\ref{prop:abstract-metric-isomap}
when Assumption~\ref{ass:embedded-submanifold}\textup{(b)} holds, applied with comparison region
$S_{\mathrm{reg}}$, gives
\[
(1-\lambda_t)d_{\mathcal M_1}(u,v)\le d_G(u,v),
\qquad u,v\in\mathcal V.
\]
For the upper bound, fix $a\in\{x,y\}$ and $v\in\mathcal V$.
The minimizing segment from $a$ to $v$ lies in $S$, as just proved.
When $a\ne v$, put $L_0:=d_{\mathcal M_1}(a,v)$ and
$m:=\lceil L_0/(\epsilon_t-2\delta_t)\rceil$.
Subdivide the segment into $m$ equal pieces and approximating the
interior points by vertices of $\mathcal{V}$ within distance of $\delta_t$.
Consecutive vertices have ambient distance at most
$L_0/m+2\delta_t\le\epsilon_t$. Then, the graph path length is at most
\[
L_0+2\delta_t(m-1)
\le\left(1+\frac{2\delta_t}{\epsilon_t-2\delta_t}\right)L_0
\le(1+\lambda_t)L_0,
\]
using $\delta_t\le\lambda_t\epsilon_t/4$ and $\lambda_t<1$.
Thus, such $\delta_{t}$ coverings ensure $E_{t}$ holds.

Finally, condition on
$\mathcal H_t:=\mathcal F_{t-1}\vee\sigma(I_t)$.
On a nontrivial step, the endpoints and $\widetilde S_t$ are fixed,
and the samples are conditionally i.i.d. from normalized volume on
$\widetilde S_t$. By Lemma~\ref{lem:sampling-for-delta-net},
\eqref{us:eq-family-mass-small}, and the covering bound, the samples form a
$\delta_t$-cover with conditional probability at least $1-\zeta_t$. On a trivial step, the singleton-graph convention
makes $E_t$ hold without sampling. Hence
\[
\mathbb P(E_t\mid\mathcal H_t)\ge1-\zeta_t.
\]
Taking conditional expectations with respect to
$\mathcal F_{t-1}\subseteq\mathcal H_t$, the tower property implies
\[
\mathbb P(E_t\mid\mathcal F_{t-1})\ge1-\zeta_t.
\]
The complementary bound follows.
\end{proof}

\subsubsection{Proof of Lemma~\ref{lem:geodesic-tube-covering-number-bound}}\label{sec:proof-of-geodesic-tube-covering-number-bound}
To bound the covering number of a geodesic tube in $\mathcal{M}_{1}$, we parameterize its normal portion by a Fermi map and its endpoint portions by exponential maps. Then, nets in the parameter domain can be pushed forward to obtain nets in the tube, only using forward Lipschitz estimates.

\begin{lemma}\label{lem:transfer-covering-numbers}
Let $\varphi:(X,d_X)\to(Y,d_Y)$ be $L_{\varphi}$-Lipschitz, with $L_{\varphi}>0$. If $K\subset X$ and $\{x_i\}_{i=1}^N$ is a $\delta$-net of $K$, then $\{\varphi(x_i)\}_{i=1}^N$ is an $L_{\varphi}\delta$-net of $\varphi(K)$. Consequently,
\[
    \mathcal N\bigl(\varphi(K),L_{\varphi}\delta\bigr)
    \leq
    \mathcal N(K,\delta).
\]
\end{lemma}

\begin{proof}
For any $y\in\varphi(K)$, choose $x\in K$ with $y=\varphi(x)$. Since $\{x_i\}_{i=1}^N$ is a $\delta$-net of $K$, there exists $i$ such that $d_X(x,x_i)\leq\delta$. Hence
\[
    d_Y\bigl(y,\varphi(x_i)\bigr)
    =
    d_Y\bigl(\varphi(x),\varphi(x_i)\bigr)
    \leq
    L_{\varphi}d_X(x,x_i)
    \leq
    L_{\varphi}\delta.
\]
\end{proof}

The next lemma supplies the required Lipschitz estimates for the Fermi map and the exponential maps centered at the endpoints.

\begin{lemma}\label{lem:lipschitz-constant-of-fermi-and-endpoint-caps}
Let $\mathcal M$ be a complete $d$-dimensional Riemannian manifold satisfying $|\operatorname{sec}_{\mathcal M}|\leq\tilde K$ for some $\tilde K\geq0$. Let $\gamma:[0,\ell]\to\mathcal M$ be a unit-speed geodesic segment, and let $E_2,\ldots,E_d$ be a parallel orthonormal frame of the normal bundle along $\gamma$. For $r>0$, let $B^m(r):=\{z\in\mathbb R^m:\|z\|\leq r\},$ and define $\Psi \colon [0, \ell] \times B^{d-1}(r) \rightarrow \mathcal{M}$ by
\[
    \Psi(t,z)
    :=
    \exp_{\gamma(t)}
    \left(
        \sum_{i=2}^d z_iE_i(t)
    \right).
\]
Define
\begin{equation}\label{eqn:lipschitz-constant-for-fermi-and-endpoint}
    L_{r}^{\tilde{K}}
    :=
    \max\left\{
        \cosh(\sqrt{\tilde K}\,r),
        \frac{\sinh(\sqrt{\tilde K}\,r)}
             {\sqrt{\tilde K}\,r}
    \right\},
\end{equation}
where the quotient is interpreted as $1$ when $\tilde K r=0$. Then:
\begin{enumerate}
    \item For all $(t,z),(s,w)\in[0,\ell]\times B^{d-1}(r)$,
    \[
        d_{\mathcal M}\bigl(\Psi(t,z),\Psi(s,w)\bigr)
        \leq
        L_{r}^{\tilde{K}}\,
        d_1\bigl((t,z),(s,w)\bigr).
    \]
    where $d_1\bigl((t,z),(s,w)\bigr) := |t-s|+\|z-w\|.$

    \item For every $p\in\mathcal M$, the restriction of $\exp_p$ to
    \[
        B_{T_p\mathcal M}(r)
        :=
        \{v\in T_p\mathcal M:\|v\|\leq r\}
    \]
    is $L_{r}^{\tilde{K}}$-Lipschitz. 
\end{enumerate}
\end{lemma}
\begin{proof}
For the rest of the proof, set
\begin{equation} \label{eqn:hyperbolic-comparison-funtions}
    C_{\widetilde K}(\tau)
    :=
    \cosh(\sqrt{\widetilde K}\,\tau)
\qquad \text{and} \qquad
    S_{\widetilde K}(\tau)
    :=
    \begin{cases}
        \dfrac{\sinh(\sqrt{\widetilde K}\,\tau)}
              {\sqrt{\widetilde K}},
        & \widetilde K>0,\\[2mm]
        \tau,
        & \widetilde K=0.
    \end{cases}
\end{equation}
We use the following standard two-sided Jacobi-field estimates
\cite[Theorem~6.5.2]{Jost2017}. If $\sigma$ is a unit-speed
geodesic on a manifold with
$-\widetilde K\leq\operatorname{sec}\leq\widetilde K$, then a
Jacobi field $J$ along $\sigma$ satisfies
\begin{equation}\label{eq:two-special-jacobi-estimates}
    \begin{aligned}
        J(0)=0
        &\quad\Longrightarrow\quad
        \|J(\tau)\|
        \leq
        S_{\widetilde K}(\tau)\|D_\tau J(0)\|,\\
        D_\tau J(0)=0
        &\quad\Longrightarrow\quad
        \|J(\tau)\|
        \leq
        C_{\widetilde K}(\tau)\|J(0)\|.
    \end{aligned}
\end{equation}

\paragraph{Lipschitz Estimate for Exponential Map.} Fix $p\in\mathcal M$ and
write $v=\rho u$, where $\rho>0$ and $\|u\|=1$. Along the radial
geodesic $\sigma(\tau)=\exp_p(\tau u)$, we have $D(\exp_p)_{\rho u}(u)=\dot\sigma(\rho)$ and $\bigl\|D(\exp_p)_{\rho u}(u)\bigr\|=1$. 

If $w\perp u$, consider the geodesic variation and its Jacobi field:
\[
    \Sigma(\tau, s)
    :=
    \exp_p\left(\tau\left(u+\frac{s}{\rho}w\right)\right)
\qquad \text{and} \qquad
    J_w(\tau):=\partial_s\Sigma(\tau,0),
\]
respectively. By construction, $J_w(0)=0, D_\tau J_w(0)=\frac{w}{\rho}, $ and $J_w(\rho)=D(\exp_p)_{\rho u}(w).$
Hence
\[
    \bigl\|D(\exp_p)_{\rho u}(w)\bigr\| = \|J_{w}(\rho)\| \leq S_{\widetilde{K}}(\rho) \| D_{\tau}(J_{w}(0))\| = \frac{S_{\widetilde K}(\rho)}{\rho}\|w\|.
\]
by application of \eqref{eq:two-special-jacobi-estimates}.

By the Gauss lemma, the images of the radial and transverse
components are orthogonal. Thus, writing
$\xi=au+w$ with $w\perp u$, we obtain
\[
    \bigl\|D(\exp_p)_{\rho u}(\xi)\bigr\|^2
    = \left\| D(\exp_{p})_{\rho u}(au) + D(\exp_{p})_{\rho u}(w)\right\|^{2} = 
    a^2+
    \bigl\|D(\exp_p)_{\rho u}(w)\bigr\|^2.
\]
Therefore, we have that for $\rho \leq r$:
\begin{equation*}
    \|D(\exp_{p})_{\rho u}(\xi)\|^{2} \leq a^{2} + \left(\frac{S_{\widetilde{K}}(\rho)}{\rho}\right)^{2} \|w\|^{2} \leq \left(L_{r}^{\widetilde{K}}\right)^{2} \|\xi\|^{2}.
\end{equation*}
The second inequality follows from $\tau \mapsto S_{\tilde{K}}(\tau)/\tau$ being non-decreasing and $1 \leq L_{r}^{\widetilde{K}}$. Therefore, $\|D(\exp_{p})_{\rho u}(\xi)\| \leq L_{r}^{\widetilde{K}}\|\xi\|$. The same conclusion holds at the origin as
$D(\exp_p)_0=\operatorname{id}$. Integrating along the line
segment between two vectors in the convex ball
$B_{T_p\mathcal M}(r)$ shows
\begin{equation}\label{eqn:exponential-map-Lipschitz-estimate}
    d_{\mathcal M}\bigl(\exp_p(v),\exp_p(w)\bigr)
    \leq
    L_{r}^{\widetilde{K}}\|v-w\|,
    \qquad
    v,w\in B_{T_p\mathcal M}(r).
\end{equation}

\paragraph{Lipschitz Estimate for $\Psi$} Extend $\gamma$ and its parallel normal
frame slightly beyond the endpoints when necessary. Fix
$(t,z)$ with $\rho:=\|z\|>0$, and put
\[
    u_a
    :=
    \rho^{-1}\sum_{i=2}^d z_iE_i(a).
\]
The field $a\mapsto u_a$ is parallel along $\gamma$. Consider the geodesic variation and associated Jacobi field
\[
    \Sigma(\tau,s)
    :=
    \exp_{\gamma(t+s)}(\tau u_{t+s})
\qquad \text{and} \qquad
    J_{\parallel}(\tau)
    :=
    \Sigma_s (\tau,0).
\]
This Jacobi field satisfies $J_{\parallel}(0)=\dot\gamma(t) $ and $D_\tau J_{\parallel}(0)=0,$ the second equality following from the parallelness of
$a\mapsto u_a$. Therefore, since $\|\dot{\gamma}\| = 1$, application of \eqref{eq:two-special-jacobi-estimates} results in
\[
    \bigl\|D\Psi_{(t,z)}(\partial_t)\bigr\|
    =
    \|J_{\parallel}(\rho)\|
    \leq
    C_{\widetilde K}(\rho) \|J_{\parallel}(0)\| = C_{\widetilde K}(\rho) 
    \leq
    L_{r}^{\widetilde{K}}.
\]
For $z=0$ the same estimate follows from
$\Psi(t,0)=\gamma(t)$. Integration in the $t$-variable gives
\begin{equation}\label{eq:transitive-estimate}
    d_{\mathcal M}\bigl(\Psi(t,z),\Psi(s,z)\bigr) \leq L_{r}^{\widetilde{K}}|t-s|.
\end{equation}

For fixed $s$, the linear map $I_{s} \colon \mathbb{R}^{d-1} \rightarrow T_{\gamma(s)}\mathcal{M}$ defined by
\[
    I_s(z)=\sum_{i=2}^d z_iE_i(s),
\]
is an isometric embedding and $\Psi(s,z)=\exp_{\gamma(s)}(I_s(z)).$ Therefore, we may apply the exponential map estimate \eqref{eqn:exponential-map-Lipschitz-estimate}. Thus, by application of triangle inequality followed by \eqref{eqn:exponential-map-Lipschitz-estimate} and \eqref{eq:transitive-estimate} we have:
\[
    \begin{aligned}
    d_{\mathcal M}\bigl(\Psi(t,z),\Psi(s,w)\bigr)
    &\leq
    d_{\mathcal M}\bigl(\Psi(t,z),\Psi(s,z)\bigr)
    +
    d_{\mathcal M}\bigl(\Psi(s,z),\Psi(s,w)\bigr) \\
    &\leq 
    L_{r}^{\widetilde{K}}
    \bigl(|t-s|+\|z-w\|\bigr).
    \end{aligned}
\]
\end{proof}

We are now in a position to prove the covering-number bound for the geodesic tube
$T_r(\gamma_{x\rightarrow y}^{\mathcal M_1})$.

\begin{proof}[Proof of Lemma~\ref{lem:geodesic-tube-covering-number-bound}]
For notational simplicity, let $\gamma := \gamma_{x\rightarrow y}^{\mathcal M_1},$ parameterized by arclength on $[0,\ell]$, and set $\tilde K := \tilde K_{\mathcal M_1},$ and $L:=L_{r}^{\tilde{K}}.$ Additionally, let $Q := [0,\ell]\times B^{d-1}(r),$ and $\mathcal C_r(\gamma):=\Psi(Q),$ and define the endpoint region $\mathcal E_r(\gamma):=\overline B_{\mathcal M_1}(x,r)\cup\overline B_{\mathcal M_1}(y,r).$ Then
\begin{equation}\label{eqn:tube-decomposition}
    T_r(\gamma)
    =
    \mathcal C_r(\gamma)
    \cup
    \mathcal E_r(\gamma).
\end{equation}
Indeed, for any $q\in T_r(\gamma)$, compactness of $\gamma([0,\ell])$ gives a closest point $\gamma(t_0)$. If $q=\gamma(t_0)$, then $q\in\mathcal C_r(\gamma)$ by taking the zero normal vector. Otherwise, if $t_0\in(0,\ell)$, let $\nu$ be a minimizing geodesic from $\gamma(t_0)$ to $q$. The first variation formula at the closest interior point gives $\dot\nu(0)\perp\dot\gamma(t_0),$ so $q\in\mathcal C_r(\gamma)$. If $t_0\in\{0,\ell\}$, then
$q\in\mathcal E_r(\gamma)$. The reverse inclusion in
\eqref{eqn:tube-decomposition} is immediate.

For the rest of the proof, we will use the standard Euclidean estimate
\begin{equation}\label{eqn:euclidean-ball-covering}
    \mathcal N\bigl(B^m(r),\eta\bigr)
    \leq
    \left(
        1+\frac{2r}{\eta}
    \right)^m,
    \qquad
    m\geq0.
\end{equation}
For $m\geq1$, this follows by taking a maximal $\eta$-separated subset: the Euclidean balls of radius $\eta/2$ centered at its points are disjoint and contained in $B^m(r+\eta/2)$. The case $m=0$ is immediate.

\paragraph{Bounding $\mathcal N(\mathcal C_r(\gamma),\delta)$.}
Equip $Q$ with the metric $d_1$ from Lemma~\ref{lem:lipschitz-constant-of-fermi-and-endpoint-caps}, and set $\eta:=\delta/(2L).$ There exist $\eta$-nets $A\subset[0,\ell]$ and $B \subseteq B^{d-1}(r)$ satisfying
\begin{equation*}
    |A| \leq 1+\frac{\ell}{\eta} = 1+\frac{2L\ell}{\delta} \qquad \text{and} \qquad |B|
    \leq
    \left(
        1+\frac{2r}{\eta}
    \right)^{d-1}
    =
    \left(
        1+\frac{4Lr}{\delta}
    \right)^{d-1}.
\end{equation*}
The first estimate $|A|$ follows from a uniform discretization of $[0,\ell]$ and the second from \eqref{eqn:euclidean-ball-covering}. 

Given these $\eta$-nets, notice that for every $(t,z)\in Q$, we may choose $a\in A$ and $b\in B$ with $|t-a|\leq\eta$ and $\|z-b\|\leq\eta.$ Then $d_1\bigl((t,z),(a,b)\bigr) \leq 2\eta = \frac{\delta}{L},$ implying $A\times B$ is a $\delta/L$-net of $Q$. By Lemmas~\ref{lem:transfer-covering-numbers} and \ref{lem:lipschitz-constant-of-fermi-and-endpoint-caps}, it follows that
$\Psi(A\times B)$ is a $\delta$-net of $\mathcal C_r(\gamma)$. Therefore,
\begin{equation}\label{eqn:covering-number-on-normal-cylinder}
    \mathcal N\bigl(\mathcal C_r(\gamma),\delta\bigr)
    \leq
    \left(
        1+\frac{2L\ell}{\delta}
    \right)
    \left(
        1+\frac{4Lr}{\delta}
    \right)^{d-1}.
\end{equation}

\paragraph{Bounding $\mathcal N(\mathcal E_r(\gamma),\delta)$.}
For each $p\in\{x,y\}$,
Lemmas~\ref{lem:transfer-covering-numbers} and
\ref{lem:lipschitz-constant-of-fermi-and-endpoint-caps}, together with
\eqref{eqn:euclidean-ball-covering}, give
\[
    \mathcal N
    \bigl(
        \overline B_{\mathcal M_1}(p,r),
        \delta
    \bigr)
    \leq
    \mathcal N
    \bigl(
        B^d(r),
        \delta/L
    \bigr)
    \leq
    \left(
        1+\frac{2Lr}{\delta}
    \right)^d.
\]
Hence
\begin{equation}\label{eqn:covering-number-of-endpoint-caps}
    \mathcal N\bigl(\mathcal E_r(\gamma),\delta\bigr)
    \leq
    2
    \left(
        1+\frac{2Lr}{\delta}
    \right)^d.
\end{equation}

\paragraph{Bounding $\mathcal{N}(T_r(\gamma), \delta)$.}
By \eqref{eqn:tube-decomposition}, subadditivity of covering numbers, and
\eqref{eqn:covering-number-on-normal-cylinder} and \eqref{eqn:covering-number-of-endpoint-caps},
\begin{align*}
    \mathcal N\bigl(T_r(\gamma),\delta\bigr)
    &\leq
    \left(
        1+\frac{2L\ell}{\delta}
    \right)
    \left(
        1+\frac{4Lr}{\delta}
    \right)^{d-1}
    +
    2
    \left(
        1+\frac{2Lr}{\delta}
    \right)^d
    \\
    &\leq
    \left(
        3+\frac{2L(\ell+2r)}{\delta}
    \right)
    \left(
        1+\frac{4Lr}{\delta}
    \right)^{d-1}.
\end{align*}
For the last inequality, set $b:=4Lr/\delta$ and use $2(1+b/2)^d\leq(2+b)(1+b)^{d-1}.$
\end{proof}
\subsubsection{Proof of Lemma~\ref{lem:sampling-mass-lower-bound}}\label{sec:proof-of-sampling-mass-lower-bound}
For notational simplicity, write
\[
D:=D_{\mathrm{reg}},\qquad
V_{\mathrm{reg}}:=\operatorname{vol}_{\mathcal M_1}(S_{\mathrm{reg}}),
\qquad
L(r):=L_r^{\widetilde K_{\mathcal M_1}},
\]
and let $C_{\log}$ be the constant in~\eqref{eqn:uniform-rauche-comparison-constant}.
By Assumption~\ref{ass:region-of-interest}, $D>0$ and
$V_{\mathrm{reg}}>0$. We first obtain a uniform lower bound on the
volume of intrinsic balls intersected with $S_{\mathrm{reg}}$.

\begin{lemma}\label{lem:region-thickness}
Assume Assumptions~\ref{ass:embedded-submanifold} and~\ref{ass:region-of-interest}
hold. Define
\[
a_{\mathrm{reg}}
:=\frac{V_{\mathrm{reg}}}
{D^d[C_{\log}L(D)]^d}>0.
\]
Then, for every $z\in S_{\mathrm{reg}}$ and every $0<s\le D$,
\begin{equation}\label{us:eq-region-thickness}
\operatorname{vol}_{\mathcal M_1}
\bigl(B_{\mathcal M_1}(z,s)\cap S_{\mathrm{reg}}\bigr)
\ge a_{\mathrm{reg}}s^d.
\end{equation}
\end{lemma}

\begin{proof}
Fix $z\in S_{\mathrm{reg}}$ and $0<s\le D$, and set $\tau:=s/D$.
Consider
\[
F_\tau(w):=
\exp_z^{\mathcal M_1}\bigl(\tau\log_z^{\mathcal M_1}(w)\bigr),
\qquad w\in S_{\mathrm{reg}}.
\]
By Assumption~\ref{ass:region-of-interest}, $S_{\mathrm{reg}}$ is within the normal neighborhood
of $z$, thus by geodesic convexity 
\begin{equation}\label{eqn:contraction-inclusion}
F_\tau(S_{\mathrm{reg}})
\subseteq S_{\mathrm{reg}}\cap B_{\mathcal M_1}(z,s).
\end{equation}

Let $v:=\log_z^{\mathcal M_1}(w)$. By the chain rule,
\begin{equation}\label{eqn:chain-rule-for-contraction}
D_wF_\tau
=D_{\tau v}\exp_z^{\mathcal M_1}
\circ(\tau I)\circ D_w\log_z^{\mathcal M_1}.
\end{equation}
The exponential-map estimate of Lemma~\ref{lem:lipschitz-constant-of-fermi-and-endpoint-caps} gives
$\|D_v\exp_z^{\mathcal M_1}\|\le L(D)$. Since logarithm and
exponential are inverse on this normal neighborhood,
$\sigma_{\min}\bigl(D_w\log_z^{\mathcal M_1}\bigr)
\ge L(D)^{-1}.$
Also, $F_\tau(w)\in S_{\mathrm{reg}}$, so
\eqref{eqn:Rauch-comparison-estimate} gives
$\|D_{F_\tau(w)}\log_z^{\mathcal M_1}\|\le C_{\log}$.
Thus,
$\sigma_{\min}\bigl(D_{\tau v}\exp_z^{\mathcal M_1}\bigr)
\ge C_{\log}^{-1}.$
Applying these estimates to \eqref{eqn:chain-rule-for-contraction} results in the following lower bound
\[
|\operatorname{Jac}F_\tau(w)|
\ge\left(\frac{\tau}{C_{\log}L(D)}\right)^d.
\]
The map $F_\tau$ is injective, since $\tau>0$ and the radial
exponential map is injective on the relevant normal neighborhood.
Thus, change of variables gives
\[
\begin{aligned}
\operatorname{vol}_{\mathcal M_1}
\bigl(F_\tau(S_{\mathrm{reg}})\bigr)
= \int_{\ROI}|\operatorname{Jac}F_{\tau}(w)|d\operatorname{vol}_{\mathcal M_{1}}(w)\ge
\left(\frac{\tau}{C_{\log}L(D)}\right)^dV_{\mathrm{reg}}
&=a_{\mathrm{reg}}s^d.
\end{aligned}
\]
The inclusion~\eqref{eqn:contraction-inclusion} finishes the proof.
\end{proof}

We next upper bound the volume of the full tube. The argument uses the tube decomposition from the proof of Lemma~\ref{lem:geodesic-tube-covering-number-bound}
and the differential bounds of
Lemma~\ref{lem:lipschitz-constant-of-fermi-and-endpoint-caps}.

\begin{lemma}\label{lem:tube-volume}
Assume Assumptions~\ref{ass:embedded-submanifold} and~\ref{ass:region-of-interest}
hold. Let $\gamma$ be a minimizing geodesic segment of length
$\ell\ge0$ in $\mathcal M_1$. For every $r>0$,
\begin{equation}\label{us:eq-tube-volume}
\operatorname{vol}_{\mathcal M_1}(T_r(\gamma))
\le L(r)^d
\bigl(\omega_{d-1}\ell r^{d-1}+2\omega_dr^d\bigr),
\end{equation}
where $\omega_j$ denotes the volume of the Euclidean unit ball in
dimension $j$, with $\omega_0:=1$. In particular, defining
\[
C_{\mathrm{tube}}
:=L(2D)^d\max\{\omega_{d-1},2\omega_d\},
\]
we have, for $0<r\le2D$,
\begin{equation}\label{us:eq-tube-volume-uniform}
\operatorname{vol}_{\mathcal M_1}(T_r(\gamma))
\le C_{\mathrm{tube}}(\ell+r)r^{d-1}.
\end{equation}
\end{lemma}

\begin{proof}
Suppose first that $\ell>0$, and parameterize $\gamma$ by arclength
on $[0,\ell]$. Let $E_2,\ldots,E_d$ be a parallel orthonormal frame
normal to $\dot\gamma$, and define
\[
\Phi(t,z):=
\exp_{\gamma(t)}^{\mathcal M_1}
\left(\sum_{i=2}^d z_iE_i(t)\right),
\qquad
(t,z)\in[0,\ell]\times B^{d-1}(r).
\]
The decomposition from the proof of Lemma~\ref{lem:geodesic-tube-covering-number-bound}
gives
\[
T_r(\gamma)
=\Phi\bigl([0,\ell]\times B^{d-1}(r)\bigr)
\cup B_{\mathcal M_1}(\gamma(0),r)
\cup B_{\mathcal M_1}(\gamma(\ell),r).
\]
By Lemma~\ref{lem:lipschitz-constant-of-fermi-and-endpoint-caps}, each coordinate derivative of $\Phi$
has norm at most $L(r)$. Hadamard's determinant inequality therefore
gives $|\operatorname{Jac}\Phi|\le L(r)^d$. The area formula implies
\[
\operatorname{vol}_{\mathcal M_1}
\bigl(\Phi([0,\ell]\times B^{d-1}(r))\bigr)
\le L(r)^d\omega_{d-1}\ell r^{d-1}.
\]

For either endpoint $p$ of $\gamma$, completeness of $\mathcal M_1$
gives
$B_{\mathcal M_1}(p,r)=
\exp_p^{\mathcal M_1}(B_{T_p\mathcal M_1}(r))$.
The exponential-map estimate of Lemma~\ref{lem:lipschitz-constant-of-fermi-and-endpoint-caps} and the
area formula give
\[
\operatorname{vol}_{\mathcal M_1}(B_{\mathcal M_1}(p,r))
\le L(r)^d\omega_dr^d.
\]
Adding these estimates proves~\eqref{us:eq-tube-volume}.
If $\ell=0$, the tube is a single ball, and the endpoint estimate
alone suffices. Finally, monotonicity of $L$ gives
\eqref{us:eq-tube-volume-uniform} when $r\le2D$.
\end{proof}

The following transfers the thickness of $S_{\mathrm{reg}}$ to its
intersection with a geodesic tube.

\begin{lemma}\label{lem:clipped-thickness}
Assume Assumptions~\ref{ass:embedded-submanifold} and~\ref{ass:region-of-interest}
hold. Let $x,y\in S_{\mathrm{reg}}$, let
$\gamma:=\gamma^{\mathcal M_1}_{x\to y}$, and let $r>0$.
For every $v\in S_r(\gamma)$ and $0<s\le\min\{r,D\}$,
\begin{equation}\label{us:eq-clipped-thickness}
\operatorname{vol}_{\mathcal M_1}
\bigl(B_{\mathcal M_1}(v,s)\cap S_r(\gamma)\bigr)
\ge\frac{a_{\mathrm{reg}}}{2^d}s^d.
\end{equation}
\end{lemma}

\begin{proof}
Choose a closest point $p\in\gamma$ to $v$, and set
$b:=d_{\mathcal M_1}(v,p)\le r$. If $b\ge s/2$, let $z$ be the
point at distance $s/2$ from $v$ along the minimizing geodesic from
$v$ to $p$. If $b<s/2$, set $z:=p$. In either case,
$z\in S_{\mathrm{reg}}$ by geodesic convexity.

We claim that
\begin{equation}\label{us:eq-thick-patch-inclusion}
B_{\mathcal M_1}(z,s/2)\cap S_{\mathrm{reg}}
\subseteq B_{\mathcal M_1}(v,s)\cap S_r(\gamma).
\end{equation}
Indeed, if $b\ge s/2$ and
$w\in B_{\mathcal M_1}(z,s/2)\cap S_{\mathrm{reg}}$, then
\[
d_{\mathcal M_1}(v,w)\le s,
\qquad
\operatorname{dist}_{\mathcal M_1}(w,\gamma)
\le d_{\mathcal M_1}(w,p)
\le s/2+(b-s/2)\le r.
\]
If $b<s/2$, the corresponding estimates are
\[
d_{\mathcal M_1}(v,w)\le b+s/2\le s,
\qquad
\operatorname{dist}_{\mathcal M_1}(w,\gamma)\le s/2\le r.
\]
Thus~\eqref{us:eq-thick-patch-inclusion} holds. Applying
Lemma~\ref{lem:region-thickness} at $z$ with radius $s/2$ proves
\eqref{us:eq-clipped-thickness}.
\end{proof}

\begin{proof}[Proof of Lemma~\ref{lem:sampling-mass-lower-bound}]
The set $S_r(\gamma)$ contains
$B_{\mathcal M_1}(x,\min\{r,D\})\cap S_{\mathrm{reg}}$.
Lemma~\ref{lem:region-thickness} therefore gives
$\operatorname{vol}_{\mathcal M_1}(S_r(\gamma))>0$.
This volume is finite because $S_r(\gamma)\subseteq S_{\mathrm{reg}}$.
Thus, $\mu_{\gamma,r}$ is well-defined.

For $0<\delta\le\min\{r,D\}$, Lemmas~\ref{lem:tube-volume}
and~\ref{lem:clipped-thickness} give, uniformly in
$v\in S_r(\gamma)$,
\[
\begin{aligned}
\mu_{\gamma,r}(B_{\mathcal M_1}(v,\delta))
&=\frac{
\operatorname{vol}_{\mathcal M_1}
(B_{\mathcal M_1}(v,\delta)\cap S_r(\gamma))}
{\operatorname{vol}_{\mathcal M_1}(S_r(\gamma))}
\\
&\ge\frac{a_{\mathrm{reg}}\delta^d}
{2^d C_{\mathrm{tube}}(\ell+r)r^{d-1}}.
\end{aligned}
\]
This proves~\eqref{eqn:lower-mass-bound-general} with
\begin{equation}\label{eqn:lower-bound-mass-constant}
c_{\mathrm{mass}}
:=\frac{a_{\mathrm{reg}}}{2^d C_{\mathrm{tube}}}>0.
\end{equation}
\end{proof}

\section{Proof of Theorem~\ref{thm:main-conv-thm-body}}\label{sec:proof-of-main-conv-thm}
We state the full convergence theorem; specialization thereof deferred to Appendix~\ref{sec:specialization-for-main-body}.
\begin{theorem}[Fully parameterized convergence guarantees]\label{thm:main-conv-thm}
Suppose Assumptions~\ref{ass:embedded-submanifold} and~\ref{ass:region-of-interest} hold.  Initialize $x_{0} \in \ROI$ and fix $T \geq 1$. Let $c_{1}, \goodProb \in (0,1)$ and define $\widetilde{c}_{1} = c_{1}\goodProb$. Choose $\beta > 1/(\SCf \widetilde{c}_{1})$, and set $\sigma_{0}^{2} := 24f(\widehat{x}_{N})$ and $\sigma_{1}^{2} = 12/\SCf$. Let 
\begin{equation} \label{eqn:sequence-of-stepsize-radius-scaling}
    \eta_{t} := \frac{\beta}{\gamma + t}, \quad 0 \leq t \leq T-1 \qquad \alpha_{t} := \frac{\beta}{(\gamma+t)^{2}}, \quad 1 \leq t \leq T,
\end{equation}
where $\gamma > \max\{1, 2(\beta \SCf \widetilde{c}_{1} - 1)\}$ is chosen so that $\beta/\gamma = \eta_{0} < \min\{c_{1}/(\smoothf \widetilde{\sigma}_{1}^{2}), 1/(\SCf c_{1})\}$, with $\widetilde{\sigma}_{0}^{2}:=\goodProb^{-1}\sigma_{0}^{2} + 4(\ROIdiam/\gamma)^{2}$ and $\widetilde{\sigma}_{1}^{2} := \goodProb^{-1}\sigma_{1}^{2}$. 

Let $c_\Psi > 0$ be the constant defined by \eqref{eqn:explicit-constant-for-lower-bound-control-of-comparison} and choose
\begin{equation} \label{eqn:sequence-of-metric-equivalence}
    0 < \lambda_{t} \leq \frac{c_{\Psi}}{8} \alpha_{t}^{2}, \qquad 1 \leq t \leq T.
\end{equation}
Under Assumption~\ref{ass:embedded-submanifold}\textup{(a)}, let
\[
0<L_{\mathrm{cmp}}
<\min\left\{
\frac{\pi}{\fundBound},
\frac{\pi}{\sqrt{(\kmaxExt)_+}}
\right\},
\qquad
H_{\mathrm{cmp}}
:=
\frac{\fundBound^2}{24}
+\frac{5\kmagExt C_J(L_{\mathrm{cmp}})}{18},
\]
where $C_J(\cdot)$ is defined by
\eqref{eqn:constants-for-Jacobi-control}, and set
\[
C_\rho
:=
\frac12\min\left\{
s_{\mathrm{loc},S_{\mathrm{reg}}},
\frac{L_{\mathrm{cmp}}}{C_{\mathrm{loc},S_{\mathrm{reg}}}},
r_{\mathrm{inj}}^{\mathcal M_1}(S_{\mathrm{reg}}),
\frac{r_{\mathrm{inj}}^{\mathcal M_2}(S_{\mathrm{reg}})}
{C_{\mathrm{loc},S_{\mathrm{reg}}}}
\right\}.
\]

Under Assumption~\ref{ass:embedded-submanifold}\textup{(b)},
instead set
\[
H_{\mathrm{cmp}}:=K_{\mathrm{met}},
\qquad
C_\rho
:=
\frac12\min\left\{
s_{\mathrm{loc},S_{\mathrm{reg}}},
\frac{r_{\mathrm{met}}}{C_{\mathrm{loc},S_{\mathrm{reg}}}},
r_{\mathrm{inj}}^{\mathcal M_1}(S_{\mathrm{reg}})
\right\}.
\]
In both branches, the local compatibility constants are those of
Lemma~\ref{lem:local-intrinsic-extrinsic-compatibility}.
Choose
\begin{align}
0<\epsilon_t\leq
\min\left\{
\frac{\sqrt{\lambda_t}}
{C_{\mathrm{loc},S_{\mathrm{reg}}}\sqrt{H_{\mathrm{cmp}}}},
\frac{\alpha_tD_{\mathrm{reg}}}
{C_{\mathrm{loc},S_{\mathrm{reg}}}},
C_\rho
\right\}.
\label{eqn:connectivity-parameter-main-theorem}
\end{align}
The first term is interpreted as $+\infty$ when
$H_{\mathrm{cmp}}=0$. Define
\begin{equation}\label{eqn:data-covering-main-theorem}
    \delta_{t} = \frac{\lambda_{t} \epsilon_{t}}{4}.
\end{equation}

Set $\widetilde r_{\max}^{(t)}
:=\alpha_tD_{\mathrm{reg}}+\rho_{S_{\mathrm{reg}}}(\epsilon_t)$,
and let $c_{\mathrm{mass}}(\ROI,\ROIdiam,\kmagInt) > 0$ be the constant in
Lemma~\ref{lem:sampling-mass-lower-bound}. Consider running
Algorithm~\ref{algo:graph-algorithm} for $T$ iterations, with
Algorithm~\ref{algo:tube-sample-oracle} drawing conditionally i.i.d. samples
from normalized Riemannian volume on $\widetilde S_t$ at every
nontrivial iteration. Suppose
\begin{align}\label{eqn:per-iteration-sample-rates-formal}
    n_t&\ge\left\lceil
\frac{2^d(D_{\mathrm{reg}}+\widetilde r_{\max}^{(t)})
(\widetilde r_{\max}^{(t)})^{d-1}}
{c_{\mathrm{mass}}(\ROI,\ROIdiam,\kmagInt)\delta_t^d}
\left(\log\mathcal N_t(\delta_t/2)
+\log\frac1\zeta\right)
\right\rceil,\\
\zeta&:=\min\left\{
\frac\beta{T(\gamma+T)},1-p_{\mathrm{good}}\right\},
\qquad 1\le t\le T. \nonumber
\end{align}
where $\mathcal{N}_{t}(s)$ is defined by \eqref{eqn:covering-number-family}. 
Then
\begin{equation}\label{eqn:expected-objective-gap-bound-formal}
    \mathbb{E}[f_{P_{N}}(x_{T}) - f_{P_{N}}(\widehat{x}_{N})] \leq \frac{\nu + \beta \ROIdiam^{2}}{\gamma + T},
\end{equation}
where $\widehat{x}_{N}$ minimizes \eqref{eqn:objective-functional-informal} in $\ROI$ and
\begin{align}\label{eqn:constant-on-objective-decrease}
    \nu \geq \max \Biggr\{ &\gamma (f_{P_{N}}(x_{0}) - f_{P_{N}}(\widehat{x}_{N})), 
    \frac{\ROIdiam^{2}(\ROIraucheC(4\beta^{2} + 3) + 2\beta)^{2}}{\SCf(1-c_{1})^{2}(\beta \SCf \tilde{c}_{1} - 1)} + \frac{\smoothf \beta^{2} \widetilde{\sigma}_{0}^{2}}{2(\beta \SCf \tilde{c}_{1} - 1)},\\
    & \frac{2\ROIdiam^{2}[\beta^{2} + \ROIraucheC(4\beta^{3} + 3\beta) + \smoothf\beta^{2}\sigma_{1}^{2}/2]}{(\beta \SCf \tilde{c}_{1} - 1)} + \frac{ \smoothf\beta^{2}  (\widetilde{\sigma}_{0}^{2}+2(\sigma_{0}^{2}+4(\ROIdiam/\gamma)^{2}))}{2(\beta \SCf \tilde{c}_{1} - 1) } \Biggr\}.  \label{eqn:constant-on-objective-decrease-2}
\end{align}
\end{theorem}

Let $t \geq 0$. For the proof of Theorem~\ref{thm:main-conv-thm} denote the random objective gap and increments by
\begin{equation*}
    \Delta_{t} := f_{P_{N}}(x_{t}) - f_{P_{N}}(\widehat x_{N}), \qquad \text{and} \qquad R_{t+1} := \Delta_{t+1} - \Delta_{t}. 
\end{equation*}
Moreover, denote the conditional probabilities of interest by
\begin{equation*}
    p_{t} := \mathbb{P}(E_{t+1} \mid \mathcal{F}_{t}), \qquad \text{and} \qquad q_{t} := \mathbb{P}(E_{t+1}^{c} \mid \mathcal{F}_{t}).
\end{equation*}
Lastly, for notational simplicity, define the following constants
\begin{gather}
    a := \beta \SCf \tilde{c}_{1}, \qquad 
    b := \frac{\smoothf  \beta^{2}\widetilde{\sigma}_{0}^{2}}{2}, \qquad
    C_{\operatorname{align}} := \ROIraucheC(4\beta^{2} + 3) + 2\beta, \label{eqn:align-constant} \\
    C_{\operatorname{stab}} := \ROIdiam^{2} \left[\beta^{2} + \ROIraucheC(4\beta^{3} + 3\beta) + \frac{\smoothf \beta^{2} \sigma_{1}^{2}}{2}\right] + \frac{\smoothf \beta^{2}}{2}\left[\sigma_{0}^{2} + 4 \left(\frac{\ROIdiam}{\gamma}\right)^{2}\right]. \label{eqn:stab-constant}
\end{gather}
Thus, the lower bounding condition of $\nu$ -- \eqref{eqn:constant-on-objective-decrease} and \eqref{eqn:constant-on-objective-decrease-2} -- may be compactly written as
\begin{equation}\label{eqn:compressed-constant-on-objective-decrease}
    \nu \geq \max \left\{\gamma \Delta_{0}, \frac{b + C_{\operatorname{align}}^{2} D_{\operatorname{reg}}^{2}/[\SCf(1-c_{1})^{2}]}{a-1}, \frac{b+2C_{\operatorname{stab}}}{a-1}\right\}.
\end{equation}
\begin{proof}[Proof of Theorem~\ref{thm:main-conv-thm}] On the successive good events, the local-metric equivalence gives the conditional estimator bounds needed for an SGD-type argument. As a result, the following indicator-weighted bound is obtained; proof deferred to Appendix~\ref{sec:proof-of-sgd-like-rate-on-good-event}.
\begin{lemma}\label{lem:sgd-like-rate-on-good-event}
    Under the assumptions of Theorem~\ref{thm:main-conv-thm}, 
    \begin{equation*}
        \mathbb{E}[\mathbf{1}_{\Omega_{T}} \Delta_{T}] \leq \frac{\nu}{\gamma+ T}. 
    \end{equation*}
\end{lemma}
The iterates $x_{t}$ are in a bounded neighborhood of $\widehat{x}_{N}$ almost surely by Assumption~\ref{ass:region-of-interest}, and thus $f_{P_{N}}$ is bounded. By sampling according to \eqref{eqn:per-iteration-sample-rates-formal}, we ensure $\mathbb{P}(\Omega_{T}^{c})$ is small enough to control the second term. The following is the resulting bound; its proof is deferred to
Appendix~\ref{sec:proof-of-control-of-bad-event-at-appropriate-error-rate}.
\begin{lemma}\label{lem:control-of-bad-event-at-appropriate-error-rate}
    Under the assumptions of Theorem~\ref{thm:main-conv-thm},
    \begin{equation*}
        \mathbb{E}\left[\mathbf{1}_{\Omega_{T}^{c}} \Delta_{T}\right] \leq \frac{\beta \ROIdiam^{2}}{\gamma + T}.
    \end{equation*}
\end{lemma}
We complete the proof by direct application of the preceding two lemmas:
\begin{equation*}
\mathbb{E}[\Delta_{T}] = \mathbb{E}\left[\mathbf{1}_{\Omega_{T}} \Delta_{T}\right] + \mathbb{E}\left[\mathbf{1}_{\Omega_{T}^{c}} \Delta_{T}\right] \leq \frac{\nu + \beta \ROIdiam^{2}}{\gamma + T}.
\end{equation*}
\end{proof}

\subsection{Specialization to Theorem~\ref{thm:main-conv-thm-body}}\label{sec:specialization-for-main-body}

\begin{proof}
We check that the choices satisfy the hypotheses
of Theorem~\ref{thm:main-conv-thm}, then pick $C_{\lambda}$ and $C_{\epsilon}$.

\paragraph{Admissibility of the tunable optimization parameters.}
Set
\[
    c_1=p_{\mathrm{good}}=\frac12,
    \qquad
    \widetilde c_1:=c_1p_{\mathrm{good}}=\frac14,
    \qquad
    \beta:=\frac{8}{\SCf },
    \qquad
    \gamma:=\frac{576\smoothf }{\SCf ^2}.
\]
Then $\beta \SCf \widetilde c_1=2,$
so
    $\beta
    > 1/(\SCf \widetilde{c}_{1}).$ 
Additionally, Lemma~\ref{lem:functional-SC-smooth-parameters} gives
    $0<\SCf \leq 1\leq \smoothf .$
As a result,
\[
    \gamma=\frac{576\smoothf }{\SCf ^2}
    \geq 576
    >
    \max\left\{
        1,\,
        2\bigl(\beta \SCf \widetilde c_1-1\bigr)
    \right\}
    =2.
\]
Since 
    $\sigma_1^2=12/\SCf,$ we have $
    \widetilde\sigma_1^2
    =
    24/\SCf .$
Hence
\[
    \eta_0
    =
    \frac{\SCf }{72\smoothf }
    <
    \frac{\SCf }{48\smoothf }
    =
    \frac{c_1}{\smoothf \widetilde\sigma_1^2}.
\]
Also, since $\SCf \leq1\leq \smoothf $,
\[
    \eta_0
    =
    \frac{\SCf }{72\smoothf }
    \leq \frac1{72}
    <
    \frac{2}{\SCf }
    =
    \frac{1}{\SCf c_1}.
\]
Thus the optimization parameters satisfy the hypotheses of Theorem~\ref{thm:main-conv-thm-body}. 

\paragraph{Uniform $C_{\lambda}$ and $C_{\epsilon}$, Failure Probability, and $C_{\operatorname{conv}}$.}
Let
\[
\overline\alpha:=\sup_{t\ge1}\alpha_t
=\frac{\beta}{(\gamma+1)^2}<1,
\qquad
C_\lambda:=\frac{c_\Psi}{8}.
\]
By Lemma~\ref{lem:comparison-with-quadratic},
$0<c_\Psi\le4/(1+\sqrt5)<2$.
Thus, $0<C_\lambda<1/4$, and
$\lambda_t:=C_\lambda\alpha_t^2$ satisfies
\eqref{eqn:sequence-of-metric-equivalence}.

With $H_{\mathrm{cmp}}$ and $C_\rho$ as in
Theorem~\ref{thm:main-conv-thm}, choose
\[
C_\epsilon:=\min\left\{
\frac{\sqrt{C_\lambda}}
{C_{\mathrm{loc},S_{\mathrm{reg}}}\sqrt{H_{\mathrm{cmp}}}},
\frac{D_{\mathrm{reg}}}{C_{\mathrm{loc},S_{\mathrm{reg}}}},
\frac{C_\rho}{\overline\alpha}
\right\}>0.
\]
Then $\epsilon_t:=C_\epsilon\alpha_t$ satisfies
\eqref{eqn:connectivity-parameter-main-theorem} for every $t\ge1$.
Now, under these choices of $C_{\lambda}$ and $C_{\epsilon}$,
\[
\zeta
=
\min\left\{
\frac{8}{m_fT(\gamma+T)},\frac12
\right\}
=:\zeta_T.
\]
Thus~\eqref{eqn:explicit-sampling-budget} is exactly the explicit sampling condition~\eqref{eqn:per-iteration-sample-rates-formal}. Lemma~\ref{lem:sampling-to-bound-probability-of-bad-event} gives the required conditional success probability at each iteration.
Lastly, by the localization assumption, Assumption~\ref{ass:region-of-interest},
\[
f_{P_N}(\widehat x_N)\leq\frac12D_{\mathrm{reg}}^2,
\qquad
f_{P_N}(x_0)-f_{P_N}(\widehat x_N)
\leq\frac12D_{\mathrm{reg}}^2.
\]
Let $\nu_\star$ be any fixed constant satisfying \eqref{eqn:constant-on-objective-decrease}. Theorem~\ref{thm:main-conv-thm} gives
\[
\mathbb E_{\mathrm{alg}}
\left[
f_{P_N}(x_T)-f_{P_N}(\widehat x_N)
\right]
\leq
\frac{\nu_\star+(8/m_f)D_{\mathrm{reg}}^2}{\gamma+T} = \frac{C_{\mathrm{conv}}}{\gamma + T}, \qquad C_{\mathrm{conv}}
:=
\nu_\star+\frac8{m_f}D_{\mathrm{reg}}^2.
\]
\end{proof}

\subsection{SGD-Like Rate Under Metric Equivalence: Proof of Lemma~\ref{lem:sgd-like-rate-on-good-event}}\label{sec:proof-of-sgd-like-rate-on-good-event}
We further decompose into two sub-events: when $x_{t}$ is close to the minimizer and otherwise. 
To make this decomposition precise, define
\begin{equation}\label{eqn:progress-split-for-good-event-iteration-t}
    \widehat{t} = \gamma + t, \qquad \rho_{t} 
    := \frac{(a-1)\nu - b}{\nu(\widehat{t}-1)} \qquad \text{and} \qquad  \theta_{t} := \frac{\rho_{t}}{2}. 
\end{equation}
By \eqref{eqn:constant-on-objective-decrease}, $(a-1)\nu - b > 0$. Also, the stepsize restriction gives $\gamma > \beta \SCf c_{1} > a$, so $0 < \theta_{t} < \rho_{t} \leq (a-1)/(\widehat{t} - 1) < 1$. Now, define the event
\begin{equation*}
    A_{t} := \left\{ \Delta_{t} \geq \frac{\theta_{t}\nu}{\widehat{t}+1} \right\}.
\end{equation*}

The event $A_t\in\mathcal F_t$ captures being sufficiently far from
the minimizer. On $A_t$, we establish the required second-moment and
alignment estimates conditional on $(\mathcal F_t,E_{t+1})$.
These one-step conditional estimates enable us to bound the corresponding decrease in the objective gap.
$\rho_t$ is chosen so that this contraction removes a $\rho_t$-fraction of the target tolerance at level $\nu/(\widehat{t}+1)$. This is the content of the following; proof is deferred to Appendix~\ref{sec:proof-of-one-step-progress-graph-biased-SGD}.

\begin{lemma}\label{lem:one-step-progress-graph-biased-SGD}
    Suppose the assumptions of Theorem~\ref{thm:main-conv-thm} hold. Let $(x_{t})_{t=0}^{T-1}$ be obtained by Algorithm~\ref{algo:graph-algorithm}. If
    \begin{equation}\label{eqn:induction-hypothesis-for-good-event}
        \mathbb{E}\left[\mathbf{1}_{\Omega_{t}} \Delta_{t} \right] \leq \frac{\nu}{\widehat{t}},
    \end{equation}
    then,
    \begin{equation}\label{eqn:enough-progress-when-SGD-holds}
        \mathbb{E}[\mathbf{1}_{\Omega_{t+1}}\mathbf{1}_{A_{t}} \Delta_{t+1}] \leq  \frac{(1 - \rho_{t})\nu}{\widehat{t}+1}.
    \end{equation}
\end{lemma}  
While we cannot verify that Assumption~\ref{assump:gradient-estimator-assumptions} holds for $g_{t}$ on $A_{t}^{c}\cap \Omega_{t+1}$, we can still control the objective gap using the metric equivalence of event $E_{t+1}$ and the fact that $\Delta_{t}$ is ``small enough.'' The following makes this formal; its proof is deferred to \Cref{sec:proof-of-stability-if-iterate-is-good-enough}.

\begin{lemma}\label{lem:stability-if-iterate-is-good-enough}
    Suppose the assumptions of Theorem~\ref{thm:main-conv-thm} hold.
    For any $0 \leq t \leq T-1$, 
    \begin{equation}
        \mathbb{E}\left[\mathbf{1}_{\Omega_{t+1}}\mathbf{1}_{A_{t}^{c}}  \Delta_{t+1}\right] \leq \frac{\rho_{t} \nu}{\widehat{t}+1}.
    \end{equation}
\end{lemma}

With these two Lemmas in place, the proof for Lemma~\ref{lem:sgd-like-rate-on-good-event} follows from an induction argument.

\begin{proof}[Lemma~\ref{lem:sgd-like-rate-on-good-event}] First, note that when $t = 0$, we have that 
    $\mathbb{E}[1_{\Omega_{0}}\Delta_{0}] = \mathbb{E}[\Delta_{0}] = \Delta_{0} \leq\nu/\gamma$
by \eqref{eqn:compressed-constant-on-objective-decrease}. Therefore, the base case holds. Now, for $0 \leq t \leq T-1$  suppose that 
\begin{equation*}
    \mathbb{E}[\mathbf{1}_{\Omega_{t}}\Delta_{t}] \leq \frac{\nu}{\widehat{t}}
\end{equation*}
Then, it follows from direct application of Lemmas~\ref{lem:one-step-progress-graph-biased-SGD} and~\ref{lem:stability-if-iterate-is-good-enough} that
\begin{align*}
    \mathbb{E}\left[\mathbf{1}_{\Omega_{t+1}}\Delta_{t+1}\right] & = \mathbb{E}[\mathbf{1}_{A_{t}} \mathbf{1}_{\Omega_{t+1}}\Delta_{t+1}] +  \mathbb{E}\left[\mathbf{1}_{A_{t}^{c}} \mathbf{1}_{\Omega_{t+1}} \Delta_{t+1}\right] \leq \frac{(1-\rho_{t})\nu}{\widehat{t}+1} + \frac{ \rho_{t} \nu}{\widehat{t}+1} = \frac{\nu}{\widehat{t}+1},
\end{align*}
completing the induction step. Therefore, by induction Lemma~\ref{lem:sgd-like-rate-on-good-event} is proven.
\end{proof}

\subsubsection{Expected One-Step Improvement: Proof of Lemma~\ref{lem:one-step-progress-graph-biased-SGD}}\label{sec:proof-of-one-step-progress-graph-biased-SGD}

Fix $0\le t\le T-1$. We establish the conditional moment bound
and, on $A_t$, the conditional alignment needed for one-step descent.
Note that the parameter restrictions imply
$\eta_t\le\eta_0<c_1p_{\mathrm{good}}m_f/(12L_f)<1$ and
$0<\alpha_{t+1}<1$.

\paragraph{Conditional second moment and alignment.}
On $E_{t+1}$, the corresponding metric equivalence conditions hold. Then, by
\eqref{eqn:sequence-of-metric-equivalence} and Lemma~\ref{lem:comparison-with-quadratic} give
\[
\frac{1+\lambda_{t+1}}{1-\lambda_{t+1}}
\le1+4\lambda_{t+1}
\le1+\frac{c_\Psi}{2}\alpha_{t+1}^2
<\widetilde\Psi(\ell_{t+1},\alpha_{t+1}\ell_{t+1}).
\]
Thus, Lemma~\ref{lem:conditions-on-metric-equivalence-for-tube-containment} implies every graph-path
vertex lies within $T_{r_{t+1}}(\gamma^{\mathcal M_1}_{x_t\to y_{t+1}})$.

Since $c_{\Psi} \leq 2$ and $\alpha_{t} \in(0,1]$, we have $\lambda_{t+1}\le1/2$. Thus, in conjunction with the
stepsize schedules we have
\begin{equation}\label{us:eq-rounding-bound}
\begin{aligned}
\mathfrak e_{t+1}
&=\frac{\epsilon_{t+1}}{2\eta_t(1-\lambda_{t+1})}
\le\frac{\epsilon_{t+1}}{\eta_t}
\\
&\le
\frac{D_{\mathrm{reg}}}{C_{\mathrm{loc},S_{\mathrm{reg}}}}
\frac{\gamma+t}{(\gamma+t+1)^2}
\le\frac{D_{\mathrm{reg}}}{\widehat{t}}
\le\frac{D_{\mathrm{reg}}}{\gamma}.
\end{aligned}
\end{equation}
Lemma~\ref{lem:sampling-to-bound-probability-of-bad-event} gives $p_t\ge p_{\mathrm{good}}$.
Hence Lemma~\ref{lem:metric-equivalence-gives-second-moment-bounds} implies
\[
\begin{aligned}
\mathbb E[\|g_t\|_{x_t}^2\mid\mathcal F_t,E_{t+1}]
&\le\frac1{p_t}
\left(24f_{P_N}(\widehat x_N)
+\frac{12}{m_f}\|\nabla_{\mathcal M_1}f_{P_N}(x_t)\|_{x_t}^2\right)
+4 \mathfrak e_{t+1}^2
\\
&\le\widetilde\sigma_0^2
+\widetilde\sigma_1^2
\|\nabla_{\mathcal M_1}f_{P_N}(x_t)\|_{x_t}^2.
\end{aligned}
\]

We next verify the alignment condition on $A_t$.
By the failure budget, $q_t\le\zeta\le\eta_t$.
Also, $D_t\le D_{\mathrm{reg}}$ and
$\eta_{t}\mathfrak{e}_{t+1}\le\epsilon_{t+1}\le\alpha_{t+1}D_{\mathrm{reg}}$.
The left-hand side of the sufficient condition in
Lemma~\ref{lem:tube-containment-to-gradient-estimator-control} is therefore bounded by
\[
\begin{aligned}
&C_{\log}
\left[\left(\frac{2\lambda_{t+1}}{1-\lambda_{t+1}}\eta_t
+2\alpha_{t+1}\right)D_t+\eta_{t}\mathfrak{e}_{t+1}\right]
+2q_t\eta_tD_t
\\
&\quad\le
C_{\log}(4\alpha_{t+1}^2\eta_t+3\alpha_{t+1})D_{\mathrm{reg}}
+2\eta_t^2D_{\mathrm{reg}}
\\
&\quad\le
\frac{\eta_tD_{\mathrm{reg}}}{\widehat{t}}
\bigl[C_{\log}(4\beta^2+3)+2\beta\bigr]
=\frac{\eta_tC_{\mathrm{align}}D_{\mathrm{reg}}}{\widehat{t}}.
\end{aligned}
\]
Here we used $\widehat{t}>1$,
$\alpha_{t+1}\le\eta_t/\widehat{t}$, and
$\alpha_{t+1}^2\le\beta^2/\widehat{t}$.
By the definition of $\theta_t$ and~\eqref{eqn:constant-on-objective-decrease},
\[
\frac{\theta_t\nu}{\widehat{t}+1}
=\frac{(a-1)\nu-b}{2(\widehat{t}^2-1)}
\ge
\frac{C_{\mathrm{align}}^2D_{\mathrm{reg}}^2}
{2m_f(1-c_1)^2\widehat{t}^2}.
\]
On $A_t$, the Polyak--\L{}ojasiewicz bound therefore gives
\[
\frac{\eta_tC_{\mathrm{align}}D_{\mathrm{reg}}}{\widehat{t}}
\le(1-c_1)\eta_t
\sqrt{\frac{2m_f\theta_t\nu}{\widehat{t}+1}}
\le(1-c_1)\eta_t\sqrt{2m_f\Delta_t}
\le(1-c_1)\eta_t
\|\nabla_{\mathcal M_1}f_{P_N}(x_t)\|_{x_t}.
\]
Lemma~\ref{lem:tube-containment-to-gradient-estimator-control} now yields, on $A_t$,
\[
\left\langle\nabla_{\mathcal M_1}f_{P_N}(x_t),
\mathbb E[g_t\mid\mathcal F_t,E_{t+1}]\right\rangle_{x_t}
\ge c_1\|\nabla_{\mathcal M_1}f_{P_N}(x_t)\|_{x_t}^2.
\]
All iterates and their minimizing connecting segments lie in
$S_{\mathrm{reg}}$, so the normal-step identities and the smoothness
inequality apply. Consequently, on the $\mathcal F_t$-measurable
event $A_t$, the preceding conditional estimates give
\begin{equation}\label{eqn:one-step-conditional-estimate}
\mathbb E[R_{t+1}\mid\mathcal F_t,E_{t+1}]
\le
-2\eta_tm_f
\left(c_1-\frac{\eta_tL_f\widetilde\sigma_1^2}{2}\right)\Delta_t
+\frac{\eta_t^2L_f\widetilde\sigma_0^2}{2}.
\end{equation}

\paragraph{Using SGD Tools.} First, we begin with the decomposition
\begin{align}
    & \mathbb{E}[\mathbf{1}_{A_{t}}\mathbf{1}_{\Omega_{t+1}}\Delta_{t+1}] \nonumber 
    = \mathbb{E}[\mathbf{1}_{A_{t}}\mathbf{1}_{\Omega_{t}}\mathbf{1}_{E_{t+1}}(\Delta_{t} + R_{t+1})] \nonumber \\
    & \qquad = \mathbb{E}[\mathbf{1}_{A_{t}}\mathbf{1}_{\Omega_{t}}\mathbf{1}_{E_{t+1}}\Delta_{t}] + \mathbb{E}[\mathbf{1}_{A_{t}}\mathbf{1}_{\Omega_{t}}\mathbf{1}_{E_{t+1}} R_{t+1}].\label{eqn:before-using-one-step-improvement}
\end{align}
We bound these two terms separately. 

To deal with the first term $\Delta_{t} \geq 0$, we have that
\begin{equation*}
    \mathbb{E}[\mathbf{1}_{A_{t}}\mathbf{1}_{\Omega_{t}}\mathbf{1}_{E_{t+1}}\Delta_{t}] \leq \mathbb{E}[\mathbf{1}_{A_{t}}\mathbf{1}_{\Omega_{t}}\Delta_{t}],
\end{equation*}
since $\Delta \geq 0$. To deal with the increment term, we have
\begin{align*}
    & \mathbb{E}[\mathbf{1}_{A_{t}} \mathbf{1}_{\Omega_{t}}\mathbf{1}_{E_{t+1}} R_{t+1}]\\
    & \qquad = \mathbb{E}[\mathbf{1}_{A_{t}} \mathbf{1}_{\Omega_{t}}\mathbb{E}[\mathbf{1}_{E_{t+1}} R_{t+1} \mid \mathcal{F}_{t}]]\\
    & \qquad = \mathbb{E}[\mathbf{1}_{A_{t}} \mathbf{1}_{\Omega_{t}}p_{t}\mathbb{E}[R_{t+1} \mid \mathcal{F}_{t}, E_{t+1}]] \\
    & \qquad \leq  - 2\eta_{t} \SCf \left(c_{1} - \frac{\eta_{t}\smoothf \widetilde{\sigma}_{1}^{2}}{2}\right)\mathbb{E}[\mathbf{1}_{A_{t}} \mathbf{1}_{\Omega_{t}}p_{t} \Delta_{t}] + \frac{\eta_{t}^{2} \smoothf \widetilde{\sigma}_{0}^{2}}{2} \mathbb{E}[\mathbf{1}_{A_{t}}\mathbf{1}_{\Omega_{t}}p_{t}]\\
    & \qquad \leq - \eta_{t} \SCf c_{1} \mathbb{E}[\mathbf{1}_{A_{t}} \mathbf{1}_{\Omega_{t}}p_{t}\Delta_{t}] + \frac{\eta_{t}^{2}\smoothf \widetilde{\sigma}_{0}^{2}}{2}. 
\end{align*}
The first equality follows from the tower property and second from definition of conditional probability. The first inequality uses \eqref{eqn:one-step-conditional-estimate}. The next uses
$\eta_tL_f\widetilde\sigma_1^2\le c_1$ and
$\mathbf1_{A_t}\mathbf1_{\Omega_t},p_t\in[0,1]$.

Returning to \eqref{eqn:before-using-one-step-improvement}, it then immediately follows from these  two bounds that: 
\begin{align*}
    & \mathbb{E}[\mathbf{1}_{A_{t}} \mathbf{1}_{\Omega_{t+1}} \Delta_{t+1}]\\
    & \qquad \leq \mathbb{E}[\mathbf{1}_{A_{t}}\mathbf{1}_{\Omega_{t}}\Delta_{t}] - \eta_{t} \SCf c_{1} \mathbb{E}[\mathbf{1}_{A_{t}} \mathbf{1}_{\Omega_{t}} p_{t} \Delta_{t}] + \frac{\eta_{t}^{2} \smoothf \widetilde{\sigma}_{0}^{2}}{2}\\
    & \qquad = \mathbb{E}[\mathbf{1}_{A_{t}} \mathbf{1}_{\Omega_{t}}(1 - \eta_{t} \SCf c_{1} p_{t}) \Delta_{t}] + \frac{\eta_{t}^{2} \smoothf \widetilde{\sigma}_{0}^{2}}{2}\\
    & \qquad \leq \left(1 - \eta_{t}\SCf\tilde{c}_{1} \right) \mathbb{E} \left[\mathbf{1}_{A_{t}} \mathbf{1}_{\Omega_{t}} \Delta_{t}\right] + \frac{\eta_{t}^{2} \smoothf  \widetilde{\sigma}_{0}^{2}}{2}
\end{align*}
where the first inequality follows from applying the previous two bounds to \eqref{eqn:before-using-one-step-improvement}, and the second from the almost-sure lower bound $p_{t} \geq \goodProb$. Since $\mathbf{1}_{A_{t}} \Delta_{t} \leq \Delta_{t}$, applying the induction hypothesis \eqref{eqn:induction-hypothesis-for-good-event} and recalling $\eta_{t} = \beta/\widehat{t}$ results in the following:
\begin{align*}
    &\left(1 - \eta_{t}\SCf\tilde{c}_{1} \right) \mathbb{E} \left[\mathbf{1}_{A_{t}} \mathbf{1}_{\Omega_{t}} \Delta_{t}\right] + \frac{\eta_{t}^{2} \smoothf  \widetilde{\sigma}_{0}^{2}}{2}\\
    & \qquad \leq \left(1 - \frac{\beta \SCf \tilde{c}_{1}}{\widehat{t}}\right)\frac{\nu}{\widehat{t}} + \frac{\beta^{2} \smoothf  \widetilde{\sigma}_{0}^{2}}{2 \widehat{t}^{2}}\\
    & \qquad = \left(1 - \frac{a}{\widehat{t}}\right) \frac{\nu}{\widehat{t}} + \frac{b}{\widehat{t}^{2}}\\
    & \qquad = \left(1 - \rho_{t}\right)\frac{\nu(\widehat{t} - 1)}{\widehat{t}^{2}}  \\
    & \qquad \leq \frac{(1 - \rho_{t}) \nu}{\widehat{t}+1}.
\end{align*}
The last inequality follows from $\widehat{t}^{2} \geq (\widehat{t} + 1)(\widehat{t} - 1)$, finishing the proof.

\subsubsection{Stability Bound: Proof of Lemma~\ref{lem:stability-if-iterate-is-good-enough}}\label{sec:proof-of-stability-if-iterate-is-good-enough}
Fix iteration $0\leq t\leq T-1$. It is enough to prove
\begin{equation}\label{eqn:inequality-to-achieve-stability-result}
    \mathbb E\left[
        \mathbf 1_{\Omega_{t+1}}\mathbf 1_{A_t^c}
        R_{t+1}
    \right]
    \leq
    \frac{(\rho_t-\theta_t)\nu}{\widehat{t}+1}.
\end{equation}
This is because on $A_t^c$ we have $\Delta_t<\theta_{t}\nu/(\widehat{t} + 1)$. Hence \eqref{eqn:inequality-to-achieve-stability-result} implies
\begin{align*}
    \mathbb E[\mathbf 1_{\Omega_{t+1}}\mathbf 1_{A_t^c}\Delta_{t+1}]
    & = \mathbb{E}[\mathbf 1_{\Omega_{t+1}} \mathbf 1_{A_{t}^{c}} \Delta_{t}] + \mathbb{E}[\mathbf 1_{\Omega_{t+1}} \mathbf 1_{A_{t}^{c}} R_{t+1}]
    \leq 
    \frac{\theta_{t} \nu}{\widehat{t}+1}
    +\frac{(\rho_t-\theta_t)\nu}{\widehat{t}+1}
    =\frac{\rho_{t}\nu}{\widehat{t}+1}.
\end{align*}
Henceforth, set $\nabla f_t:=\nabla_{\mathcal M_1}f(x_t)$. 

Before we begin, notice that by Assumption~\ref{ass:region-of-interest}, we have
\begin{equation}\label{eqn:diameter-gradient-bound-stability}
    \|\nabla f_t\| = \left\|- \frac{1}{N} \sum_{r=1}^{N} \log_{x_{t}}(z_{r})\right\|
    \leq
    \frac1N\sum_{r=1}^{N}d_{\mathcal M_1}(x_t,z_r)
    \leq\ROIdiam.
\end{equation}
Additionally, since $f$ satisfies local smoothness on $\ROI$ by Lemma~\ref{lem:functional-SC-smooth-parameters}, we have
\begin{align*}
    R_{t+1} &\leq \langle \nabla f_{t}, \log_{x_{t}}^{\mathcal{M}_{1}}(x_{t+1}) \rangle + \frac{\smoothf}{2}\|\log_{x_{t}}^{\mathcal{M}_{1}}(x_{t+1})\|^{2}
    = - \eta_{t} \langle \nabla f_{t}, g_{t} \rangle + \frac{\smoothf \eta_{t}^{2}}{2}\|g_{t}\|^{2},
\end{align*}
since $\log_{x_{t}}^{\mathcal{M}_{1}}(x_{t+1}) = -\eta_{t} g_{t}$. Therefore, it immediately follows that
\begin{align}
    &\mathbb E\left[
        \mathbf 1_{\Omega_{t+1}}\mathbf 1_{A_t^c}
        R_{t+1}
    \right]
    \leq
    \underbracket{\mathbb E\left[
        \mathbf 1_{\Omega_{t+1}}\mathbf 1_{A_t^c}
        \langle\nabla f_t,-\eta_tg_t\rangle
    \right]}_{(\mathrm{A})}
    +
    \underbracket{\frac{\smoothf\eta_t^2}{2}
    \mathbb E\left[
        \mathbf 1_{\Omega_{t+1}}\mathbf 1_{A_t^c}
        \|g_t\|^2
    \right]}_{(\mathrm{B})}.
    \label{eqn:pre-bound-handling-stability}
\end{align}

\paragraph{Bounding $(\mathrm{A})$.} First, inserting the unbiased estimator $\tilde{g}_{t}$, we have 
\begin{align}
    &\mathbb E [\mathbf{1}_{\Omega_{t+1}} \mathbf{1}_{A_{t}^{c}}\langle \nabla f_{t}, -\eta_{t} g_{t} \rangle] = \mathbb{E}\left[\mathbf{1}_{\Omega_{t+1}} \mathbf{1}_{A_{t}^{c}} \langle \nabla f_{t}, -\eta_{t} \tilde{g}_{t}\rangle\right] + \mathbb{E} \left[\mathbf{1}_{\Omega_{t+1}} \mathbf{1}_{A_{t}^{c}} \langle \nabla f_{t}, -\eta_{t}(g_{t} - \tilde{g}_{t})\rangle\right]\nonumber  \\
    & \quad = \mathbb{E}\left[ \mathbf{1}_{\Omega_{t}} \mathbf{1}_{A_{t}^{c}}\left(-\eta_{t} \left\langle \nabla f_{t}, \mathbb{E}\left[\mathbf{1}_{E_{t+1}} \tilde{g}_{t} \mid \mathcal{F}_{t} \right]\right \rangle \right) \right]\nonumber \\
    & \qquad \qquad 
    + 
    \mathbb{E}\left[\mathbf{1}_{\Omega_{t}} \mathbf{1}_{A_{t}^{c}} \left \langle \nabla f_{t}, -\eta_{t} \mathbb{E}\left[ \mathbf{1}_{E_{t+1}} (g_{t} - \tilde{g}_{t}) \mid \mathcal{F}_{t}\right] \right\rangle \right]. \label{eqn:intermediate-checkpoint-for-stability-bound}
\end{align}
The second equality follows from the tower property and from noting that $\tilde{g}_{t}$, $g_{t}$, and $\mathbf{1}_{E_{t+1}}$ are the only random variables that are not $\mathcal{F}_{t}$-measurable. For the first term, note that
\begin{equation*}
    \nabla f_{t} = \mathbb{E}[\tilde{g}_{t} \mid \mathcal{F}_{t}] = \mathbb{E}[\mathbf 1_{E_{t+1}} \tilde g_{t} \mid \mathcal{F}_{t}] + \mathbb{E}[\mathbf 1_{E_{t+1}^{c}} \tilde g_{t} \mid \mathcal{F}_{t}]
\end{equation*}
since $\tilde{g}_{t}$ is an unbiased estimator of $\nabla f_{t}$. Therefore, it directly follows that 
\begin{align}
    &\mathbb{E}\left[ \mathbf{1}_{\Omega_{t}} \mathbf{1}_{A_{t}^{c}}\left(-\eta_{t} \left\langle \nabla f_{t}, \mathbb{E}\left[\mathbf{1}_{E_{t+1}} \tilde{g}_{t} \mid \mathcal{F}_{t} \right]\right \rangle \right) \right] \nonumber\\
    & \qquad = 
    \mathbb{E}\left[\mathbf{1}_{\Omega_{t}} \mathbf{1}_{A_{t}^{c}}\left(-\eta_{t} \|\nabla f_{t}\|^{2} + \eta_{t} \left\langle \nabla f_{t}, \mathbb{E}\left[\mathbf{1}_{E_{t+1}^{c}} \tilde g_t \mid \mathcal F_{t}\right]\right\rangle  \right]\right] \nonumber \\
    & \qquad \leq \mathbb{E}\left[\mathbf{1}_{\Omega_{t}} \mathbf{1}_{A_{t}^{c}} \eta_{t} \left \langle \nabla f_{t}, \mathbb{E} \left[ \mathbf{1}_{E_{t+1}^{c}} \tilde{g}_{t} \mid \mathcal{F}_{t} \right] \right \rangle \right] \nonumber \\
    & \qquad \leq \eta_{t} \mathbb{E}\left[\|\nabla f_{t} \| \cdot \left \| \mathbb{E} [ \mathbf{1}_{E_{t+1}^{c}} \tilde{g}_{t} \mid \mathcal{F}_{t} ] \right\| \right] \nonumber \\
    & \qquad \leq \eta_{t} \mathbb{E}\left[\|\nabla f_{t} \| \cdot  \mathbb{E} [ \mathbf{1}_{E_{t+1}^{c}} \|\tilde{g}_{t}\| \mid \mathcal{F}_{t} ]  \right] \nonumber\\
    & \qquad \leq \eta_{t} \mathbb{E} \left[ q_{t} D_{t}\|\nabla f_{t} \| \right] \nonumber \\
    & \qquad \leq \frac{\beta^{2}\ROIdiam^{2}}{\widehat{t}^{2}}. \label{eqn:first-term-in-bounding-linear-term-for-stability}
\end{align}
The first inequality follows from dropping the negative term, the second from Cauchy-Schwarz then dropping the indicators, the third Jensen's, the fourth from definition of conditional probability and $\|\tilde{g}_{t}\| = \|\log_{x_{t}}^{\mathcal{M}_{1}}(y_{t+1})\| = d_{\mathcal{M}_{1}}(x_{t}, y_{t+1}) \leq D_{t}$, and the last from Assumption~\ref{ass:region-of-interest} with \eqref{eqn:diameter-gradient-bound-stability} and $q_{t} \leq \beta/(T(\gamma+T)) \leq \eta_{t} = \beta/\widehat{t}$ by Lemma~\ref{lem:sampling-to-bound-probability-of-bad-event} and \eqref{eqn:per-iteration-sample-rates-formal}.

For the second term, first note by Corollary~\ref{cor:implicit-gradient-discrepancy}, with $\lambda_{t+1} = \lambda_{1}^{(t+1)} = \lambda_{2}^{(t+1)}$ and $D_{t} \leq \ROIdiam$,
\begin{equation}
    \|\widetilde g_t - g_{t} \| \leq \frac{C_{\log}}{\eta_{t}} 
\Bigg[
    \left(
    \frac{
    2\lambda_{t+1}
    }{
    1-\lambda_{t+1}} \eta_{t}
    +
    2\alpha_{t+1}
    \right)
    \ROIdiam
    +
    \eta_{t}\mathfrak{e}_{t+1}
\Bigg] \leq  \frac{\ROIraucheC \ROIdiam}{\eta_{t}}(4\alpha_{t+1}^2\eta_t+3\alpha_{t+1}). \label{eqn:norm-bound-on-discrepancy}
\end{equation}
The second inequality follows from
\[
\eta_{t}\mathfrak{e}_{t+1}
=\frac{\epsilon_{t+1}}{2(1-\lambda_{t+1})}
\le\epsilon_{t+1}
\le\alpha_{t+1}D_{\mathrm{reg}},
\]
and $\lambda_{t+1}\le \frac{c_{\Psi}}{8}\alpha_{t+1}^{2} \leq \min\{\alpha_{t+1}^2,1/2\}$. Therefore, 
\begin{align}
    &\mathbb{E}\left[\mathbf{1}_{\Omega_{t}} \mathbf{1}_{A_{t}^{c}} \left \langle \nabla f_{t}, -\eta_{t} \mathbb{E}\left[ \mathbf{1}_{E_{t+1}} (g_{t} - \tilde{g}_{t} )\mid \mathcal{F}_{t}\right] \right\rangle \right]\nonumber\\
    &\qquad\leq \mathbb{E}\left[\mathbf{1}_{\Omega_{t}} \mathbf{1}_{A_{t}^{c}} \eta_{t}\| \nabla f_{t} \| \|\mathbb{E}[\mathbf{1}_{E_{t+1}} (g_{t} - \tilde{g}_{t}) \mid \mathcal{F}_{t}] \|\right]\nonumber\\
    &\qquad\leq \eta_{t} \mathbb{E}\left[ \| \nabla f_{t} \| \|\mathbb{E}[\mathbf{1}_{E_{t+1}} (g_{t} - \tilde{g}_{t}) \mid \mathcal{F}_{t}] \|\right]\nonumber\\
    & \qquad = \eta_{t} \mathbb{E}\left[ p_{t} \| \nabla f_{t} \| \|\mathbb{E}[(g_{t} - \tilde{g}_{t}) \mid \mathcal{F}_{t}, E_{t+1}] \|\right] \nonumber \\
    &\qquad\leq
    \ROIraucheC\ROIdiam^2
    \left(
        4\alpha_{t+1}^{2} \eta_{t} + 3\alpha_{t+1}
    \right)\nonumber\\
    &\qquad\leq
    \frac{\ROIraucheC\ROIdiam^2(4\beta^3+3\beta)}{\widehat{t}^2}.
    \label{eqn:second-term-bound-for-first-term-of-stability}
\end{align}
The first inequality follows from Cauchy-Schwarz, the second from dropping the indicators, the third from direct application of \eqref{eqn:norm-bound-on-discrepancy} after Jensen's and \eqref{eqn:diameter-gradient-bound-stability} giving $\|\nabla f(x_{t})\| \leq \ROIdiam$, and the last from definition of $\alpha_{t+1}$ and $\eta_{t}$, and $\widehat{t} \geq 1$. Applying \eqref{eqn:first-term-in-bounding-linear-term-for-stability} and \eqref{eqn:second-term-bound-for-first-term-of-stability} to \eqref{eqn:intermediate-checkpoint-for-stability-bound} gives
\begin{equation}\label{eqn:bound-on-linear-term-stability}
    \mathbb E [\mathbf{1}_{\Omega_{t+1}} \mathbf{1}_{A_{t}^{c}}\langle \nabla f_{t}, -\eta_{t} g_{t} \rangle] \leq \frac{\ROIdiam^{2}\big(\beta^{2} + \ROIraucheC (4\beta^{3} + 3\beta)\big)}{\widehat{t}^{2}}.
\end{equation}

\paragraph{Bounding (B).} Bounding this term is straightforward:
\begin{align}
    &\frac{\smoothf \eta_{t}^{2}}{2}\mathbb E \left[ \mathbf{1}_{\Omega_{t+1}} \mathbf 1_{A_{t}^{c}} \|g_{t}\|^{2}\right] \leq \frac{L_{f}\eta_{t}^{2}}{2} \mathbb{E}\left[\mathbf{1}_{E_{t+1}} \|g_{t}\|^{2}\right]\nonumber \\
    &\qquad = 
    \frac{L_{f}\eta_{t}^{2}}{2} \mathbb{E}\left[ p_{t}\mathbb{E}[\|g_{t}\|^{2} \mid \mathcal{F}_{t}, E_{t+1}]\right] \nonumber\\
    & \qquad \leq \frac{L_{f} \eta_{t}^{2}}{2}\mathbb{E}\left[24f_{P_N}(\widehat x_N)
 +\frac{12}{m_f}\|\nabla_{\mathcal M_1}f_{P_N}(x_t)\|^2 + 4p_{t}\mathfrak{e}_{t+1}^{2}\right] \nonumber \\
    & \qquad \leq \frac{L_{f} \eta_{t}^{2}}{2}\left(\sigma_{0}^{2} + \ROIdiam^{2}\left(\sigma_{1}^{2} + \frac{4}{\gamma^{2}}\right) \right) \notag\\
    & \qquad \leq \frac{L_{f} \beta^{2}}{2\widehat{t}^{2}}\left(\sigma_{0}^{2} + \ROIdiam^{2}\left(\sigma_{1}^{2} + \frac{4}{\gamma^{2}}\right)\right)\label{eqn:bound-on-quadratic-term-stability-bound}
\end{align}
The equality follows from the tower-property. The first inequality follows from dropping indicators, the second inequality from application of Lemma~\ref{lem:metric-equivalence-gives-second-moment-bounds} and $\lambda_{t+1} \leq 1/2$. The third from definition of $\sigma_{0}^{2}$ and $\sigma_{1}^{2}$, alongside using \eqref{eqn:diameter-gradient-bound-stability} giving $\|\nabla f(x_{t})\| \leq \ROIdiam$ and $p_t\mathfrak e_{t+1}^2
\le \mathfrak e_{t+1}^2
\le \ROIdiam^2/\gamma^2$. 

Applying \eqref{eqn:bound-on-linear-term-stability} and \eqref{eqn:bound-on-quadratic-term-stability-bound} to \eqref{eqn:pre-bound-handling-stability} yields
\begin{align*}
    & \mathbb E\left[
        \mathbf 1_{\Omega_{t+1}}\mathbf 1_{A_t^c}
        R_{t+1}
    \right]
    \leq \frac{\ROIdiam^{2}\big(2\beta^{2} + 2\ROIraucheC (4\beta^{3} + 3\beta)+ \smoothf\beta^{2}\sigma_{1}^{2}\big) + \smoothf \beta^{2}(\sigma_{0}^{2}+ 4(\ROIdiam/\gamma)^{2})}{2\widehat{t}^{2}} = \frac{C_{\operatorname{stab}}}{\widehat{t}^{2}}.
\end{align*}
Finally, by the condition on $\nu$ \eqref{eqn:compressed-constant-on-objective-decrease} we have 
\begin{equation*}
    \frac{C_{\operatorname{stab}}}{\widehat{t}^{2}} \leq \frac{(a-1)\nu - b}{2\widehat{t}^{2}} \leq \frac{(a-1)\nu - b}{2(\widehat{t}^{2} - 1)} = \frac{(\rho_{t} - \theta_{t})\nu}{\widehat{t}+1}.
\end{equation*}
The second inequality follows from $\widehat{t} > 1$. This establishes \eqref{eqn:inequality-to-achieve-stability-result}, completing the proof.

\subsection{Bounding Under Metric Equivalence Failure: Proof of Lemma~\ref{lem:control-of-bad-event-at-appropriate-error-rate}}\label{sec:proof-of-control-of-bad-event-at-appropriate-error-rate}
    By Assumption~\ref{ass:region-of-interest}, every iterate belongs to $\ROI$. Thus, $f(x) \leq \ROIdiam^{2}$ almost surely. Therefore, it immediately follows that 
    \begin{equation}\label{eqn:control-under-metric-equivalence-failure}
        \mathbb{E}\left[\mathbf{1}_{\Omega_{T}^{c}}\Delta_{T}\right] = \mathbb{P}(\Omega_{T}^{c}) \mathbb{E}\left[\Delta_{T} \mid \Omega_{T}^{c}\right] \leq \mathbb{P}(\Omega_{T}^{c}) \ROIdiam^{2}. 
    \end{equation}
    To continue, we note we have the following chain of inequalities:
    \begin{align*}
        &\mathbb{P}(\Omega_{T}^{c}) 
        = \sum_{t=0}^{T-1} \mathbb{E}[\mathbf{1}_{\Omega_{t}}q_{t}]
        \leq \sum_{t=0}^{T-1} \mathbb{E}[q_{t}] 
        \leq T \zeta
        \leq \frac{\beta}{T+\gamma}.
    \end{align*}
    The first equality follows from $\Omega_{T}^{c} = \bigsqcup_{t=0}^{T-1}(\Omega_{t} \cap E_{t+1}^{c})$ and  $\mathbb{P}(\Omega_{t} \cap E_{t+1}^{c}) = \mathbb{E}[\mathbf{1}_{\Omega_{t}}\mathbb{E}[\mathbf{1}_{E_{t+1}^{c}} \mid \mathcal{F}_{t}] ] = \mathbb{E}[\mathbf{1}_{\Omega_{t}} q_{t}]$. The first inequality follows from dropping the indicator, where we use the Tower Property with $\Omega_{t} \in \mathcal{F}_{t}$. The second and third follows from the bounds almost surely bounds on $q_{t}$ from \eqref{eqn:per-iteration-sample-rates-formal} and Lemma~\ref{lem:sampling-to-bound-probability-of-bad-event}. Applying this bound to \eqref{eqn:control-under-metric-equivalence-failure} recovers the desired bound, completing the proof.

\section{Additional Experimental Results}\label{sec:additional-numerical-results}
\begin{figure}
    \centering
    \includegraphics[width=1\linewidth]{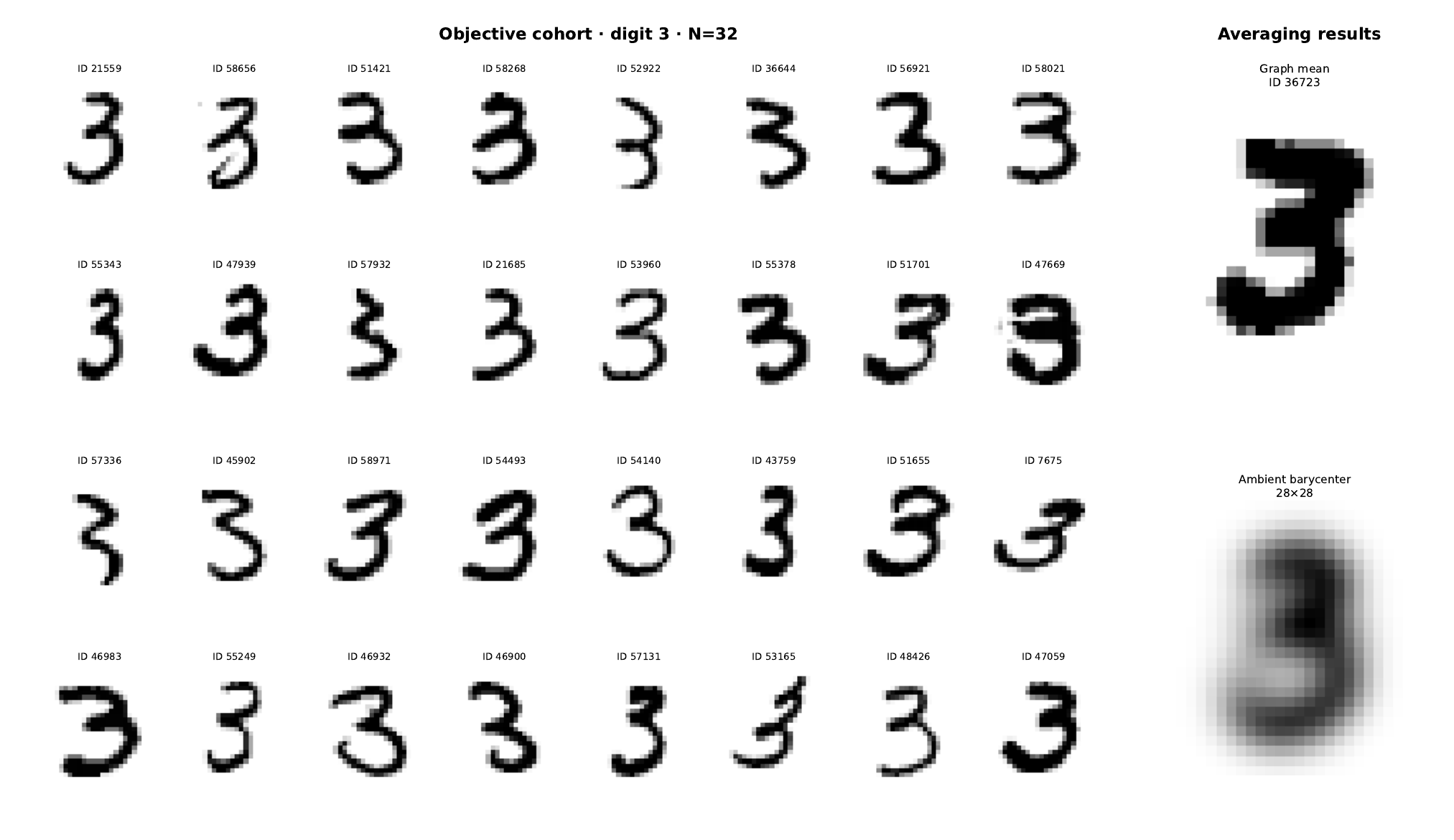}
    \caption{Single-digit averaging for MNIST digit $3$.
Left: the $N=32$ reference images defining the empirical
Fr\'echet objective, labeled by source ID.
Right: the computed graph-based mean and the ambient
Wasserstein barycenter on the $28\times28$ pixel grid.}
    \label{fig:single-3-appendix}
\end{figure}
\begin{figure}
    \centering
    \includegraphics[width=1\linewidth]{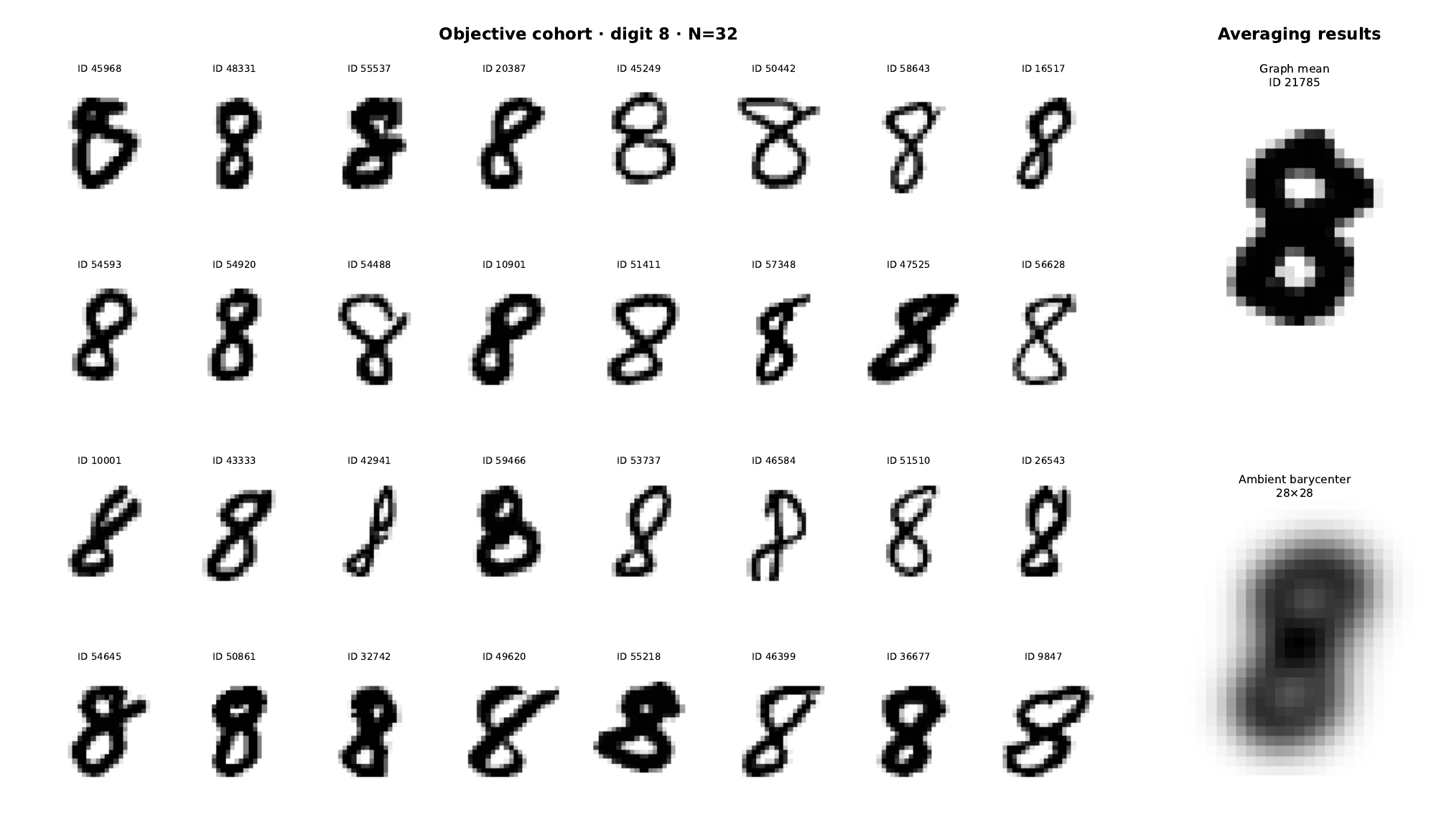}
    \caption{Single-digit averaging for MNIST digit $8$.
Left: the $N=32$ reference images defining the empirical
Fr\'echet objective, labeled by source ID.
Right: the computed graph-based mean and the ambient
Wasserstein barycenter on the $28\times28$ pixel grid.}
    \label{fig:single-8-appendix}
\end{figure}
\begin{figure}
    \centering
    \includegraphics[width=1\linewidth]{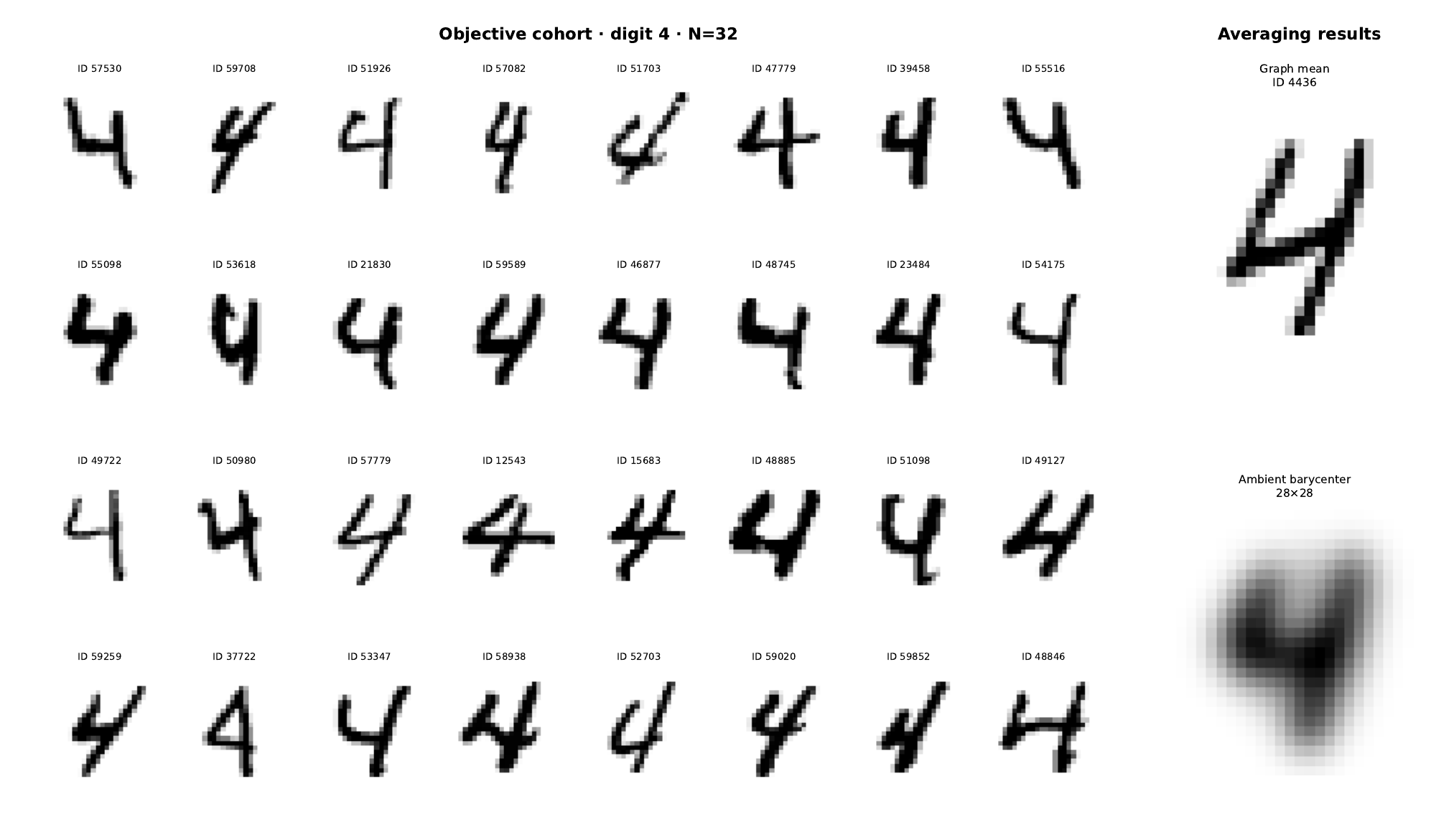}
    \caption{Single-digit averaging for MNIST digit $4$.
Left: the $N=32$ reference images defining the empirical
Fr\'echet objective, labeled by source ID.
Right: the computed graph-based mean and the ambient
Wasserstein barycenter on the $28\times28$ pixel grid.}
    \label{fig:single-4-appendix}
\end{figure}
\begin{figure}
    \centering
    \includegraphics[width=1\linewidth]{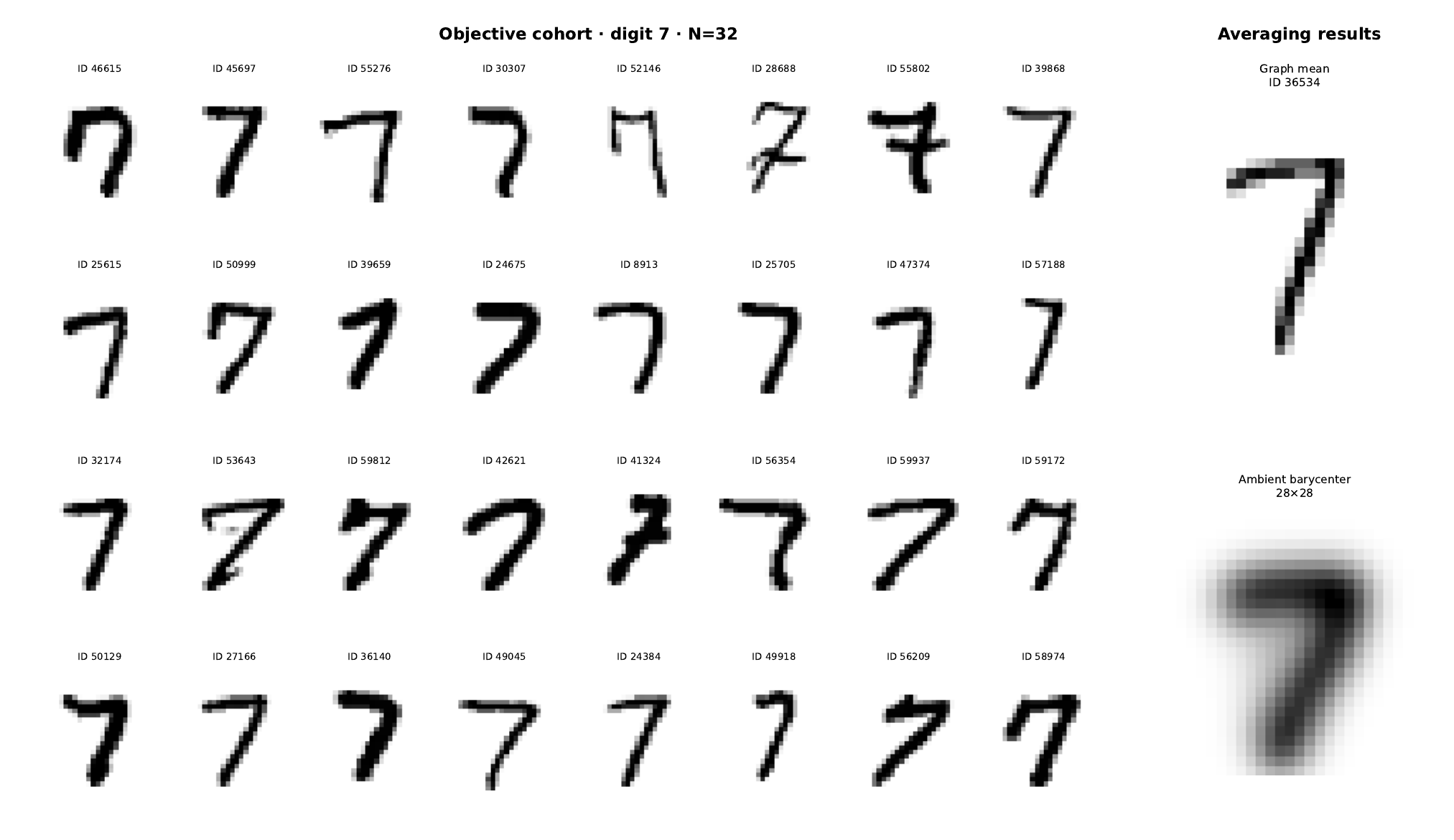}
    \caption{Single-digit averaging for MNIST digit $7$.
Left: the $N=32$ reference images defining the empirical
Fr\'echet objective, labeled by source ID.
Right: the computed graph-based mean and the ambient
Wasserstein barycenter on the $28\times28$ pixel grid.}
    \label{fig:single-7-appendix}
\end{figure}

\begin{figure}
    \centering
    \includegraphics[width=1\linewidth]{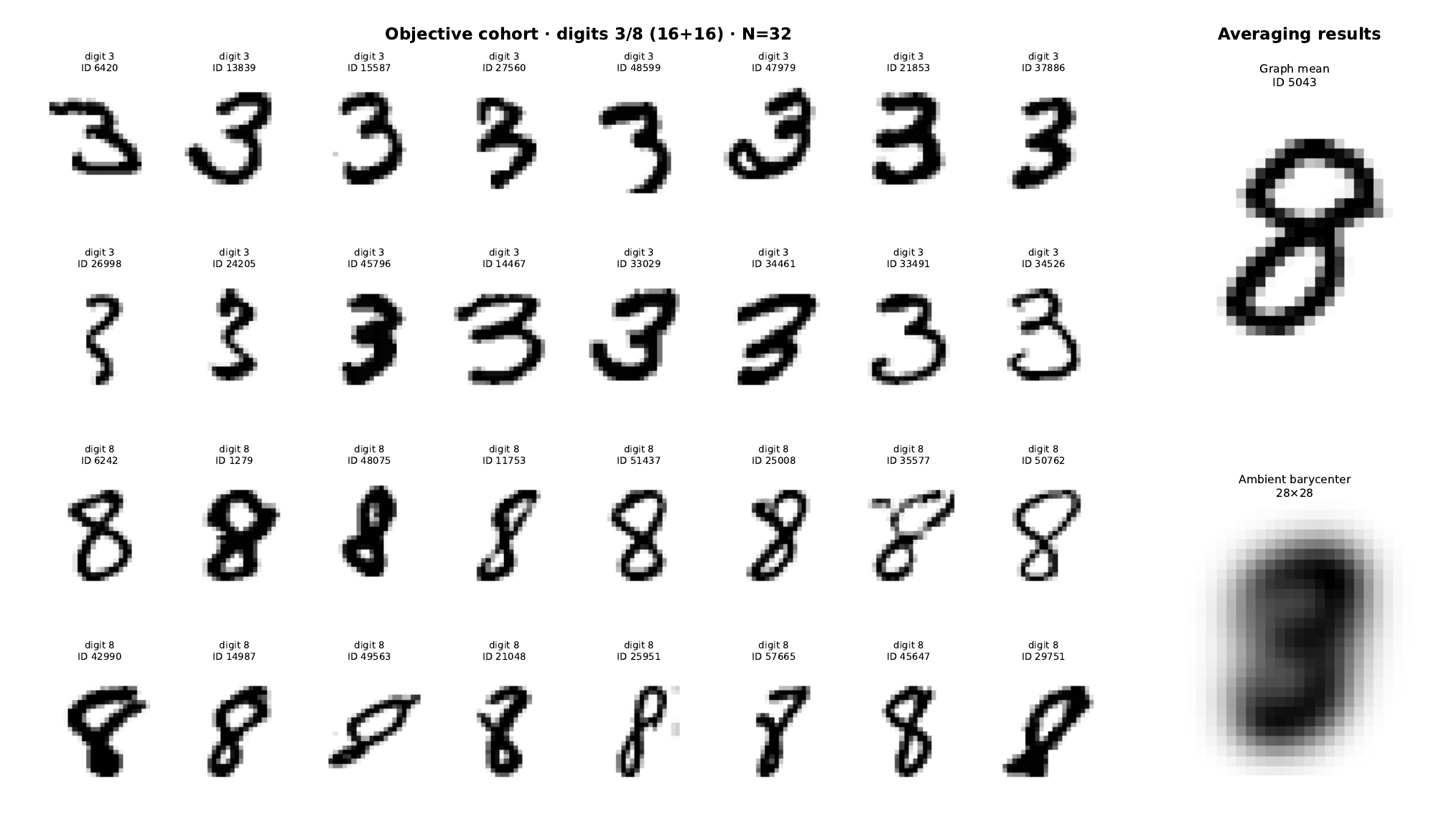}
    \caption{Balanced two-digit averaging for MNIST digits
$3$ and $8$. Left: the $N=32$ reference images, consisting
of $16$ observations from each digit.
Right: the computed graph-based mean and the ambient
Wasserstein barycenter on the $28\times28$ pixel grid.}
    \label{fig:mixed-3-8-appendix}
\end{figure}
\begin{figure}
    \centering
    \includegraphics[width=1\linewidth]{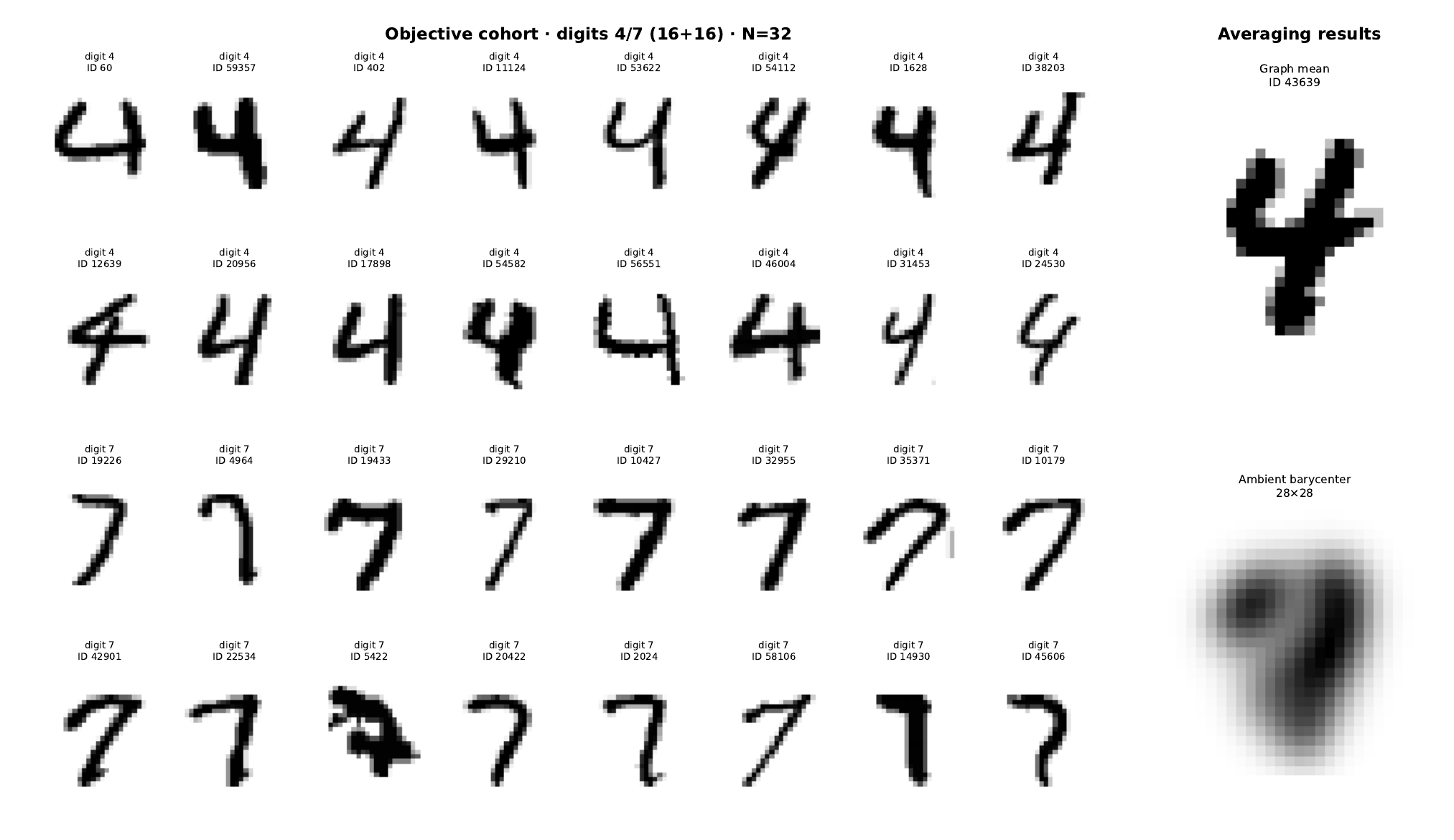}
    \caption{Balanced two-digit averaging for MNIST digits
$4$ and $7$. Left: the $N=32$ reference images, consisting
of $16$ observations from each digit.
Right: the computed graph-based mean and the ambient
Wasserstein barycenter on the $28\times28$ pixel grid.}
    \label{fig:mixed-4-7-appendix}
\end{figure}
\begin{figure}
    \centering
    \includegraphics[angle=90, width=0.86\linewidth]{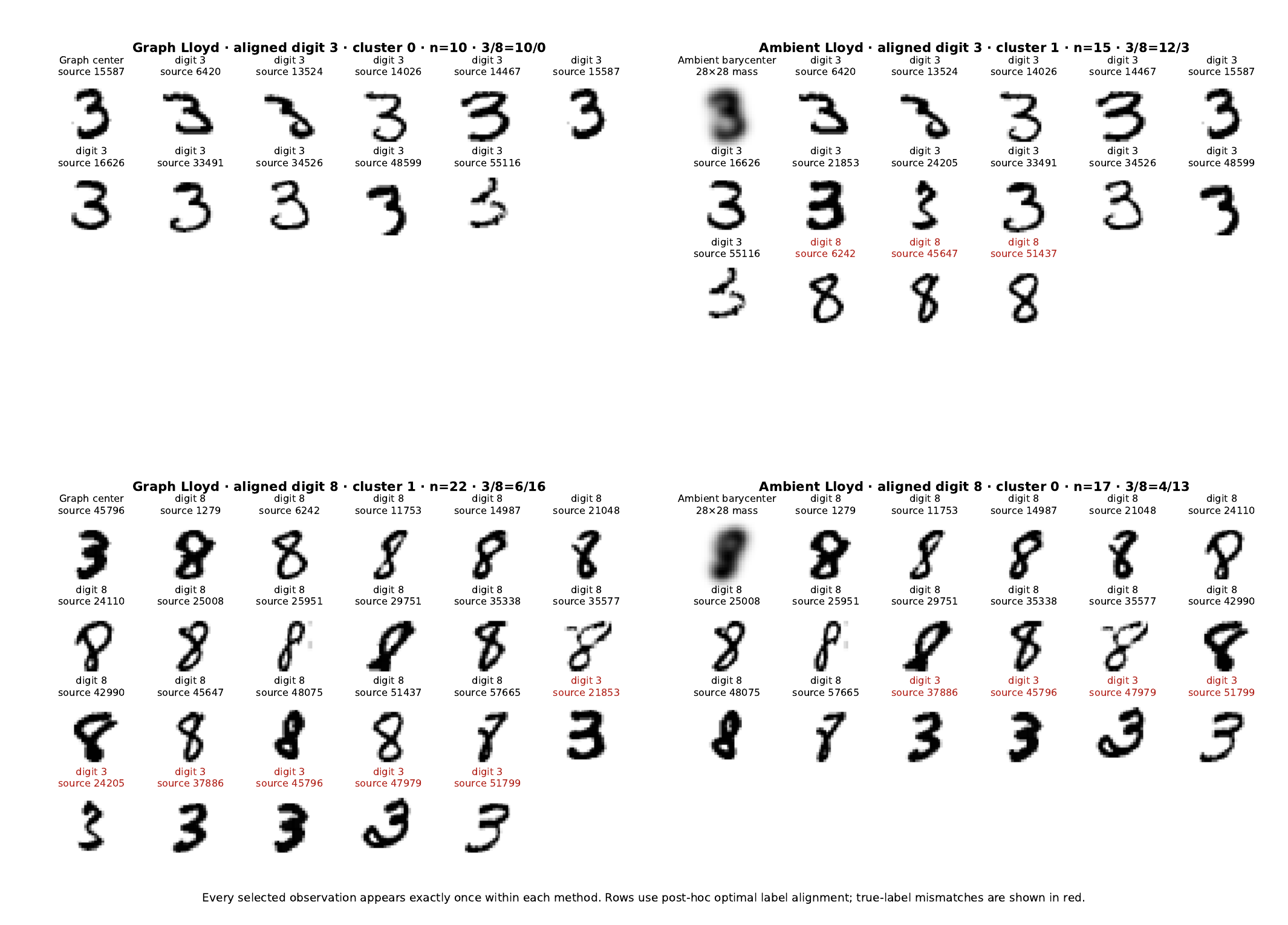}
    \caption{Graph-based and ambient Lloyd clustering for
a balanced MNIST $3/8$ cohort ($N=32$).
Each method displays its computed cluster centers and
the observations assigned to each cluster.
Clusters are aligned post hoc with the true digit labels,
and incorrectly assigned observations are marked in red.
Each observation appears exactly once within each method.
For this example, graph Lloyd correctly assigns $26/32$
observations, compared with $25/32$ for ambient Lloyd.}

    \label{fig:38-appendix}
\end{figure}
\begin{figure}
    \centering
    \includegraphics[angle=90,width=0.86\linewidth]{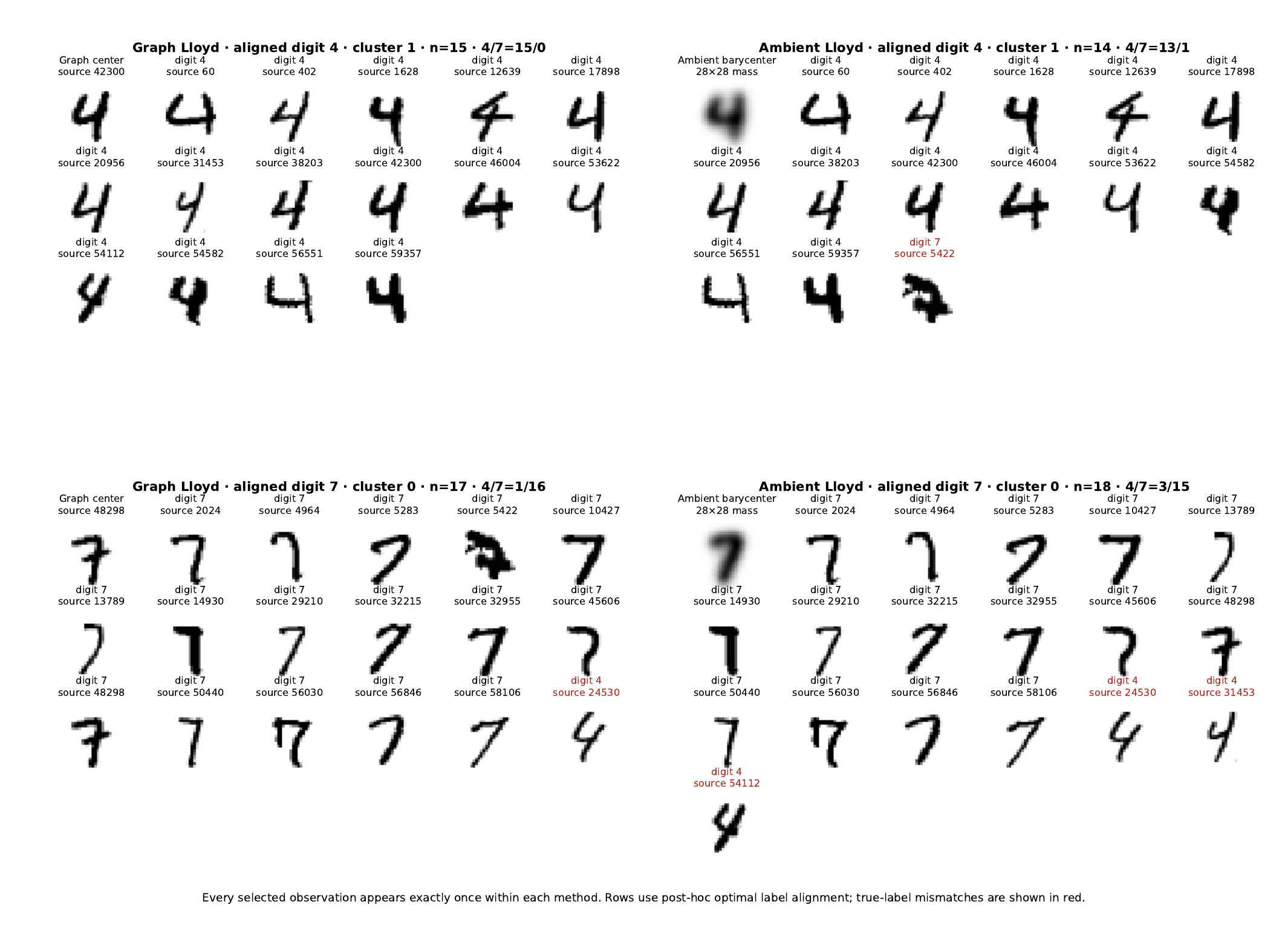}
    \caption{Graph-based and ambient Lloyd clustering for
a balanced MNIST $4/7$ cohort ($N=32$).
Each method displays its computed cluster centers and
the observations assigned to each cluster.
Clusters are aligned post hoc with the true digit labels,
and incorrectly assigned observations are marked in red.
Each observation appears exactly once within each method.
For this example, graph Lloyd correctly assigns $31/32$
observations, compared with $28/32$ for ambient Lloyd.
}
    \label{fig:47-appendix}
\end{figure}

\end{appendix}
\end{document}